\documentclass[12pt,fancyChapter, fancyPart, squeezeCommittee]{these-ubl2}

\bibliography{Biblio}	

\usepackage[nottoc, notlof, notlot]{tocbibind}
\usepackage{enumitem}

\newtheorem{thm}{Théorème}
\numberwithin{thm}{section}
\newtheorem{thmx}{Théorème}

\newtheorem{cor}[thm]{Corollaire}
\newtheorem{lem}[thm]{Lemme}
\newtheorem{prop}[thm]{Proposition}
\newtheorem{df}[thm]{Définition}
\theoremstyle{remark}
\newtheorem{rk}[thm]{Remarque}

\newtheorem{ex}[thm]{Exemple}

\newcommand{\w}{\omega}

\renewcommand{\S}{\Sigma}

\newcommand{\g}{\gamma}

\renewcommand{\L}{\Lambda}

\newcommand{\M}{\mathcal{M}}

\newcommand{\ra}{\rightarrow}

\newcommand{\moduleadh}{\overline{\mathcal{M}}_S (J) \backslash \mathcal{M}_S (J)}

\newcommand{\moduleadhcont}{\overline{\mathcal{M}}_S (J;p) \backslash \mathcal{M}_S (J;p)}

\newcommand{\ev}{\mathrm{ev}}

\newcommand{\ulim}{u_{\infty }}

\DeclareMathOperator{\codim}{codim}

\DeclareMathOperator{\virdim}{virdim}

\DeclareMathOperator{\ind}{ind}

\DeclareMathOperator{\Ima}{Im}

\DeclareMathOperator{\Card}{Card}

\usepackage{times} 

\logouniversite{./Couverture-these/MathSTIC/logo-mathSTIC} 
\scalelogouniversite{1.0} 
\logoetablissement{./Couverture-these/MathSTIC/logo-etablissements/logo-UN_noir}
\scalelogoetablissement{0.5}
\vspace{2.5cm}
\unite{DE L'UNIVERSITÉ DE NANTES}
\ubl{COMUE UNIVERSITÉ BRETAGNE LOIRE}
\ecoledoc{Ecole Doctorale N° 601}
\nomecole{Math\'{e}matiques et Sciences et Technologies de l'Information et de la Communication}
\spec{Sp\'{e}cialit\'{e} : \textit{Mathématiques et leurs Interactions}}
\usepackage{eso-pic}
               \newcommand\BackgroundIm{
               \put(0,0){
               \parbox[b][43cm]{48cm}{%
               \vfill
               \includegraphics[height=21cm,width=33cm,
               keepaspectratio]{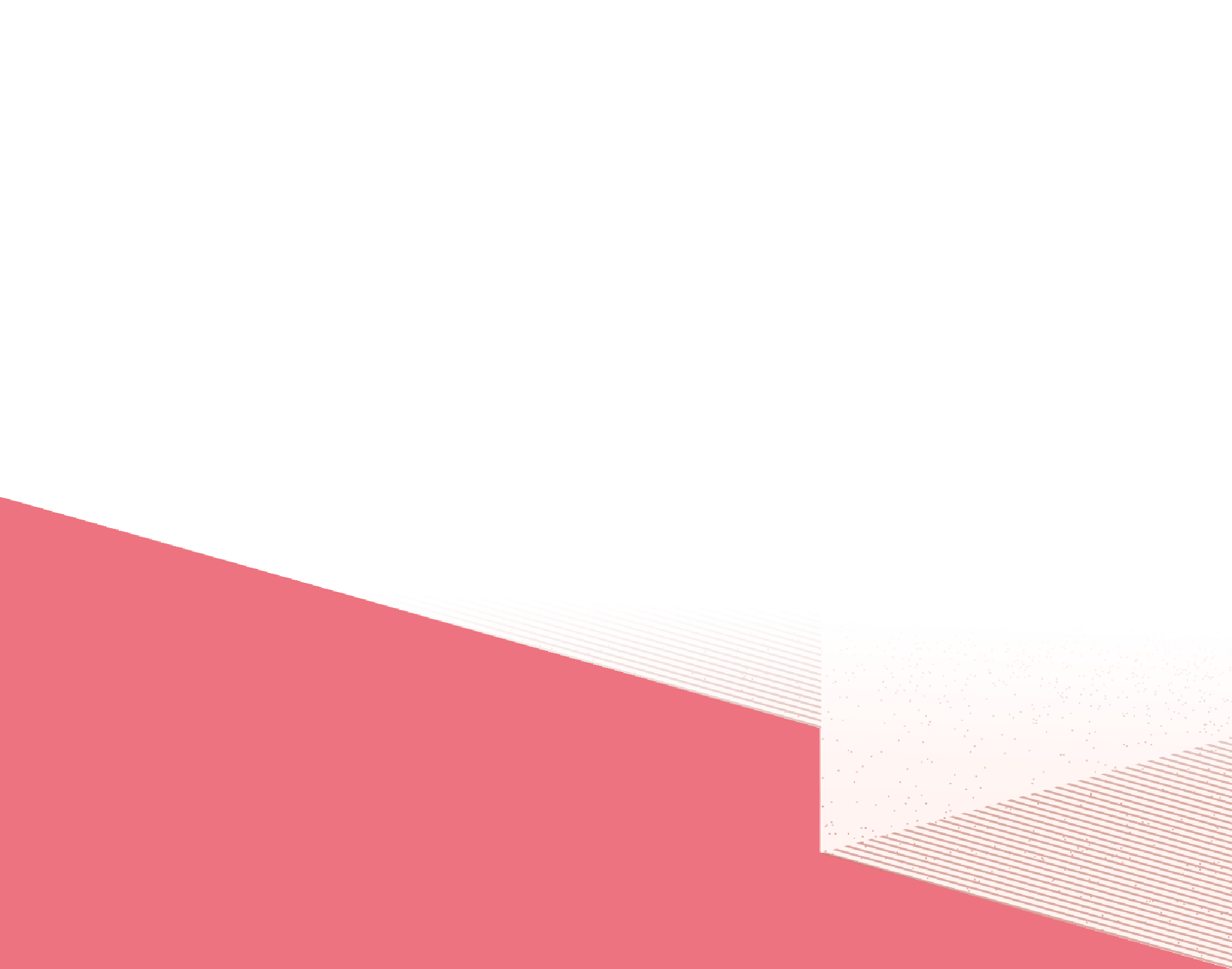}
               \vfill
               }}}
               
\begin{document}
\AddToShipoutPicture{\BackgroundIm}
\author{Fabien KÜTLE}
\title{Courbes symplectiques de haute auto-intersection dans les surfaces symplectiques}
\lieu{Université de Nantes}
\date{\date \textcolor}
\uniterecherche{Unit\'{e} de recherche : Laboratoire de Mathématiques Jean Leray (LMJL)}
  \raggedbottom{
\rapporteur{
{\setlength{\baselineskip}{0.5\baselineskip} 
\begin{description}
\item \hspace{3.33cm}\textbf{Jean-François \textsc{Barraud}}, Professeur, Université de Toulouse
\item \hspace{3.33cm}\textbf{Samuel \textsc{Lisi}}, Associate Professor, Université du Mississippi
\end{description}
}
}
\jury{
{\setlength{\baselineskip}{0.5\baselineskip} 
\begin{description}
\item Pr\'{e}sident : \hspace{1.55cm} \textbf{Klaus \textsc{Niederkrüger}}, Professeur, Université de Lyon 1
\item Examinateurs : \hspace{1.0cm}\textbf{Penka \textsc{Georgieva}}, Maître de conférence, Sorbonne Université
\item \hspace{3.33cm} \textbf{Agnès \textsc{Gadbled}}, Lectrice Hadamard, Université Paris-Saclay
\item \hspace{3.33cm} \textbf{Jean-François \textsc{Barraud}}, Professeur, Université de Toulouse 
\item \hspace{3.33cm} \textbf{Paolo \textsc{Ghiggini}}, Chargé de recherche CNRS, Université de Nantes
\item \hspace{3.33cm} \textbf{Gilles \textsc{Carron}}, Professeur, Université de Nantes
\item Directeurs de thèse : \hspace{0.005cm} \textbf{Marco \textsc{Golla}}, Chargé de recherche CNRS, Université de Nantes
\item \hspace{3.335cm} \textbf{Vincent \textsc{Colin}}, Professeur, Université de Nantes
\end{description}
}
}
}

\date{28 Octobre 2021} 

\maketitle

\chapter*{} 
\ClearShipoutPicture
\thispagestyle{empty}
\clearemptydoublepage
\chapter*{Remerciements}
	Le long chemin semé d'embûches que représente une thèse ne se parcourt pas en étant livré à soi-même. De nombreuses personnes ont contribué à la création du présent manuscrit. Je tiens, comme le veut la tradition nantaise, à leur accorder quelques pages pour les remercier. L'encadrant d'une thèse joue un rôle clé dans cette traversée, et je ne pouvais pas rêver mieux que d'avoir Marco à mes côtés pour me guider et accompagner mes premiers pas dans la recherche. Merci de m'avoir fait découvrir de belles mathématiques, diverses et variées, aux cours de ces trois ans ; de m'avoir laissé la possibilité de choisir parmi plusieurs pistes de recherches pendant les premiers mois de ma thèse ; et plus généralement, de m'avoir fait confiance plus tard vis-à-vis des pistes de recherches que je voulais explorer. Je te remercie également pour ta patience, pour ta disponibilité (à partir de 10h, bien entendu), pour les cafés assortis des petits pains de la boulangerie, et surtout, je te suis très reconnaissant d'avoir été à l'écoute dans les moments les plus difficiles et d'avoir adapté notre manière de travailler en fonction. Merci aussi pour les « bien fait ! » les plus gratifiants que j'ai pu entendre ! Je réalise la chance que j'ai eu d'avoir un encadrant aussi bienveillant et impliqué. Merci Vincent, pour nous avoir suivi régulièrement aux cours des trois dernières années et pour nous avoir conseillé lors des moments charnières. Merci à Jean-François Barraud et Samuel Lisi, qui ont rapporté cette thèse, pour leur lecture attentive, leurs remarques pertinentes qui m'ont permis d'améliorer le manuscrit, ainsi que pour avoir partagé mon enthousiasme vis-à-vis du sujet de cette thèse. Merci aussi aux autres personnes ayant accepté de faire partie de mon jury : Penka Georgieva, Agnès Gadbled, Klaus Niederkrüger, Paolo et Gilles.

 	Je remercie également Chris Wendl pour s'être déplacé à Nantes afin de répondre à mes nombreuses questions. Les discussions qu'on a eues ensemble m'ont été très précieuses pour mener à bien mes raisonnements. Merci d'avoir écrit le livre « Holomorphic Curves in Low Dimensions ». Il m'a servi d'ouvrage de référence du tout début à la toute fin de ma thèse (et je crois d'ailleurs pouvoir me vanter d'être en possession de l'exemplaire à la fois le plus lu et le plus usé à ce jour). L'influence de ce livre sur ce manuscrit est manifeste.

	Je suis également reconnaissant envers de nombreux autres chercheurs du LMJL. À commencer par Paolo, qui a su répondre à mes questions les plus alambiquées, m'orienter vers des références pertinentes, mais aussi repérer très tôt des erreurs dans une version embryonnaire de la démonstration du résultat principal. Merci Samuel également, à la fois pour le cours de Géométrie différentielle que tu as donné quand j'étais en M2, mais aussi, et surtout, pour m'avoir orienté vers Marco et Vincent lors de ma recherche d'un stage en Topologie. Sans ton intervention, cette thèse n'aurait probablement jamais vue le jour ! Merci Friedrich pour ton accueil, ta sympathie et ta disponibilité (merci notamment pour nous avoir aidé Anh et moi à résoudre un exercice ardu du livre de Milnor et Stasheff sur les classes caractéristiques). Merci François pour tous les sages et précieux conseils dont tu nous fais profiter, ainsi que pour la considération et l'estime que tu portes aux jeunes du laboratoire. Merci Erwan pour l'intérêt authentique que tu as porté envers le sujet de thèse et les discussions intéressantes qui ont suivi. La liste est bien trop longue pour mentionner tout le monde individuellement, mais je remercie plus généralement toutes les personnes qui ont dispensés les nombreux enseignements de Master 2, aussi passionnants que divers, auxquels j'ai pu assister pendant les quatre dernières années (je crois que si je pouvais me le permettre, je pourrai passer une majeure partie de mon temps à suivre ce genre de cours). Merci également à toutes les personnes avec qui j'ai eu le plaisir de travailler dans le cadre de mes enseignements : Christoph (qui a fourni un travail titanesque pour assurer des cours de qualité aux étudiants de Licence 1 pendant la crise sanitaire), Xavier, Matthias, Nicolas et Vincent.

	Un laboratoire de mathématiques tomberait très rapidement dans le chaos le plus total sans le travail indispensable de ses secrétaires. J'ai eu le plaisir au cours ma thèse de travailler auprès d'un certain nombre d'entre elles, toutes aussi compétentes qu'amicales et rassurantes : merci Brigitte, Stéphanie, Anaïs, Annick, Béatrice, Alexandra et Caroline. Sans vous je n'aurai pas pu voyager pour assister aux conférences, organiser correctement le séminaire des doctorants ou préparer ma soutenance. Je remercie aussi Eric et Saïd pour leur expertise informatique, Anh et Claude pour nous permettre l'accès à la connaissance et Clémentine et ses collègues pour nous permettre de travailler dans des locaux propres. Merci à Bertrand de la cafet, pour sa sympathie et ses anecdotes quotidiennes lors des pauses café (sa mémoire des noms m'impressionnera toujours). Merci aussi à Marine (qui dispose toujours d'une anecdote de biologie à la fois rigolote et instructive sous le coude), Jenny, Joel et tous les autres membres du Potagex : les séances de jardinage collectif m'ont fait beaucoup de bien, me permettant de me changer les idées, de décompresser après une journée passée le nez dans les livres et les articles, et aussi de cuisiner de bons petits plats avec les légumes récoltés !

	Il est maintenant temps de remercier mes différents compagnons doctorants, postdocs, ou ATER avec qui j'ai eu le plaisir de traverser cette période de ma vie. Sans eux, la vie au labo n'aurait pas été aussi agréable et stimulante. Comme on a pu en faire l'expérience avec la crise du Covid19, une thèse ne se vit pas de la même manière sans les pauses cafés et les repas au RU passés en bonne compagnie ! Je vais tenter de leur rendre hommage dans l'ordre (approximatif) où je les ai rencontrés. 

	Merci tout d'abord Anh, pour les innombrables situations hilarantes que tu créées au quotidien en étant tout simplement toi-même. Les mathématiques n'étant jamais aussi intéressantes et stimulantes que lorsqu'elles sont partagées, je te suis reconnaissant pour les nombreuses séances d'apprentissage en groupe concernant diverses notions d’algèbre, de topologie ou de géométrie. Mais les moments passés hors du labo ne sont pas moins importants à mes yeux. Je te remercie donc aussi pour les sorties, pour l'exploit d'avoir réussi à me faire prononcer (un peu) quelques sonorités vietnamiennes et pour les séances de cuisine (qui ont porté leurs fruits : c'est grâce à toi que je sais manger et cuisiner avec des baguettes, confectionner des nems et des rouleaux de printemps ou encore préparer des plats à base d’herbes douteuses récoltées près de Centrale !). Merci Hélène, tu es la première doctorante du labo que j'ai rencontrée. Ton accueil sympathique a d'emblée réussi à me mettre à l'aise et à m'aider à m'intégrer au groupe. Merci pour les soirées jeux de société passées en ta compagnie et celle de Laura. Merci Matthieu de m'avoir fait confiance en me léguant la responsabilité d'organiser le séminaire des doctorants et de m'avoir initié au Potagex. J'espère que tu continues à concocter ces jeux de mots dont toi seul détient le secret. Je ne doute pas du fait que tu continues à montrer l'exemple avec ton engagement militant et à semer des graines de coriandre partout sur ton chemin. Merci Solène pour m'avoir fait découvrir tout un tas de jeux de par ta passion (ou devrait-on parler d'addiction ?) pour les jeux de société en tout genre et d'avoir essayé sans relâche de m'initier à la taroinche (même si je suis irrécupérable malheureusement !). Merci Caro pour tous les moments de rigolade passés en ta compagnie. Je crois que je n'ai jamais vu une personne aussi décontractée pendant sa troisième année de thèse ! Merci Zeinab d'avoir toujours le mot pour rire et d'avoir su me rappeler régulièrement à quel point tu avais apprécié ma salade quinoa champignon clémentine ! Merci Thomas pour les astuces et conseils d'ancien que tu partages aux doctorants. Il est certain qu'on se serait souvent senti bien déboussolés sans toi ! Merci Hala pour tes conseils concernant la cuisine libanaise et les bons moments passés pendant les rencontres doctorales. Merci Matilde, pour ta sociabilité extraordinaire et ta capacité à répandre le rire et la bonne humeur autour de toi peu importe la situation, que ce soit grâce à tes histoires ou tes maladresses. Merci Germain, pour les nombreuses discussions passionnantes, sérieuses ou ridicules, qu'on a pu avoir autour d'un éventail gigantesque de sujets. Ta capacité à rire de tout, ton habileté à rassurer les doctorants plus jeunes et ta passion pour la pédagogie sont très inspirantes. Il est aussi  toujours rassurant de savoir qu'il existe des gens aussi étourdis que soi, donc merci pour ça également ! Merci Côme pour tous nous déculpabiliser lorsqu'on pense avoir un rythme de vie chaotique et pour partager avec moi cette passion pour la randonnée. On dit souvent que les mathématiciens transforment le café en théorèmes, mais il semblerait que dans ton cas les cookies fonctionnent aussi bien ! Merci Maha pour ton amabilité sans faille, pour la sérénité que tu inspires et pour m'avoir régulièrement fait goûter de bons thés et d'excellentes pâtisseries libanaises. Merci Claire pour les anecdotes loufoques (ma préférée concernant le mystère des morceaux de fromage poussant dans ton jardin) et pour les séances de piscine passées ensemble. Merci Trung pour tes « ba !» en réponse à mes « hi ! ». J'espère que tu garderas toujours ce côté décalé qui fait ton charme. Merci Meissa pour ta vitalité et tes anecdotes à base d'arnaques en ligne et de course-poursuites pour des sombres histoires de terrain de foot à Nantes. Merci Mohamad pour ton calme olympien, ton humour et pour les tapes amicales sur l'épaule. Merci Mael pour ton enthousiasme lors de discussions mathématiques en tout genre. Merci Alex, fidèle compagnon de bureau, pour ta gentillesse inégalable, pour être capable d'animer d'une main de maître les diverses soirées jeux de société ou jeux de rôles, et pour avoir eu le courage de mener une offensive contre l'école doctorale lorsque celle-ci perdait la tête. Merci aussi pour les petites énigmes en tout genre qui pimentent le quotidien. C'était non seulement un véritable plaisir de travailler, de m'amuser et de rire à tes côtés, mais aussi de se soutenir mutuellement dans les moments difficiles. \c Ca me manque déjà ! Merci Amiel pour ton ouverture d'esprit, ton esprit militant et pour avoir toujours des astuces sur à peu près tout. Je suis nostalgique de l'époque où l'armoire du bureau était remplie de jeux de société coopératifs et où la vague de désordre sur ton bureau menaçait d'envahir le mien à chaque instant ! Tu m'as aidé à apprendre à décompresser un peu plus et à éviter de prendre la thèse trop à cœur. Merci Azeddine pour ta quiétude et pour avoir toujours partagé tes friandises quand on se croisait dans le bureau. Merci Arthur pour les discussions toujours enrichissantes qu'on a pu avoir, qu'elles tournent autour de l'écologie ou du sens du métier de chercheur. Merci Karzan pour ta gentillesse et pour les conversations qu'on avait parfois lors de la pause café du midi. Le porte-clé que tu m'as offert est toujours accroché à mon trousseau. Merci Adrian pour ton talent à se poser de bonnes questions en maths et à mobiliser plusieurs personnes pendant de nombreuses heures sur un même problème. C'était un plaisir de réfléchir avec toi au labo, autour d'une bière à Anvers ou en mangeant des falafels (sans oublier de payer bien évidemment !). L'incroyable bestiaire de boules unités sorties de l'imagination sans limite de tes étudiants reste un des plus grands moments de fou rire que j'ai pu avoir au labo ! Merci Antoine de m'avoir accompagné lors de nombreuses séances de natation et d'avoir repris le séminaire des doctorants d'une main de maître ! Ta capacité à comprendre profondément des exposés sur un nombre aussi vaste de sujets m'impressionne véritablement ! Pense de temps en temps à la fameuse fenêtre métabolique quand il faut se faire plaisir autour d'un bon repas ! Merci Samuel d'avoir toujours des curiosités mathématiques amusantes à partager. J'espère que tu tortures un peu moins tes étudiants maintenant ! Merci Anthony pour ton engagement auprès des doctorants au moment où ils en avaient le plus besoin et pour tes discours passionnants autour de sujets concernant aussi bien la politique que l'intelligence artificielle. Merci Adrien d'avoir su remettre de l'ordre dans le bureau (c'était pas gagné !), d'incarner un élément moteur de la vie de laboratoire, d'être toujours partant pour aller boire un verre et de m'avoir appris l'existence du hockey subaquatique. J'espère que tes efforts pour ressusciter le dossier partagé des doctorants porteront un jour leurs fruits ! Merci Jérôme pour ton humour et ton talent inégalé pour mixer les proverbes. Merci Silvère pour la sympathie naturelle que tu dégages, ta bonne humeur et pour les parties passées sur Among Us pendant le deuxième confinement. Merci également à Fakrhi, Agnès, Lisa, Khaled, Mazen, Hai, Lucas, Jean et Ludovic pour toutes les interactions sympathiques et enrichissantes qu'on a pu avoir.

	Je suis également reconnaissant à mes amis doctorants de laboratoires extérieurs, qui m'ont aidé à faire vivre le temps d'une année le séminaire des doctorants, armés de leurs exposés tous aussi passionnants les uns que les autres : merci Octave, Chloé, Victor, Alexandre, Cyril, Théo et bien d'autres encore. Merci Fabski pour ton humour si particulier, pour les nombreux accords, pour les désaccords encore plus nombreux, pour la découverte de l'aquarelle ou de l'harmonica, pour les bons moments passés à l’Artichaut autour de cafés de hipster, ainsi que pour tout ce qu'on a pu partager depuis qu'on se connaît. Merci également Shella, cuisinière des plus talentueuses que je connaisse, de m'avoir fait découvrir les endroits les plus merveilleux de l'Indonésie. Merci Mélo, pour les soirées séries après le boulot, les bêtises avec la fontaine de chocolat ou encore les goûters-bières réconfortants pendant la période du couvre-feu. Merci Céline, Gwen et Cookie pour toutes les bonnes soirées passées ensembles. Merci Guillaume pour être toujours capable de me surprendre en m'embarquant régulièrement dans des activités nouvelles. Tu as réussi à me faire découvrir (et apprécier !) à la fois les escape game, l'accrobranche et plus récemment le théâtre. Merci Loulou pour m'avoir accompagné et soutenu tout au long de cette interminable traversée du désert qu'incarne la période de rédaction en crise sanitaire. Sans nos soirées visio régulières, toujours remplies de discussions sans fin, passionnantes, à des heures souvent complètement déraisonnables, sur un panel absolument incroyable de sujets, je ne sais pas comment je m'en serai sorti. Merci Zuchemin, compagnon de grandes randonnées, pour tes coups de téléphones réguliers qui remontent le moral, qui changent les idées et qui permettent toujours d'approfondir ses réflexions politiques, et merci Emmanuelle pour ton accueil quand je passais parfois presque à l'improviste. Merci Valentin, Eleonore, Benjamin et Richard pour les agréables ballades dans Nantes en votre compagnie. Merci Maman et Mamie pour m'avoir aidé à organiser le pot, et pour n'avoir jamais douté en ma capacité à réussir dans mes études. Merci Ploplo de partager tant de points communs avec moi, d'être capable si souvent de me faire rire aux éclats et d'être une si grande source de fierté à mes yeux. Enfin, pour avoir grandi avec moi et avoir hautement contribué à façonner la personne que je suis, merci Aelma.

\clearemptydoublepage
\frontmatter 
\renewcommand{\contentsname}{Table des matières}
\tableofcontents 
\renewcommand{\listfigurename}{Table des figures}%
\listoffigures

\chapter*{Introduction}
\addcontentsline{toc}{chapter}{Introduction}
\chaptermark{Introduction}
À l'instar de la géométrie complexe, la géométrie symplectique est une géométrie de dimension paire. Mais cette dernière n'est pas aussi rigide que la première et elle présente plutôt des propriétés à mi-chemin entre flexibilité et rigidité. La notion de structure symplectique a été initialement introduite dans le cadre de la mécanique classique, mais elle constitue désormais un objet d'étude intéressant en soit dans le cadre de la topologie différentielle.

Une structure symplectique sur une variété différentielle $M$ de dimension $2n$ est une $2$--forme différentielle $\omega$ (c'est-à-dire une forme bilinéaire antisymétrique sur chaque espace tangent) fermée (c'est-à-dire qui vérifie $d \omega =0$) et non dégénérée (c'est-à-dire que pour tout $u \in TM$, l'application $\omega (u, \cdot)$ n'est pas nulle). On peut notamment y penser comme un moyen de mesurer les aires des surfaces plongées dans $M$. En effet, dans $\mathbb{R}^{2n}$, la forme symplectique standard est donnée par $\omega_{st} = d x_1 \wedge dx_2 + \dots + dx_{2n-1} \wedge dx_{2n}$. Pour une surface $S$ plongée dans $M$, la quantité $\int_S \omega$ représente alors la somme des aires des projections de $S$ sur les $n$ plans respectivement engendrés par les vecteurs $\frac{\partial}{\partial x_{2k-1}}$,$\frac{\partial}{\partial x_{2k}}$.
Dans la même direction, la condition de non dégénérescence assure que la $2n$--forme différentielle $\omega^{ n}$ définit une forme volume sur $M$.

Un théorème central en topologie symplectique, le théorème de stabilité de Moser assure qu'en déformant continument une forme symplectique $\omega$ à travers un chemin de formes symplectiques dans la même classe d'homologie que $\omega$, on ne change pas la structure symplectique sous-jacente (pour être plus précis les deux formes symplectiques aux extrémités du chemin sont isotopes). Ce résultat, qui traduit une certaine forme de flexibilité, a des conséquences retentissantes sur la nature même de ce que désigne la notion de structure symplectique. Parmi les conséquences notables, le fait le plus marquant est qu'à la différence des variétés riemanniennes, on ne peut pas distinguer deux variétés symplectiques à l'aide de mesures locales, comme la courbure par exemple. En effet, le théorème de Darboux affirme que la structure symplectique au voisinage d'un point peut toujours être identifiée à la restriction de la forme symplectique standard $\omega_{st}$ à un ouvert de $\mathbb{R}^{2n}$. Ce sont donc les propriétés globales des variétés symplectiques qui sont intéressantes et on parle plutôt en ce sens de topologie symplectique.

On peut étudier les variétés symplectiques à plusieurs relations d'équivalence près, comme celle donnée par les symplectomorphismes (qui sont des difféomorphismes préservant les formes symplectiques) ou celle donnée par les déformations symplectiques (qui correspondent à une déformation de la forme symplectique à travers un chemin de formes symplectiques). Le volume constitue un invariant d'une variété symplectique considérée à symplectomorphisme près. On pourrait alors légitimement se demander si la notion de symplectomorphisme coïncide avec celle de difféomorphisme préservant le volume.
Si c'est bien le cas en dimension $2$, la situation commence à différer à partir de la dimension $4$. En effet, bien qu'il soit toujours possible d'étirer une boule par un difféomorphisme afin qu'elle soit arbitrairement longue et fine tout en préservant son volume, ce n'est plus forcément le cas dans la catégorie symplectique. Pour mieux le constater, intéressons nous aux formes possibles que peuvent avoir la boule unité dans $(\mathbb{R}^4, \omega_{st})$
$$B^4 = \{(x_1,x_2,x_3,x_4) \in \mathbb{R}^4 \mid x_1^2+x_2^2+x_3^2+x_4^2 \leq 1 \}$$
sous l'action des symplectomorphismes.
Il se trouve qu'il est possible d'aplatir symplectiquement la boule dans certaines directions afin qu'elle soit longue et fine via des symplectomorphismes de la forme
$$\varphi_k : (x_1,x_2,x_3,x_4) \rightarrow \left( 2^{k} x_1, \frac{1}{2^{k}} x_2, 2^{k} x_3, \frac{1}{2^k} x_4 \right),$$
mais des difféomorphismes définis de manières similaires, aplatissant la boule dans d'autres directions, comme les applications suivantes
$$\psi_k : (x_1,x_2,x_3,x_4) \rightarrow \left( \frac{1}{2^k} x_1, \frac{1}{2^k} x_2, 2^k x_3, 2^k x_4 \right),$$
ne sont pas des symplectomorphismes pour des entiers $k$ non nuls (on peut voir aisément qu'ils ne préservent pas l'aire de la projection sur le plan engendré par $\frac{\partial}{\partial x_{1}}$,$\frac{\partial}{\partial x_{2}}$ par exemple). 
Plus généralement, le théorème de non compression de Gromov montre qu'il n'est pas possible de comprimer symplectiquement la boule unité afin qu'elle rentre entièrement dans un cylindre de rayon $r$ strictement plus petit que $1$
$$Z(r) = D^2(r) \times \mathbb{R}^2 = \{(x_1,x_2,x_3,x_4) \in \mathbb{R}^4 \mid x_1^2+x_2^2\leq r^2 \}.$$

\begin{thm}[Théorème de non compression de Gromov,~\cite{Gromov}]
Si $r <1$, il n'existe pas de symplectomorphisme tel que $\varphi (B^4) \subset Z(r)$.
\end{thm}
Ce théorème est symptomatique d'un forme de rigidité symplectique. Un grand intérêt dans le domaine de la topologie symplectique consiste à essayer de tracer la limite entre souplesse et rigidité.

Il existe bien d'autres invariants pour les variétés symplectiques, mais ceux-ci sont en général assez difficiles à définir et proviennent des théories sophistiquées faisant notamment appel aux équations de Seiberg--Witten ou bien à des techniques liées aux courbes pseudoholomorphes.

Dans ce manuscrit, on s'intéresse essentiellement aux variétés symplectiques de dimension $4$, qu'on appellera des \emph{surfaces symplectiques}, ainsi qu'aux sous-variétés symplectiques de dimension $2$ qu'elles contiennent, qu'on appellera des \emph{courbes symplectiques}.

\emph{Dans tous les énoncés évoqués, sauf mention du contraire, les variétés considérées sont différentielles, compactes, connexes, orientées et sans bord.}

Le résultat principal de cette thèse tire son inspiration d'un théorème montré par Hartshorne~\cite{hartshorne} dans le cadre de la géométrie algébrique complexe et d'un théorème de topologie symplectique dû à McDuff~\cite{McDuff}. Il constitue une manifestation de plus de la rigidité des variétés symplectiques et un pas de plus dans la direction des similarités entre structures complexes et symplectiques. On montre plus précisément qu'une courbe symplectiquement plongée dans une surface symplectique détermine à la fois le type difféomorphisme de la surface symplectique et le plongement dans la surface symplectique si son nombre d'auto-intersection est suffisamment grand par rapport à son genre.

\subsection*{Intersection des surfaces dans les variétés de dimension $4$}

%

La topologie de basse dimension est le domaine qui s'intéresse aux variétés de dimension inférieure ou égale à $4$. Les outils utilisés pour les étudier sont de nature différente comparés à ceux utilisés en dimension supérieure. D'un côté, l'astuce de Whitney permet d'affirmer que les cobordismes jouent un rôle important en dimension supérieure ou égale à $5$. De l'autre côté, la topologie des variétés de dimensions $1$, $2$ et $3$ est quant à elle gouvernée par la géométrie. Si la compréhension de la topologie des variétés de dimension $3$ est maintenant considérablement avancée depuis la démonstration de la conjecture de Poincaré et de la conjecture de géométrisation de Thurston par Grigori Perelman en 2003 grâce à des méthodes de géométrie riemannienne, les variétés de dimension $4$ demeurent assez mystérieuses.

Un invariant qui revêt un rôle central dans l'étude des variétés de dimension $4$ est donné par l'étude des intersections entre les surfaces que contient cette variété. Plus précisément, étant donné des surfaces $S_1$ et $S_2$ plongées dans une variété $M$ de dimension $4$, on peut toujours les perturber afin qu'elles s'intersectent transversalement. Comme $2+2=4$, ces surfaces s'intersectent en des points isolés, qui sont en nombre fini par compacité. En comptant ces points en leur attribuant un signe en fonction des orientations en jeu (comme il est souvent l'usage le faire pour définir des invariants en topologie), on obtient un entier, noté $S_1 \cdot S_2$, appelé le nombre d'intersection entre les deux surfaces, qui ne dépend en réalité que des classes d'homologies des deux surfaces :

$$S_1 \cdot S_2 = \sum\limits_{x \in S_1 \cap S_2} \varepsilon (x),$$
où $$\varepsilon(x) = 
\left \{
\begin{array}{rl}
+1 & \text{si} ~ T_x S_1 \oplus T_x S_2 ~ \text{a la même orientation que} ~ T_x M. \\ 
-1 &\text{sinon.}
\end{array}
\right.$$

En prenant deux copies d'une même surface $S$ plongée dans $M$ et en les perturbant pour qu'elles s'intersectent transversalement, on définit de manière similaire le nombre d'auto-intersection $S^2 = S \cdot S$ de la surface $S$.

Comme ces nombres d'intersection ne dépendent que des classes d'homologie des surfaces, on peut s'en servir pour définir une forme $\mathbb{Z}$--bilinéaire sur la partie sans torsion du second groupe d'homologie $H_2(M; \mathbb{Z})$, dont la classe d'isomorphisme constitue alors un invariant topologique de $M$. Cette forme permet d'encoder toutes les informations sur les nombres d'intersection entre les surfaces dans $M$ à l'aide d'une quantité finie d'information. En effet, en prenant une base de la partie sans torsion de $H_2 (M; \mathbb{Z})$ et en considérant tous les nombres d'intersection possibles entre les éléments de cette base, on obtient une matrice carrée à coefficients entiers de taille $b_2(M) \times b_2(M)$. 

Cette forme d'intersection joue un rôle central lorsqu'on étudie les variétés lisses de dimension $4$. Parmi les nombreux résultats notables liés à cette forme d'intersection, on peut citer le théorème de Freedman~\cite{freedman} qui ramène le problème de classification des variété lisses de dimension $4$ simplement connexes à homéomorphisme près à la classification des formes $\mathbb{Z}$--bilinéaires à isomorphisme près ; ou encore le théorème de Donaldson~\cite{DON83,DON87} dont une des conséquences est l'existence de certaines variétés topologiques qui n'admettent pas de structure lisse (comme la variété $E_8$ par exemple). De très nombreuses questions restent ouvertes dans ce domaine, mais la plus célèbre d'entre elles concerne l'existence, ou non, de variétés lisses homéomorphes mais non difféomorphes à la sphère $\mathcal{S}^4$. 

La forme d'intersection ne donne pas assez d'information pour distinguer plus finement les structures lisses (à difféomorphisme près), mais les nombres d'intersections regagnent de l'intérêt lorsqu'on rajoute des structures additionnelles, comme une structure complexe ou une structure symplectique par exemple.

\subsection*{La topologie symplectique, un sujet presque complexe}

En 1939, Weyl introduit l'adjectif symplectique, qui signifie complexe en grec, pour qualifier le groupe $Sp(2n)$, le groupe des automorphismes linéaires réels de $\mathbb{C}^{n}$ conjuguant la multiplication par $i$ à elle-même. Son objectif est de lever une ambigüité. En effet, le nom qui précédait son usage était \og groupe du linéaire complexe \fg{}, ce qui pouvait prêter à confusion avec le groupe des automorphismes complexes. Le nom est particulièrement bien choisi puisque les variétés symplectiques se révèlent être très similaires aux variétés complexes.

Dans les surfaces complexes, la rigidité se manifeste géométriquement par la propriété de positivité d'intersection (entre autres). Deux courbes complexes sans composante commune s'intersectent toujours (localement) positivement, ce qui impose de nombreuses restrictions. Cette propriété de positivité d'intersection n'est plus valable pour les courbes symplectiques dans les surfaces symplectiques, mais l'introduction d'une structure auxiliaire dite presque complexe permet néanmoins de tracer de nombreux parallèles entre les deux catégories.

En effet, toute variété symplectique $(M, \omega)$ peut être munie d'une structure presque complexe, c'est-à-dire d'une application linéaire $J_x$ sur chaque espace tangent en $x \in M$, dont le carré est égal à la multiplication par $-1$ et qui dépend de manière lisse de $x$. Cette structure presque complexe peut même être choisie de sorte à être compatible avec la forme symplectique (c'est-à-dire de manière à ce que le $2$--tenseur $\omega(\cdot , J \cdot)$ définisse un métrique riemannienne) et de manière à préserver les espaces tangents d'une courbe symplectique préalablement choisie. Lorsqu'on peut choisir des coordonnées locales telles que l'application $x \mapsto J_x$ est constante, on dit que la structure presque complexe est intégrable. Dans ce cas, cette notion coïncide avec celle de variété complexe.

En 1985, Gromov introduit la notion de courbes pseudoholomorphes, ce qui marque le début de l'étude moderne de la géométrie symplectique. Celles-ci sont définies comme des applications $u$ d'une surface de Riemann $(\Sigma, j )$ vers une variété presque complexe $(M,J)$ qui vérifient l'équation de Cauchy--Riemann non linéaire 
$$du \circ j =  J \circ du.$$
Les courbes pseudoholomorphes entretiennent un lien étroit avec les courbes symplectiques. En effet, l'image de toute courbe pseudoholomorphe est une courbe symplectique singulière, c'est-à-dire une courbe symplectique dont les singularités sont modelées sur les singularités des courbes complexes ; inversement toute courbe symplectique singulière est l'image d'une courbe pseudoholomorphe. Les courbes pseudoholomorphes dans les surfaces symplectiques manifestent beaucoup de propriétés similaires au courbes complexes dans les surfaces complexes. En effet, elles vérifient la propriété de positivité d'intersection et la formule d'adjonction (qui permet de relier le genre d'une courbe à sa classe d'homologie) et enfin, de telles courbes peuvent être regroupées de manière à constituer des espaces de dimension finie, qu'on appelle espaces de modules.

Les courbes pseudoholomorphes ont donc une utilité toute particulière pour prouver des théorèmes de rigidité et c'est notamment grâce à leurs propriétés que Gromov a démontré son théorème de non compression. La grande majorité (si ce n'est la totalité) des résultats présentés dans ce manuscrit se reposent sur les propriétés de ces courbes. Elles permettent de tracer de nombreux parallèles entre les surfaces complexes et les surfaces symplectiques et ce, bien qu'il existe des variétés symplectiques n'admettant pas de structure complexe intégrable (en revanche les variétés complexes projectives sont toujours munies de formes symplectiques).

\subsection*{Surfaces réglées, surfaces rationnelles}

L'existence des techniques pseudoholomorphes incite à transposer dans le cadre symplectique certains résultats valables pour les surfaces complexes.

En géométrie complexe, l'existence de certaines courbes remplissant des conditions particulières peuvent donner des informations sur la topologie globale des surfaces qui les contiennent. On peut notamment citer le critère d'Enriques pour les surfaces réglées (voir~\cite[Corollary VI.18]{beauville1996complex}). Celui-ci affirme qu'une surface complexe non singulière $V$ est réglée (c'est-à-dire birationnellement équivalente au produit cartésien de la droite complexe projective $\mathbb{C}P^1$ avec une courbe non singulière ; ou encore de dimension de Kodaira égale à $-\infty$) si et seulement si elle contient une courbe complexe irréductible $X$ qui n'est pas un diviseur exceptionnel (c'est-à-dire une courbe rationnelle non singulière d'auto-intersection $-1$) telle que la classe canonique $K_V$ de $V$ vérifie $K_V \cdot X <0 $. Si la courbe $X$ est non singulière, la formule d'adjonction permet de reformuler cette condition en $X^2 > 2g-2$, où $g$ désigne le genre de $X$. Un théorème de Hartshorne va encore plus loin et assure qu'en renforçant la condition sur l'auto-intersection de $X$, on peut aussi connaître la manière dont $X$ est plongée dans la surface. Avant d'énoncer ce théorème, on introduit un peu de terminologie. \'Etant donné deux surfaces complexes $V_1$ et $V_2$, ainsi qu'une courbe complexe $X$, on dit que deux plongements $X \rightarrow V_1$ et $X \rightarrow V_2$ sont \emph{équivalents} s'il existe une application birationnelle $f : V_1 \rightarrow V_2$ qui est un isomorphisme sur un voisinage ouvert de $X$ dans $V_1$ et qui induit l'application identité sur $X$. Une surface est dite géométriquement réglée si c'est l'espace total d'un fibré holomorphe en $\mathbb{C}P^1$ au-dessus d'une courbe complexe non singulière.

\begin{thm}[{\cite[Theorem~4.1]{hartshorne}}]
Soit $X$ une courbe algébrique complexe non singulière de genre $g$ dans une surface algébrique complexe non singulière $V$. On suppose que $X^2 > 4g +5$. Alors $V$ est une surface réglée et le plongement de $X$ dans $V$ est équivalent à une section d'une surface géométriquement réglée. De plus, l'énoncé est également vrai si $X^2 \geq 4g +5$, à l'exception de la cubique non singulière de $\mathbb{C}P^2$ (ou d'un plongement équivalent).
\end{thm}

Du côté symplectique, les techniques pseudoholomorphes ont permis d'obtenir certains résultats de cette nature. Le premier d'entre eux a été démontré par Gromov en 1985. Une surface symplectique est dite minimale si elle ne contient pas de diviseur exceptionnel symplectique (c'est-à-dire une courbe de genre $0$ symplectiquement plongée d'auto-intersection $-1$).

\begin{thm}[\cite{Gromov}]
Soit $(M, \omega)$ une surface symplectique minimale et $S \hookrightarrow M$ une courbe rationnelle c'est-à-dire de genre égal à $0$) symplectiquement plongée dans $(M, \omega)$ telle que $S ^2 = 1$, alors $M$ est difféomorphe à $\mathbb{C}P^2$.
\end{thm}

Ce résultat a ensuite été étendu par McDuff avec l'étude des courbes rationnelles symplectiquement plongées d'auto-intersection positive.

\begin{thm}[\cite{McDuff}]
Soit $(M, \omega)$ une surface symplectique minimale, et $S \hookrightarrow M$ une courbe rationnelle symplectiquement plongée dans $(M, \omega)$ telle que $S \cdot S \geq 0$, alors $M$ est difféomorphe à $\mathbb{C}P^2$ ou à un fibré en sphères symplectique au-dessus d'une surface.
\end{thm}

Puis McDuff a de nouveau élargi son ampleur grâce à l'étude des courbes rationnelles symplectiquement positivement immergées (c'est-à-dire dont tous les points multiples sont des points doubles transverses positifs).

\begin{thm}[Théorème de McDuff pour les sphères symplectiquement positivement immergées, \cite{McDuffImmersed}]
Soit $S$ une courbe rationnelle symplectiquement positivement immergée dans une surface symplectique $(M, \omega)$ telle que $c_1([S]) \geq 2$. Alors $(M, \omega)$ contient également une courbe rationnelle symplectiquement plongée d'auto-intersection positive. En particulier $(M, \omega)$ est une surface symplectiquement réglée.
\end{thm}

Les trois derniers théorèmes ont été démontrés grâce à une étude approfondie des propriétés des espaces de modules de courbes pseudoholomorphes. Le dernier d'entre eux, celui de McDuff sur les courbes rationnelles symplectiquement positivement immergées, permet notamment d'obtenir une caractérisation des surfaces symplectiquement réglées (et des surfaces symplectiques rationnelles) par la présence de certaines courbes, celle-ci donnant suffisamment d'informations pour classifier les surfaces symplectiquement réglées (et les surfaces symplectiques rationnelles) à difféomorphisme près.

D'autres théorèmes similaires ont été obtenus grâce à l'utilisation de la théorie de Seiberg--Witten, comme le théorème suivant par exemple.

\begin{thm}[\cite{Taubes95,Taubes96a,Taubes96b}]
Soit $(M, \omega )$ une surface symplectique et $S$ une courbe symplectiquement plongée dans $(M, \omega)$ de genre $g$. Si $[S]^2 > 2g-2$ et $S$ n'est pas un diviseur exceptionnel, alors $(M,\omega)$ contient une courbe rationnelle symplectiquement plongée d'auto-intersection positive. En particulier $(M, \omega)$ est une surface symplectiquement réglée.
\end{thm}

Si ces derniers résultats sont plus forts, leurs preuves ont le désavantage de ne pas être aussi géométriques et on peut espérer trouver des preuves de tels théorèmes ne faisant appel qu'à des techniques symplectiques et pseudoholomorphes.


Le théorème qui suit est le résultat principal de cette thèse. Il consiste à adapter le théorème de Hartshorne au cadre symplectique. Une surface symplectique est dite relativement minimale par rapport à un sous-ensemble si elle ne contient pas de diviseur exceptionnel symplectique disjoint de ce sous-ensemble.

\begin{thmx}\label{t:A}
Soit $(M, \omega )$ une surface symplectique et $S$ une courbe symplectiquement plongée dans $(M, \omega)$ de genre $g$ telle que $(M, \w)$ est relativement minimale par rapport à $S$. Si $[S]^2 > 4g +5$, alors
la surface symplectique $(M,\omega)$ est un fibré en sphères symplectique au-dessus d'une surface de genre $g$ et $S$ est l'image d'une section de ce fibré.
De plus, l'énoncé est également vrai si on suppose seulement que $[S]^2 \geq 4g +5$, à l'unique exception près des cubiques symplectiques non singulières dans le plan projectif complexe.
\end{thmx}

On donne deux preuves de ce théorème. La première d'entre elles, dont on esquisse seulement les grandes lignes, utilise le résultat de Taubes (et donc la théorie de Seiberg-Witten). La seconde démonstration est plus intéressante et géométrique. Elle ne fait pas appel à la théorie de Seiberg-Witten et repose essentiellement sur des techniques pseudoholomorphes. On y démontre notamment des propriétés de non compacité de certains espaces de modules, assurant l'existence de courbes nodales dont on se servira pour démontrer une version faible du théorème de Taubes (remplaçant l'hypothèse $[S]^2 > 2g-2$ par $[S]^2 \geq 4g$) sans utiliser la théorie de Seiberg--Witten. Cette preuve sera également l'occasion d'étudier les courbes symplectiques singulières (c'est-à-dire dont les singularités sont modelées sur celles des courbes complexes) dans les surfaces réglées.

On présente ensuite une conséquence directe du Théorème~\ref{t:A} concernant la classification de certains remplissages symplectiques.

\subsection*{Remplissages symplectiques}

Les variétés de contact sont les analogues de dimension impaire des variétés symplectiques. Une forme de contact sur une variété différentielle $X$ de dimension $2n+1$ est une $1$--forme différentielle $\alpha$ vérifiant $\alpha \wedge (d \alpha)^n \neq 0$. Une structure de contact sur $X$ est une distribution d'hyperplans $\xi \subset T X$ qui est localement définie comme le noyau d'une forme de contact. Tout comme les variétés symplectiques, les variétés de contact sont naturellement orientées (quand la structure de contact provient d'une forme de contact définie globalement) et ne possèdent pas d'invariant local.

Les variétés de contact apparaissent naturellement comme bords des variétés symplectiques. En effet, étant donné une variété symplectique à bord $(M, \omega)$, si on dispose d'un champ de vecteur de Liouville (c'est-à-dire vérifiant $\mathcal{L}_v \omega = \omega$) défini au voisinage du bord d'une variété symplectique et transverse au bord, alors le bord $\partial M$ est naturellement muni d'une structure de contact. Selon que le champ de vecteur pointe vers l'intérieur ou l'extérieur de $M$, on dit que le bord est concave ou convexe. Lorsqu'une variété de contact $(X, \xi)$ est le bord concave d'une variété symplectique $(M, \omega)$, on dit que $(M, \omega)$ est un bouchon symplectique de $(X, \xi)$. Si en revanche la variété de contact $(X, \xi)$ est le bord convexe de $(M, \omega)$, on parle plutôt de remplissage symplectique fort. Si une variété de contact peut posséder beaucoup de bouchons symplectiques (par exemple, le complémentaire d'une boule de n'importe quelle surface symplectique fermée est un bouchon symplectique de la sphère $\mathcal{S}^3$ munie de sa structure de contact standard), c'est loin d'être le cas pour les remplissages symplectiques. Une question importante en topologie de contact consiste donc à se demander quels peuvent être les remplissages d'une variété de contact donnée, à éclatements et déformations symplectiques près.

Dans ce manuscrit on s'intéresse plus particulièrement au cas des variétés de contact de dimension $3$. Obtenir une classification des remplissages symplectiques d'une variété de contact arbitraire de dimension $3$ semble hors d'atteinte puisque certaines d'entre elles admettent une infinité de remplissages forts deux à deux non homéomorphes ou non difféomorphes~\cite{smith01, OS03, AEMS08}. En revanche on dispose également de résultats de finitude ou d'unicité (à éclatements et déformations symplectiques près). Les résultats allant dans ce sens sont très nombreux. Les premiers d'entre eux concernent $\mathcal{S}^3$~\cite{Gromov, eliashberg91}, $\mathcal{S}^1 \times \mathcal{S}^2$~\cite{eliashberg91}, le fibré cotangent unitaire de $\mathcal{S}^2$~\cite{Hind00} ou les espaces lenticulaires~\cite{McDuff, Hind03, Lisca} par exemple. Plus généralement, Wendl a démontré un résultat de classification des remplissages des structures de contact planaires via des fibrations de Lefschetz~\cite{Wendl10}.

Dans~\cite{Gay11}, Gay et Mark ont montré qu'à un type de courbe $C$ donnée (c'est-à-dire à un genre et un nombre d'auto-intersection) d'auto-intersection strictement positive, on peut associer une variété de contact $(X_C, \xi_C)$ correspondant au bord d'un voisinage tubulaire de $C$ (voir Chapitre~4 pour plus de détails). Le Théorème~\ref{t:A} permet alors d'obtenir le résultat suivant.

\begin{thmx}\label{t:B}
Soit $S$ une courbe de genre $g$ symplectiquement plongée dans une surface symplectique. Si on a $S^2 > 4g+5$ ou bien $S^2 \geq 4g+ 5$ et $g \neq 1$, alors la variété de contact $(X_S, \xi_S)$ admet un unique remplissage symplectique fort minimal à déformation symplectique près, qui est un fibré en disques au-dessus d'une surface de genre $g$. Si en revanche $g=1$ et $S^2 = 9$, la variété de contact $(X_S, \xi_S)$ admet deux remplissages symplectiques forts minimaux à déformation symplectique près.
\end{thmx}

Le Théorème~\ref{t:B} apporte ainsi un point de vue complémentaire au résultat, énoncé dans~\cite{BMV}, selon lequel la variété de contact associée à une courbe symplectiquement plongée de genre $g \geq 2$ et d'auto-intersection $s$ telle que $0 < s \leq 2g -4$ admet beaucoup de remplissages minimaux.

\subsection*{Le problème d'isotopie symplectique}

Dans l'optique de redémontrer complètement le théorème de Taubes sans utiliser la théorie de Seiberg-Witten, la continuité logique consiste à s'intéresser aux courbes symplectiquement plongées de genre $g$ et d'auto-intersection $s$ dans les surfaces symplectiques qui vérifient $2g-2 <  s < 4g$. En particulier, on cherche à savoir si les espaces de modules associés à ces courbes permettent de trouver des courbes nodales. Le premier cas dans l'ordre lexicographique est celui des courbes elliptiques plongées d'auto-intersection $1$ ($g=s=1$). On montre cependant (dans le cas où la surface symplectique considérée est minimale) un résultat de compacité pour les espaces de modules associés qui va à l'encontre de l'objectif recherché. Mais si ce résultat ne nous est d'aucune aide pour prouver l'existence de courbes nodales, il est en revanche idéal pour résoudre des problèmes d'isotopie symplectique.

Le problème d'isotopie symplectique consiste à se demander sous quelles conditions des courbes symplectiques (possiblement singulières, avec des hypothèses sur le genre, le nombre d'auto-intersection, le type et le nombre de singularités) dans une surface kählérienne (c'est-à-dire une surface symplectique munie d'une structure complexe intégrable compatible avec la forme symplectique) sont symplectiquement isotopes à des courbes complexes (par isotopie, on désigne ici une famille continue de courbes symplectiques avec les mêmes types de singularité). Ce genre de problème repose sur des techniques liées au courbes pseudoholomorphes : on choisit une structure presque complexe $J$ qui rend la courbe $C$ en question pseudoholomorphe et on considère un chemin générique de structure presque complexe $(J_t)_{t \in [0,1]}$ qui relie $J$ à la structure complexe intégrable. L'objectif est alors de prouver, pour tout $t \in [0,1]$, l'existence d'une courbe $J_t$--holomorphe $C_t$ qui remplit les conditions souhaitées.

Le premier résultat de cette nature a été démontré par Gromov dans~\cite{Gromov}, il concerne les courbes de degré $1$ et $2$ dans $\mathbb{C}P^2$. Il a ensuite été étendu aux courbes planes de degré $3$ par Sikorav~\cite{Sikorav}, puis aux courbes planes de degré inférieur ou égal à $6$ par Shevchishin~\cite{Shevchishin} et enfin aux courbes planes de degré inférieur ou égal à $17$ par Siebert et Tian~\cite{SiebertTian}.

Du côté des courbes singulières, le problème a été résolu pour la cubique cuspidale par Ohta et Ono~\cite{OO99}. Barraud~\cite{barraud} et Shevchishin~\cite{she04} ont résolu le problème pour les courbes nodales de $\mathbb{C}P^2$ (c'est-à-dire dont les seules singularités sont des points doubles positifs) de genre inférieur ou égal à $4$. Francisco~\cite{Francisco} s'est occupée du cas où les courbes ne possèdent que des points doubles ou des singularités cuspidales simples (c'est-à-dire localement données par l'équation $y^2=x^3$). Enfin, Golla et Starkston~\cite{GS} ont étudié les courbes rationnelles cuspidales de $\mathbb{C}P^2$ jusqu'au degré $5$ et Golla et l'auteur~\cite{GK} ont étendu ce résultat aux courbes rationnelles cuspidales $\mathbb{C}P^2$ de degré $6$ et $7$. L'article~\cite{GS} présente également des configurations de courbes symplectiques qui ne sont pas équisingulièrement isotopes à des configurations de courbes complexes comme la fausse configuration de Pappus, qui est un arrangement de neuf droites symplectiques dans $\mathbb{C}P^2$, ou bien l'exemple donné par Orevkov d'une courbe symplectique singulière irréductible de degré $8$ dans $\mathbb{C}P^2$.

Pour des résultats hors du cadre du plan projectif complexe, Auroux,  Kulikov et  Shevchishin~\cite{AKS07} ont résolu le problème pour les sections des surfaces de Hirzebruch (c'est-à-dire les surfaces géométriquement réglées au dessus d'une courbe rationnelle) et certaines courbes singulières dans les surfaces de Hirzebruch appelées courbes de Hurwitz. Siebert et Tian~\cite{SiebertTian} se sont occupés de toutes les courbes symplectiques dans les surfaces de Hirzebruch pour lesquelles l'application de projection se retreint à une application de degré inférieur ou égal à $7$. Dans~\cite{HI10}, Hind et Ivrii ont résolu le problème localement pour les sections symplectiques des surfaces géométriquement réglées. Plus précisément, ils ont montré que toutes les sections dans une même classe d'homologie qui sont disjointes d'une même section sont symplectiquement deux à deux isotopes. Ce résultat ne permet en revanche pas de conclure en toute généralité. En effet, les courbes symplectiques ne vérifient pas la propriété de positivité d'intersection et par conséquent, rien n'assure l'existence d'une section disjointe à toutes les sections dans une classe d'homologie donnée.

Dans l'autre direction, Smith~\cite{smith} a montré, pour des genres arbitrairement grands, qu'il existe des familles infinies de courbes symplectiques homologues et deux à deux non symplectiquement isotopes dans des surfaces symplectiques non simplement connexes. Il a également montré que tout éclatement d'une surface projective simplement connexe contient une courbe symplectique (connexe) qui n'est isotope à aucune courbe complexe.

On résout dans cette thèse le problème d'isotopie symplectique pour les sections des surfaces géométriquement réglées au-dessus de courbes elliptiques (c'est-à-dire non singulières, de genre $1$). Celui-ci repose essentiellement sur la compacité de l'espace de modules des courbes pseudoholomorphes elliptiques d'auto-intersection $1$ pour une structure presque complexe générique, sur un théorème de transversalité automatique pour les courbes pseudoholomorphes et sur des techniques birationnelles introduites dans~\cite{GS}.

\begin{thmx}\label{t:C}
Soit $T$ une courbe elliptique symplectiquement plongée dans une surface géométriquement réglée au-dessus d'une courbe elliptique. Si $T$ est une section alors elle est symplectiquement isotope à une section complexe. 
\end{thmx}

Notons que contrairement au Théorème~\ref{t:A}, on fait ici des hypothèses globales sur les surfaces qu'on considère.

\subsection*{Organisation du manuscrit}

Le Chapitre 1 introduit les notions fondamentales dont on fera usage dans ce manuscrit. Il couvre de brefs rappels de topologie symplectique, la définition des nombres d'intersection et de la forme d'intersection, l'étude d'intersections locales de fonctions holomorphes. On y présente ensuite les principales propriétés des courbes pseudoholomorphes, en allant de la positivité d'intersection et de la formule d'adjonction, aux théorèmes de transversalité automatique et en passant par des résultats de généricité sur divers espaces de modules. On termine ce chapitre en définissant les opérations d'éclatement et de contraction ainsi que les pinceaux et les fibrations de Lefschetz topologiques.

Le Chapitre 2 expose des résultats centraux concernant les courbes symplectiques rationnelles. On y introduit les notions de surfaces symplectiquement réglées et de surfaces symplectiques rationnelles. On présente également un résumé des preuves du théorème de McDuff et Gromov sur les courbes rationnelles symplectiquement plongée d'auto-intersection positive et du théorème de McDuff sur les sphères positivement symplectiquement immergées avec première classe de Chern plus grande que $2$. On explique aussi comment ces théorèmes permettent de caractériser les surfaces symplectiquement réglées et les surfaces symplectiques rationnelles. La compréhension de ces résultats et de leurs preuves sera utile pour le Chapitre 3. On met l'accent sur les spécificités liées au caractère rationnel des courbes dans les preuves.

Le Chapitre 3 concerne l'étude des courbes symplectiques de genre arbitraire. Il s'ouvre sur une brève explication du résultat et des arguments utilisés par Hartshorne dans~\cite{hartshorne}. Celle-ci est suivie de deux démonstrations du Théorème~\ref{t:A}. La première démonstration, qu'on esquisse seulement, fait appel au théorème de Taubes, qui utilise la théorie de Seiberg--Witten. Le seconde plus géométrique et de nature plus symplectique, se base sur les propriétés de certains espaces de modules de courbes pseudoholomorphes. On termine ce chapitre en montrant l'optimalité de la borne dans le Théorème~\ref{t:A} et en discutant d'éventuelles pistes d'amélioration pour certains résultats intermédiaires.

Dans le Chapitre 4, on donne une présentation rapide de la notion de remplissages symplectiques. On y démontre le Théorème~\ref{t:B}, qui est une application directe du Théorème~\ref{t:A}.

Enfin le Chapitre 5 se focalise sur les courbes elliptiques symplectiquement plongées d'auto-intersection $1$. Cette fois on démontre, sans utiliser le Théorème de Taubes, un résultat de compacité pour les espaces de modules associés (quand la surface symplectique considérée est minimale). Celui-ci a pour conséquence de résoudre le problème d'isotopie symplectique pour les sections des surfaces géométriquement réglées au-dessus de courbes elliptiques, ce qui fait l'objet du Théorème~\ref{t:C}. 
\clearemptydoublepage
\mainmatter 

\chapter{Surfaces symplectiques et courbes pseudoholomorphes}
L'objectif de ce Chapitre est d'introduire les prérequis nécessaires pour comprendre les démonstrations présentées dans les chapitres suivants. On présente notamment de brefs rappels de topologie symplectique dans la Section~\ref{s:varsymp}. On introduit ensuite dans la Section~\ref{s:nbint} la définition des nombres d'intersection entre deux surfaces dans une variété de dimension $4$, ainsi que la définition de la forme d'intersection d'un telle variété. On y présente aussi rapidement quelques notions élémentaires concernant les singularités de courbes algébriques. La Section~\ref{s:Jholo} expose quant à elle les propriétés des courbes pseudholomorphes dans les surfaces symplectiques. En particulier, on y discute de la positivité d'intersection, de la formule d'adjonction, ainsi que de la topologie de divers espaces de modules. Dans la Section~\ref{s:eclat}, on explique les opérations d'éclatement et de contraction symplectiques. Ces opérations sont inspirées de la géométrie algébrique complexe et  permettent notamment de définir la notion de transformation birationnelle entre deux surfaces symplectiques. Enfin dans la Section~\ref{s:pinceauLefschetz}, on introduit les notions de fibration de Lefschetz symplectique et de pinceau de Lefschetz symplectique, qui permettent de décrire de manière efficace la topologie de certaines surfaces symplectiques à l'aide des courbes qu'elles contiennent.

\section{Variétés symplectiques}\label{s:varsymp}
Soit $M$ une variété différentielle de dimension $2n$ munie d'une $2$--forme différentielle $\omega$. On dit que $(M, \omega)$ est une \emph{variété symplectique} si $\omega$ est fermée (c'est-à-dire $d \omega =0$) et non dégénérée (ou de manière équivalente $\omega^n \neq 0$), on appelle alors $\omega$ la \emph{forme symplectique} de $(M, \omega)$.

\begin{rk}
La condition sur la parité de la dimension provient du fait qu'une forme bilinéaire antisymétrique sur un espace vectoriel de dimension impair est nécessairement dégénérée. Notons qu'une variété symplectique est toujours naturellement orientée puisque $\omega^n$ définit une forme volume sur $M$. Enfin, la forme symplectique $\omega$ définit une classe de cohomologie non nulle $[\omega] \in H^2(M; \mathbb{R})$, donc $b_2 (M) \geq 1$. Par exemple, pour $n \geq 2$, $\mathcal{S}^{2n}$ la sphère de dimension $2n$ n'admet pas de structure symplectique.
\end{rk}

\begin{ex}
\leavevmode
\begin{enumerate}
\item Pour l'espace vectoriel $\mathbb{R}^{2n}$ muni des coordonnées $(x_1, y_1, \dots, x_n, y_n)$ , la $2$--forme différentielle $\omega_{st} = \sum_{i_=1}^n x_i \wedge y_i$ est une forme symplectique, appelée la forme symplectique standard sur $\mathbb{R}^{2n}$.
\item On considère l'espace vectoriel $\mathbb{C}^{n+1}$, muni des coordonnées $(z_1, \dots,z_n)$, $z_j = p_j +i q_j$, pour $j \in \{1, \dots, n \}$. La forme symplectique standard $\omega_{st}$ passe au quotient pour la relation d'équivalence sur $\mathbb{C}^{n+1} \backslash \{0 \}$ donnée par la multiplication par un scalaire complexe non nul. On obtient alors une forme symplectique $\omega_{FS}$ sur le quotient $\mathbb{C}P^n$, qu'on appelle la forme de Fubini--Study.
\item Une forme volume sur une surface orientée est une forme symplectique.
\item Étant donné deux variétés symplectiques $(M_1,\omega_1)$ et $(M_2, \omega_2)$, la variété $M_1 \times M_2$ munie de la $2$--forme différentielle $\omega_1 \oplus \omega_2 =\pi_1^*\omega_1 + \pi_2^*\omega_2$ est une variété symplectique, où $\pi_1$ et $\pi_2$ désignent respectivement la projection sur le premier facteur et la projection sur le second facteur.
\end{enumerate}
\end{ex}

Deux variétés symplectiques $(M, \omega)$ et $(M', \omega')$ sont dites \emph{symplectomorphes} s'il existe un difféomorphisme $\Phi : M \rightarrow M'$ tel que $\Phi^*\omega' = \omega$, on dit alors que $\Phi$ est un \textit{symplectomorphisme}. Cette notion donne lieu à une relation d'équivalence sur l'ensemble des variétés symplectiques. Une autre relation d'équivalence qu'on va considérer dans ce manuscrit est donnée par les déformations symplectiques. On dit qu'une variété symplectique $(M', \omega')$ est une \emph{déformation symplectique} de $(M, \omega)$ s'il existe un difféomorphisme $\Phi : M \rightarrow M'$, et $(\omega_t)_{t \in [0,1]}$ un chemin de formes symplectiques sur $M$ satisfaisant $\omega_0 = \omega$, et $\omega_1 = \Phi^* \omega'$.

Une grande différence entre les variétés symplectiques et les variétés riemanniennes par exemple, est l'absence d'invariant local (comme la courbure). L'étude des variétés symplectiques se rapproche donc plus de celle des variétés différentielles et relève alors plutôt du domaine de la topologie. Cette propriété est une conséquence des théorèmes suivants.

\begin{thm}[Théorème de stabilité de Moser]\label{t:Moser}
Soit $(\omega_t)_{t \in [0,1]}$ une famille de formes symplectiques sur une variété sans bord $M$ telle que pour tout $t \in [0,1]$, $[\omega_t] = [\omega_0]$. Alors il existe une isotopie $(\varphi_t)_{t \in [0,1]}$ telle que $\varphi_0 = \mathrm{id}_M$, et pour tout $t \in [0,1]$, $\varphi_t^* \omega_t = \omega_0$.
\end{thm}

Le théorème de stabilité de Moser permet en effet de démontrer le théorème de Darboux, qui affirme que toute variété symplectique est localement symplectomorphe à l'espace euclidien muni de sa forme standard.

\begin{thm}[Théorème de Darboux]
Soit $(M, \omega)$ une variété symplectique. Pour tout point $p \in M$, il existe un voisinage ouvert $U$ de $p$ dans $M$ tel que $(U, \omega_{\vert U})$ est symplectomorphe à un ouvert $V \subset \mathbb{R}^{2n}$ muni de la forme symplectique standard $\omega_{st\vert V}$.
\end{thm}

\begin{rk}
On pourrait définir de manière équivalente les variétés symplectiques en demandant que $(M, \omega)$ soit localement symplectomorphe à $(\mathbb{R}^{2n}, \omega_{st})$ et que les changements de cartes soient des symplectomorphismes.
\end{rk}

On dit qu'une sous-variété $S$ d'une variété symplectique $(M, \omega)$ est symplectique si $\omega_{\vert S}$ est non dégénérée. Par conséquent, $(S, \omega_{\vert S})$ est également une variété symplectique et, si $M$ et $S$ sont compactes, possède un volume $\int_S \omega$ non nul. Ainsi, la classe d'homologie de $S$ ne peut être nulle, ce qui traduit une certaine forme de rigidité des sous-variétés symplectiques comparées aux sous-variétés différentielles (une sous-variété symplectique compacte d'une variété compacte ne peut être incluse dans une boule par exemple). On dit également qu'une variété $S$ plongée dans $M$ est \emph{symplectiquement plongée} dans $(M, \omega)$ lorsque le tiré en arrière de la forme symplectique $\omega$ par le plongement définit une forme symplectique sur $S$. En revanche, comme la condition de non dégénérescence de $\omega_{\vert S}$ est une condition ouverte, une perturbation $\mathcal{C}^\infty$ arbitrairement petite d'une sous-variété symplectique reste une sous-variété symplectique, ce qui dénote une certaine forme de flexibilité.

Dans ce manuscrit, on s'intéresse tout particulièrement aux variétés symplectiques de dimension $2$ et $4$. Dans tout ce qui suit, afin d'être en phase avec la terminologie concernant les variétés complexes,  on appellera \emph{surface symplectique} toute variété symplectique de dimension $4$. De même, le terme \emph{courbe symplectique} désignera toute variété symplectique de dimension $2$. \emph{Sauf mention du contraire, toutes les variétés considérées seront différentielles, compactes, connexes, orientées, sans bord et tous les difféomorphismes préserveront les orientations.} On conclut cette section en présentant un théorème de Weinstein qui nous assure que la structure symplectique est entièrement déterminée au voisinage d'une sous-variété symplectique (ici énoncé pour les surfaces symplectiques).

\begin{thm}[Théorème du voisinage symplectique, \cite{WEINSTEIN1971329}]
Soit $(M, \omega)$, $(M', \omega')$ deux surfaces symplectiques et $S \subset M$, $S' \subset M'$ deux courbes symplectiquement plongées. On suppose qu'il existe un isomorphisme $\nu S \rightarrow \nu S'$ entre les fibrés normaux de $S$ et $S'$, qui se restreint en un symplectomorphisme $f : (S, \omega_{\vert S}) \rightarrow (S', \omega'_{\vert S'})$. Alors $f$ s'étend en un symplectomorphisme entre des voisinages tubulaires des courbes $S$ et $S'$.
\end{thm}

Par conséquent, la structure symplectique au voisinage d'une courbe symplectiquement plongée $S$ est complètement déterminée par son aire symplectique $\int_S \omega$ et la classe d'isomorphisme de son fibré normal (orienté).

Guadagni a également démontré une version du théorème du voisinage symplectique pour des configurations de sous-variétés symplectiques qui s'intersectent transversalement positivement en des points deux à deux distincts. On l'énonce ici dans le cadre des surfaces symplectiques.

\begin{thm}[Théorème du voisinage symplectique pour les plombages de courbes, \cite{Guadagni}]\label{t:sympnghdconf}
Soit $(M, \omega)$, $(M', \omega')$ deux surfaces symplectiques et des familles finies de courbes symplectiques $S_1, ..., S_\ell \subset M$, $S_1', \dots, S_\ell ' \subset M'$ qui s'intersectent transversalement positivement en des points deux à deux distincts. On suppose qu'il existe un difféomorphisme entre des voisinages de $\bigcup_{i=1}^\ell S_i$ et de $\bigcup_{i=1}^\ell S_i'$. Alors il existe un voisinage de $\bigcup_{i=1}^\ell S_i$ qui est une déformation symplectique d'un voisinage de $\bigcup_{i=1}^\ell S_i'$.
\end{thm}


Pour une introduction détaillée à la topologie symplectique avec une présentation centrée autour des basses dimensions, on renvoit à~\cite{ozbagci2004surgery}.


\section{Nombres d'intersection}\label{s:nbint}

\subsection{Forme d'intersection d'une variété de dimension 4}

Une particularité en dimension $4$ réside dans le fait que, génériquement, les surfaces lisses dans une variété de dimension $4$ s'intersectent transversalement en un nombre fini de points (à cause de la compacité). On peut alors compter ces points avec signe, en tenant compte des différentes orientations en jeu, et obtenir ainsi des invariants topologiques très utiles. Dans cette section, on rappelle la définition de la forme d'intersection d'une variété orientée de dimension $4$ et des principaux invariants associés. Le livre~\cite{scorpan} donne un aperçu global de l'utilité de ces invariants en topologie de basse dimension.

\begin{df}
Soit $M$ une variété orientée de dimension $4$. La forme d'intersection de $M$ est la forme $\mathbb{Z}$--bilinéaire symétrique définie par 
$$\begin{array}{ccccc}
\mathcal{Q}_M & : & H_2^{free} (M; \mathbb{Z}) \times H_2^{free}(M; \mathbb{Z}) & \longrightarrow & \mathbb{Z} \\
 & & (A,B) & \longmapsto & \langle PD^{-1}(A) \smile PD^{-1}(B), [M] \rangle, \\
\end{array}$$
où $H_2^{free} (M; \mathbb{Z})$ désigne le quotient du groupe d'homologie $H_2 (M; \mathbb{Z})$ par son sous-groupe de torsion, et 
$$\begin{array}{ccccc}
P D & : & H^2(M; \mathbb{Z}) & \longrightarrow & H_2(M; \mathbb{Z}) \\
 & & \alpha & \longmapsto & [M] \frown \alpha
 \end{array}$$
désigne la dualité de Poincaré.
\end{df}

La classe d'isomorphisme de la forme d'intersection constitue un invariant topologique de $M$. Par linéarité, les classes de torsion sont dans le noyau de la forme d'intersection. On supposera donc implicitement sans perte de généralité qu'on se restreint à la partie sans torsion de $H_2(M; \mathbb{Z})$. La forme d'intersection d'une variété (orientée, compacte et sans bord) est unimodulaire (son déterminant est égal à $\pm 1$), et son rang est égal au nombre de Betti $b_2 (M)$.

Pour $A,B \in H_2(M; \mathbb{Z})$, l'entier $Q_M(A,B)$ est appelé le nombre d'intersection entre $A$ et $B$, parfois noté $A \cdot B$ (on note aussi $A^2 = A \cdot A$). Si la définition est valable pour les variétés topologiques, ce nombre peut être interprété de manière géométrique dans le cas des variétés différentielles. En effet, on peut dans ce cas représenter les classes d'homologies $A$ et $B$ par des surfaces orientées $S_A$ et $S_B$ plongées dans $M$ (voir~\cite{scorpan}). Quitte à les perturber de manière $\mathcal{C}^\infty$, on peut supposer que $S_A$ et $S_B$ s'intersectent transversalement en des points isolés. Par compacité, il y a un nombre fini de tels points $p$, qu'on peut compter en leur attribuant un signe $\varepsilon (p)$ de la manière suivante :
$$\varepsilon(p) = 
\left \{
\begin{array}{rl}
+1 & \text{si} ~ T_p S_A \oplus T_p S_B ~ \text{a la même orientation que} ~ T_p M, \\ 
-1 &\text{sinon.}
\end{array}
\right.$$

On a alors l'égalité suivante : $$S_A \cdot S_B = \sum\limits_{p \in S_A \cap S_B} \varepsilon (p).$$

D'après la loi d'inertie de Sylvester, on peut diagonaliser la forme d'intersection $Q_M$ sur $\mathbb{R}$ et écrire $H_2(M; \mathbb{R}) = H_2^+ (M; \mathbb{R}) \oplus H_2^-(M; \mathbb{R})$, où $H_2^+ (M; \mathbb{R})$ et $H_2^- (M; \mathbb{R})$ sont des espaces vectoriels réels de dimensions maximales pour lesquels la restriction de $Q_M$ à $H_2^+ (M; \mathbb{R})$ est définie positive (c'est-à-dire pour toute classe $A$ non nulle, $Q_M (A,A) > 0$) et la restriction de $Q_M$ à $H_2^- (M; \mathbb{R})$ est définie négative (c'est-à-dire pour toute classe $A$ non nulle, $Q_M (A,A) < 0$). On peut alors définir la \emph{signature} de $Q_M$ par $\sigma(M) = \dim H_2^+ (M; \mathbb{R}) - \dim H_2^- (M; \mathbb{R})$. On notera également $b_2^+(M) = \dim H_2^+ (M; \mathbb{R})$, et $b_2^-(M) = \dim H_2^- (M; \mathbb{R})$. Remarquons qu'on a alors $$b_2(M) = b_2^+ (M) + b_2^- (M) \quad \text{et} \quad \sigma (M) = b_2^+ (M) - b_2^- (M).$$

\begin{rk}
Si $(M, \omega)$ est une surface symplectique, on a $[\omega \wedge \omega] = [\omega] \smile [\omega]$, et puisque $\omega \wedge \omega$ est  une forme volume sur $M$, le réel $\int_M \omega \wedge \omega = \langle [\omega] \smile [\omega], [M] \rangle$ est strictement positif. Ainsi, on a nécessairement $b_2^+ (M) \geq 1$.
\end{rk}

Un autre invariant obtenu à partir de la forme d'intersection est la parité. On dit que $Q_M$ est \emph{paire} si pour toute classe $A \in H_2(M; \mathbb{Z})$, $Q_M (A,A)$ est pair, sinon on dit que $Q_M$ est \emph{impaire}.

\begin{ex}
\leavevmode
\begin{enumerate}
\item Le groupe $H_2 (\mathbb{C}P^2; \mathbb{Z})$ est engendré par la classe d'homologie $h$ d'une droite complexe. La forme d'intersection $\mathcal{Q}_{\mathbb{C}P^2}$ est donc représentée par la matrice $\begin{pmatrix} 1 \end{pmatrix}$. On en déduit que la forme d'intersection de $\mathbb{C}P^2$ est impaire et on a $b_2^+(\mathbb{C}P^2) =1$ et $b_2^-(\mathbb{C}P^2) =0$.

\item Le groupe $H_2 (\mathcal{S}^2 \times \mathcal{S}^2; \mathbb{Z})$ est engendré par les classes d'homologie $[\mathcal{S}^2 \times \{ * \}]$ et $[\{ * \} \times \mathcal{S}^2]$. La forme d'intersection $\mathcal{Q}_{\mathcal{S}^2 \times \mathcal{S}^2}$ est donc représentée dans cette base par la matrice 
$\begin{pmatrix}
0 & 1 \\
1 & 0
\end{pmatrix}$. 
On en déduit que la forme d'intersection de $\mathcal{S}^2 \times \mathcal{S}^2$ est paire et on a $b_2^+(\mathcal{S}^2 \times \mathcal{S}^2) =1$ et $b_2^-(\mathcal{S}^2 \times \mathcal{S}^2) =1$.
\end{enumerate}
\end{ex}

\subsection{Intersection locale de courbes holomorphes}

Dans cette partie, on s'intéresse aux propriétés locales des courbes complexes, possiblement singulières, dans les surfaces complexes. Pour cela, il suffit d'étudier les courbes dans des cartes affines, c'est-à-dire les zéros de polynômes de $\mathbb{C}[x,y]$. On s'intéresse donc principalement aux courbes définies dans un voisinage de l'origine $O$ de $\mathbb{C}^2$. Pour cette étude locale, il n'est pas important de faire la différence entre le fait que l'équation définissant la courbe soit définie par un polynôme ou par une fonction analytique. Le livre \cite{wall2004singular} fournit une exposition détaillée sur ce sujet.

Une approche pour étudier localement les courbes algébriques dans les surfaces consiste à s'intéresser à leur paramétrisation. 

\begin{df}
Soit $C$ une courbe algébrique dans une surface complexe $X$. Une application $\varphi$ définie au voisinage de $0$ dans $\mathbb{C}$ et à valeurs dans $C$ est une \emph{paramétrisation locale} de la courbe $C$ en $O$ s'il existe un voisinage $U$ de $0$ et un voisinage $V$ de $O$ tels que l'application $\varphi : U \rightarrow V \cap C$ est surjective.

On dit de plus que c'est une \emph{bonne paramétrisation} si $\varphi$ est bijective.
\end{df}

On peut alors compter, avec multiplicité, le nombre de points d'intersection entre deux courbes définies localement au voisinage de $O$.

\begin{df}
Soit $C_1$ une courbe algébrique plane définie au voisinage de $0$ par l'équation analytique $g(x,y) = 0$ et $C_2$ une courbe algébrique plane définie au voisinage de $0$ par une bonne paramétrisation $\varphi$. Alors le \emph{nombre d'intersection local} entre $C_1$ et $C_2$ en $O$ est défini comme étant l'ordre de la fonction analytique $g \circ \varphi$. On le note $(C_1 \cdot C_2)_O$.
\end{df}

On peut montrer que le nombre d'intersection est symétrique et qu'il ne dépend pas du choix des cartes (voir \cite{wall2004singular}). Cette définition coïncide avec celle donnée dans la partie précédente pour des intersections transverses, et on peut montrer que pour deux courbes sans composante commune $C_1$ et $C_2$ dans une surface complexe, on a $[C_1] \cdot[C_2] = \sum_{p \in C_1 \cap C_2} (C_1 \cdot C_2)_p$. On remarque également que le nombre d'intersection local entre deux germes de courbes complexes distincts est positif : deux courbes complexes sans composante commune dans une surface s'intersectent toujours localement positivement. En particulier, toujours sous l'hypothèse que les deux courbes n'ont pas de composante commune, on a $[C_1] \cdot[C_2] = 0$ si et seulement si $C_1$ et $C_2$ sont disjointes. De même pour $p \in C_1 \cap C_2$, on a $(C_1 \cdot C_2)_p = 1$ si et seulement si $p$ est un point non singulier pour $C_1$ et $C_2$ en lequel les deux courbes s'intersectent transversalement.

\begin{ex}
Soit $C_1$ et $C_2$ deux courbes complexes dans $\mathbb{C}P^2$ de degrés respectifs $d_1$ et $d_2$. Comme le groupe $H_2 (\mathbb{C}P^2; \mathbb{Z})$ est engendré par la classe d'homologie $h$ d'une droite complexe et puisqu'on a $[C_1] \cdot h =d_1$ et $[C_2] \cdot h =d_2$, les courbes $C_1$ et $C_2$ ont pour classes d'homologie respectives $d_1 h$ et $d_2 h$. On retrouve alors le théorème de Bézout pour les courbes algébriques complexes : 
$$\sum_{p \in C_1 \cap C_2} (C_1 \cdot C_2)_p = d_1 h \cdot d_2 h = d_1 d_2.$$ 
\end{ex}

\begin{ex}
On considère la singularité complexe $p$ donnée localement au voisinage de $O$ par la courbe $C$ d'équation $y^2-x^3=0$ (parfois appelée singularité cuspidale simple). La courbe $C_y$ paramétrée localement par $t \mapsto (0,t)$ au voisinage de $O$ vérifie $(C \cdot C_y)_O = 2$. Elle représente localement l'intersection en $p$ entre $C$ et une courbe non singulière qui n'est pas tangente à $C$ en $p$. La courbe $C_x$ paramétrée localement par $t \mapsto (t,0)$ au voisinage de $O$ vérifie $(C \cdot C_x)_O = 3$. Elle représente localement l'intersection en $p$ entre $C$ et une courbe non singulière qui est tangente à $C$ en $p$. 
\end{ex}

Soit $C$ une courbe algébrique complexe dans une surface complexe $X$ et $p$ un point appartenant à $C$. À partir des nombres d'intersection locaux, on peut définir la \emph{multiplicité} de $p$ en tant que point de $C$ comme l'entier naturel non nul 
$$\min \{(C \cdot C')_p, ~C' \text{ germe de courbe complexe en } p \}.$$ 
La multiplicité de $p$ est égale à $1$ si et seulement $p$ est un point non singulier de $C$.

\begin{ex}
Les singularités cuspidales simples sont des points singuliers de multiplicité $2$.
\end{ex}

\section{Courbes pseudoholomorphes} \label{s:Jholo}

Les techniques pseudoholomorphes, introduites par Gromov~\cite{Gromov}, se sont révélées être très efficaces en topologie symplectique, notamment en ce qui concerne l'étude des surfaces symplectiques. Elles permettent de faire de nombreux parallèles entre les variétés symplectiques et les variétés complexes. Dans cette section, on rappelle les définitions et les résultats qui posent le cadre propice à la compréhension des raisonnements exposés dans les chapitres suivants. Cette partie est très inspirée de l'excellente référence~\cite{Wendl}. Par soucis de clarté, on reprendra les mêmes notations dans ce manuscrit. On renvoie donc le lecteur à cette référence, ainsi qu'à~\cite{wendl2014lectures}, pour de plus amples détails sur le sujet.

\subsection{Structures presque complexes}

Pour un fibré vectoriel $E$ de dimension finie, une \emph{structure complexe} sur $E$ est une application $J : E \rightarrow E$ linéaire fibre à fibre, satisfaisant $J \circ J = -\mathrm{id}_{E}$. Pour une variété différentielle $M$ de dimension $2n$, on appelle \emph{structure presque complexe} $J$ sur $M$ toute structure complexe sur $TM$, c'est-à-dire toute application entre fibrés vectoriels $J : TM \rightarrow TM$, linéaire fibre à fibre, satisfaisant $J \circ J = -\mathrm{id}_{TM}$. Remarquons qu'une variété complexe est toujours naturellement munie d'une structure presque complexe canonique, donnée par la multiplication par le nombre imaginaire $i$ sur chaque espace tangent. On dit que $J$ est \emph{intégrable} lorsqu'elle provient d'une structure complexe sur $M$. De même, on dit que $J$ est \emph{intégrable au voisinage d'un point} $p \in M$ si $J$ est intégrable sur un ouvert $U \subset M$ contenant $p$. Notons qu'une structure presque complexe sur une surface est toujours intégrable, on parle dans ce cas de \emph{surface de Riemann}. Enfin on appelle $(M,J)$ une \emph{variété presque complexe} et lorsque $M$ est de dimension $4$, on dit que $(M,J)$ est une \emph{surface presque complexe}.

\begin{rk}
\leavevmode
\begin{enumerate}
\item En fixant un point $p \in M$, et en appliquant le déterminant à l'égalité $J_p \circ J_p = -\mathrm{id}_{T_pM}$, on remarque que la définition n'a effectivement du sens que lorsque $M$ est de dimension paire. 
\item Une variété presque complexe $(M,J)$ est naturellement orientée par sa structure presque complexe $J$. En effet, en tout point $p \in M$, on définit une base de $T_p M$ de la forme $(v_1, J v_1, \dots, v_n, J v_n)$ comme étant une base positive.
\end{enumerate}
\end{rk}

Cette notion se marie remarquablement bien avec la notion de variété symplectique. En effet, on dit qu'une structure presque complexe $J$ sur une variété symplectique $M$ est \emph{dominée} par $\omega$ si, pour tout $u \in TM$ non nul, $\omega(u , Ju) > 0$. Si de plus, le $2$--tenseur défini pour $u,v \in TM$, par $g_J (u,v) = \omega(u, Ju)$ est symétrique, on dit que $J$ est \emph{compatible} avec $\omega$. Dans ce dernier cas, $g_J$ défini une métrique riemannienne sur $M$.

Dans la suite, on notera $\mathcal{J}(M)$ l'ensemble des structures presque complexes lisses sur $M$, $\mathcal{J}_\tau (M, \omega)$ l'ensemble des structures presque complexes lisses dominées par $\omega$, et $\mathcal{J}(M, \omega)$ l'ensemble des structures presque complexes lisses compatibles avec $\omega$. Ces ensembles, munis de leur topologie $\mathcal{C}^\infty$ naturelle, sont des espaces complets métrisables.

\begin{prop}
Pour une variété symplectique $(M, \omega)$, les espaces $\mathcal{J}(M, \omega)$ et $\mathcal{J}_\tau (M, \omega)$ sont non vides et contractiles.
\end{prop}

Les structures presque complexes se révèlent être très adaptées à l'étude des courbes symplectiques dans les variétés symplectiques, comme le montre la proposition suivante.

\begin{prop}
Soit $(M, \omega)$ une variété symplectique et $S \subset M$ une sous-variété lisse de dimension $2$. Alors $S$ est une sous-variété symplectique si et seulement si il existe une structure presque complexe $J$ sur $M$ dominée par $\omega$ qui préserve $TS$.
\end{prop}

\subsection{Courbes simples et revêtements multiples}\label{s:revmult}
Soit $(\Sigma , j)$ une surface de Riemann et $(M, J)$ une variété presque complexe de dimension (réelle) $2n$. On dit qu'une application $u : \Sigma \rightarrow M$ de classe $\mathcal{C}^\infty$ est une \emph{courbe pseudoholomorphe} (on emploiera aussi le terme \emph{courbe $J$--holomorphe}) si elle vérifie l'équation de Cauchy--Riemann non linéaire suivante (ici $du$ désigne la différentielle de $u$)
$$du \circ j =  J \circ du.$$

Dans le cas où $J$ est intégrable, cette équation nous indique juste que la différentielle de $u$ est $\mathbb{C}$--linéaire, ce qui signifie simplement que $u$ est une application holomorphe au sens usuel. 

Lorsqu'on étudie les courbes pseudoholomorphes, on est souvent amené à s'intéresser aux applications holomorphes entre des surfaces de Riemann, ne serait-ce que pour des questions de reparamétrisation. Le degré d'une telle application nous donne beaucoup d'informations sur l'application elle-même.

\begin{prop}
Soit $(\Sigma , j)$, $(\Sigma' , j')$ deux surfaces de Riemann et $ \varphi : (\Sigma , j) \rightarrow (\Sigma' , j')$ une application holomorphe. Alors le degré $\deg \varphi \in \mathbb{Z}$ de l'application $\varphi$ est positif, et une des situations suivantes est vérifiée :
\begin{enumerate}
\item $\deg \varphi =0$, auquel cas $\varphi$ est une application constante ;
\item $\deg \varphi =1$, auquel cas $\varphi$ est un biholomorphisme ;
\item $\deg \varphi \geq 2$, auquel cas $\varphi$ est un revêtement ramifié.
\end{enumerate}
\end{prop}

On introduit désormais le vocabulaire qu'on utilisera fréquemment en ce qui concerne les courbes pseudoholomorphes. Si $ \varphi : (\Sigma , j) \rightarrow (\Sigma' , j')$ est une application holomorphe, et $u' : ( \Sigma',j') \rightarrow (M, J)$ est une courbe $J$--holomorphe, leur composition $u = u' \circ \varphi : ( \Sigma, j) \rightarrow (M,J)$ définit également une courbe $J$--holomorphe. Dans le cas où $\deg \varphi =1$, on dit que $u$ est une \emph{reparamétrisation} de $u'$. Si $k = \deg \varphi \geq 2$, on dit que $u$ est un \emph{revêtement ramifié} de degré $k$ de $u'$ (on dira parfois que $u$ est une courbe \emph{multiplement revêtue}). On dit qu'une courbe $J$--holomorphe non constante est \emph{simple} lorsqu'elle n'est le revêtement ramifié d'aucune autre courbe $J$--holomorphe.

On dit qu'une application $u : \Sigma \rightarrow M$ est \emph{injective quelque part} s'il existe un point $z \in \Sigma$, qu'on appelle \emph{point injectif} de $u$, tel que l'application $d_z u : T_z \Sigma \rightarrow T_{u(z)} M$ est injective et $u^{-1}(u(z)) = \{z \}$. On appelle \emph{point non immergé} tout point $z$ en lequel l'application $d_z u$ n'est pas injective. Dans le cas où $u$ est une courbe $J$--holomorphe, $z$ est un point non immergé si et seulement si $d_z u = 0$. Les courbes multiplement revêtues ne possèdent évidemment aucun point injectif. Le résultat suivant montre que la réciproque est également vraie.

\begin{prop}\label{p:revetement}
Soit $u : (\Sigma, j) \rightarrow (M,J)$ une courbe $J$--holomorphe. Les propositions suivantes sont équivalentes : 
\begin{enumerate}
\item $u$ est injective quelque part, 
\item $u$ est simple,
\item $u$ possède au plus un nombre fini de points d'auto-intersection et de points non immergés.
\end{enumerate}
De plus, si la courbe $u$ n'est ni simple, ni constante, elle nécessaire un revêtement multiple de degré $k \geq 2$ d'une courbe simple.
\end{prop}

On regroupe désormais quelques propriétés classiques des courbes pseudoholomorphes, qu'on utilisera fréquemment par la suite.

On définit la classe d'homologie d'une courbe $J$--holomorphe $u : \Sigma \rightarrow M$ par $[u]:= u_*[\Sigma] \in H_2(M;\mathbb{Z})$. Si la courbe $u$ est un revêtement multiple de degré $k \in \mathbb{N}^*$ d'une autre courbe $J$--holomorphe $u'$, remarquons que $[u] = k [u']$. Enfin, dès que $J$ est dominée par une structure symplectique $\omega$, les courbes $J$--holomorphes constantes sont complètement caractérisées par leurs classes d'homologie. En effet, toute courbe $J$--holomorphe doit alors satisfaire l'inégalité suivante
$$\int u^* \omega \geq 0,$$
avec égalité si et seulement si $u$ est constante. Ainsi, une telle courbe $J$--holomorphe est constante si et seulement si elle est homologue à $0 \in H_2(M; \mathbb{Z})$.

On appelle \emph{genre géométrique} d'une courbe $J$--holomorphe simple $u : \Sigma \rightarrow M$ le genre de la surface $\Sigma$. Lorsque le genre géométrique d'une courbe pseudoholomorphe est nul, on dit que la courbe est \emph{rationnelle} et on dit qu'elle est \emph{irrationnelle} sinon. Pour terminer, on présente quelques propositions classiques concernant le genre des surfaces de Riemann et le genre géométrique des courbes holomorphes.

\begin{prop}[Formule de Riemann--Hurwitz]
Soit $\varphi : (\Sigma,j) \rightarrow (\Sigma',j')$ une application holomorphe de degré $k \geq 1$ entre deux surfaces de Riemann. Notons $Z(d \varphi) \in \mathbb{Z}$ la somme des ordres des points critiques de $\varphi$ (définis comme l'ordre des zéros de $d \varphi$). On a alors $$- \chi (\Sigma) + k \chi ( \Sigma') = Z (d \varphi) \geq 0.$$
\end{prop}

En particulier, pour une courbe $J$--holomorphe $u : ( \Sigma,j) \rightarrow (M, J)$ non constante de genre géométrique $g$ telle qu'il existe une courbe $u' : ( \Sigma',j') \rightarrow (M, J)$ de genre géométrique $g'$ et une application holomorphe $ \varphi : (\Sigma , j) \rightarrow (\Sigma' , j')$  de degré $k \geq 1$ vérifiant $u = u' \circ \varphi$, alors on a $2g -2 + k(2-2g') = Z(d \varphi)$. Ainsi, on a $g'=0$ ou $2(g-g')=(k-1)(2g'-2) + Z(d \varphi) \geq 0$, ce qui montre que $g \geq g'$. En fait, ce dernier fait est un peu plus général, comme le montre la proposition suivante.

\begin{prop}\label{p:ineqgenresurfaces}
Soit $\Sigma_g$, $\Sigma_h$ deux surfaces de genres respectifs $g$ et $h$, et $\varphi : \Sigma_g \rightarrow \Sigma_h$ une application continue de degré non nul. Alors $g \geq h$.
\end{prop}

\begin{proof}
Supposons $h >g$. L'application $\varphi^* : H^1 ( \Sigma_h; \mathbb{Z}) \rightarrow H^1 (\Sigma_g; \mathbb{Z})$ est morphisme d'un groupe abélien libre de rang $h$ vers un groupe abélien libre de rang strictement plus petit $g$, donc son noyau est non nul. Prenons désormais $\alpha \in (\mathrm{ker} \varphi^* )\backslash \{ 0 \}$. Par la dualité de Poincaré, on a $[\Sigma_h] \frown \alpha \neq 0$. On obtient alors 
$$0 = \varphi_* ([\Sigma_g] \frown \varphi^* (\alpha)) = \varphi^* ([\Sigma_g]) \frown \alpha = \deg \varphi \times [\Sigma_h] \frown \alpha ,$$
ce qui implique $\deg \varphi = 0$.
\end{proof}

\subsection{Positivité d'intersection, formule d'adjonction}

Un des grands avantages des courbes pseudoholomorphes sur les courbes symplectiques réside dans leur propriété de positivité d'intersection. Cette dernière assure, à l'instar des courbes holomorphes, que deux courbes $J$--holomorphes (connexes) d'images distinctes s'intersectent localement positivement.

\begin{thm}[Positivité d'intersection]
Soit $(M,J)$ une surface presque complexe et $u: \Sigma \rightarrow M$, $v : \Sigma \rightarrow M$ deux courbes $J$--holomorphes (connexes) d'images distinctes. Alors les courbes $u$ et $v$ s'intersectent en un nombre fini de points, et on a
$$[u] \cdot [v] \geq \Card \left\{ (z,z') \in \Sigma \times \Sigma' ~ | ~ u(z) = v(z') \right\},$$
avec égalité si et seulement si toutes les intersections sont transverses. En particulier, on a $[u] \cdot [v] = 0$ si et seulement si les images de $u$ et $v$ sont disjointes, et $[u] \cdot [v] = 1$ si et seulement si $u$ et $v$ s'intersectent exactement une fois, de manière transverse, en un point qui est un point injectif pour les deux courbes.
\end{thm}

\begin{rk}
En fait, le théorème est vrai plus généralement pour deux courbes $J$--holomorphes sans composante commune (voir la Section~\ref{s:CompactificationGromov} pour la définition de composante d'une courbe).
\end{rk}

\begin{ex}
La positivité d'intersection se montre aisément dans le cas où un point $p \in M$ est une intersection transverse entre deux courbes $J$--holomorphes $u$ et $v$. En effet, étant donné un vecteur non nul $X$ dans l'espace tangent $T_p u$ de la courbe $u$ en $p$, le couple $(X, JX)$ définit une base positive $T_p u$. De même, pour un vecteur non nul $Y$ dans l'espace tangent $T_p v$ de la courbe $v$ en $p$, le couple $(Y, JY)$ définit une base positive $T_p v$. Le fait que le quadruplet $(X, JX,Y, JY)$ forme une base positive de $T_p M$ permet de conclure.
\end{ex}

En fait, un résultat plus fort dû à McDuff, nous assure qu'en réalité, du point de vue topologique, les singularités des courbes pseudoholomorphes dans les surfaces presque complexes sont localement modelées sur les singularités des courbes complexes dans les surfaces complexes. On peut alors évaluer les contributions locales de chaque point d'intersection entre deux courbes dans le calcul du nombre d'intersection de la même manière que pour les singularités complexes. Mentionnons tout de même que Micallef et White ont démontré un résultat similaire dans~\cite{MicallefWhite}, dans le cas où la structure presque complexe est dominée par une forme symplectique (en effet, dans ce cas, les courbes pseudoholomorphes sont des courbes d'aire minimale dans leur classe d'homotopie, ce qui correspond bien aux hypothèses du théorème de Micallef et White).

\begin{thm}[\cite{McDuff-Jhol}] \label{t:McDuff92} 
Soit $(M,J)$ une variété presque complexe et $C$ l'image d'une courbe $J$--holomorphe. Alors $C$ possède un nombre fini de points singuliers (c'est-à-dire de points qui ne sont pas injectifs pour une paramétrisation simple de $C$). De plus, il existe une structure presque complexe $J_0$ arbitrairement proche de $J$ pour la topologie $\mathcal{C}^0$ et une courbe $J_0$--holomorphe qui est $\mathcal{C}^1$--proche de $C$, d'image $C_0$, telle que $J_0$ est intégrable au voisinage de chacun des points singuliers de $C_0$ et telle que la paire $(V,C_0)$ est homéomorphe à $(V,C)$.
\end{thm}

En particulier, lorsque qu'une courbe pseudoholomorphe $u : (\Sigma, j) \rightarrow (M,J)$ simple est immergée, elle possède un nombre fini de points d'auto-intersection, qu'on peut compter avec signe et multiplicité. Par positivité d'intersection, tous ces points d'auto-intersection sont comptés positivement. On peut alors définir l'entier naturel
$$\delta (u) = \frac{1}{2} \sum\limits_{u(z) = u(z')} i(z,z'),$$
où la somme se fait sur les paires $(z,z')$ de points distincts de $\Sigma$ pour lesquelles $u(z) = u(z')$, et $i(z,z') \in \mathbb{N}^*$ est le nombre algébrique d'intersection local entre des images par $u$ de voisinages suffisamment petit de $z$ et $z'$. Puisque pour toutes telles paires $(z,z')$, on a $i(z,z') >0$, on obtient que $\delta(u) \geq 0$, avec égalité si et seulement si $u$ est plongée.

Lorsque $u$ est simple, mais n'est pas immergée, on définit $\delta(u)$ comme étant égal à $\delta (u')$, où $u'$ désigne une perturbation immergée $\mathcal{C}^\infty$--proche de $u$. Le choix de $u'$ peut être effectué de manière à ce que le nombre $\delta(u')$ ne dépende pas du choix de $u'$ effectué (plus précisément en effectuant la perturbation localement de manière complexe au voisinage de chacun des points singuliers de $u$ grâce au Théorème~\ref{t:McDuff92}). De plus si $u$ possède un point de non immersion, l'entier $\delta(u)$ est strictement positif.

Pour les courbes pseudolohomorphes simples, on dispose d'une formule reliant la classe d'homologie d'une courbe et son genre. Avant de l'énoncer, rappelons brièvement quelques propriétés concernant la première classe de Chern (une présentation plus détaillée peut être trouvée dans \cite{milnor1974characteristic}, \cite{bott1995differential} ou \cite{husemoller1994fibre}). Étant donné un fibré vectoriel complexe $E$ sur une variété $X$, on note $c_1 (E) \in H^2(X; \mathbb{Z})$ sa première classe de Chern. Dans le cas d'une variété presque complexe $(M,J)$, on peut alors considérer la première classe de Chern de son fibré tangent $c_1(TM,J)$, qu'on appelle la première classe de Chern de $(M,J)$ (notons que pour une variété symplectique $(M, \omega)$ et $J \in \mathcal{J}_\tau (M, \omega)$, la première classe de Chern $c_1(TM,J)$ ne dépend pas de la structure presque complexe $J$ choisie). Dans ce manuscrit, on utilise l'abus de notation suivant : pour $A \in H_2(M; \mathbb{Z})$, on note $c_1(A) := \langle c_1 (TM,J), A \rangle$. On rappelle que dans le cas où $E$ est un fibré en droite complexe au-dessus d'une surface $X$, l'entier $\langle c_1 (E), [X] \rangle$ correspond au nombre de zéros d'une section générique de $E$ (c'est-à-dire transverse à la section nulle), comptés avec signe. Rappelons également que les fibrés en droites complexes sur une variété $X$ sont complètement classifiés par leurs premières classes de Chern et que pour deux fibrés vectoriels complexes $E_1$ et $E_2$ sur une même variété $X$, on a $c_1 (E_1 \oplus E_2) = c_1(E_1) + c_1 (E_2)$.

\begin{thm}[Formule d'adjonction]\label{t:adjformula}
Soit $(M,J)$ une surface presque complexe et $u : \Sigma \rightarrow M$ une courbe $J$--holomorphe injective quelque part. Alors 
$$c_1([u])=  \chi ( \Sigma) + [u] \cdot [u] - 2 \delta(u),$$
où $\delta (u) \geq 0$, avec égalité si et seulement si $u$ est plongée.
\end{thm}

\begin{proof}
Commençons par donner une preuve dans le cas où $u$ est plongée, d'image $S \subset M$. On a alors $TM_{|S}= TS \oplus \nu S$, où $\nu S$ désigne le fibré normal de $S$, donc $c_1([u])= \langle c_1(TS), [S] \rangle + \langle c_1(\nu S), [S] \rangle$. Comme $TS$ et $\nu S$ sont des fibrés en droites complexes, les entiers $\langle c_1(TS), [S] \rangle$ et $\langle c_1(\nu S), [S] \rangle$ correspondent respectivement aux nombres de zéros de sections génériques des fibrés en droites correspondants, comptés avec signe. Le théorème de Poincaré-Hopf nous assure alors que $\langle c_1(TS), [S] \rangle = \chi (S)$, et grâce au Théorème du voisinage tubulaire, on obtient $\langle c_1(\nu), [S] \rangle = [S]^2$, ce qui nous donne bien $c_1([u]) = \chi ( \Sigma) + [u] \cdot [u] $.

Lorsque $u$ est immergée, on peut quand même définir un fibré normal $N_u$ sur $\Sigma$ (en utilisant le fait que les restrictions de $u$ à des ouverts suffisamment petits sont des plongements). On peut également supposer sans perte de généralité que tous les points d'auto-intersection de $u$ sont des points doubles qui sont transverses et positifs. On applique alors le même raisonnement en remarquant que $u^*TM = T \Sigma \oplus N_u$. Cette fois-ci, calculons le nombre d'auto-intersection de $u$ en prenant une petite perturbation générique $u'$. Chaque point double positif apporte une contribution de $2$ et les autres intersections entre $u$ et $u'$ correspondent au nombre de zéros d'une section générique de $N_u$ (voir Figure~\ref{f:adjunction}), ainsi on a $c_1 (N_u) = [u] \cdot [u] - 2 \delta(u)$, ce qui conclut la preuve.
\begin{figure}[h]
	\centering
	\includegraphics[scale=0.5]{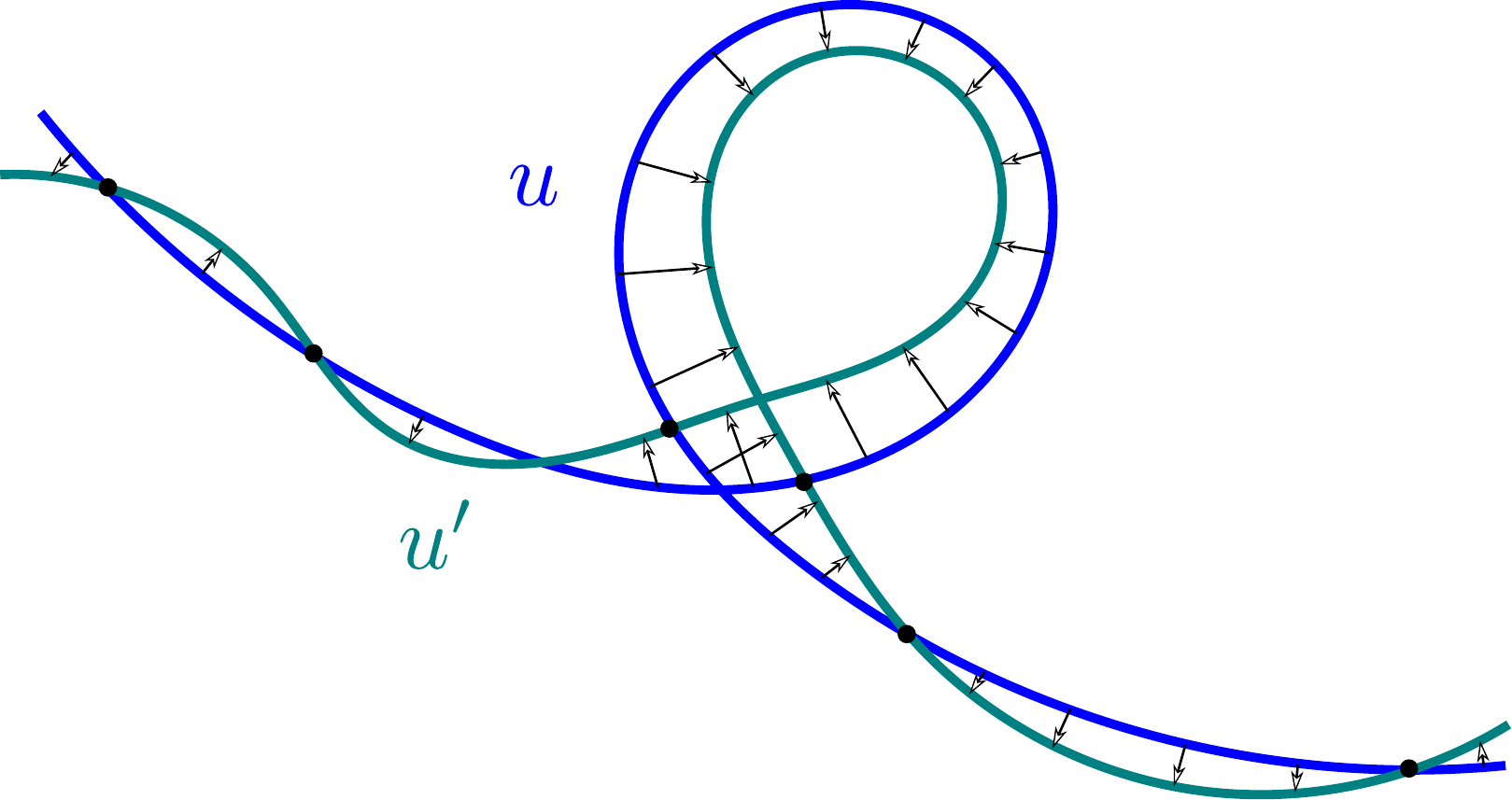}\\
	\caption{Calcul du nombre d'auto-intersection de $u$ à l'aide d'une perturbation générique $u'$.}
	\label{f:adjunction}
\end{figure}

\end{proof}

\begin{rk}\label{r:adjonctionplongée}
Dans le cas où $u$ est une courbe plongée de genre $g$, la formule d'adjonction nous indique que toutes les courbes simples de genre $g$ homologues à $u$ sont également plongées.
\end{rk}

On appelle \emph{genre arithmétique} d'une courbe pseudoholomorphe $u$ simple l'entier naturel $p_a(u)= g + \delta(u)$. D'après la formule d'adjonction, on a $p_a(u) = \frac{1}{2} ([u] \cdot [u] - c_1([u])) +1$, ce qui correspond au genre des courbes pseudohomorphes plongées homologues à $u$ (ainsi pour les courbes plongées, les notions de genre géométrique et genre arithmétique coïncident). Le genre arithmétique de $u$ peut également être vu comme le genre de la surface obtenue en prenant d'abord une perturbation immergée de $u$, pour laquelle tous les points d'auto-intersection sont des points doubles transverses positifs, puis en lissant tous ces points doubles.

Le lissage d'un point double désigne l'opération suivante. Dans un premier temps, on trouve des coordonnées complexes $(z_1,z_2)$ au voisinage d'un point double positif, dans lesquelles l'intersection est localement décrite par l'équation $z_1 z_2 =0$ (deux droites complexes qui se croisent), et on remplace les portions de courbes dans ce voisinage par une perturbation arbitrairement petite donnée localement par l'équation $z_1 z_2 = \varepsilon$. On utilise ensuite des fonctions plateaux pour raccorder les extrémités du germe de courbe ainsi obtenu au reste de la courbe.

\emph{Lorsqu'on ne précise pas de quel genre on parle, on considère dans ce manuscrit qu'on parle par défaut du genre géométrique.}

\begin{ex}
Les cubiques singulières de $\mathbb{C}P^2$ (comme la cubique nodale donnée par l'équation $y^2z =x^2(x+z)$ ou la cubique cuspidale donnée par l'équation $y^2z=x^3$) sont paramétrées par $\mathbb{C}P^1$, ce sont donc des courbes de genre géométrique $0$ (autrement dit ce sont des courbes rationnelles). Quand on réalise des perturbations arbitrairement petites des équations qui les définissent (en rajoutant un terme de la forme $\varepsilon z^3$ pour $\varepsilon>0$ arbitrairement proche de $0$ par exemple), on obtient des cubiques lisses. D'après la formule du genre pour les courbes algébriques planes, le genre d'une cubique lisse est égal à $\frac{1}{2} (3-1)(2-1) = 1$. Ainsi, les cubiques singulières de $\mathbb{C}P^2$ sont de genre arithmétique $1$.
\end{ex}

Il sera parfois pratique d'utiliser la terminologie suivante, introduite par Golla et Starskton dans \cite[Section~2.2]{GS}.
Celle-ci permet de décrire une classe plus générale de courbes symplectiques, dans laquelle les courbes sont autorisées à posséder des points singuliers localement modelés sur des singularités complexes.

\begin{df}
Une \emph{courbe symplectique singulière} dans une surface symplectique $(M, \omega)$ est un sous-ensemble $S \subset  X$ satisfaisant les deux propriétés suivantes : 
\begin{enumerate}
\item il existe un nombre fini de points $p_1, \dots, p_m$ de $S$ tels que $S \backslash \{p_1, \dots, p_m \}$ est une sous-variété symplectique (possiblement non compacte) de dimension $2$,
\item pour tout $i \in \{1, \dots, n \}$, il existe un voisinage ouvert $U \subset M$ de $p_i$, $C \subset \mathbb{C}^2$ un germe de singularité holomorphe défini dans un voisinage $V$ de l'origine $O$, et un symplectomorphisme $(U,S\cap U,p_i) \to (V, C\cap V,O)$.
\end{enumerate}
Les points $p_1, \dots, p_m$ sont appelés les singularités (ou les points singuliers) de $S$. On dit de plus qu'une courbe symplectique singulière $S$ est \emph{réductible} si elle peut s'écrire comme l'union de deux courbes symplectiques singulières (non vides) distinctes. On dit que $S$ est \emph{irréductible} sinon.
\end{df}

Cette définition est motivée par les courbes pseudoholomorphes, comme le montre la caractérisation suivante qui prolonge le Théorème \ref{t:McDuff92}.

\begin{thm}[{\cite[Lemma~3.4]{GS}},\cite{McDuff-Jhol},\cite{ MicallefWhite}] \label{t:SingSymp}
Soit $(M, \omega)$ une surface symplectique. Si $S$ est une courbe symplectique singulière dans $(M, \omega)$, alors il existe un ensemble contractile de structures presque complexes $J$ compatibles avec $\omega$ sur $M$ tel que $S$ est l'image d'une courbe $J$--holomorphe.

Réciproquement, si $J$ est une structure presque complexe sur $M$ compatible avec $\omega$ et $u$ est une courbe $J$--holomorphe, alors l'image de $u$ est une courbe symplectique singulière.
\end{thm}

Les définitions de genre géométrique et de genre arithmétique peuvent s'étendre aux courbes symplectiques singulières irréductibles, et la formule d'adjonction énoncée précédemment est également valable dans ce cas.

\subsection{Espaces de modules, transversalité}

Soit $(M,J)$ une variété complexe de dimension $2n$, $g$ et $m$ deux entiers naturels et $A \in H_2 (M; \mathbb{Z})$ une classe d'homologie. L'\emph{espace de modules des courbes $J$--holomorphes non paramétrées} dans $(M,J)$, homologues à $A$, de genre géométrique $g$ et avec $m$ points marqués est défini comme l'ensemble quotient
$$\mathcal{M}_{g,m} (A;J) = \{ ( \Sigma, j ,u , ( \zeta_1, \dots, \zeta_m)) \} \diagup \sim$$
où
\begin{enumerate}
\item $(\Sigma,j)$ est une surface de Riemann de genre $g$,
\item $u : (\Sigma ,j ) \rightarrow (M,J)$ est une courbe $J$--holomorphe vérifiant $[u] =A$,
\item $(\zeta_1, \dots, \zeta_m)$ est un $m$--uplet de points deux à deux distincts de $\Sigma$,
\item $(\Sigma, j ,u , ( \zeta_1, \dots, \zeta_m)) \sim (\Sigma', j' ,u' , ( \zeta_1', \dots, \zeta_m'))$ s'il existe un biholomorphisme $\varphi : (\Sigma, j)  \rightarrow (\Sigma',j')$ tel que $u = u' \circ \varphi$ et qui vérifie pour tout $i \in \{1, \dots, m \}$, $\varphi(\zeta_i) = \zeta_i'$.
\end{enumerate}

On notera parfois simplement $\mathcal{M}_g (A;J) = \mathcal{M}_{g,0}(A;J)$. Notons qu'on dispose toujours d'une application d'oubli $\mathcal{M}_{g,m} (A;J) \rightarrow \mathcal{M}_g (A;J)$ qui à une classe d'équivalence $[(\Sigma, j ,u , ( \zeta_1, \dots, \zeta_m))]$ associe la classe d'équivalence $[(\Sigma, j ,u)]$.

L'espace de modules admet une topologie métrisable, avec la notion de convergence suivante. Une suite de classes d'équivalences $[(\Sigma_k, j_k ,u_k , ( \zeta_1^k, \dots, \zeta_m^k))]$ converge vers $[(\Sigma, j ,u , ( \zeta_1, \dots, \zeta_m))] \in \mathcal{M}_{g,m} (A;J)$ si et seulement si on peut choisir des représentants de la forme 
$$(\Sigma_k, j_k' ,u_k' , ( \zeta_1^k, \dots, \zeta_m^k)) \sim (\Sigma_k, j_k ,u_k , ( \zeta_1^k, \dots, \zeta_m^k))$$
tels que les suites $(u_k')$ et $(j_k')$ convergent respectivement vers $u$ et $j$ pour la topologie $\mathcal{C}^\infty$. Lorsque le contexte ne laissera place à aucune ambigüité, on notera tout simplement $u$ à la place de $[(\Sigma, j ,u , ( \zeta_1, \dots, \zeta_m))]$.

On note $\mathcal{M}_{g,m}^* (A;J)$ le sous-ensemble de $\mathcal{M}_{g,m} (A;J)$ constitué des courbes injectives quelque part. C'est un sous-ensemble ouvert. D'après la formule d'adjonction, on remarque que pour toutes courbes $u, u' \in \mathcal{M}_{g,m}^* (A;J)$, on a $\delta(u) = \delta(u')$. En particulier, si $u$ est une courbe plongée, toutes les autres courbes de l'espace de modules $\mathcal{M}_{g,m}^* (A;J)$ sont également plongées.

On appelle \emph{indice} d'une courbe $u : \Sigma \rightarrow M$ de genre $g$ l'entier 
$$\ind (u) = (n-3)\chi(\Sigma) + 2 c_1([u]),$$
et \emph{dimension virtuelle} de l'espace de modules des courbes injectives quelque part correspondant l'entier $$\virdim \mathcal{M}_{g,m}^* (A;J) = \ind (u) +2m.$$
Dans la suite, on utilisera ces formules seulement pour $n=2$, ce qui nous donne $$\ind (u) = -\chi(\Sigma) + 2 c_1([u]),$$
ou encore, dans le cas des courbes injectives quelque part, d'après la formule d'adjonction 
\begin{equation}\label{eq:indexformula}
\ind(u) = \chi(\Sigma) + 2 [u]^2 - 4 \delta (u).
\end{equation}

La proposition suivante est une simple reformulation en termes d'indices de courbes pseudoholomorphes de la formule de Riemann--Hurwitz vue précédemment dans la Section~\ref{s:revmult}.

\begin{prop}\label{p:RHindiceversion}
Soit $(M,J)$ une surface presque complexe, $u$ une courbe $J$--holomorphe non constante dans $(M,J)$ et $\tilde{u} = u \circ \varphi$ un revêtement multiple de $u$ de degré $k \geq 1$. Notons $Z(d \varphi) \in \mathbb{Z}$ la somme des ordres des points critiques de $\varphi$. On a alors 
$$\ind (\tilde{u}) = k \ind (u) + Z(d \varphi).$$
En particulier, on a $\ind (\tilde{u}) \geq k \ind (u)$, avec égalité si et seulement si le revêtement ne possède aucun point de ramification.
\end{prop}

Enfin, pour toute courbe $J$--holomorphe $( \Sigma, j ,u , ( \zeta_1, \dots, \zeta_m))$, on appelle \emph{groupe d'automorphisme} de $u$, noté $\mathrm{Aut}(u)$, l'ensemble des biholomorphismes $\varphi: (\Sigma,j) \rightarrow (\Sigma,j)$ tels que $u \circ \varphi = u$ et pour tout $i \in \{1, \dots, m \}$, $u ( \zeta_i) = \zeta_i$. On remarque que lorsque $u$ n'est pas constante, $\mathrm{Aut}(u)$ est fini, et lorsque $u$ est injective quelque part, $\mathrm{Aut}(u)$ est trivial.

On présente désormais quelques résultats de régularité sur les espaces de modules. En effet, pour une structure presque complexe $J$ générique, les espaces de modules de courbes injectives quelque part peuvent être munis de structures de variétés lisses (pas nécessairement compactes) de dimension finie, égale à leur dimension virtuelle. Pour comprendre certaines idées importantes derrière ce genre de résultats, on esquisse très rapidement le cadre dans lequel ils se situent.  Le sujet des courbes pseudoholomorphes étant très vaste, il ne peut être abordé dans le détail de manière exhaustive sans trop s'éloigner du but de ce manuscrit. Les détails techniques qui ne seront pas précisés ici sont couverts dans \cite{Wendl} et \cite{wendl2014lectures}, ou encore dans \cite{mcduff2004j}. Pour une variété presque complexe $(M,J)$, on appelle opérateur de Cauchy--Riemann non linéaire l'opérateur 
$$\begin{array}{ccccc}
\overline{\partial}_J & : & \mathcal{B} & \longrightarrow & \mathcal{E} \\
					&	& (j,u) & \longmapsto  & d u + J \circ du \circ j ,
\end{array}$$
où $\mathcal{B}$ désigne un espace de Banach de dimension infinie constitué des paires convenables $(j,u)$, avec $u : \Sigma \rightarrow M$ une application d'une certaine régularité et $j$ une structure complexe sur $\Sigma$ ; et $\mathcal{E}$ désigne un fibré vectoriel de Banach de dimension infini de base $\mathcal{B}$ (dont chaque fibre au-dessus de $(j,u)$ est un espace de Banach de sections du fibré au-dessus de $\Sigma$ constitué des morphismes anti-linéaires complexes $(T\Sigma, j) \rightarrow (u^*TM,J)$). La condition de $J$--holomorphicité pour une application lisse $u  : (\Sigma,j) \rightarrow (M,J)$ peut alors être reformulée par l'équation $\overline{\partial}_J u = 0$. Les espaces de modules de courbes paramétrées (c'est-à-dire avant le quotient par la relation d'équivalence décrite au début de cette section) peuvent être vus comme des sous-ensembles de $\overline{\partial}_J^{-1} \{ 0 \}$. On remarque que $\overline{\partial}_J$ peut être vu comme une section du fibré $\mathcal{E}$. On peut également voir les espaces de modules de courbes paramétrées comme des sous-ensembles de l'intersection entre la section $\overline{\partial}_J$ et la section nulle du fibré.

En chaque point $(j,u)$, on peut définir le linéarisé de l'opérateur de Cauchy--Riemann $D_{(j,u)} \overline{\partial}_J : T_{(j,u)} \mathcal{B} \mapsto \mathcal{E}_{(j,u)}$, qui est un opérateur borné entre deux espaces de Banach. L'opérateur $D_{(j,u)} \overline{\partial}_J$ est un \emph{opérateur de Fredholm}, c'est-à-dire que la dimension de son noyau et la codimension de son image sont des espaces vectoriels de dimension finie. L'\emph{indice de Fredholm} de $D_{(j,u)} \overline{\partial}_J$ est alors défini par l'entier 
$$\ind D_{(j,u)} \overline{\partial}_J = \dim \ker D_{(j,u)} \overline{\partial}_J - \codim \Ima D_{(j,u)} \overline{\partial}_J,$$ 
et correspond à la quantité $\ind (u)$ introduite précédemment. Lorsque le linéarisé en une courbe $J$--holomorphe $u : (\Sigma ,j) \rightarrow (M,J)$ de l'opérateur de Cauchy--Riemann $D_{(j,u)} \overline{\partial}_J$ est  surjectif, il est possible d'appliquer le théorème des fonctions implicites pour les variétés de Banach de dimension infinie pour munir l'espace de modules correspondant d'une structure lisse au voisinage de $u$. On dit dans ce cas que la courbe pseudoholomorphe $u$ est \emph{Fredholm régulière} (on dira parfois aussi que la courbe $u$ est transverse car elle peut être vue comme un point d'intersection transverse entre $\overline{\partial}_J (\mathcal{B}) \subset \mathcal{E}$ et la section nulle de $\mathcal{E}$). La propriété de régularité Fredholm permet non seulement de s'assurer que l'espace de modules est lisse au voisinage de $u$, mais permet également d'affirmer que $u$ survit à toute perturbation arbitrairement petite $J'$ de $J$ (c'est-à-dire qu'on peut trouver une courbe $J'$--holomorphe $\mathcal{C}^\infty$--proche de $u$).

On note $\mathcal{M}_{g,m}^{\mathrm{reg}} (A;J)$ l'ensemble des courbes de $\mathcal{M}_{g,m} (A;J)$ Fredholm régulières dont le groupe d'automorphisme est trivial. L'espace de modules $\mathcal{M}_{g,m}^{\mathrm{reg}} (A;J)$ est nécessairement un sous-ensemble ouvert de l'espace de modules $\mathcal{M}_{g,m} (A;J)$.

\begin{thm}
L'ensemble $\mathcal{M}_{g,m}^{\mathrm{reg}} (A;J)$ admet naturellement une structure de variété lisse orientée (non nécessairement compacte) avec 
$$\dim \mathcal{M}_{g,m}^{\mathrm{reg}} (A;J)=\virdim \mathcal{M}_{g,m}^{\mathrm{reg}} (A;J).$$
\end{thm}

\begin{rk}
Si on lève la condition sur le groupe d'automorphisme, on obtient que l'espace de modules correspondant est muni d'une structure d'orbifold.
\end{rk}

On peut s'assurer que les courbes $J$--holomorphes injectives quelque part sont Fredholm régulières en réalisant une perturbation arbitrairement petite de la structure presque complexe $J$, de telle sorte que pour toute courbe $J$--holomorphe $u : (\Sigma,j) \rightarrow (M,J)$, le linéarisé en $(j,u)$ de l'opérateur de Cauchy--Riemann $D_{(j,u)} \overline{\partial}_J$ est surjectif. Autrement dit, pour $J$ générique, toute courbe injective quelque part est transverse (Fredholm régulière).

Pour des entiers naturels $g,n$ fixés, on note $\mathcal{J}_\tau^{\mathrm{reg}}(M, \omega)$ l'ensemble des structures presque complexes sur $M$ dominées par $\omega$ telles que pour toute classe d'homologie $A \in H_2(M ;\mathbb{Z})$, toute courbe $u \in \mathcal{M}_{g,m}(A;J)$ avec un point injectif est Fredholm régulière. 

\begin{thm}\label{t:généricité}
Soit $(M, \omega)$ une variété symplectique et $g, m$ deux entiers naturels. L'ensemble $\mathcal{J}_\tau^{\mathrm{reg}}(M, \omega)$ est un sous-ensemble comaigre de $\mathcal{J}_\tau (M, \omega)$. En particulier tout $J \in \mathcal{J}_\tau (M, \omega)$ admet une perturbation $\mathcal{C}^\infty$ arbitrairement petite $J'$ dominée par $\omega$ telle que toutes les courbes $J'$--holomorphes injectives quelque part sont Fredholm régulières.
\end{thm}

\begin{rk}\label{rk:comaigre}
On rappelle que pour un espace topologique $X$, un sous-ensemble $Y$ de $X$ est dit \emph{comaigre} s'il contient une intersection dénombrable d'ouvert dense ; qu'une intersection dénombrable de sous-ensembles comaigres est un sous-ensemble comaigre ; que si $X$ est un espace métrique complet, le théorème de Baire nous indique tout sous-ensemble comaigre de $X$ est dense dans $X$. Lorsqu'on emploie le terme \og générique \fg{} pour un objet mathématique par la suite, cela sous-entendra qu'un résultat est vrai pour un ensemble comaigre de tels objets.
\end{rk}

\begin{rk}\label{rk:généricitéhorsrelcompact}
En fait, dans l'énoncé du Théorème~\ref{t:généricité}, on peut se contenter de réaliser la perturbation de la structure presque complexe $J$ en $J'$ dans un ouvert relativement compact $\mathcal{U} \subset M$ (et de laisser ainsi $J$ inchangée dans $M \backslash \mathcal{U}$). La conclusion reste alors vraie pour toutes les courbes $J'$--holomorphes qui possèdent un point injectif dans $\mathcal{U}$ (voir~\cite[Theorem~2.12]{Wendl}).
\end{rk}

Puisque toute courbe injective possède un groupe d'automorphisme trivial, on obtient le résultat suivant.

\begin{cor}\label{c:smoothmoduliinj}
Soit $(M, \omega)$ une variété symplectique. Pour une structure presque complexe générique $J$ dominée par $\omega$, l'espace de modules $\mathcal{M}^*_{g,m}(A;J)$ est une variété lisse orientée (non nécessairement compacte) de dimension égale à la dimension virtuelle de $\mathcal{M}_{g,m}(A;J)$.
\end{cor}

D'après la Remarque~\ref{rk:comaigre}, en prenant une intersection dénombrable sur tous les choix de $g$, $m$ et $A$ possibles, on peut  s'assurer que pour $J$ générique, tous les espaces de modules de la forme $\mathcal{M}^*_{g,m}(A;J)$ sont lisses (autrement dit que toutes les courbes injectives quelque part sont transverses).

\begin{ex}\label{e:toreautoint0}
Dans la variété $\mathcal{T}^2 \times \mathcal{S}^2$ munie d'une structure complexe scindée, les courbes complexes de la forme $\mathcal{T}^2 \times \{ * \}$ sont plongées, de genre $1$ et d'auto-intersection $0$. Leur indice est donc égal à $2 \times 1 -2 +2 \times 0 = 0$. Mais l'ensemble de ces courbes forme une variété de dimension réelle $2$. La structure presque complexe de départ n'est donc pas générique, et lorsqu'on la perturbe en une structure presque complexe générique sur $\mathcal{T}^2 \times \mathcal{S}^2$, seul un nombre fini de telles courbes pseudoholomorphes survivent à la perturbation.
\end{ex}

\subsubsection{Espaces de modules à paramètre}

On dispose de résultats similaires pour des familles de structures presque complexes dépendant de manière lisse d'un paramètre réel, ce qui se révèle être d'une grande utilité lorsqu'on considère des problèmes d'isotopie symplectique par exemple. Soit $(\omega_s)_{s \in [0,1]}$ un chemin lisse de formes symplectiques sur $M$ et $\mathcal{J}_{\tau}(M, (\omega_s)_{s \in [0,1]})$ l'espace des chemins lisses de structures presque complexes $(J_s)_{s \in [0,1]}$ tels que pour tout $s \in [0,1]$, $J_s \in \mathcal{J}_\tau (M,\omega_s)$. L'espace $\mathcal{J}_{\tau}(M, (\omega_s)_{s \in [0,1]})$ est non vide, contractile et muni d'une structure lisse naturelle. Pour deux structures presque complexes $J \in \mathcal{J}_\tau (M, \omega_0)$ et $J' \in \mathcal{J}_\tau (M, \omega_1)$, l'ensemble 
$$\mathcal{J}_\tau(M, (\omega_s)_{s \in [0,1]}; J, J') = \{ (J_s)_{s \in [0,1]} \in \mathcal{J}_{\tau}(M, (\omega_s)_{s \in [0,1]}) \mid J_0 = J ~ \mathrm{et}   ~ J_1 = J' \}$$ possède les mêmes propriétés.

Étant donné une famille $(J_s)_{s \in [0,1]} \in \mathcal{J}_{\tau}(M, (\omega_s)_{s \in [0,1]})$ de structures presque complexes, on peut définir un espace de modules à paramètre
$$ \mathcal{M}_{g,m}(A; (J_s)_{s \in [0,1]}) = \{ (s,u) \mid s \in [0,1] ~ \mathrm{et} ~ u \in \mathcal{M}_{g,m} (A; J_s) \}.$$
Cet espace est muni d'une topologie naturelle pour laquelle une suite $((s_k, u_k))$ converge vers une courbe $(s,u)$ si et seulement si $(s_k)$ converge vers $s$ et $(u_k)$ converge vers $u$. 

On note $\mathcal{M}_{g,m}^*(A; (J_s)_{s \in [0,1]})$ le sous-ensemble de $\mathcal{M}_{g,m}(A; (J_s)_{s \in [0,1]})$ constitué des paires $(s,u)$ pour lesquelles $u$ est injective quelque part, et $\mathcal{M}_{g,m}^\mathrm{reg}(A; (J_s)_{s \in [0,1]})$ le sous-ensemble de $\mathcal{M}_{g,m}(A; (J_s)_{s \in [0,1]})$ constitué des paires $(s,u)$ pour lesquelles $u$ est Fredholm régulière et possède un groupe d'automorphisme trivial. Notons que ces deux sous-ensembles sont ouverts.

On présente désormais des résultats de régularité pour ces espaces de modules (tout d'abord pour l'ensemble des courbes Fredholm régulières avec groupe d'automorphisme trivial, puis pour l'ensemble des courbes injectives quelque part pour un chemin de structures presque complexes générique). 

\begin{thm}\label{t:moduliparamreg}
Soit $(J_s)_{s \in [0,1]} \in \mathcal{J}_{\tau}(M, (\omega_s)_{s \in [0,1]})$. L'espace $\mathcal{M}_{g,m}^\mathrm{reg}(A; (J_s)_{s \in [0,1]})$ est naturellement muni d'une structure de variété lisse (non nécessairement compacte) orientée à bord avec 
$$\partial \mathcal{M}_{g,m}^\mathrm{reg}(A; (J_s)_{s \in [0,1]}) = -(\{ 0 \} \times \mathcal{M}_{g,m}^\mathrm{reg}(A; J_0)) \cup (\{ 1 \} \times \mathcal{M}_{g,m}^\mathrm{reg}(A; J_1)).$$
De plus, la projection
$$
\begin{array}{ccc}
\mathcal{M}_{g,m}^\mathrm{reg}(A; (J_s)_{s \in [0,1]})  & \longrightarrow & [0,1] \\
(s,u) & \longmapsto & s 
\end{array}
$$
est une submersion.
\end{thm}

\begin{thm}\label{t:moduliparaminj}
Soit $J \in \mathcal{J}_\tau^{\mathrm{reg}} (M, \omega_0)$ et $J' \in \mathcal{J}_\tau^\mathrm{reg} (M, \omega_1)$. Alors l'espace $J \in \mathcal{J}_\tau (M, (\omega_s)_{s \in [0,1]}; J, J')$ contient un ensemble comaigre $\mathcal{J}_\tau^\mathrm{reg} (M, (\omega_s)_{s \in [0,1]}; J, J')$ tel que pour tout chemin de structures presque complexes $(J_s)_{ s \in [0,1]} \in \mathcal{J}_\tau^\mathrm{reg} (M, (\omega_s)_{s \in [0,1]}; J, J')$, l'espace de modules $\mathcal{M}_{g,m}^*(A ; (J_s)_{s \in [0,1]})$ est naturellement muni d'une structure de variété lisse (non nécessairement compacte) orientée à bord avec 
$$\partial \mathcal{M}_{g,m}^*(A; (J_s)_{s \in [0,1]}) = -(\{ 0 \} \times \mathcal{M}_{g,m}^*(A; J_0)) \cup (\{ 1 \} \times \mathcal{M}_{g,m}^*(A; J_1)).$$
De plus, toutes les valeurs critiques de la projection
$$
\begin{array}{ccc}
\mathcal{M}_{g,m}^*(A; (J_s)_{s \in [0,1]})  & \longrightarrow & [0,1] \\
(s,u) & \longmapsto & s 
\end{array}
$$
appartiennent à $]0,1[$.
\begin{figure}[h]
	\centering
	\includegraphics[scale=0.3]{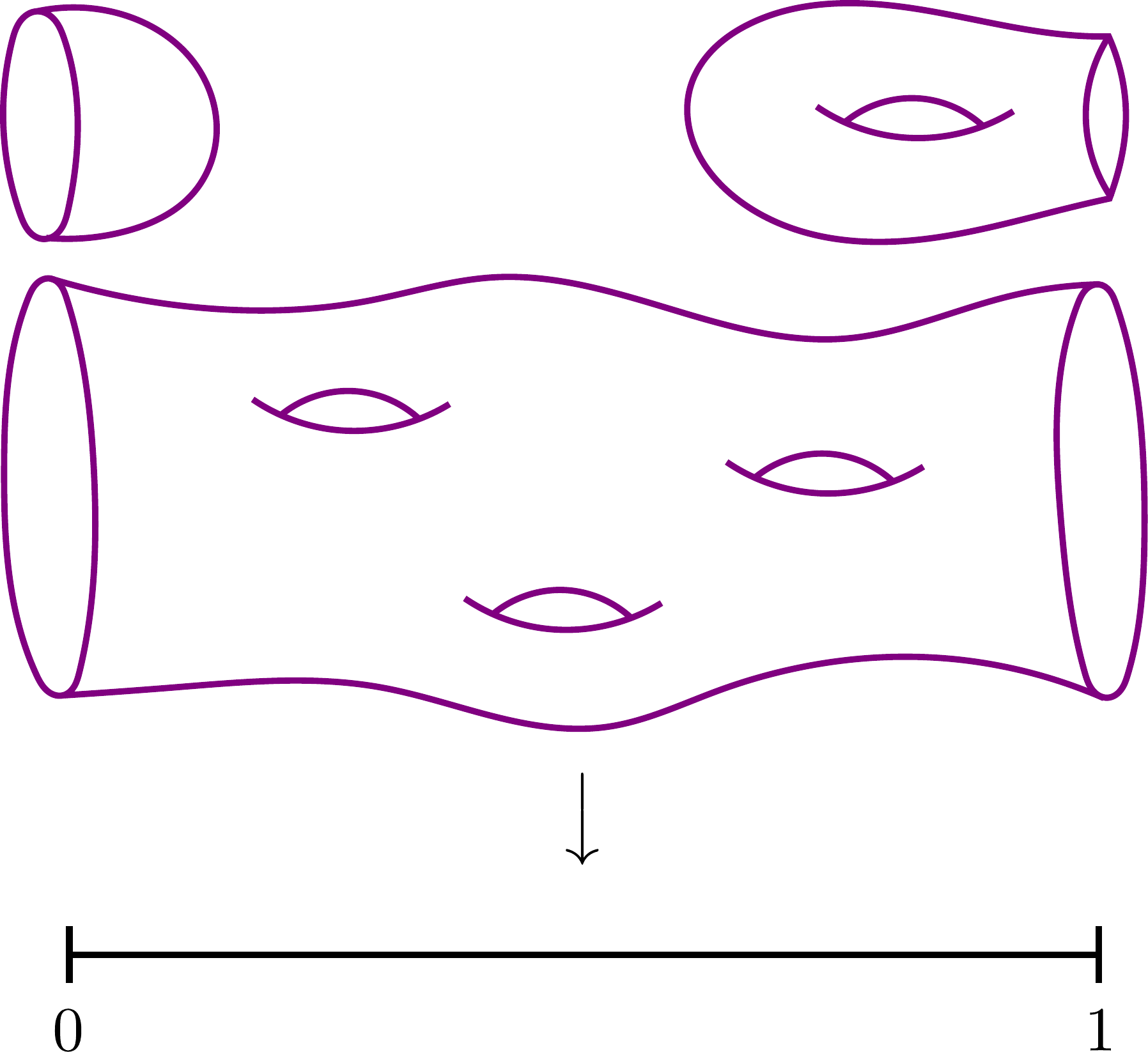}\\
	\caption{Représentation schématique d'un espace de modules de la forme $\mathcal{M}_{g,m}^*(A ; (J_s)_{s \in [0,1]})$, où $(J_s)_{s \in [0,1]}$ est un chemin générique de structures presque complexes dominées par un chemin de formes symplectiques $(\omega_s)_{s \in [0,1]}$.}
	\label{f:moduliparam}
\end{figure}
\end{thm}

Il est important de remarquer que pour un chemin générique de structures presque complexes $(J_s)_{ s \in [0,1]} \in \mathcal{J}_\tau^\mathrm{reg} (M, (\omega_s)_{s \in [0,1]})$, les structures presque complexes $J_s$ ne sont, quant à elles, pas nécessairement toutes génériques. Ainsi, même en choisissant un chemin de structures presque complexes $(J_s)_{ s \in [0,1]}$ générique, certaines courbes de l'espace de modules $\mathcal{M}_{g,m}^*(A; (J_s)_{s \in [0,1]})$ peuvent ne pas être Fredholm régulières. C'est la raison de l'existence possible de valeurs critiques dans le Théorème~\ref{t:moduliparaminj}. Par conséquent, la topologie de l'espace de modules $\mathcal{M}_{g,m}^*(A; J_s)$ peut dépendre de la valeur du paramètre $s$.

Ces différents résultats, avec ou sans paramètre, imposent des contraintes sur les indices des courbes $J$--holomorphes injectives quelque part dans le cas où $J$ est générique.

\begin{cor}
Soit $(M, \omega)$ une variété symplectique et $J \in \mathcal{J}_\tau (M, \omega)$. Après une perturbation générique de $J$ dominée par $\omega$, toutes les courbes $J$--holomorphes $u$ injectives quelque part vérifient $\ind (u) \geq 0$.
\end{cor}

\begin{proof}
Notons $g$ le genre d'une courbe $J$--holomorphe injective quelque part $u$. D'après le Corollaire \ref{c:smoothmoduliinj}, l'espace de modules $\mathcal{M}_g^*([u];J)$ est une variété lisse non vide de dimension $\ind (u)$, d'où $\ind(u) \geq 0$.
\end{proof}

\begin{cor}\label{c:indgeq-1}
Soit $M$ une variété et $( \omega_s)_{s \in [0,1]}$ un chemin de formes symplectiques sur $M$. Alors pour une famille de structures presque complexes générique $(J_s)_{s \in [0,1]} \in \mathcal{J}_\tau (M, (\omega_s)_{ s \in [0,1]})$, toute courbe $J_s$--holomorphe $u$ injective quelque part vérifie $\ind (u) \geq -1$.
\end{cor}

\begin{proof}
Notons $g$ le genre d'une courbe $J_s$--holomorphe injective quelque part $u$. D'après le Théorème~\ref{t:moduliparaminj}, l'espace de modules $\mathcal{M}_{g,m}^*(A ; (J_s)_{s \in [0,1]})$ est une variété lisse non vide de dimension $\ind (u)+1$, d'où $\ind(u) \geq -1$.
\end{proof}

\begin{rk}\label{r:indgeq-1à0}
Dans le cadre décrit dans ce manuscrit, les indices des courbes pseudoholomorphes sont toujours pairs. Ainsi, la conclusion du Corollaire~\ref{c:indgeq-1} peut être renforcée en $\ind (u) \geq 0$. Ce n'est en revanche pas forcément le cas dans un cadre plus général (lorsqu'on considère des courbes pseudoholomorphes épointées par exemple, voir~\cite[Section 8.3]{Wendl}).
\end{rk}

\subsection{Contraintes ponctuelles}

Pour chaque espace de modules $\mathcal{M}_{g,m} (A;J)$, on dispose d'une \emph{application d'évaluation}
$$\mathrm{ev} : \mathcal{M}_{g,m}(A;J) \rightarrow M^m.$$

Grâce à cette application d'évaluation, on peut construire des espaces de modules de courbes avec des contraintes ponctuelles de la manière suivante. Soit $Z \subset M^m$ une sous-variété. Pour tous $J \in \mathcal{J}_\tau (M, \omega)$, $(J_s)_{s \in [0,1]} \in \mathcal{J}_\tau (M, (\omega_s)_{s \in [0,1]})$, on définit les \emph{espaces de modules contraints} par
$$\mathcal{M}_{g,m}(A;J;Z) = \mathrm{ev}^{-1}(Z) \subset \mathcal{M}_{g,m}(A;J),$$
$$\mathcal{M}_{g,m}(A;(J_s)_{s \in [0,1]};Z) = \{ (s,u) \mid s \in [0,1], u \in \mathcal{M}_{g,m}(A;J_s;Z) \}.$$
Autrement dit, les éléments de $\mathcal{M}_{g,m}(A;J;Z)$ peuvent être paramétrés par des courbes des $J$--holomorphes $u : \Sigma \rightarrow M$ avec des points marqués $\zeta_1, \dots , \zeta_m \in \Sigma$ satisfaisants $(u(\zeta_1), \dots, u(\zeta_m)) \in Z$.

Dans la suite, on utilisera principalement ces espaces de modules contraints dans le cas où $Z$ est un point $(p_1, \dots, p_m) \in M^m$, et on notera $\mathcal{M}_{g,m}(A;J;p_1, \dots, p_m) = \mathcal{M}_{g,m}(A;J;Z)$. Dans ce cas, les éléments de $\mathcal{M}_{g,m}(A;J;Z)$ peuvent être paramétrés par des courbes des $J$--holomorphes $u : \Sigma \rightarrow M$ avec des points marqués $\zeta_1, \dots , \zeta_m \in \Sigma$ satisfaisants $(u(\zeta_1), \dots, u(\zeta_m))=(p_1, \dots, p_m)$. 

On dit qu'une courbe $u \in \mathcal{M}_{g,m} (A;J; Z)$ est \emph{Fredholm régulière pour le problème avec contraintes} si c'est une courbe Fredholm régulière telle que l'application d'évaluation $\mathrm{ev} :\mathcal{M}_{g,m}(A;J) \rightarrow M^m$ est transverse à la sous-variété $Z$ (dans le cas où $Z = (p_1, \dots, p_m)$, cela revient à demander que $(p_1, \dots, p_m)$ est une valeur régulière de $\mathrm{ev}$). On note $\mathcal{M}^\mathrm{reg}_{g,m} (A;J; Z)$ l'ensemble des courbes Fredholm régulières pour le problème avec contraintes. L'espace $\mathcal{M}^\mathrm{reg}_{g,m} (A;J; Z)$ est alors une sous-variété lisse de $\mathcal{M}_{g,m}^{\mathrm{reg}} (A;J)$ de codimension $\codim (Z)$ et un sous-ensemble ouvert de $\mathcal{M}_{g,m} (A;J; Z)$.

Ceci motive à définir la \emph{dimension virtuelle de l'espace de modules avec contraintes} par
$$\virdim \mathcal{M}_{g,m}(A;J;Z)
 = \virdim \mathcal{M}_{g,m}(A;J) - \codim (Z),$$
qu'on appellera parfois l'\emph{indice contraint} d'une courbe $u \in \mathcal{M}_{g,m}(A;J;Z)$.
Lorsque $(M,J)$ est une surface presque complexe, on a alors
\begin{align*}
\virdim \mathcal{M}_{g,m}(A;J;p_1,\dots, p_m )
& = \virdim \mathcal{M}_{g,m}(A;J) - 4m \\
& = -(2-2g) + 2 c_1 (A) +2m -4m \\
& = \ind (u) - 2m.
\end{align*}
Autrement dit chaque contrainte ponctuelle diminue de $2$ la dimension virtuelle de l'espace de modules.

De la même manière que précédemment, on note $\mathcal{M}^*_{g,m}(A;J;Z) \subset \mathcal{M}_{g,m}(A;J;Z)$ et $\mathcal{M}^*_{g,m}(A;(J_s)_{s \in [0,1]} ;Z) \subset \mathcal{M}_{g,m}(A;(J_s)_{s \in [0,1]} ;Z)$ les sous-ensembles constitués des courbes injectives quelque part. On note également $\mathcal{J}_\tau^{\mathrm{reg}}(M, \omega;Z)$ l'ensemble des structures presque complexes $J \in \mathcal{J}_\tau (M, \omega)$ tel que toute courbe $J$--holomorphe injective quelque part est Fredholm régulière pour le problème avec contraintes $Z$. On présente désormais quelques résultats de régularité pour les espaces de modules contraints.

\begin{thm}
Soit $(M, \omega)$ une variété symplectique, $g$ et $m$ deux entiers naturels, $A \in H_2 (M; \mathbb{Z})$ et $Z$ une sous-variété de $M^m$. L'ensemble $\mathcal{J}_\tau^{\mathrm{reg}}(M, \omega;Z)$ est un sous-ensemble comaigre de $\mathcal{J}_\tau (M, \omega )$. 
\end{thm}

\begin{rk}
Il est important ici de souligner que l'espace $\mathcal{J}_\tau^{\mathrm{reg}}(M, \omega;Z)$ dépend du choix de $Z$ (on dira parfois que $J$ est \og générique pour $Z$ \fg{}). En effet, on ne peut en général s'attendre à ce qu'une structure presque complexe soit générique pour une quantité indénombrable de contraintes (comme tous les $m$--uplets de $M^m$ par exemple). Alternativement pour $J \in \mathcal{J}_\tau^\mathrm{reg} (M, \omega)$, on peut pour un choix générique d'une sous-variété $Z \subset M^m$, s'assurer que l'espace de modules $\mathcal{M}^*_{g,m}(A;J;Z)$ est une variété lisse.
\end{rk}

\begin{ex}\label{e:cubiquespar9pts}
L'ensemble des polynômes homogènes à trois variables avec coefficients complexes forme un $\mathbb{C}$--espace vectoriel de dimension $10$. Donc l'ensemble $\mathcal{C}_3$ des courbes complexes de degré $3$ dans l'espace projectif $\mathbb{C}P^2$ forme un espace projectif complexe de dimension $9$ (notons que l'ensemble des cubiques non singulières est un sous-ensemble ouvert de $\mathcal{C}_3$).

Fixons désormais $p_1, \dots, p_8 \in \mathbb{C}P^2$ des points génériques. L'ensemble des cubiques passant par ces huit points forme un sous-espace projectif $\mathcal{C}_3^{p_1, \dots, p_8}$ de $\mathcal{C}_3$ de dimension complexe $1$. En effet, chaque condition d'appartenance d'un $p_i$ à une courbe $C$ définit une équation $\mathbb{C}$--linéaire en les coefficients du polynôme $P$ qui définit $C$. Comme les points $p_i$ sont génériques, les diverses contraintes définissent un système linéaire de rang $8$ en les $10$ coefficients de $P$. L'ensemble $\mathcal{C}_3^{p_1, \dots, p_8}$ est par conséquent un sous-espace projectif complexe de $\mathcal{C}_3$ de dimension $1$. Prenons $C_1$ et $C_2$ deux cubiques distinctes passant par $p_1, \dots, p_8$ d'équations respectives $P_1(x,y,z) = 0$ et $P_2(x,y,z) = 0$. L'ensemble $\mathcal{C}_3^{p_1, \dots, p_8}$ est alors exactement l'ensemble des courbes ayant une équation de la forme $ \lambda_1 P_1(x,y,z) + \lambda_2 P_2(x,y,z) = 0$, avec $\lambda_1, \lambda_2 \in \mathbb{C}$.

Par le théorème de Bézout, les courbes $C_1$ et $C_2$ s'intersectent en un neuvième point $p_9$. Mais la relation obtenue au paragraphe précédent nous indique alors que toutes les cubiques passant par $p_1, \dots, p_8$ passent également par $p_9$. L'espace de modules des courbes complexes de $\mathbb{C}P^2$ de genre $1$, homologues à $3h$ (où $h$ désigne la classe de la droite) passant par $p_1, \dots, p_9$ forme donc une variété (non compacte) de dimension réelle $2$. Or en utilisant la formule d'adjonction, on constate que la dimension virtuelle de cet espace de modules est égale à $0 +2(0 +3^2) -2 \times 9 =0$. Cette différence est due au fait que la contrainte $\{ p_1, \dots, p_9 \}$ n'est pas générique pour la structure complexe standard sur $\mathbb{C}P^2$.
\end{ex}

\begin{thm}
Soit $M$ une variété et $(\omega_s)_{s \in [0,1]}$ un chemin lisse de structures symplectiques sur $M$, $J \in \mathcal{J}_\tau^{\mathrm{reg}}(M, \omega_0;Z), J' \in \mathcal{J}_\tau^{\mathrm{reg}}(M, \omega_1;Z)$, $g$ et $m$ deux entiers naturels et $A \in H_2 (M; \mathbb{Z})$.
Il existe un sous-ensemble comaigre
$$\mathcal{J}_\tau^{\mathrm{reg}}(M, (\omega_s)_{s \in [0,1]};Z;J,J') \subset \mathcal{J}_\tau(M, (\omega_s)_{s \in [0,1]};Z;J,J')$$
tel que pour tout chemin $(J_s)_{s \in [0,1]} \in \mathcal{J}_\tau^{\mathrm{reg}}(M, (\omega_s)_{s \in [0,1]};Z;J,J')$, l'espace de modules $ \mathcal{M}_{g,m}(A; (J_s)_{s \in [0,1]} ; Z)$ est une variété lisse (non nécessairement compacte) à bord, avec
$$\partial \mathcal{M}_{g,m}^*(A; (J_s)_{s \in [0,1]};Z) = -(\{ 0 \} \times \mathcal{M}_{g,m}^*(A; J_0;Z)) \cup (\{ 1 \} \times \mathcal{M}_{g,m}^*(A; J_1;Z)).$$
De plus, toutes les valeurs critiques de la projection
$$
\begin{array}{ccc}
\mathcal{M}_{g,m}^*(A; (J_s)_{s \in [0,1]};Z)  & \longrightarrow & [0,1] \\
(s,u) & \longmapsto & s 
\end{array}
$$
appartiennent à $]0,1[$.
\end{thm}

\subsection{Contraintes sur les dérivées} \label{s:contraintesdérivées}

Soit $(M,J)$ une variété presque complexe de dimension réelle $2n$. Les points marqués permettent également d'extraire des informations sur les dérivées des courbes pseudoholomorphes, ce qui donne lieu à de nombreux résultats de généricité particulièrement utiles lorsqu'on s'intéresse à des configurations de courbes. Dans ce manuscrit, on se restreint aux contraintes à l'ordre $1$ sur les dérivées (voir~\cite[Section 2.1.5]{Wendl} pour une présentation plus générale). Afin de pouvoir réaliser de telles contraintes, on introduit l'espace des $1$--jets
$$\mathrm{Jet}^1_J (M) = \{ (p, \Phi ) \mid p \in M ~ \text{et} ~ \Phi : ( \mathbb{C}, i) \rightarrow (T_p M, J) ~ \text{est}  ~ \mathbb{C}\text{--linéaire} \}.$$
On peut voir l'espace $\mathrm{Jet}^1_J (M)$ comme l'ensemble des classes d'équivalence d'applications $J$--holomorphes définies sur un voisinage de l'origine de $\mathbb{C}$ à valeurs dans $(M,J)$, où deux applications sont équivalentes si leurs dérivées coïncident en l'origine. Ici, dans le couple $(p, \Phi)$, $p$ sert à représenter la valeur d'une fonction $J$--holomorphe en l'origine, et $\Phi$ sa dérivée en l'origine. L'espace $\mathrm{Jet}^1_J (M)$ est une variété non compacte de dimension $\dim \mathrm{Jet}^1_J (M) =4n$.

Étant donné un espace de modules $\mathcal{M}_{g,m} (A;J)$, un choix coordonnées complexes au voisinage de chacun des points marqués de chacune des courbes de $\mathcal{M}_{g,m} (A;J)$ permet de définir une application d'évaluation sur l'espace des $1$--jets
$$\mathrm{ev} : \mathcal{M}_{g,m} (A;J) \rightarrow (\mathrm{Jet}^1_J (M))^m$$
qui nous permet d'imposer des contraintes sur les dérivées des courbes de $\mathcal{M}_{g,m} (A;J)$ en les points marqués.
Pour en tirer des informations pertinentes, il est intéressant de définir de telles sous-variétés $Z \subset (\mathrm{Jet}^1_J (M))^m$ de façon à ce que l'application d'évaluation ne dépende pas du choix de coordonnées effectué. 
On peut par exemple considérer 
$$Z_{\mathrm{crit}} = \left\{ (p,0) \in \mathrm{Jet}^1_J (M) \mid p \in M \right\}.$$
Lorsque $m=1$, l'ensemble $\mathcal{M}_{g,\mathrm{crit}} (A;J) =\mathrm{ev}^{-1} (Z_{\mathrm{crit}}) \subset \mathcal{M}_{g,1} (A;J)$ est alors constitué des courbes de $\mathcal{M}_{g,1} (A;J)$ qui admettent un point critique en leur point marqué. Remarquons que $\codim Z_{\mathrm{crit}} = 4n - \dim Z_{\mathrm{crit}} = 2n$.
On peut également combiner cette contrainte sur les dérivées avec des contraintes ponctuelles $Z \subset M^m$. On considère par exemple la sous-variété
$$Z \times Z_\mathrm{crit} \subset M^m \times \mathrm{Jet}_J^1(M),$$
ainsi que l'application d'évaluation correspondante
$$\mathrm{ev} : \mathcal{M}_{g,m+1} (A;J) \rightarrow M^m \times \mathrm{Jet}^1_J (M).$$
On note alors
$$\mathcal{M}_{g,m, \mathrm{crit}} (A;J;Z) = \mathrm{ev}^{-1}(Z \times Z_\mathrm{crit}) \subset \mathcal{M}_{g,m+1} (A;J; Z),$$
l'espace de modules constitué des courbes $[(\Sigma, j ,u , ( \zeta_1, \dots, \zeta_{m+1}))]$ de $\mathcal{M}_{g,m+1}$ vérifiant $(u( \zeta_1), \dots, u( \zeta_m)) \in Z$ et $d_{\zeta_{m+1}} u =0$.
On définit la dimension virtuelle d'un tel espace de modules de la manière suivante 
\begin{align*}
\virdim \mathcal{M}_{g,m,\mathrm{crit}} (A;J;Z)
& = \virdim \mathcal{M}_{g,m+1} (A;J) - \codim Z - 2n\\
& = \virdim \mathcal{M}_{g,m}(A;J;Z) -2(n-1).
\end{align*}

On présente désormais un théorème de régularité pour ces espaces de modules. Comme précédemment, on note $\mathcal{M}_{g,m, \mathrm{crit}}^* (A;J;Z) \subset \mathcal{M}_{g,m, \mathrm{crit}} (A;J;Z)$ le sous-ensemble ouvert constitué des courbes injectives quelque part.

\begin{thm}\label{t:crit}
Soit $Z \subset M^m$ une sous-variété. Pour $J \in \mathcal{J}_\tau^\mathrm{reg} (M, \omega;Z)$ l'espace $\mathcal{M}_{g,m, \mathrm{crit}}^* (A;J;Z)$ est une variété lisse (non nécessairement compacte) de dimension égale à sa dimension virtuelle.
\end{thm}

On peut également considérer pour un entier $k \in \{1, \dots, m\}$, l'espace de modules $\mathcal{M}_{g,m, \mathrm{crit}_k}^* (A;J;Z) \subset \mathcal{M}_{g, m}^* (A;J;Z)$ constitué des courbes ayant un point critique en leur $k$--ième point marqué. En suivant le même procédé que précédemment, on obtient le résultat suivant.

\begin{thm}\label{t:critk}
Soit $Z \subset M^m$ une sous-variété. Pour $J \in \mathcal{J}_\tau (M, \omega)$ générique l'espace $\mathcal{M}_{g,m, \mathrm{crit_k}}^* (A;J;Z)$ est une variété lisse (non nécessairement compacte) de dimension égale à sa dimension virtuelle $$\virdim \mathcal{M}_{g,m, \mathrm{crit_k}} (A;J;Z) = \virdim \mathcal{M}_{g,m} (A;J;Z) - 2n.$$
\end{thm}

On obtient alors le corollaire suivant en appliquant le théorème de Sard aux applications d'oubli $\mathcal{M}_{g,m, \mathrm{crit}}^* (A;J;Z) \rightarrow \mathcal{M}_{g,m}^* (A;J;Z)$ et $\mathcal{M}_{g,m, \mathrm{crit}_k}^* (A;J;Z) \rightarrow \mathcal{M}_{g,m}^* (A;J;Z)$.

\begin{cor}
Soit $(M, \omega)$ une variété symplectique de dimension $2n \geq 4$ et $Z \subset M^m$ une sous-variété. Pour $J \in \mathcal{J}_\tau (M, \omega)$ générique et $A \in H_2(M; \mathbb{Z})$ le sous-ensemble des courbes immergées de $\mathcal{M}_{g, m}^* (A;J;Z)$ est un ouvert dense.
\end{cor}

\begin{rk}\label{r:pasdecritind0}
En conséquence du Théorème~\ref{t:crit} et du Théorème~\ref{t:critk}, toute courbe $u$ d'un espace de modules de la forme $\mathcal{M}_{g, m}^* (A;J;Z)$, avec $J$ générique par rapport à $Z$, dont l'indice contraint est nul par rapport à $Z$ est immergée.
\end{rk}

En considérant l'ensemble 
$$Z_\mathrm{tan} = \left\{ ((p, \Phi), (p, \Psi)) \in (\mathrm{Jet}^1_J)^2 \mid p \in M  ~ \text{et} ~ \exists c \in \mathbb{C}, ~ \Psi = c \Phi \neq 0 \right\},$$
et en appliquant le même procédé (remarquons que $\codim Z_\mathrm{tan} = 4n-2$), on obtient la proposition suivante.

\begin{prop}
Soit $(M, \omega)$ une variété symplectique de dimension $2n \geq 4$ et $Z \subset M^m$ une sous-variété. Pour $J \in \mathcal{J}_\tau (M, \omega)$ générique le sous-ensemble des courbes sans point d'auto-intersection tangentiel de $\mathcal{M}_{g, m}^* (A;J;Z)$ est un ouvert dense.
\end{prop}

Plus généralement, on peut généraliser ce procédé à des produits cartésiens d'espaces de modules et obtenir de nombreux résultats de généricité de cette manière. On obtient alors les principes généraux suivants. Soit $(M, \omega)$ une surface symplectique, $J \in \mathcal{J}_\tau (M ,\omega)$ et $u_1, \dots, u_\ell$ un nombre fini de courbes $J$--holomorphes simples Fredholm régulières dans $M$. Il existe $J'$ une perturbation générique $\mathcal{C}^\infty$--proche de $J$, et $u_1', \dots, u_\ell'$ des courbes $J'$--holomorphes $\mathcal{C}^\infty$--proches de $u_1, \dots, u_\ell$ telles que :
\begin{itemize}
\item pour tout $i \in \{1, \dots, \ell \}$, la courbe $u_i$ est immergée, et tous ses points d'auto-intersection sont des points doubles transverses,
\item pour tous $i,j \in \{1, \dots, \ell \}$ distincts, les courbes $u_i$ et $u_j$ ne s'intersectent qu'en des points injectifs, et ces intersections sont transverses,
\item pour tous $i,j,k \in \{1, \dots, \ell \}$ deux à deux distincts, il n'existe pas d'intersection commune aux trois courbes $u_i$, $u_j$ et $u_k$.
\end{itemize}
Ainsi les seuls types de singularités pouvant apparaître dans la configuration donnée par les courbes Fredholm régulières $u_1', \dots, u_\ell'$ sont des points doubles transverses positifs.

\subsection{Transversalité automatique}

Dans le cas des surfaces presque complexes, certaines conditions numériques permettent d'assurer la régularité Fredholm des courbes pseudoholomorphes immergées, et ce, indépendamment de la structure presque complexe choisie. On désigne généralement ce genre de théorèmes sous l'appellation de \og théorèmes de transversalité automatique \fg{}. C'est un des rares moyens de vérifier explicitement qu'une courbe est Fredholm régulière et c'est également une propriété incontournable lorsqu'on s'intéresse à des problèmes d'isotopie symplectique.

\begin{thm}[{\cite[Theorem 2.44]{Wendl}}]\label{t:transauto1} 
Soit $(M,J)$ une surface presque complexe et $u$ une courbe $J$--holomorphe immergée de genre $g$. Si $\ind (u) > 2g-2$, alors la courbe $u$ est Fredholm régulière.
\end{thm}

Cela signifie que les courbes $J$--holomorphes rationnelles immergées sont automatiquement transverses dès que leur indice est positif. Comme le montre le théorème suivant, ceci est également vrai pour l'indice contraint lorsqu'on s'intéresse au problème avec des contraintes ponctuelles pour des courbes rationnelles. 

\begin{thm}[{\cite[Theorem~2.46]{Wendl}}] \label{t:transauto2} 
Soit $(M,J)$ une surface presque complexe, $p_1, \dots , p_m \in M$ et $u \in \mathcal{M}_{g,m}(J;p_1, \dots , p_m)$ une courbe rationnelle immergée. Si $\ind (u) - 2m \geq 0$, alors u est Fredholm régulière pour le problème avec contraintes ponctuelles $p_1, \dots, p_m$.
\end{thm}

Plus généralement, Sikorav a montré le résultat suivant. Mais avant de l'énoncer, je tiens à remercier Chris Wendl qui, lors de sa visite au LMJL au début de l'année 2020, alors que je n'avais pas encore connaissance de l'existence de ce résultat, l'a redémontré indépendamment après une des discussions qu'on a eues et a proposé de fournir une annexe pour ce manuscrit.

\begin{thm}[{\cite[Proposition~1]{Sikorav}}] \label{t:transautoSiko}
Soit $(M,J)$ une surface presque complexe, $p_1, \dots , p_m \in M$ et $u \in \mathcal{M}_{g,m}(J ; p_1, \dots , p_m)$ une courbe immergée. Si $\ind (u)-2m  > 2g-2$, alors $u$ est Fredholm régulière pour le problème avec contraintes $p_1, \dots, p_m$.
\end{thm}

\begin{rk}
En utilisant la formule de l'indice, la condition $\ind (u) > 2g-2 + 2m$ peut être reformulée en $c_1([u]) >2m$. D'après la formule d'adjonction, si la courbe $u$ est plongée, cette condition est aussi équivalente à $[u]^2 >2g -2 +m$.
\end{rk}

\begin{ex}
\leavevmode
\begin{enumerate}
\item Soit $J$ une structure presque complexe sur $\mathbb{C}P^2$ et $S$ l'image d'une courbe $J$--holomorphe plongée de degré $d \in \mathbb{N}^*$ dans $(\mathbb{C}P^2, J)$ (c'est à dire homologue à $dh$, où $h$ désigne la classe d'homologie d'une droite complexe). La courbe $S$ vérifie la formule du genre (qui est une conséquence directe de la formule d'adjonction) 
$$g = \frac{(d-1)(d-2)}{2}.$$
Donc on a 
$$[S]^2 = d^2 = 2g -2 +3d > 2g - 2.$$
Ainsi, toutes les courbes pseudholomorphes plongées dans $\mathbb{C}P^2$ sont automatiquement Fredholm régulières d'après le Théorème~\ref{t:transautoSiko}.
\item Les cubiques lisses passant par $p_1, \dots, p_8$ dans l'Exemple~\ref{e:cubiquespar9pts} sont automatiquement Fredholm régulières pour le problème avec contraintes $p_1, \dots,p_8$ d'après le Théorème~\ref{t:transautoSiko}. Par contre, les hypothèses du Théorème~\ref{t:transautoSiko} ne sont plus satisfaites pour le problème avec contraintes $p_1, \dots, p_9$. Cet exemple montre en particulier que les conditions sur les indices contraints dans le Théorème~\ref{t:transauto2} et le Théorème~\ref{t:transautoSiko} sont optimales.
\item On considère la variété $\mathcal{S}^2 \times \mathcal{S}^2$ munie sa structure complexe standard, qu'on note ici $J_{st}$. Les courbes complexes de la forme $\mathcal{S}^2 \times \{ * \}$ et $\{ * \} \times \mathcal{S}^2$ sont des courbes plongées de genre $0$ et d'auto-intersection $0$. Elle satisfont donc les hypothèses du Théorème~\ref{t:transauto1}. Donc les espaces de modules de la forme $\mathcal{M}_{g,m}^*([\mathcal{S}^2 \times \{ * \}];J_{st})$ ou $\mathcal{M}_{g,m}^*([\{ * \} \times \mathcal{S}^2];J_{st})$ sont des variétés lisses de dimension réelle $2+ 2 \times 0 =2$ (en fait on peut montrer qu'ils sont difféomorphes à $\mathcal{S}^2$ et ce indépendamment de la structure presque complexe choisie).
\item Dans l'Exemple~\ref{e:toreautoint0} concernant la variété $\mathcal{T}^2 \times \mathcal{S}^2$ munie d'une structure complexe scindée, les courbes de la forme $\mathcal{T}^2 \times \{ * \}$ sont plongées d'indice $0$ mais elle forment une variété de dimension réelle $2$. En effet, ces courbes sont de genre $1$ et d'auto-intersection $0$ et ne satisfont donc pas aux hypothèses du Théorème~\ref{t:transauto1}. Cet exemple montre en particulier que la condition sur l'indice du Théorème~\ref{t:transauto1} est optimale.
\end{enumerate}
\end{ex}

\subsection{Compactification de Gromov}\label{s:CompactificationGromov}

L'espace de modules $\mathcal{M}_{g,m}(A;J)$ admet une compactification naturelle, qu'on appelle la \emph{compactification de Gromov}, notée $\overline{\mathcal{M}}_{g,m}(A;J)$, de sorte que $\mathcal{M}_{g,m}(A;J)$ est un sous-ensemble ouvert de $\overline{\mathcal{M}}_{g,m}(A;J)$. L'objectif de cette section est d'expliquer en quoi consiste ce procédé de compactification.

On définit l'espace de modules des \emph{courbes nodales stables $J$--holomorphes non paramétrées}, de genre $g$, avec $m$ points marqués et homologues à $A \in H_2 (M)$ comme l'ensemble de classes d'équivalence suivant 
$$\overline{\mathcal{M}}_{g,m}(A;J) = \{(\Sigma,j,u, ( \zeta_1, \dots , \zeta_m), \Delta) \} \diagup \sim$$
où
\begin{enumerate}
\item $(\Sigma,j)$ est une union disjointe de surfaces de Riemann $\Sigma_1, \dots, \Sigma_\ell$,
\item $u : (\Sigma ,j ) \rightarrow (M,J)$ est une courbe $J$--holomorphe vérifiant 
$$[u]:= \sum\limits_{i=1}^\ell u_* ([\Sigma_i]) =A,$$
\item $(\zeta_1, \dots, \zeta_m)$ est un $m$--uplet de points deux à deux distincts de $\Sigma$,
\item $\Delta = \{ \{ \hat{z}_1, \check{z}_1 \}, \dots, \{ \hat{z}_r, \check{z}_r \} \}$ est un ensemble fini de paires non ordonnées de points de $\Sigma$ deux à deux distincts, et tous distincts des points marqués $\zeta_i$, qu'on appelle \emph{points nodaux}, tel que pour tout $i \in \{1, \dots, r\}$, $u \left( \hat{z}_i \right) = u\left( \check{z}_i \right)$,
\item $(\Sigma, j ,u , ( \zeta_1, \dots, \zeta_m), \Delta) \sim (\Sigma', j' ,u' , ( \zeta_1', \dots, \zeta_m'), \Delta' )$ s'il existe un biholomorphisme $\varphi : (\Sigma, j)  \rightarrow (\Sigma',j')$ tel que 
\begin{itemize}
\item $u = u' \circ \varphi$,
\item pour tout $i \in \{1, \dots, m \}$, $\varphi(\zeta_i) = \zeta_i'$,
\item on peut écrire $\Delta = \{ \{ \hat{z}_1, \check{z}_1 \}, \dots, \{ \hat{z}_r, \check{z}_r \} \}$, $\Delta' = \{ \{ \hat{z}_1', \check{z}_1' \}, \dots, \{ \hat{z}_r', \check{z}_r' \} \}$ de telle sorte que pour tout $i \in \{1, \dots, r \}$, $\varphi (\hat{z}_i) = \varphi (\hat{z}_i')$ et $\varphi (\check{z}_i) = \varphi (\check{z}_i')$,
\end{itemize}
\item la condition sur le genre $g$ se traduit de la manière suivante : pour chaque paire de points nodaux $\left\{ \hat{z}_i , \check{z}_i \right\} \in \Delta$, on retire des disques arbitrairement petits centrés en $\hat{z}_i$ et  $\check{z}_i$, et on colle ensemble les deux bords ainsi créés, pour obtenir au final une surface de Riemann connexe, compacte, sans bord $\hat{\Sigma}$ de genre $g$,
\item la condition de stabilité est donnée par la Définition \ref{d:stable}.
\end{enumerate}

La condition de stabilité a été introduite par Kontsevich afin de définir une topologie séparée sur $\overline{\mathcal{M}}_{g,m}$, permettant ainsi d'assurer l'unicité des limites de suites de courbes pseudoholomorphes convergentes. 

\begin{df}\label{d:stable}
Soit $(\Sigma ,j,u, ( \zeta_1 , \dots, \zeta_m ), \Delta )$ une courbe nodale. On note $\dot{S}$ la surface épointée obtenue en retirant de $S$ tous les points marqués $\zeta_1 , \dots, \zeta_m$ et tous les points nodaux de $\Delta$. On dit que $(\Sigma ,j,u, ( \zeta_1 , \dots, \zeta_m ), \Delta )$ est \emph{stable} lorsque chaque composante connexe de $\dot{S}$ sur laquelle $u$ est constante possède une caractéristique d'Euler strictement négative.
\end{df}

\begin{ex}
Une composante constante de genre $0$ est stable dès qu'elle possède au moins $3$ points marqués ou nodaux. Ces composantes constantes stables de genre $0$ sont appelées des \emph{bulles fantômes}. Une composante constante de genre $1$ est stable dès qu'elle possède au moins $1$ point marqué ou un point nodal. Les composantes de genre supérieur ou égal à $2$ sont toujours stables. 
\begin{figure}[h]
	\centering
	\includegraphics[scale=0.5]{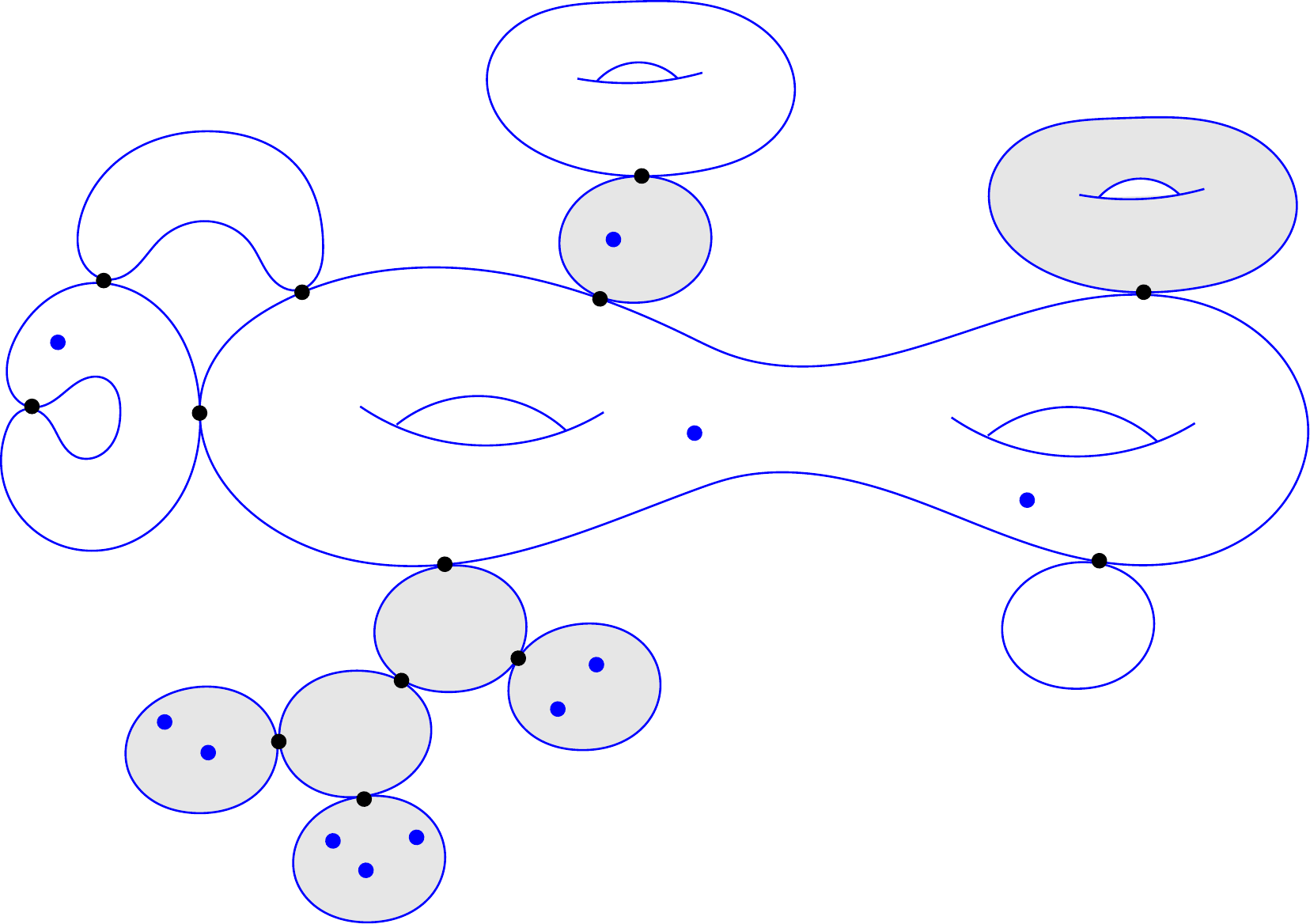}\\
	\caption{Représentation possible d'une courbe nodale stable de genre $6$ avec $11$ points marqués (en bleu) et $13$ paires de points nodaux (en noir). Les composantes constantes sont grisées.}
	\label{f:nodalcurve}
\end{figure}
\end{ex}

Le compactifié de Gromov $\overline{\mathcal{M}}_{g,m}$ est un espace topologique métrisable, qui contient naturellement $\mathcal{M}_{g,m}$ comme sous-ensemble ouvert (c'est le sous-ensemble des courbes nodales de la forme $[\Sigma,j,u, (\zeta_1, \dots, \zeta_m),\emptyset]$, qu'on appellera \emph{courbes lisses} par la suite). Plus généralement, pour une courbe nodale $[\Sigma,j,u, (\zeta_1, \dots, \zeta_m),\Delta]$ (qu'on notera parfois simplement $u$), on peut restreindre $u$, $j$ et les points marqués $\zeta_1, \dots, \zeta_m$ à chacune des composantes connexes de $\Sigma$ et obtenir de cette façon les \emph{composantes lisses} de $[\Sigma,j,u, (\zeta_1, \dots, \zeta_m),\Delta]$, qui sont des éléments d'espaces de modules de la forme $\mathcal{M}_{h,k} (B;J)$.


On s'intéresse désormais à l'énoncé du théorème de compacité de Gromov. Pour une variété symplectique $(M, \omega)$, munie d'une structure presque complexe $J$ dominée par $\omega$, on définit l'\emph{énergie} d'une courbe $J$--holomorphe $u$ comme l'aire
$$E_\omega (u) = \int_\Sigma u^* \omega.$$
C'est un réel positif, qui est nul si et seulement si $u$ est constante. Remarquons que cette quantité ne dépend que des classes $[u] \in H_2 (M; \mathbb{Z})$ et $[\omega] \in H^2_{dR} (M; \mathbb{R})$. En particulier toutes les courbes d'un espace de modules de la forme $\mathcal{M}_{g,m} (A;J)$ possèdent la même énergie.

\begin{thm}[Théorème de compacité de Gromov]
Soit $M$ une variété et $A \in H_2(M; \mathbb{Z})$. On considère $(\omega_k)$ une suite de formes symplectiques sur $M$ qui converge de manière $\mathcal{C}^\infty$ vers une forme symplectique $\omega$, $(J_k)$ une suite de structures presque complexes telle que pour tout $k$, $J_k \in \mathcal{J}_\tau (\omega_k)$, qui converge vers $J \in \mathcal{J}_\tau (\omega)$. Soit une suite $(u_k)$ de courbes pseudoholomorphes d'énergie uniformément bornée telle que pour tout $k$, $u_k \in \mathcal{M}_{g,m} (A;J_k)$. Alors on peut extraire de la suite $(u_k)$ une sous-suite qui converge vers une courbe nodale stable de $\overline{\mathcal{M}}_{g,m} (A;J)$.
\end{thm}
On renvoit à~\cite{Wendl} pour plus de détails sur la manière dont la convergence s'effectue dans le compactifié de Gromov.

\begin{rk}\label{r:limitenonplongée}
Soit $(M, \omega)$ une surface symplectique, $J \in \mathcal{J}_\tau (M, \omega)$ et $u : \Sigma \rightarrow M$ une courbe $J$--holomorphe plongée de genre $g$. Lorsque $m \leq 1$, toute courbe nodale $[\Sigma_\infty,j_\infty,u_\infty, (\zeta_1, \dots, \zeta_m),\Delta] \in \overline{\mathcal{M}}_{g,m} ([u];J) \backslash \mathcal{M}_{g,m} ([u];J)$ possède au moins deux composantes non constantes ou au moins une composante non constante qui n'est pas plongée. Pour le montrer, raisonnons par l'absurde et supposons que $u$ ne possède qu'une composante non constante et que cette composante non constante $u_0 : \Sigma_0 \rightarrow M$ est une courbe pseudoholomorphe plongée de genre $g_0$ (en particulier il n'y pas deux points nodaux d'une même paire $\{ \hat{z}_i, \check{z}_i \}$ sur $\Sigma_0$). Par la formule d'adjonction pour les courbes pseudoholomorphes plongées, puisque $[u_0]=[u]$, on a $g=g_0$. La condition sur le genre assure alors que $u$ ne possède pas de composante constante de genre strictement positif. On considère maintenant le graphe $\Gamma$ dont les sommets sont les composantes connexes de $\Sigma$ et dont les arêtes correspondent aux paires de points nodaux de $\Delta$. D'après la condition sur le genre, on a $b_1(\Gamma)=0$, autrement dit $\Gamma$ est un arbre. Chaque bulle fantôme qui représente une feuille de $\Gamma$ possède exactement un point nodal. Puisque $m \leq 1$, ces bulles fantômes ne peuvent pas vérifier la condition de stabilité. Par conséquent $u_\infty$ n'admet pas de composante constante. Au final, on a $\Delta = \emptyset$, donc $u$ appartient à $\mathcal{M}_{g,m} ([u];J)$, ce qui contredit les hypothèses.
\end{rk}

\begin{rk}\label{r:cuspisgenus}
Soit $(M, \omega)$ une surface symplectique, $J \in \mathcal{J}_\tau (M, \omega)$, $u : \Sigma \rightarrow M$ une courbe $J$--holomorphe simple de genre $g$ et $[\Sigma_\infty,j_\infty,u_\infty,\Delta] \in \overline{\mathcal{M}}_{g} ([u];J)$ une courbe nodale qui ne possède qu'une composante non constante $u_0 : \Sigma_0 \rightarrow M$. On suppose de plus que $u_0$ est simple. Si $u_\infty$ possède une composante constante de genre strictement positif, on a $g_0 < g$ d'après la condition sur le genre. Or par la formule d'adjonction, on obtient $g_0 + \delta(u_0) = g + \delta (u)$. Par conséquent on a  $\delta(u_0) > \delta (u)$. La composante $u_0$ possède donc une singularité en l'image de chaque composante constante de genre strictement positif (car la convergence dans le compactifié de Gromov s'effectue localement de manière $\mathcal{C}^\infty$ en dehors des points nodaux : les nouvelles singularités ne peuvent apparaître qu'aux images des paires de points nodaux).
\end{rk}

\begin{ex}
On considère $\mathbb{C}P^2$ muni de sa structure complexe standard $J_{st}$ et on note $h$ la classe d'homologie d'une droite complexe.  L'espace de modules $\mathcal{M}_{1} (3h;J_{st})$ comprend l'ensemble des cubiques non singulières de $\mathbb{C}P^2$ et l'ensemble des courbes de genre $1$ qui sont des revêtements triples des droites de $\mathbb{C}P^2$. L'image d'une courbe de $\overline{\mathcal{M}}_{1} (3h;J_{st}) \backslash \mathcal{M}_{1} (3h;J_{st})$ peut être (voir la Figure~\ref{f:nodalcubics}) :
\begin{enumerate}[label=(\alph*)]
\item une cubique cuspidale (une composante simple de genre $0$ reliée à une composante constante de genre $1$ par une paire de points nodaux ; la composante constante a pour image la singularité cuspidale),
\item une cubique nodale (une composante simple de genre $0$ sur laquelle se trouve une paire de points nodaux),
\item l'union d'une droite et d'une conique qui s'intersectent transversalement (deux composantes plongées de genre $0$ d'images distinctes, reliées entre elles par deux paires de points nodaux),
\item l'union d'une droite et d'une conique qui s'intersectent tangentiellement (une composante constante de genre $1$ à la laquelle sont reliées, chacune par une paire de points nodaux, deux composantes plongées de genre $0$ d'images distinctes),
\item l'union de trois droites qui s'intersectent deux à deux transversalement en des points distincts (trois composantes plongées de genre $0$, d'images deux à deux distinctes, reliées deux à deux par trois paires de points nodaux),
\item l'union de trois droites concourantes (une composante constante de genre $1$ à la laquelle sont reliées, chacune par une paire de points nodaux, trois composantes plongées de genre $0$ d'images deux à deux distinctes),
\item l'union d'une droite comptée avec multiplicité $2$ et d'une autre droite distincte (deux possibilités : une composante constante de genre $1$ à la laquelle sont reliées, chacune par une paire de points nodaux, trois composantes plongées de genre $0$ dont exactement deux d'entre elles ont la même image ; ou bien une composante de genre $1$, qui est un revêtement ramifié double, reliée à une composante plongée de genre $0$ par une paire de points nodaux),
\item une droite comptée avec multiplicité $3$ (deux possibilités : une composante constante de genre $1$ à la laquelle sont reliées, chacune par une paire de points nodaux, trois composantes plongées de genre $0$ qui ont la même image ; ou bien une composante de genre $1$, qui est un revêtement ramifié double, reliée par une paire de points nodaux à une composante plongée de genre $0$ qui possède la même image).
\end{enumerate}
On remarque que le genre de chaque composante constante correspond au genre de la singularité sur laquelle elle est envoyée.
\begin{figure}[htbp]
	\centering
	\subfloat[Cubique cuspidale.]{\label{f:cubica}\includegraphics[height=0.20\linewidth]{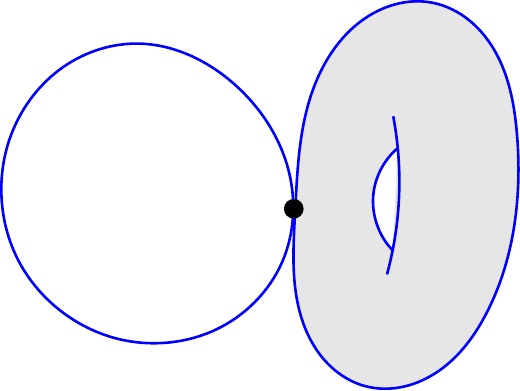}}\qquad
	\subfloat[Cubique nodale.]{\label{f:cubicb}\includegraphics[height=0.20\linewidth]{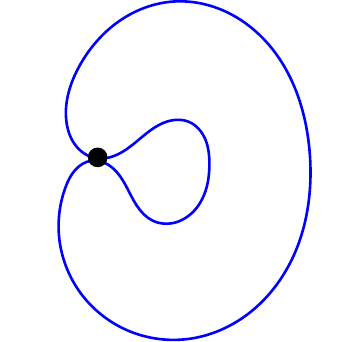}}\qquad
	\subfloat[Union d'une droite et d'une conique qui s'intersectent transversalement.]{\label{f:cubicc}\includegraphics[height=0.20\linewidth]{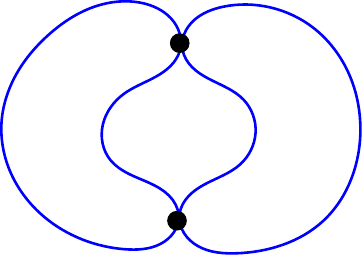}}\\
	\subfloat[Union d'une droite et d'une conique qui s'intersectent tangentiellement.]{\label{f:cubicd}\includegraphics[height=0.20\linewidth]{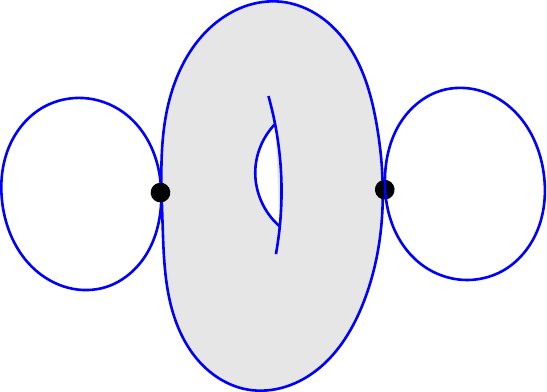}}\qquad
	\subfloat[Union de trois droites génériques.]{\label{f:cubice}\includegraphics[height=0.20\linewidth]{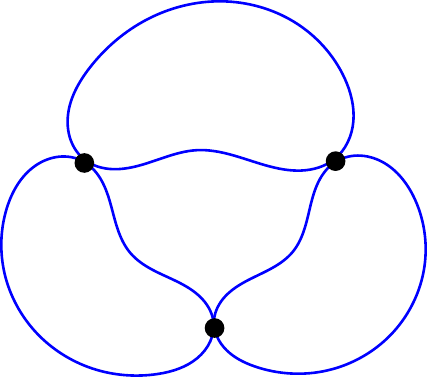}}\qquad
	\subfloat[Union de trois droites concourantes.]{\label{f:cubicf}\includegraphics[height=0.20\linewidth]{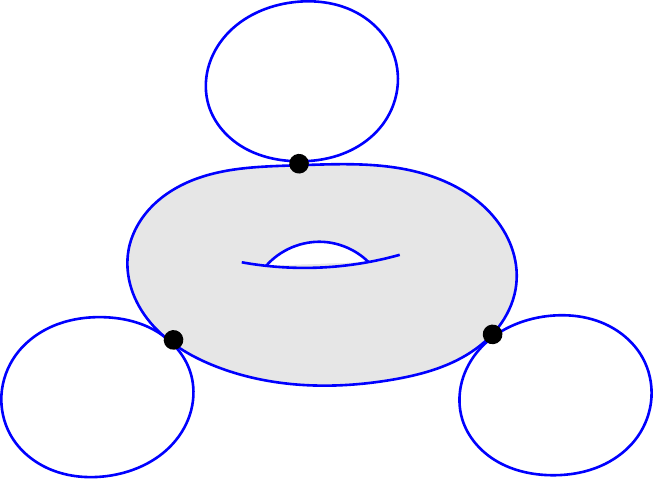}}\\	
	\subfloat[Union d'une droite comptée avec multiplicité $2$ et d'une autre droite distincte.]{\label{f:cubicg}\includegraphics[height=0.20\linewidth]{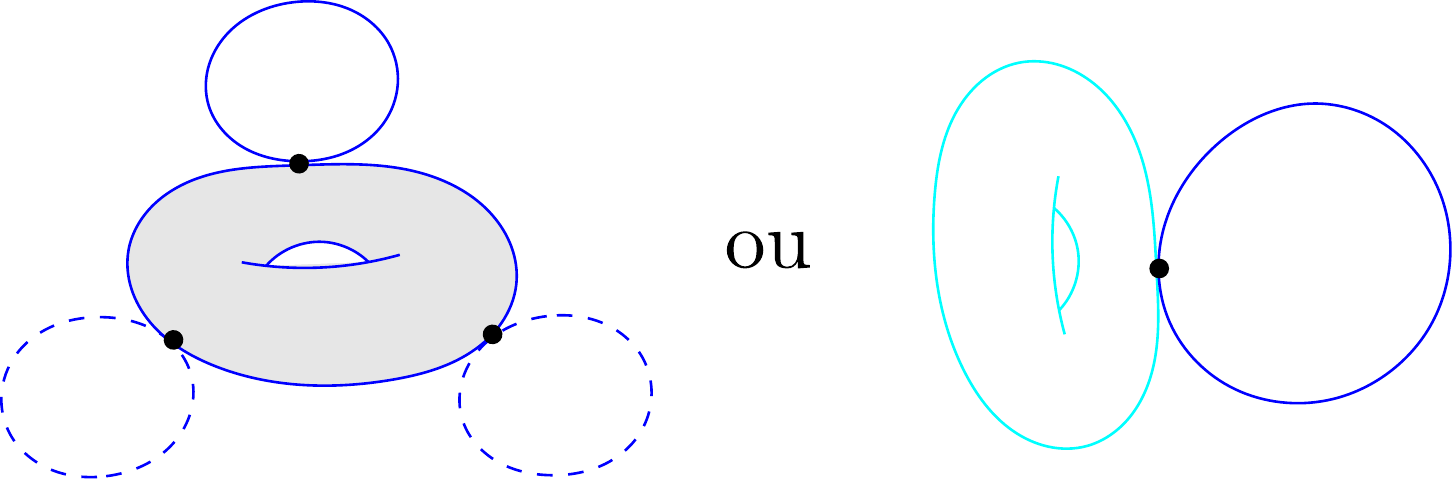}}\\
	\subfloat[Une droite comptée avec multiplicité $3$.]{\label{f:cubich}\includegraphics[height=0.20\linewidth]{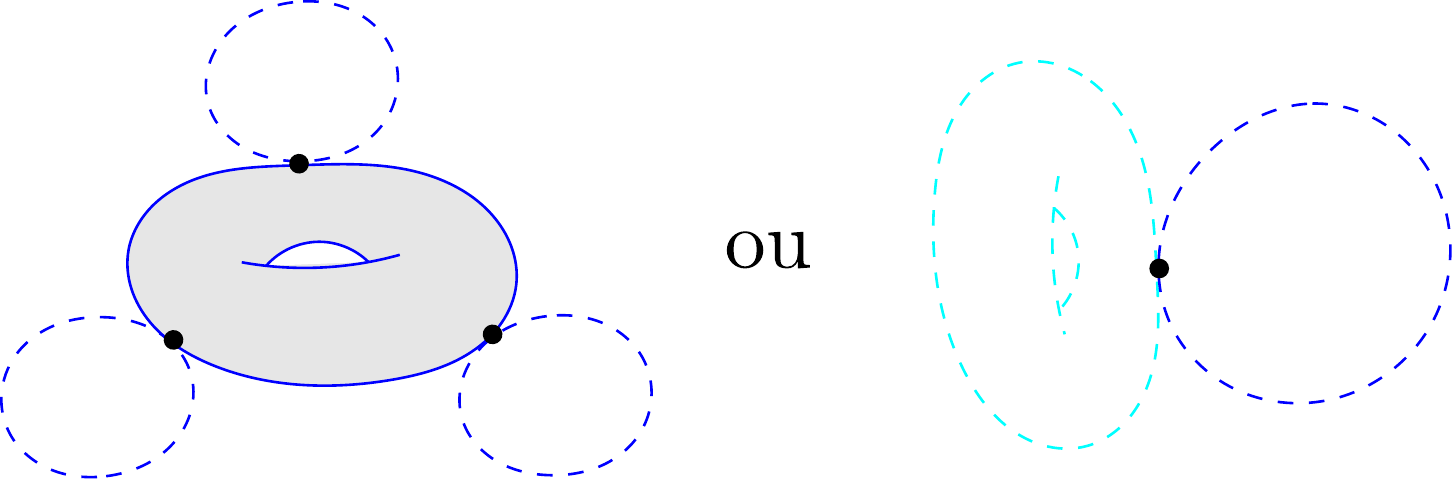}}
\caption{Les différentes possibilités pour les courbes nodales dans $\overline{\mathcal{M}}_{1} (3h;J_{st}) \backslash \mathcal{M}_{1} (3h;J_{st})$. Les composantes constantes sont grisées, les composantes qui sont des revêtements ramifiés doubles sont en cyan, les autres composantes sont simples. Les composantes en pointillé ont la même image.}\label{f:nodalcubics}
\end{figure}
\end{ex}

\begin{rk}\label{r:lissageptdoubleJholo}
Mentionnons qu'il est toujours possible de lisser de manière $J$--holomorphe un point double transverse $p$ d'une courbe $J$--holomorphe $u$ de genre $g$, simple et Fredholm régulière (voir~\cite{Sikorav}). C'est-à-dire qu'on peut trouver une courbe $J$--holomorphe $v$ de genre $g+1$, simple et Fredholm régulière, $\mathcal{C}^0$--proche de $u$ (et $\mathcal{C}^\infty$--proche de $u$ en dehors d'un voisinage de $p$). De même, il est également possible de lisser un point double transverse entre deux courbes $J$--holomorphes simples, Fredolm régulières $u_1$ et $u_2$ de genres respectifs $g_1$ et $g_2$. C'est-à-dire qu'on peut trouver une courbe $J$--holomorphe $v$ de genre $g_1+g_2$, simple et Fredholm régulière, $\mathcal{C}^0$--proche de la courbe nodale formée de l'union de $u_1$ et $u_2$ (et $\mathcal{C}^\infty$--proche de cette même courbe nodale en dehors d'un voisinage de $p$).
\end{rk}

\section{Éclatements, contractions}\label{s:eclat}

Dans cette section, on présente les opérations de contraction et d'éclatement pour les variétés symplectiques de dimension $4$. Ces opérations, inspirées de la géométrie algébrique complexe, jouent un rôle central en ce qui concerne l'étude des applications birationnelles entre deux surfaces complexes. Dans le cadre de la topologie symplectique, elles permettent de réduire l'étude des surfaces symplectiques à celles qui sont dites minimales.

\subsection{Éclatements et contractions complexes}\label{ss:eclatcomp}

On commence par présenter le modèle local dans le cas complexe. On note
$$\tilde{\mathbb{C}}^2 = \left\{ (\ell,v) \in \mathbb{C}P^1 \times \mathbb{C}^2  \mid v \in \ell \right\},$$
l'espace total du \emph{fibré tautologique} au-dessus de $\mathbb{C}P^1$, dont l'application de projection est donnée par 
$$\begin{array}{ccccc}
\pi &:	& \tilde{\mathbb{C}}^2 	& \longrightarrow 	& \mathbb{C}P^1 \\
	&	& (\ell, v) 			& \longmapsto		& \ell 
\end{array}.$$
La fibre au dessus de chaque point $\ell \in \mathbb{C}P^1$ correspond à l'ensemble des points de $\ell$, vue comme droite vectorielle de $\mathbb{C}^2$.
L'espace $\tilde{\mathbb{C}}^2$ admet une autre application de projection naturelle 
$$\begin{array}{ccccc}
\beta 	&:	& \tilde{\mathbb{C}}^2 	& \longrightarrow 	& \mathbb{C}^2 \\
		&	& (\ell, v) 			& \longmapsto		& v
\end{array},$$
qu'on appelle \emph{application de contraction}. On dit que l'espace $\tilde{\mathbb{C}}^2$ est l'\emph{éclatement} de $\mathbb{C}^2$ en l'origine. On appelle \emph{diviseur exceptionnel} (ou parfois \emph{sphère exceptionnelle}), la courbe complexe définie par $E = \pi^{-1} ( \{ 0 \}) \subset \tilde{\mathbb{C}}^2$ et on dit que $\beta$ \emph{contracte} $E$. Il est aisé de voir que $E$ est biholomorphe à $\mathbb{C}P^1$ et que $\beta$ se restreint à un biholomorphisme de $\tilde{\mathbb{C}}^2 \backslash E$ vers $\mathbb{C}^2 \backslash \{0 \}$. Autrement dit, on passe de $\mathbb{C}^2$ à $\tilde{\mathbb{C}}^2$ en remplaçant l'origine de $\mathbb{C}^2$ par l'ensemble $E \simeq \mathbb{C}P^1$ des directions complexes en l'origine. Remarquons enfin que le diviseur exceptionnel $E$ peut être vu comme la section nulle du fibré tautologique, qui est un fibré vectoriel en droites complexes de degré $-1$, d'où $$[E] \cdot [E] = \langle c_1 ( \tilde{\mathbb{C}}^2) , [\mathbb{C}P^1] \rangle = -1.$$
Ce nombre d'auto-intersection traduit une forme de rigidité du diviseur exceptionnel. En effet, on ne peut pas perturber un diviseur exceptionnel de manière complexe, sinon cela contredirait la propriété de positivité d'intersection.

\`A l'aide de ce modèle local, on peut définir plus généralement l'\emph{éclatement} d'une surface complexe $X$ en un point $z \in X$ de la manière suivante. On identifie un voisinage $U_z$ de $z$ avec un voisinage de l'origine dans $\mathbb{C}^2$ via une carte holomorphe, puis on se sert du biholomorphisme $\beta : \tilde{\mathbb{C}}^2 \backslash E \tilde\longrightarrow \mathbb{C}^2 \backslash \{0 \}$ pour munir l'ensemble 
$$\tilde{X} = \left(X \backslash \{ z \} \right) \cup \mathbb{P} ( T_z X)$$ 
d'une structure complexe (avec cette notation, on a identifié le diviseur exceptionnel $E$ avec le projectivisé complexe de $T_z X$). L'opération consistant à passer de $X$ à $\tilde{X}$ est appelée l'\emph{éclatement} de $X$ en $z$ (on dit également qu'on a \emph{éclaté} $X$ en $z$). Comme précédemment, on peut interpréter l'opération d'éclatement comme le remplacement d'un point dans une surface complexe par l'ensemble $\mathbb{P} ( T_z X)$ des directions complexes en ce point. On dispose d'une application naturelle $\beta : \tilde{X} \rightarrow X$, l'\emph{application de contraction}, qui envoie le diviseur exceptionnel $E = \mathbb{P} ( T_z X)$ sur $z$ (on dit que $\beta$ \emph{contracte} $E$) et qui se restreint en un biholomorphisme de $\tilde{X} \backslash E $ vers $X \backslash \{ z \}$. Comme le nombre d'auto-intersection de $E$ peut être calculé dans un voisinage de $E$, on a également $[E] \cdot [E] = -1$.

\begin{rk}\label{r:blowuptop}
Une autre manière possible de comprendre ce nombre d'auto-intersection consiste à remarquer que $\tilde{X}$ est difféomorphe à $X \# \overline{\mathbb{C}P}^2$, où $\overline{\mathbb{C}P}^2$ désigne le plan projectif complexe muni de l'orientation opposée à l'orientation usuelle. Cette description permet également de constater l'effet de l'opération d'éclatement sur certains invariants topologiques de $X$. En effet, on a $\pi_1 ( \tilde X) \simeq \pi_1(X)$, $b_2^+(\tilde X) = b_2^+ (X)$, $b_2^-(\tilde X) = b_2^- (X) +1 $ et pour tout $i \in \{0,1,3,4 \}$, $b_i(\tilde X) = b_i (X)$. Au niveau des formes d'intersection, on a $H_2(\tilde X;\mathbb{Z}) \simeq H_2( X;\mathbb{Z}) \oplus \langle E \rangle$ et $\mathcal{Q}_{ \tilde X} = \mathcal{Q}_ X \oplus (-1)$.
\end{rk}

Inversement, si $E$ est l'image d'une courbe rationnelle plongée holomorphe dans une surface complexe $\tilde X$ qui vérifie $[E] \cdot [E] = -1$, alors on peut identifier biholomorphiquement un voisinage de $E$ dans $\tilde X$ avec un voisinage de la section nulle du fibré tautologique (c'est le critère de contractibilité de Castelnuovo, voir~\cite[Theorem II.17]{beauville1996complex}). En se servant du modèle local, on peut alors remplacer le voisinage de $E$ dans $\tilde{X}$ par un voisinage de l'origine dans $\mathbb{C}^2$. On appelle cette opération la \emph{contraction} du \emph{diviseur exceptionnel} $E$ dans $\tilde{X}$. 

\begin{rk}
Puisque le $b_2^-$ d'une surface complexe est toujours fini, il n'est possible de réaliser des contractions successives qu'un nombre fini de fois.
\end{rk}

Enfin pour deux surfaces complexes $X$ et $Y$, on dit que $Y$ est un \emph{éclaté} (resp. un \emph{contracté}) de $X$ si $Y$ peut être obtenue en effectuant un nombre fini d'éclatements (resp. de contractions) à partir de $X$.

\begin{df}
Soit $X$ une surface complexe. On dit que $X$ est \emph{minimale} si elle ne contient aucun diviseur exceptionnel, c'est-à-dire si elle ne contient pas de courbe complexe rationnelle plongée d'auto-intersection $-1$.

Si $C \subset X$ est une courbe complexe, on dit que $X$ \emph{relativement minimale} par rapport à $C$ si $X \backslash C$ ne contient pas de diviseur exceptionnel.
\end{df}

Toute surface complexe peut être obtenue en éclatant un nombre fini de points sur surface minimale. Pour cette raison, il est souvent possible de se ramener au cas minimal lorsqu'on étudie les surfaces complexes.

\begin{rk}\label{r:birationalequiv}
Mentionnons aussi l'importance des opérations d'éclatements et des contractions dans l'étude des applications birationnelles entre deux surfaces complexes (on rappelle qu'une application rationnelle est un morphisme défini en dehors d'un fermé de Zariski et qu'une application birationnelle est une application rationnelle qui possède un inverse qui est aussi une application rationnelle). En effet, pour deux surfaces complexes $X$ et $Y$, et une application birationnelle $f : X \dashrightarrow Y$, il existe $\pi_X : \tilde X \rightarrow X$ et $\pi_Y : \tilde{Y} \rightarrow Y$ des éclatés respectifs de $X$ et $Y$ et un isomorphisme $\tilde f : \tilde X \xrightarrow{\sim} \tilde Y$ tels que $\pi_Y \circ \tilde f =  f \circ \pi_X$ (voir~\cite[Corollary II.12]{beauville1996complex}). Autrement dit, l'application $f$ peut être décomposée en plusieurs étapes : on éclate d'abord un certain nombre de fois $X$, puis on applique l'isomorphisme $\tilde f$ et enfin on effectue les contractions adéquates vers $Y$. 
\end{rk}

On s'intéresse désormais à l'effet des éclatements sur les courbes complexes (possiblement singulières) dans une surface complexe $X$. Soit $C \subset X$ une courbe complexe et $z$ un point de $C$. On note $\tilde X$ l'éclaté de $X$ en $z$ et $\beta : \tilde{X} \rightarrow X$ l'application de contraction correspondante. L'adhérence $\tilde C$ de l'ensemble $\beta^{-1} \left( C \backslash \{ z \} \right)$ est une courbe complexe de $\tilde X$, appelée la \emph{transformée propre} (ou \emph{transformée stricte}) de $C$, tandis que l'ensemble $\beta^{-1} (C)$ est appelée la \emph{transformée totale} de $C$. La transformée totale est l'union de la transformée stricte et du diviseur exceptionnel. Il est alors possible d'écrire $[C] = [\tilde C] + k [E]$, avec $k \in \mathbb{Z}$. Une étude minutieuse dans des coordonnées locales révèle que l'entier $k = [\tilde C] \cdot [E]$ est égal à la multiplicité de $z$ en tant point de $C$ (voir~\cite{wall2004singular}). Il est alors possible de calculer l'auto-intersection de $\tilde C$ de la manière suivante $[\tilde C]^2 = [C]^2 - k^2$. En particulier, si $z$ est un point non singulier, on a $[ \tilde C ] = [C] + [E]$ et $[\tilde C]^2 = [C]^2 - 1$. 

Comme les courbes $C$ et $\tilde C$ admettent des bonnes paramétrisations de même genre (pour une bonne paramétrisation $\tilde u$ de $\tilde C$, l'application $u = \beta \circ \tilde u$ est une bonne paramétrisation de $C$), les genres géométriques des deux courbes sont égaux. Enfin, on a $p_a (C) = g(C) + \delta (C)$ et $p_a (C') = g(C') + \delta (C')$. D'après la formule d'adjonction et la relation entre les nombres d'auto-intersection obtenue précédemment, on a $\delta(C') = \delta (C) - \frac{(k-1)(k-2)}{2}$, on en déduit $p_a (\tilde C) = p_a (C) - \frac{(k-1)(k-2)}{2}$.

\begin{rk}
En fait, l'opération d'éclatement permet de diminuer la \og complexité \fg{} des singularités de la courbe $C \subset X$. Au bout d'un nombre fini d'éclatements, la transformée propre de $C$ est une courbe non singulière. On appelle ce procédé la \emph{résolution des singularités} de $C$. La résolution des singularités utilisant le moins d'éclatements possible telle que la transformée totale de $C$ est l'union de courbes complexes non singulières qui s'intersectent transversalement en des points deux à deux distincts est appelée la \emph{résolution à croisements normaux}. Ce procédé est détaillé dans le livre~\cite{wall2004singular}.
\end{rk}

\subsection{Éclatements et contractions symplectiques}

On peut aussi définir les opérations d'éclatement et de contraction pour les surfaces symplectiques. Les détails techniques sont couverts dans~\cite[Chapter~3]{Wendl}. On reprend les notations de la sous-section précédente et pour tout réel $r >0$, on note $B^4(r)$ la boule de centre l'origine et de rayon $r$ dans $\mathbb{C}^2$. Commençons par remarquer que pour tout $R >0$, la $2$--forme 
$$\omega_R = \beta^* \omega_{st} + R^2 \pi^* \omega_{FS}$$
définit une forme de kählerienne sur $\tilde{\mathbb{C}}^2$.
Le diviseur exceptionnel $E$ dans $\tilde{\mathbb{C}}^2$ est alors une courbe symplectique et on peut montrer qu'un voisinage $\beta^{-1} (B^4(r))$ de $E$ dans $(\tilde{\mathbb{C}}^2, \omega_R)$ privé de $E$ est symplectomorphe à une région annulaire centrée en $0$ de $(\mathbb{C}^2, \omega_{st})$ de la forme $B^4(r) \backslash B^4 (R)$, avec $r > R$.
Étant donné une surface symplectique $(M, \omega)$ et $z \in M$, l'opération d'\emph{éclatement symplectique} (de poids $R$) de $M$ en $z$ consiste alors à retirer de $(M, \omega)$ l'image d'un plongement symplectique $(B^4(R), \omega_{st}) \hookrightarrow (M, \omega)$ centré en $z$ et à coller une copie de $(\beta^{-1} (B^4(r)), \omega_R)$ à la place. Cette opération dépend d'un certain nombre de choix (dont $r$ et $R$ font partie), mais elle est en réalité bien définie à déformation symplectique près.

\begin{df}
On appelle \emph{diviseur exceptionnel symplectique} toute courbe rationnelle $E$ symplectiquement plongée dans une surface symplectique telle que $[E]^2 =-1$.
\end{df}

Soit $( \tilde M, \tilde \omega)$ une surface symplectique et $E$ un diviseur exceptionnel symplectique dans $( \tilde M, \tilde \omega)$. Le théorème du voisinage symplectique nous assure l'existence d'un réel $R>0$ tel qu'un voisinage de $E$ dans $( \tilde M, \tilde \omega)$ est symplectomorphe à un voisinage de la section nulle dans $(\tilde{\mathbb{C}}^2, \omega_R)$. L'opération de \emph{contraction symplectique} consiste alors à remplacer le voisinage de $E$ dans $( \tilde M, \tilde \omega)$ par $B^4(R)$. Cette opération, bien définie à déformation symplectique près, ne dépend que de la classe d'isotopie symplectique du diviseur exceptionnel $E$ dans $( \tilde M, \tilde \omega)$.

\begin{rk}
Soit $(M, \omega)$ une surface symplectique munie d'une structure presque complexe $J \in \mathcal{J}_\tau (M, \omega)$ non intégrable au voisinage d'un point $z \in M$. On ne sait pas définir de manière canonique une structure presque complexe $\tilde J \in \mathcal{J}_\tau ( \tilde M, \tilde \omega)$ sur un éclaté symplectique $( \tilde M, \tilde \omega)$ de $(M, \omega)$ en $z$ de façon à ce qu'elle coïncide avec $J$ en dehors d'un voisinage du diviseur exceptionnel.

De même, soit $(\tilde M, \tilde \omega)$ une surface symplectique munie d'une structure presque complexe $\tilde J \in \mathcal{J}_\tau (\tilde M, \tilde \omega)$ non intégrable au voisinage d'un diviseur exceptionnel symplectique (ou même $\tilde J$--holomorphe) $E \subset \tilde M$. On ne sait pas non plus définir de manière canonique une structure presque complexe $J \in \mathcal{J}_\tau (M, \omega)$ sur un contracté $(M, \omega)$ de $( \tilde M, \tilde \omega)$ selon $E$ de façon à ce qu'elle coïncide avec $\tilde J$ en dehors d'un voisinage de l'image du diviseur exceptionnel par la contraction. 
\end{rk}

Les effets des éclatements symplectiques sur la topologie de $M$ sont bien évidemment les mêmes que dans le cas complexe, la Remarque~\ref{r:blowuptop} est donc également valable pour cette partie.

\begin{df}
On dit qu'une surface symplectique $(M, \omega)$ est \emph{minimale} si elle ne contient pas de diviseur exceptionnel symplectique.

On dit aussi que $(M, \omega)$ est \emph{relativement minimale} par rapport à une courbe symplectiquement plongée $S \subset (M, \omega)$ si $M \backslash S$ ne contient pas de diviseur exceptionnel symplectique.
\end{df}

Une surface symplectique $(M, \omega)$ obtenue à partir d'une autre surface symplectique $(\tilde M, \tilde \omega)$ en contractant une famille maximale $\{ E_1, \dots, E_\ell \}$ de diviseurs exceptionnels deux à deux disjoints est nécessairement minimale. En effet, sinon on peut trouver un diviseur exceptionnel symplectique $E$ dans $(M, \omega)$. En perturbant de manière symplectique $E$ en dehors des images respectives $p_1, \dots, p_\ell$ des diviseurs exceptionnels $E_1, \dots, E_\ell$ par les contractions, on contredit la maximalité de la famille $\{ E_1, \dots, E_\ell \}$.

Le paragraphe précédent illustre bien le caractère plus flexible des diviseurs exceptionnels symplectiques comparées à leurs analogues complexes. Cette flexibilité peut induire une certaine ambigüité lors des opérations de contractions, surtout quand on s'intéresse à l'effet de ces contractions sur des configurations de courbes symplectiques singulières explicites dans une surface symplectique donnée. De même, les choix des différents paramètres pour éclater un point sur une telle configuration de courbes, ne permettent pas d'observer précisément les effets de telles opérations.

Afin de contourner ces difficultés en ce qui concerne l'opération d'éclatement, une astuce consiste à choisir au préalable des coordonnées complexes intégrables compatibles avec la forme symplectique au voisinage du point à éclater. D'après le Théorème \ref{t:McDuff92}, on peut toujours faire ce choix de sorte que les courbes de la configuration sont complexes dans les voisinages considérés. On effectue alors l'opération d'éclatement de manière complexe, et l'opération ainsi effectuée peut-être interprétée comme un éclatement symplectique (voir~\cite[Theorem~3.13]{Wendl}).

Pour l'opération de contraction, on veut trouver une structure complexe \emph{intégrable} compatible avec la forme symplectique au voisinage du diviseur exceptionnel $E$ qu'on veut contracter, qui rend holomorphe la portion de la configuration de courbes dans ce voisinage. On procède en plusieurs étapes, qu'on décrit brièvement. On commence par éclater les points singuliers de la configuration sur $E$ jusqu'à obtenir la résolution à croisements normaux. Puis on utilise le théorème du voisinage symplectique pour les  plombages de courbes symplectiques (Théorème~\ref{t:sympnghdconf}) pour identifier un voisinage de la configuration de courbes à croisements normaux avec la même configuration dans un éclaté de $\mathbb{C}P^2$, via un symplectomorphisme et une déformation symplectique. On utilise ensuite ce symplectomorphisme pour tirer en arrière la structure complexe, et ainsi effectuer de manière complexe les contractions qui inversent la résolution des points singuliers de la configuration qui sont sur $E$. On obtient ainsi la structure complexe voulue au voisinage de $E$, et on peut contracter $E$ de manière complexe. Cette opération peut également être interprétée comme une contraction symplectique (voir~\cite[Theorem~3.14]{Wendl}).

Pour clore cette section, on présente un résultat d'isotopie symplectique pour les diviseurs exceptionnels symplectiques. 

\begin{thm}\label{t:divexisotopy2}
Soit $M$ une variété de dimension $4$, $(\omega_s)_{s \in [0,1]}$ un chemin lisse de formes symplectiques sur $M$ et $J \in \mathcal{J}_\tau (\omega_0)$. Alors pour une famille générique $(J_s)_{s \in [0,1]} \in \mathcal {J}_\tau (M, (\omega_s)_{s \in [0,1]})$ avec $J_0 = J$, toute courbe rationnelle $J_0$--holomorphe plongée $u_0$ dont l'image est un diviseur exceptionnel s'étend en une famille lisse $(u_s)_{s \in [0,1]}$ de diviseurs exceptionnels $J_s$--holomorphes.
\end{thm}

En particulier, la notion de minimalité ne dépend que de la classe de déformation symplectique de la variété ambiante. En effet, considérons deux structures symplectiques $\omega_0$ et $\omega_1$ sur une variété $M$ de dimension $4$. Si $\omega_0$ et $\omega_1$ sont reliées par un chemin de formes symplectiques, le Théorème \ref{t:divexisotopy2} nous assure que $(M, \omega_0)$ est minimale si et seulement si $(M, \omega_1)$ est minimale. Ce Théorème fait l'objet du chapitre 5 dans~\cite{Wendl}.

\begin{proof}[Esquisse de preuve du Théorème~\ref{t:divexisotopy2}]
Quitte à considérer seulement une composante connexe de l'espace de modules $\mathcal{M}_{0}([u_0];(J_s)_{s \in [0,1]})$, on peut supposer que l'ensemble $\mathcal{M}_{0}([u_0];(J_s)_{s \in [0,1]})$ est connexe. On commence par calculer l'indice de la courbe $u_0$ :
$$\ind (u_0) = \chi (\mathcal{S}^2) + 2 [u_0]^2 = 2 -2 =0.$$
On peut montrer que l'espace de modules à paramètre $\mathcal{M}_{0}([u_0];(J_s)_{s \in [0,1]})$ (qui est non vide puisqu'il contient $u_0$) est compact en utilisant la généricité du chemin $(J_s)_{s \in [0,1]}$, ainsi que le Corollaire~\ref{c:indgeq-1} pour trouver des contraintes sur les indices d'éventuelles courbes pouvant apparaître comme composantes de courbes nodales (voir~\cite[Chapter 4, Chapter 5]{Wendl}). Puisque $\ind (u_0) = 0$, le théorème de transversalité automatique nous assure que toutes les courbes de l'espace de modules $\mathcal{M}_{0}([u_0];(J_s)_{s \in [0,1]})$ sont automatiquement transverses (elle survivent à une perturbation arbitrairement petite de la structure presque complexe). Ainsi, d'après le Théorème \ref{t:moduliparamreg}, l'application de projection $ \pi : \mathcal{M}_{0}([u_0];(J_s)_{s \in [0,1]}) \rightarrow [0,1]$ est submersion. 
En regroupant les informations ci-dessus, on obtient que l'image de $\pi$ dans l'intervalle $[0,1]$ est non vide, ouverte et fermée. La connexité de l'ensemble $[0,1]$ permet alors de conclure.
\end{proof}

\begin{rk}
Puisque $[u_0] \cdot [u_0] = -1$, la propriété de positivité d'intersection nous indique que pour tout $s \in [0,1]$, il existe au plus un diviseur exceptionnel $J_s$--holomorphe homologue à $[u_0]$. L'application de projection $\pi$ dans la preuve du Théorème \ref{t:divexisotopy2} est donc en réalité un difféomorphisme.
\end{rk}

\begin{cor}\label{c:divexisotopy}
Soit $(M, \omega)$ une surface symplectique et $J \in \mathcal{J}_\tau (M, \omega)$ une structure presque complexe générique. Alors tout diviseur exceptionnel est symplectiquement isotope à l'image d'une unique courbe rationnelle $J$--holomorphe simple (qui est plongée).
\end{cor}

Pour finir, la Remarque~\ref{r:birationalequiv} motive la définition suivante.

\begin{df}
On dit que deux surfaces symplectiques $(M_1, \omega_1)$ et $(M_2, \omega_2)$ sont \emph{birationnellement équivalentes} si on peut obtenir $(M_2, \omega_2)$ à partir de $(M_1, \omega_1)$ par des successions d'éclatements, de contractions et de déformations symplectiques.
\end{df}

L'équivalence birationnelle définit une relation d'équivalence sur l'ensemble des surfaces symplectiques.

\begin{rk}
On peut définir de la même manière que pour le cas complexe, les notions de transformée propre et de transformée totale pour les courbes symplectiques singulières. On attire cependant l'attention sur le fait qu'il n'y a pas unicité de la transformée propre dans le cadre symplectique : sa définition dépend d'un choix de coordonnées complexes au voisinage du point éclaté. En revanche, la classe d'isotopie symplectique de la transformée propre (équisingulière dans le cas où la transformée propre possède des singularités, voir le Chapitre~5 pour les définitions) ne dépend d'aucun choix (voir~\cite[Section~3.1]{GS}). On rappelle que par définition, les singularités des courbes symplectiques singulières dans les surfaces symplectiques sont localement modelées sur celles des courbes complexes dans les surfaces complexes. Ainsi, la discussion sur la résolution des singularités à la fin de la Sous-section~\ref{ss:eclatcomp} s'étend de manière similaire aux courbes symplectiques singulières.
\end{rk}

\section{Pinceaux et fibrations de Lefschetz symplectiques}\label{s:pinceauLefschetz}

Un moyen facile de construire des exemples de variétés de dimension $4$ consiste à considérer les produits cartésiens de deux surfaces. On peut généraliser cette construction en s'intéressant aux fibrés en surfaces sur des surfaces. Les fibrations et pinceaux de Lefschetz, initialement introduits pour étudier la topologie des surfaces complexes, constituent une étape supplémentaire dans cette direction et permettent d'englober une quantité bien plus vaste d'exemples. Le but de cette section est de définir ces objets, d'abord dans le cadre des variétés lisses de dimension $4$, puis d'illustrer leur importance dans l'étude des surfaces symplectiques. 

\begin{df}
Soit $M$ une variété de dimension $4$ et $\Sigma$ une surface. Une \emph{fibration de Lefschetz} sur $M$ au-dessus de $\Sigma$ est une application lisse $\pi : M \rightarrow \Sigma$ avec un nombre fini de points critiques, telle qu'au voisinage de tout point critique $p$, il existe des coordonnées complexes $(z_1, z_2)$ au voisinage de $p$ et des coordonnées complexes au voisinage de $\pi (p)$, respectant les orientations respectives de $M$ et $\Sigma$, dans lesquelles l'application $\pi$ est localement de la forme 
$$\pi(z_1, z_2) = z_1 z_2.$$
\end{df}

\begin{rk}
La condition locale au voisinage d'un point critique peut également être écrite sous la forme $\pi (z_1, z_1) = z_1^2 + z_2^2$. Elle traduit le fait que les points critiques de $\pi$ sont non dégénérés. Cette condition peut alors être vue comme un analogue complexe des points critiques de Morse.
\end{rk}

Soit $z$ un élément de $\Sigma$, l'ensemble $\pi^{-1}(z)$ est appelé la \emph{fibre} au-dessus de $z$. Si $z$ est une valeur régulière de $\pi$, l'ensemble $\pi^{-1}(z)$ est une union de surfaces plongées dans $M$. On dit dans ce cas que c'est une \emph{fibre régulière} (notons que les fibres régulières sont deux à deux difféomorphes). Sinon l'ensemble $\pi^{-1} ( \{z\})$ est une union de surfaces immergées dans $M$, dont les points d'auto-intersection sont des points doubles transverses positifs, qui s'intersectent deux à deux transversalement de manière positive. Dans ce cas, on dit que $\pi^{-1}(\{z\})$ est une \emph{fibre singulière}. Une sous-variété $S$ de dimension $2$ qui intersecte exactement une fois, en des points réguliers et de manière transverse, chacune des fibres de la fibration de Lefschetz est appelée une \emph{section}. L'application de projection restreinte à $S$ définit alors un difféomorphisme entre $S$ et $\Sigma$.

\begin{rk}
Si $\pi : M \rightarrow \Sigma$ est une fibration de Lefschetz dont les fibres ne sont pas connexes, on dispose d'une surface $\Sigma'$ (connexe), d'un revêtement ramifié $\varphi : \Sigma' \rightarrow \Sigma$ de degré fini supérieur ou égal à $2$ et d'une fibration de Lefschetz $\pi' : M \rightarrow \Sigma'$ avec fibres connexes tels que $\pi = \pi' \circ \varphi$. Par conséquent, on supposera implicitement que toutes les fibres des fibrations de Lefschetz considérées par la suite sont connexes.
\end{rk}

\begin{df}
Soit $M$ une variété de dimension $4$ et $B$ un ensemble fini de points de $M$. Un \emph{pinceau de Lefschetz} sur $M$ de \emph{base} $B$ est la donnée d'une fibration de Lefschetz 
$$\pi : M \backslash B \rightarrow \mathbb{C}P^1,$$
telle que pour tout \emph{point base} $p \in B$, il existe des coordonnées complexes $(z_1,z_1)$ au voisinage de $p$ respectant l'orientation de $M$ dans lesquelles l'application $\pi$ est localement de la forme
$$ \pi (z_1, z_2) = [ z_1 : z_2].$$
\end{df}

Soit $z$ un élément de $\mathbb{C}P^1$, l'adhérence de l'ensemble $\pi^{-1}( \{z \})$ est appelé la \emph{fibre} au-dessus de $z$. Chaque fibre est une union de surfaces plongées dans $M$ et contient l'ensemble des points bases. On définit de manière analogue les notions de fibres régulières et fibres singulières. Le modèle local autour de chacun des points bases nous assure qu'en chaque point base, les fibres s'intersectent deux à deux de manière transverse et positive. On remarque qu'une fibre d'un pinceau de Lefschetz dont la base est non vide ne représente jamais une classe de torsion dans $H_2 (M ; \mathbb{Z})$. En effet, une telle fibre $F$ doit satisfaire l'égalité $[F] \cdot [F] = \Card B$. 

\begin{ex}
\leavevmode
\begin{enumerate}
\item L'application
$$\begin{array}{ccccc}
\pi &:& \mathbb{C}P^2  \backslash \{ [0:0:1]\} & \longrightarrow & \mathbb{C}P^1\\
& &{[x:y:z]} & \longmapsto & [x:y]
\end{array}$$ définit un pinceau de droites complexes sur $\mathbb{C}P^2$. Pour tout $[x:y] \in \mathbb{C}P^1$, l'ensemble $\overline{\pi^{-1}([x:y])}$ est la droite projective complexe de $\mathbb{C}P^2$ donnée par l'équation $yX=xY$ d'inconnues $X,Y,Z$.
\item Plus généralement, on considère deux polynômes homogènes $P, Q \in \mathbb{C}[X,Y,Z]$ distincts et de même degré $d \in \mathbb{N}^*$. Quitte à perturber génériquement $P$ et $Q$, on peut supposer que les courbes complexes projectives $C_P$ et $C_Q$ d'équations respectives $P(X,Y,Z) =0$ et $Q(X,Y,Z) =0$ sont non singulières et s'intersectent transversalement. Le théorème de Bézout nous assure que les courbes $C_P$ et $C_Q$ s'intersectent en $d^2$ points. Notons $I = C_P \cap C_Q$. L'application 
$$\begin{array}{ccccc}
\pi &:& \mathbb{C}P^2  \backslash I & \longrightarrow & \mathbb{C}P^1\\
 & & {[x:y:z]} & \longmapsto & [P(x,y,z):Q(x,y,z)]
\end{array}$$
définit un pinceau de Lefschetz avec $d^2$ points bases dont les fibres non singulières sont de genre $\frac{(d-1)(d-2)}{2}$. L'Exemple~\ref{e:cubiquespar9pts} décrit le cas $d=3$.
\begin{figure}[htbp]
	\centering
	\subfloat[Pinceau de droites passant par un point sur $\mathbb{C}P^2$.]{\label{f:pencilbasepoint}\includegraphics[height=0.30\linewidth]{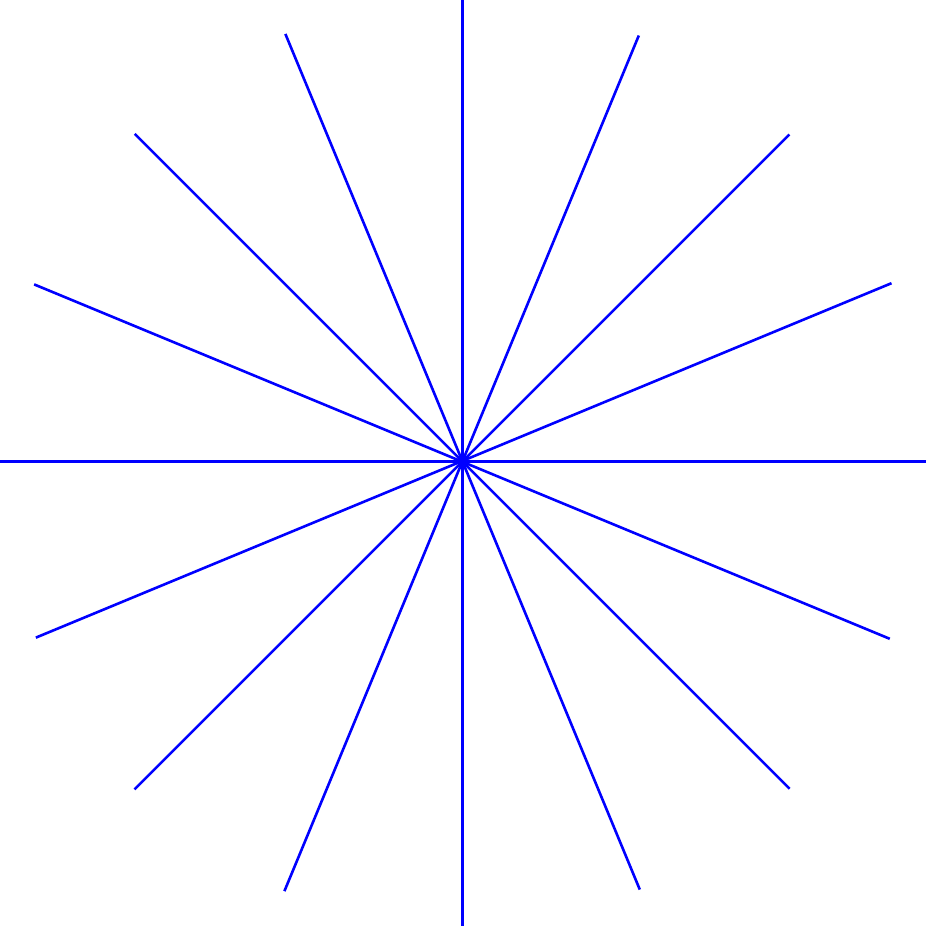}}\qquad
	\subfloat[Pinceau de coniques passant par quatre points génériques sur $\mathbb{C}P^2$.]{\label{f:pencilofconics}\includegraphics[height=0.30\linewidth]{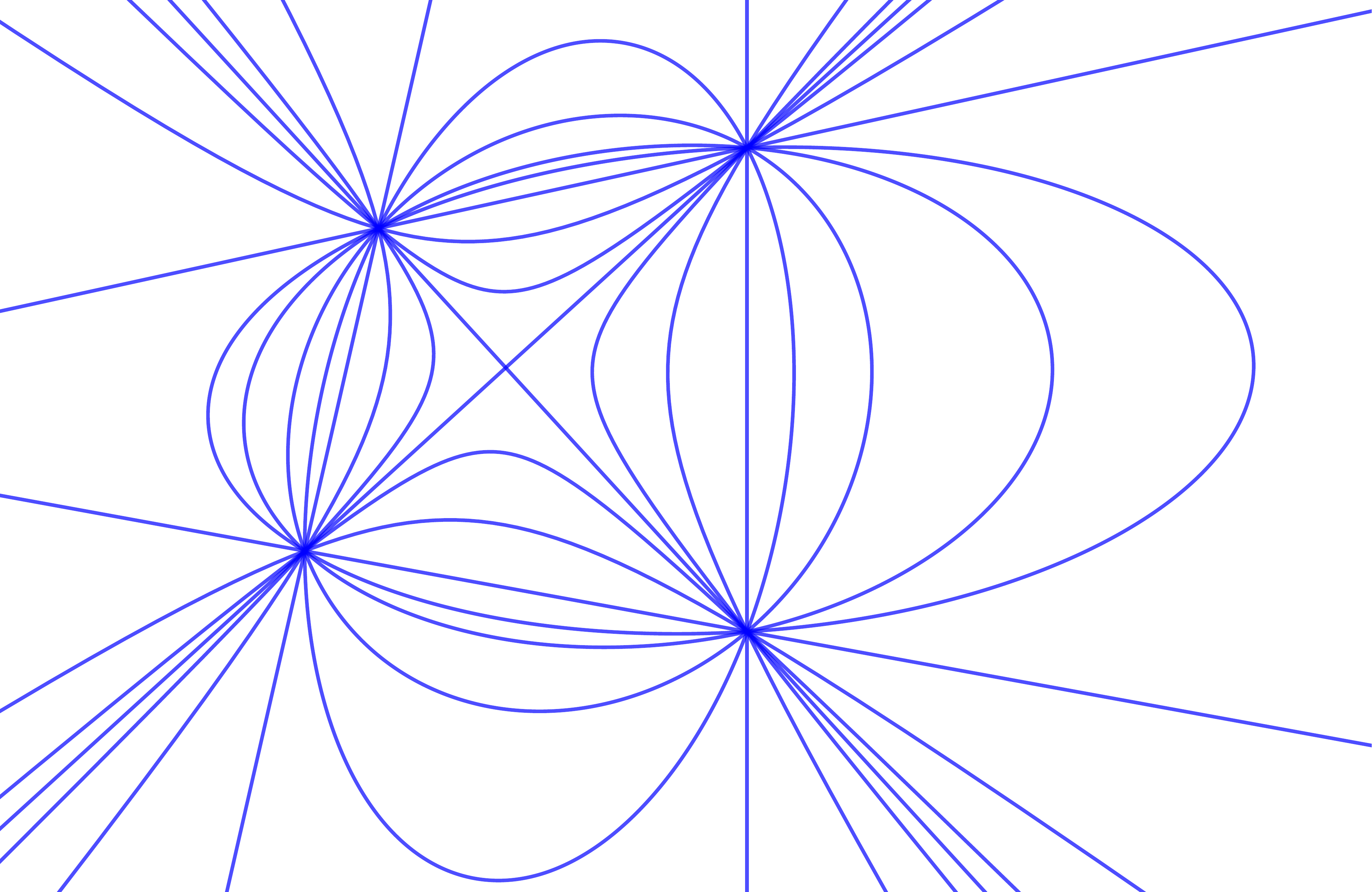}}
\caption{Exemples de pinceaux de Lefschetz.}\label{f:pencilsplane}
\end{figure}
\end{enumerate}
\end{ex}

\begin{rk}
\leavevmode
\begin{enumerate}
\item Lorsqu'on éclate tous les points bases d'un pinceau de Lefschetz, on obtient une fibration de Lefschetz au-dessus de $\mathbb{C}P^1$. Les diviseurs exceptionnels qui apparaissent suite à ces éclatements sont alors des sections de la fibration de Lefschetz ainsi obtenue.
\begin{figure}[h]
	\centering
	\includegraphics[scale=0.5]{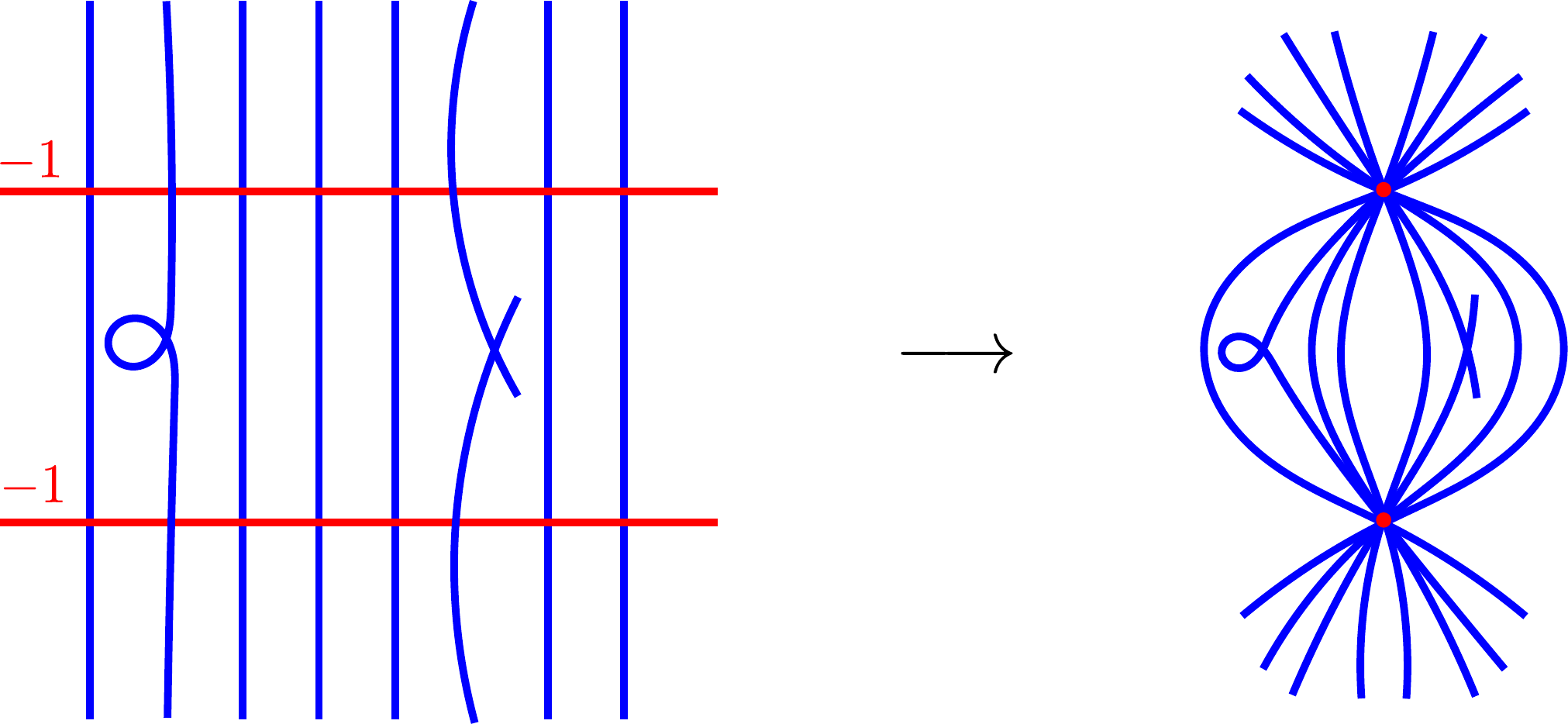}\\
	\caption{\'Eclatement des points bases d'un pinceau de Lefschetz.}
	\label{f:blownuppencil}
\end{figure}
\item Lorsqu'on éclate un point $p$ d'un pinceau (resp. fibration) de Lefschetz $\pi : M \backslash B \rightarrow \Sigma$ qui n'est pas un point base (à l'aide de coordonnées locales complexes qui rendent les fibres complexes), on obtient un pinceau (resp. fibration) $\pi : M \# \overline{\mathbb{C}P}^2 \backslash B \rightarrow \Sigma$, où la transformée totale de la fibre $F$ passant par $p$ est une fibre singulière, constituée de l'union de la transformée propre $\tilde F$ de $F$ et du diviseur exceptionnel $E$ provenant de l'éclatement en $p$ (voir la Figure~\ref{f:blownupfibration}).
\begin{figure}[h]
	\centering
	\includegraphics[scale=0.5]{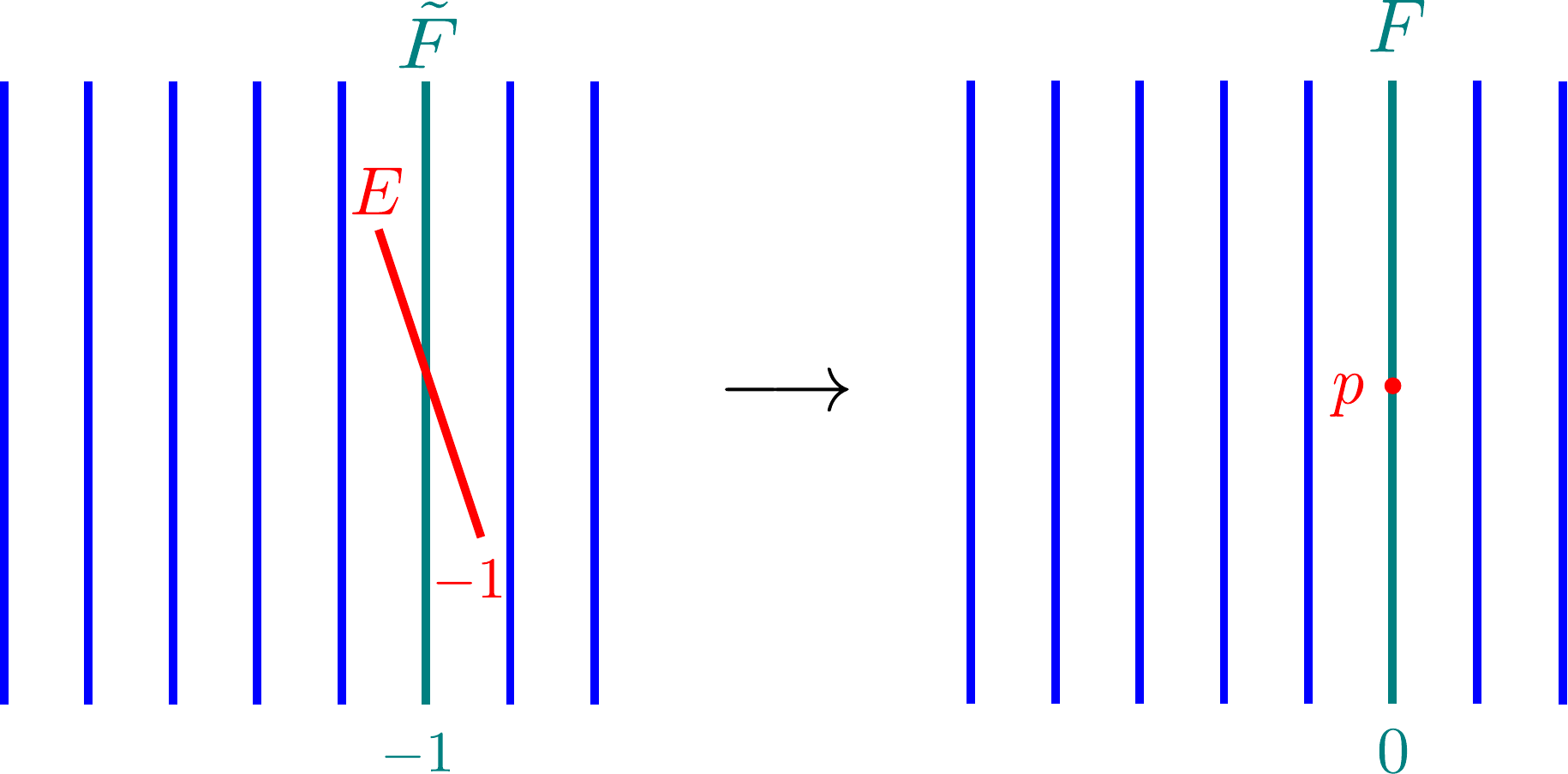}\\
	\caption{\'Eclatement d'un point régulier sur une fibre d'un fibré en surfaces au-dessus d'une surface.}
	\label{f:blownupfibration}
\end{figure}
\end{enumerate}
\end{rk}

\begin{rk}
\'Etant donné un pinceau (resp. fibration) de Lefschetz $\pi : M \backslash B \rightarrow \Sigma$, on peut toujours réaliser une perturbation $\mathcal{C}^\infty$ arbitrairement petite de $\pi$ de façon à ce que chaque fibre singulière possède exactement un point singulier. On parle dans ce cas de \emph{pinceau} (resp. \emph{fibration}) \emph{de Lefschetz générique}.
\end{rk}

On présente deux propositions concernant la topologie des pinceaux de Lefschetz qui nous seront utiles par la suite.

\begin{prop} \label{p:pinceauCP2}
Si une variété $M$ de dimension $4$ est munie d'un pinceau de Lefschetz possédant exactement un point base et dont toutes les fibres sont régulières, alors le pinceau sur $M$ est isomorphe à un pinceau de droites sur $\mathbb{C}P^2$ (en particulier $M$ est difféomorphe à $\mathbb{C}P^2$).
\end{prop}

\begin{proof}[Esquisse de preuve]
On note $\pi :  M \backslash \{p \} \rightarrow \mathbb{C}P^1$ le pinceau. On choisit une boule $B \subset M$ centrée autour du point base $p$ de $\pi :  M \backslash \{p \} \rightarrow \mathbb{C}P^1$. L'application $\pi$ restreinte à $\partial B \simeq \mathcal{S}^3$ munit $\partial B$ de la fibration de Hopf. La restriction de $\pi$ à $M \backslash B$ est donc isomorphe à un fibré en disques orienté dont le bord est la fibration de Hopf munie de l'orientation opposée. En particulier, la première classe de Chern du fibré en disques $ \pi : M \backslash B \rightarrow \mathbb{C}P^1$ est complètement déterminée. Comme les fibrés en disques orientés au-dessus d'une même base sont classifiés à isomorphisme près par leur première classe de Chern, on obtient un isomorphisme entre $\pi : M \backslash B \rightarrow \mathbb{C}P^1$ et le fibré en disques défini par la restriction d'un pinceau de droites sur $\mathbb{C}P^2$ au complémentaire d'une boule centrée autour de son point base. Cet isomorphisme peut ensuite s'étendre de manière lisse sur la boule $B$ en un isomorphisme de pinceau de $M$ vers $\mathbb{C}P^2$ grâce aux modèles locaux au voisinage des points bases.
\end{proof}

\begin{prop} \label{p:pinceau2ptbases}
Tout pinceau de Lefschetz avec au moins deux points bases admet au moins une fibre singulière.
\end{prop}

\begin{proof}[Esquisse de preuve]
On suppose qu'un pinceau de Lefschetz  $\pi : M \backslash B \rightarrow \mathbb{C}P^1$ admet au moins deux points bases et aucune fibre singulière. On choisit un point $p \in B$ et on éclate tous les points bases différents de $p$ de façon à obtenir un pinceau de Lefschetz $\tilde \pi : \tilde M \backslash \{p \} \rightarrow \mathbb{C}P^1$ avec un unique point base et aucune fibre singulière. D'après la Proposition~\ref{p:pinceauCP2}, $\tilde M$ est difféomorphe à $\mathbb{C}P^2$. Mais comme on a effectué au moins un éclatement, $\mathbb{C}P^2$ contient une sphère plongée d'auto-intersection $-1$. Puisque la forme d'intersection de $\mathbb{C}P^2$ est définie positive, on aboutit à une contradiction.
\end{proof}

On explique désormais le lien entre les surfaces symplectiques et les pinceaux et fibrations de Lefschetz.

\begin{df}
Soit $(M, \omega)$ une surface symplectique, et $\pi : M \backslash B \rightarrow \Sigma$ un pinceau (resp. fibration) de Lefschetz. On dit que $\pi : M \backslash B \rightarrow \Sigma$ est un \emph{pinceau} (resp. \emph{fibration}) \emph{de Lefschetz symplectique} si les conditions suivantes sont remplies :
\begin{enumerate}
\item la partie lisse de chaque fibre (c'est-à-dire le sous-ensemble d'une fibre constitué des points qui ne sont ni des points bases ni des points critiques) est une sous-variété symplectique,
\item pour toute structure presque complexe $J$ définie au voisinage d'un point base ou d'un point critique $p$ qui se restreint à une structure complexe sur les parties lisses de chacune des fibres, $J$ est dominée par $\omega$ en $p$.
\end{enumerate}
On dit également que $\omega$ est une \emph{forme symplectique compatible} avec $\pi$.
\end{df}

Le théorème suivant, dû à Gompf et Thurston, montre qu'une fois les exceptions évidentes mises de côté, les fibrations et pinceaux de Lefschetz peuvent être munis de formes symplectiques compatibles avec la projection.

\begin{thm}[\cite{thurston96}, \cite{gompf2002}]\label{t:GompfThurston} 
Soit $\pi : M \backslash B \rightarrow \Sigma$ une fibration ou pinceau de Lefschetz dont la fibre ne représente pas une classe de torsion dans $H_2 (M ; \mathbb{Z})$. Alors $M$ admet une structure symplectique compatible avec $\pi$. De plus deux telles structures symplectiques peuvent être reliées par un chemin lisse de structures symplectiques compatibles avec $\pi$.
\end{thm}

\begin{prop}[{\cite[Proposition~3.34]{Wendl}}]
Si la fibre d'une fibration ou pinceau de Lefschetz $\pi : M \backslash B \rightarrow \Sigma$ représente une classe de torsion dans $H_2(M ; \mathbb{Z})$, alors $B = \emptyset$ et la fibre est un tore.
\end{prop}

\begin{ex}
Le produit de la fibration de Hopf avec le fibré en cercles trivial au-dessus de $\mathcal{S}^2$ munit la variété $\mathcal{S}^3 \times \mathcal{S}^1$ d'une structure de fibré en tores au-dessus d'une sphère. Mais $\mathcal{S}^3 \times \mathcal{S}^1$ n'admet pas de structure symplectique car $H_2(\mathcal{S}^3 \times \mathcal{S}^1 ; \mathbb{Z}) =0$.
\end{ex}

De manière plus surprenante, à l'aide de techniques approximativement pseudoholomorphes, Donaldson a montré que toutes les surfaces symplectiques peuvent être munies de pinceaux de Lefschetz symplectiques (mentionnons que son résultat est en fait plus général et concerne aussi les variétés symplectiques de dimensions plus grandes). 

\begin{thm}[\cite{Donaldson96}]
Soit $(M, \omega)$ une surface symplectique telle que $[\omega]$ appartient à $H^2 (M; \mathbb{Z})$. Pour un entier $k$ suffisamment grand, on peut munir $(M, \omega)$ d'un pinceau de Lefschetz symplectique dont les fibres sont homologues au dual de Poincaré de $k[\omega]$.
\end{thm}

\begin{rk}
L'hypothèse sur la forme symplectique n'est pas restrictive. En effet, quitte à réaliser une perturbation arbitrairement petite, on peut toujours supposer que la classe de cohomologie d'une forme symplectique est à coefficients rationnels. On peut alors multiplier la forme symplectique par un entier suffisamment grand pour que sa classe d'homologie soit à coefficients entiers.
\end{rk}

Ce théorème fondamental permet notamment de ramener l'étude des surfaces symplectiques à des données combinatoires (nombres de points bases du pinceau et monodromie autour des fibres singulières), créant ainsi un pont entre la topologie symplectique de basse dimension et la théorie des groupes.

\clearemptydoublepage

\chapter{Courbes symplectiques rationnelles de haute auto-intersection dans les surfaces symplectiques}
On rappelle que dans tout ce qui suit, on appelle \emph{surface symplectique} toute variété symplectique compacte, connexe, sans bord de dimension $4$. De même, le terme \emph{courbe symplectique} désigne toute variété symplectique compacte, connexe, sans bord de dimension $2$. On s'intéresse dans ce chapitre à certains types de surfaces symplectiques : les \emph{surfaces symplectiquement réglées} et les \emph{surfaces symplectiques rationnelles}. On verra qu'elles sont caractérisées par la présence de certaines courbes rationnelles satisfaisant certaines conditions sur leur auto-intersection (ou alternativement sur leur première classe de Chern).

On commence par introduire et manipuler dans la Section \ref{s:ssr,ssr} des exemples de surfaces symplectiques rationnelles et de surfaces symplectiquement réglées. Puis dans la Section \ref{s:GromovMcDuff}, on présente un théorème de Gromov et McDuff sur les surfaces contenant des courbes rationnelles symplectiquement plongées d'auto-intersection positive. On décrit les étapes principales de la démonstration tout en soulignant les points qui relèvent spécifiquement de l'utilisation des courbes symplectiques rationnelles. Enfin, on présente dans la Section \ref{s:spheresimmergees} une amélioration de ce théorème due à McDuff qui étend le résultat aux courbes symplectiques singulières rationnelles avec première classe de Chern plus grande que $2$ dans les surfaces symplectiques. C'est cette amélioration qui permet notamment d'obtenir la caractérisation des surfaces symplectiquement réglées et celle des surfaces symplectiques rationnelles.

\section{Surfaces symplectiques rationnelles, surfaces symplectiquement réglées} \label{s:ssr,ssr}

Définissons sans plus tarder les objets centraux de ce chapitre.

\begin{df}\label{d:ruled}
On appelle \emph{surface symplectiquement réglée} toute surface symplectique birationnellement équivalente à l'espace total d'un fibré en sphères symplectique. De même, on désignera par \emph{surface rationnelle symplectique} toute surface symplectique birationnellement équivalente à $(\mathbb{C}P^2,\omega_{FS})$.
\end{df}

\begin{rk}
On met en garde la lectrice ou le lecteur sur le fait que la terminologie employée ici diffère de celle utilisée dans~\cite{Wendl}. En effet dans~\cite{Wendl}, le terme \og symplectically ruled surface \fg{} désigne seulement les fibrés en sphères symplectiques sur des surfaces et le terme \og blown-up symplectically ruled surface \fg{} est utilisé pour rendre compte des surfaces symplectiques birationnellement équivalentes à l'espace total d'un fibré en sphères symplectique qui ne sont pas des déformations symplectiques de $(\mathbb{C}P^2,\omega_{FS})$ (voir le Théorème~\ref{t:classificationsurfacesreglees} pour la justification derrière ce terme). La Définition~\ref{d:ruled} est plus générale et se rapproche de celle utilisée dans le cas complexe par~\cite{beauville1996complex} et~\cite{hartshorne1977algebraic} par exemple.
\end{rk}

Puisque $(\mathbb{C}P^2,\omega_{FS})$ éclaté en un point est un fibré en sphères symplectique, toute surface rationnelle symplectique est une surface symplectiquement réglée.

Les surfaces symplectiques munies de structures de fibrés en sphères symplectiques sont les premiers exemples de surfaces symplectiquement réglées. On commence par étudier les types de difféomorphismes possibles pour de telles surfaces symplectiques en rappelant une manière de construire les fibrés en surfaces sur des surfaces. Cette construction repose sur la compréhension de $\mathrm{Diff}_+ \left(\mathcal{S}^2 \right)$, le groupe des difféomorphismes préservant l'orientation de la fibre $\mathcal{S}^2$. Plus précisément, ce qui nous intéresse ici est le fait que $\mathrm{Diff}_+ \left(\mathcal{S}^2 \right)$ possède le même type d'homotopie que $SO(3)$. En particulier, $\mathrm{Diff}_+ \left(\mathcal{S}^2 \right)$ est connexe et son groupe fondamental est isomorphe à $\mathbb{Z} / 2 \mathbb{Z}$.

On rappelle tout d'abord qu'un fibré dont la base $B$ est contractile est globalement trivial. Lorsque la base $\Sigma$ d'un fibré en sphères $\pi : E \rightarrow \Sigma$ est une surface, on peut décomposer l'espace total $E$ du fibré en l'union de l'espace total $E_1$ du fibré en sphères trivial au-dessus du disque $\mathbb{D}^2$ et d'un fibré $E_2$ au dessus de $\Sigma \backslash \mathbb{D}^2$. Comme $\Sigma \backslash \mathbb{D}^2$ se rétracte par déformation sur un bouquet de cercles, la connexité de $\mathrm{Diff}_+ \left(\mathcal{S}^2 \right)$ nous assure que le fibré $E_2$ est également trivial. Le type de difféomorphisme de $E$ ne dépend alors que la manière dont on recolle les deux fibrés $E_1$ et $E_2$ selon leurs bords, c'est-à-dire que de la classe d'homotopie d'un chemin dans $\mathrm{Diff}_+(\mathcal{S}^2)$. Puisque $\mathrm{Diff}_+(\mathcal{S}^2)$ est difféomorphe à $\mathbb{Z} / 2 \mathbb{Z}$, il y a donc deux types de difféomorphisme possibles pour $E$ : celui qui correspond au fibré trivial $\Sigma \times \mathcal{S}^2$ et l'autre qu'on notera $\Sigma \, \tilde \times \, \mathcal{S}^2$. 

\begin{rk}\label{r:paritéréglée}
Grâce au théorème de Künneth, il est aisé de constater que le groupe $H_2(\Sigma \times \mathcal{S}^2 ; \mathbb{Z})$ est engendré par $[\S \times \{ * \}]$ et $[\{ * \} \times \mathcal{S}^2]$. La forme d'intersection de $\Sigma \times \mathcal{S}^2$ est donc paire puisque pour tous entier $k$ et $\ell$, on a 
$$(k[\S \times \{ * \}]+ \ell [\{ * \} \times \mathcal{S}^2])^2 = 2k \ell.$$

Maintenant on considère un modèle complexe de la forme $C \times \mathbb{C}P^1$ de $\Sigma \times \mathcal{S}^2$ et on éclate un point $p \in C \times \mathbb{C}P^1$. La transformée propre de la fibre passant par $p$ est une courbe rationnelle complexe $E$ d'auto-intersection $-1$, la transformée propre de la section complexe de la forme $C \times \{ * \}$ passant par $p$ est également d'auto-intersection $-1$ et les transformées propres des sections complexes de la forme $C \times \{ * \}$ ne passant pas par $p$ sont des courbes d'auto-intersection $0$ qui intersectent chacune $E$ en un point de manière transverse. On contracte ensuite $E$ afin d'obtenir un autre fibré holomorphe en $\mathbb{C}P^1$ au-dessus de $C$, qui contient cette fois-ci une section d'auto-intersection $-1$ et des sections d'auto-intersection $+1$ (voir la Figure~\ref{f:birationalbasic2}). Comme la forme d'intersection de $\Sigma \times \mathcal{S}^2$ est paire, le fibré en sphères ainsi obtenu est nécessairement difféomorphe à $\Sigma \, \tilde \times \, \mathcal{S}^2$. On en déduit que la forme d'intersection de $\Sigma \, \tilde \times \, \mathcal{S}^2$ est impaire.
\begin{figure}[h]
	\centering
	\includegraphics[scale=0.4]{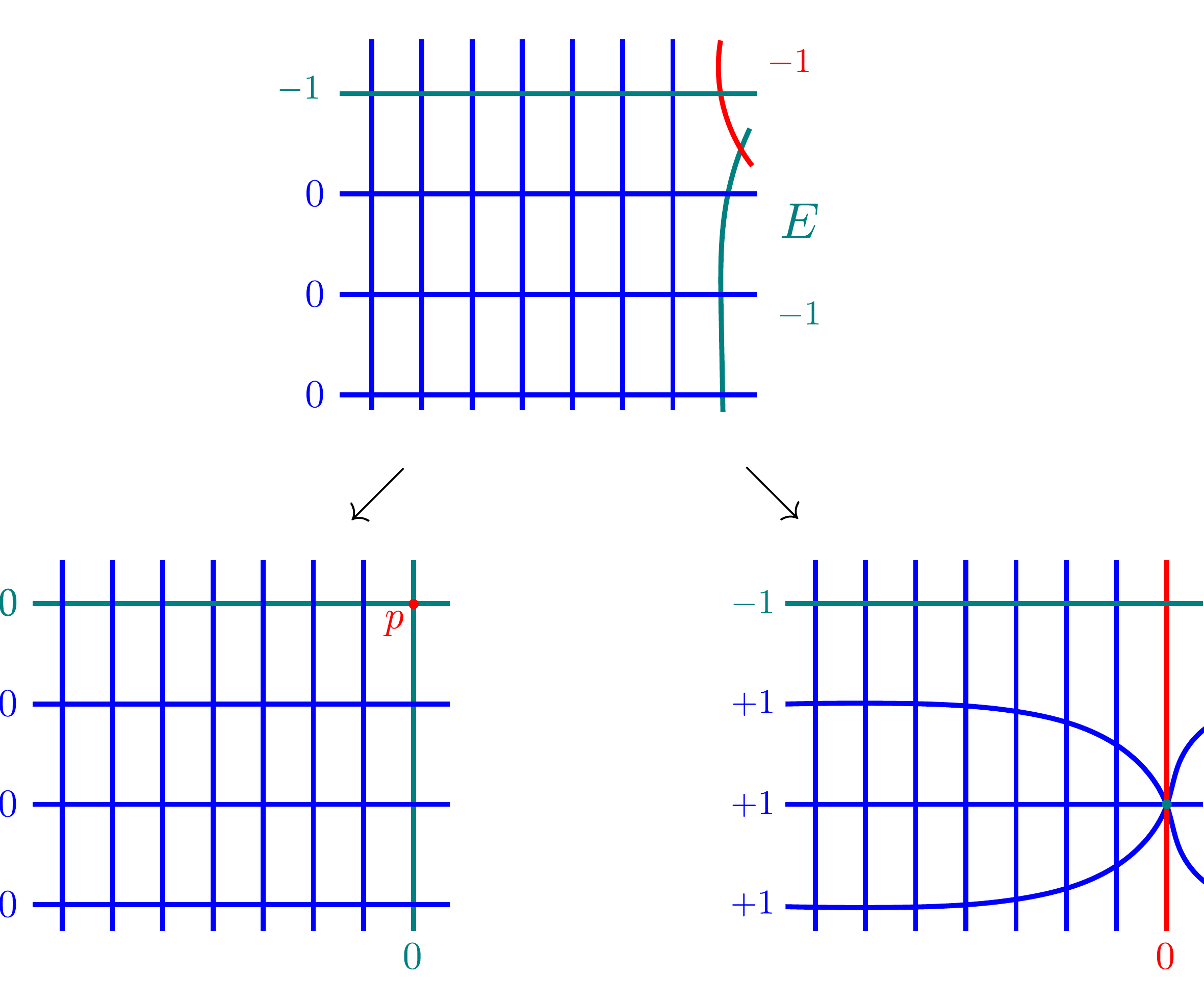}\\
	\caption{Transformation birationnelle entre $C \times \mathbb{C}P^1$ et un autre fibré holomorphe en $\mathbb{C}P^1$ au-dessus de $C$.}
	\label{f:birationalbasic2}
\end{figure}

La parité de la forme d'intersection est donc un bon moyen de déterminer les types de difféomorphisme des fibrés en sphères au-dessus d'une surface donnée. On a également montré que tous les fibrés en sphères au-dessus de surfaces admettent des sections lisses.
\end{rk}

\begin{figure}[h]
	\centering
	\includegraphics[scale=0.4]{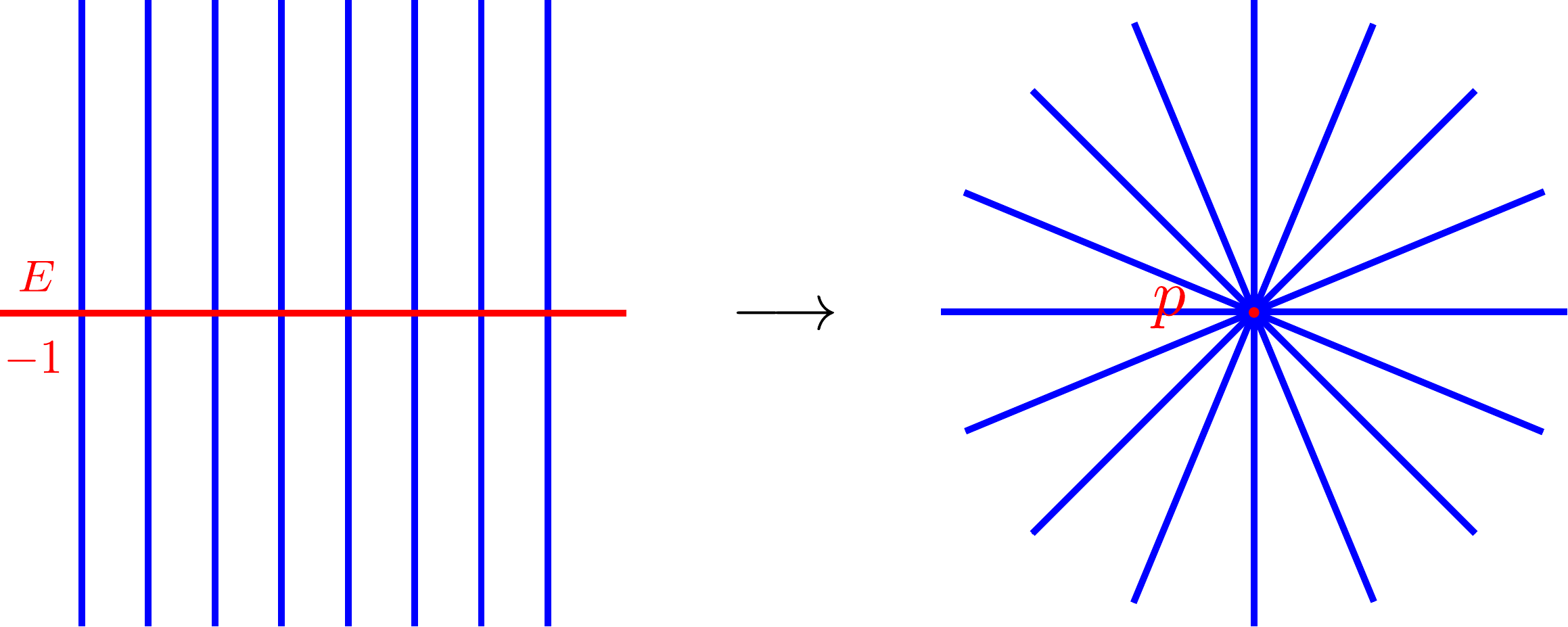}\\
	\caption{\'Eclatement du point base $p$ d'un pinceau de droites sur $\mathbb{C}P^2$.}
	\label{f:blownuppencilbasepoint}
\end{figure}

On considère désormais le point $p = [0:0:1] \in \mathbb{C}P^2$, ainsi que le pinceau des droites complexes passant par ce point. L'ensemble des droites du pinceaux est naturellement isomorphe à $\mathbb{C}P^1$ l'ensemble des droites vectorielles de $\mathbb{C}^2$. Lorsqu'on éclate en $p$, les transformées propres des droites du pinceau sont des courbes symplectiques deux à disjointes, dont l'union est égale à l'espace tout entier, c'est-à-dire à $\mathbb{C}P^2 \# \overline{\mathbb{C}P}^2$. L'application $ \pi : \mathbb{C}P^2 \rightarrow \mathbb{C}P^1$, qui à chaque point $q \in \mathbb{C}P^2 \backslash \{ p \}$ associe l'élément de $\mathbb{C}P^1$ correspondant à la fibre du pinceau passant par $q$, s'étend à l'espace $\mathbb{C}P^2 \# \overline{\mathbb{C}P}^2$ tout entier, ce qui permet alors de le munir d'une structure de fibré en sphères symplectique sur $\mathbb{C}P^1$, dont le diviseur exceptionnel est une section (voir la Figure~\ref{f:blownuppencilbasepoint}). Comme le fibré trivial $\mathcal{S}^2 \times \mathcal{S}^2$ ne contient pas de sphère d'auto-intersection $-1$ (car sa forme d'intersection est paire), $\mathbb{C}P^2 \# \overline{\mathbb{C}P}^2$ est alors difféomorphe à $\mathcal{S}^2 \, \tilde \times \, \mathcal{S}^2$.

Les exemples de surfaces symplectiquement réglées qui viennent ensuite sont les éclatés des fibrés en sphères symplectiques. Lorsqu'on éclate une surface symplectique $(M, \omega)$, munie d'une structure de fibré en sphères, en un point $p$, la transformée totale de la fibre passant par $p$ est l'union de deux sphères symplectiquement plongées d'auto-intersection $-1$ (l'une étant le diviseur exceptionnel et l'autre la transformée propre de la fibre) qui s'intersectent un point, de manière transverse et positive. Les autres fibres sont disjointes de l'éclatement. Ainsi, on remarque qu'une surface symplectique $(\tilde M, \tilde \omega)$ obtenue en éclatant $(M, \omega)$ un nombre fini de fois est munie d'une structure de fibration de Lefschetz symplectique.

Réciproquement, on considère $\pi : (M', \omega') \rightarrow \Sigma$ une fibration de Lefschetz symplectique dont les fibres régulières sont des sphères. Quitte à perturber $\pi$ légèrement, on peut supposer qu'il y a exactement un point singulier par fibre singulière. Notons $F$ une fibre régulière et $F'$ une fibre singulière. Pour des raison de genre, la fibre singulière $F'$ est constituée de deux sphères symplectiquement plongées $F_1'$ et $F_2'$ qui s'intersectent en un point, de manière transverse et positive. Comme toutes les fibres sont homologues, on a pour $i \in \{1, 2\}$,
$$0 = [F] \cdot [F_i ']  = ([F_1'] + [F_2']) \cdot [F_i'] = 1 +[F_i']^2,$$
d'où $[F_i']^2 = -1$. Ainsi, chaque fibre singulière est constituée de l'union de deux diviseurs exceptionnels. En contractant un des deux diviseurs exceptionnels contenus dans chacune des fibres singulières, on obtient un fibré en sphères symplectique au-dessus d'une surface. On a alors montré que les fibrations de Lefschetz symplectiques dont les fibres régulières sont des sphères sont des surfaces symplectiquement réglées. En fait dans les sections qui suivent, on verra que les exemples présentés ci-dessus constituent la totalité des types de difféomorphismes possibles pour les surfaces réglées.

La proposition suivante regroupe quelques propriétés topologiques des fibrés en sphères au-dessus de surfaces. Les démonstrations peuvent être trouvées dans~\cite[Chapter 7]{Wendl}.

\begin{prop}\label{p:topfibrationLefschetz}
Soit $\Sigma$ une surface, $\pi : M \rightarrow \Sigma$ une fibration de Lefschetz topologique dont la fibre générique est $\mathcal{S}^2$ telle que chaque fibre singulière possède exactement un point critique, $S \subset M$ une section  et $F$ une fibre régulière. On note $N$ le nombre de fibres singulières et $E_1, \dots, E_N$ les composantes irréductibles des fibres qui n'intersectent pas $S$. Alors
\begin{enumerate}
\item pour tout corps $\mathbb{K}$, la projection $\pi : M \rightarrow \Sigma$ et l'injection $ \iota : S \hookrightarrow M$ induisent des isomorphismes $\pi_* : \pi_1 (M) \rightarrow \pi_1 (\Sigma)$, $\iota_* : \pi_1 (S) \rightarrow \pi_1 (M)$, $\pi^* : H^1 (\Sigma; \mathbb{K}) \rightarrow H^1(M; \mathbb{K})$ et $ \iota^* : H^1(M; \mathbb{K}) \rightarrow H^1(S; \mathbb{K})$,
\item le groupe abélien $H_2(M; \mathbb{Z})$ est librement engendré par $[S]$, $[F]$ et $[E_1], \dots , [E_N]$,
\item $b_2^+ (M) = 1.$
\end{enumerate}
\end{prop}

Notons que la fibration de Lefschetz $\pi : M \rightarrow \Sigma$ de la Proposition~\ref{p:topfibrationLefschetz} admet toujours une section. En effet, on a vu lors de la Remarque~\ref{r:paritéréglée} que tous les fibrés en sphères sur des surfaces admettent des sections. Comme $\pi : M \rightarrow \Sigma$ est obtenue en éclatant un nombre fini de points sur un fibré en sphères au-dessus d'une surface, elle admet également des sections.


\begin{rk}\label{r:S2XS2birCP2}
Lorsqu'on éclate $\mathcal{S}^2 \times \mathcal{S}^2$ en un point $p$, les transformées propres des deux courbes rationnelles complexes d'auto-intersection $0$ passant par $p$ sont des diviseurs exceptionnels disjoints $E_1$ et $E_2$. Les courbes $E_1$ et $E_2$ intersectent chacune le diviseur exceptionnel $E$ donné par l'éclatement en $p$. Quand on contracte $E_1$, on obtient un pinceau de Lefschetz avec une unique fibre singulière et un unique point base : les fibres régulières correspondent aux images de courbes rationnelles complexes d'auto-intersection $0$ qui intersectent $E_1$ par la contraction, la fibre singulière est constituée de l'union des images de $E$ et $E_2$ par la contraction et le point base est l'image de $E_1$. Quand on contracte ensuite $E_2$, on obtient un pinceau de Lefschetz sans fibre singulière avec un unique point base (voir la Figure~\ref{f:birationalbasic}). D'après la Proposition \ref{p:pinceauCP2}, ce pinceau de Lefschetz est isomorphe à un pinceau de droites sur $\mathbb{C}P^2$.
\begin{figure}[h]
	\centering
	\includegraphics[scale=0.4]{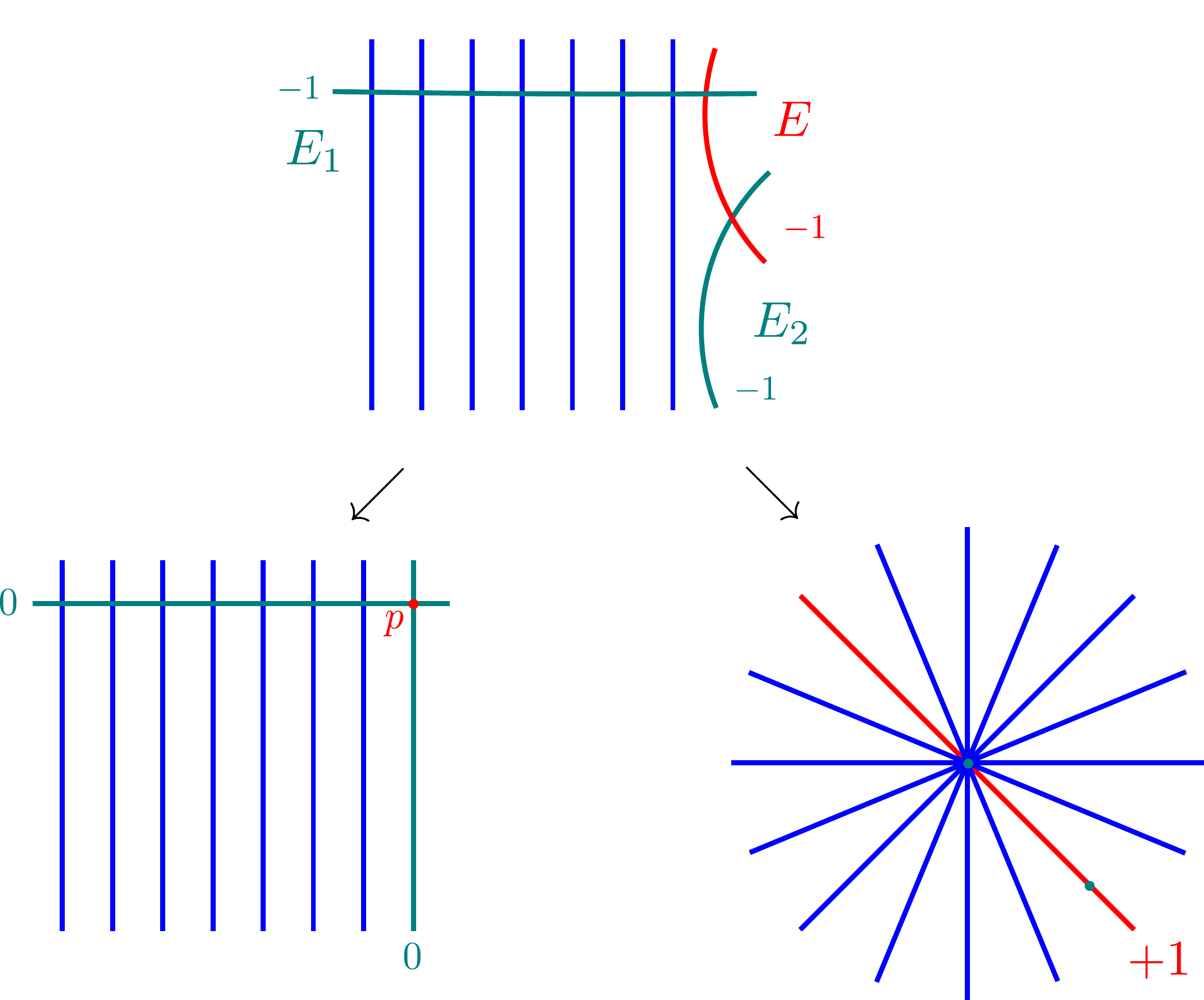}\\
	\caption{Transformation birationnelle entre $\mathcal{S}^2 \times \mathcal{S}^2$ et $\mathbb{C}P^2$.}
	\label{f:birationalbasic}
\end{figure}

En utilisant le Théorème~\ref{t:GompfThurston},
 on obtient que $\mathbb{C}P^2$ et $\mathcal{S}^2 \times \mathcal{S}^2$ sont (symplectiquement) birationnellement équivalentes. On a aussi exhibé deux modèles minimaux distincts pour $\mathcal{S}^2 \times \mathcal{S}^2 \# \overline{\mathbb{C}P}^2$ : on verra dans la Section \ref{s:spheresimmergees} que ce phénomène est spécifique aux surfaces réglées.
\end{rk} 

\section{Courbes rationnelles symplectiquement plongées d'auto-intersection positive}\label{s:GromovMcDuff}

Dans cette section, on explique les grandes lignes de la démonstration du théorème suivant, dû à McDuff~\cite{McDuff} et à Gromov~\cite{Gromov} (Gromov a d'abord traité le cas des courbes rationnelles d'auto-intersection $1$ symplectiquement plongées dans les surfaces symplectiques et McDuff s'est ensuite occupée de tous les autres cas). 

\begin{thm}[Théorème de McDuff et Gromov] \label{t:McDuff}
Soit $S$ une courbe rationnelle symplectiquement plongée dans une surface symplectique $(M, \omega)$. On suppose que $[S]^2 \geq 0$. Alors $M$ est une surface symplectiquement réglée. En particulier lorsque $(M, \omega)$ est minimale :
\begin{enumerate}
\item \label{i:1} Si $[S]^2 = 0$, alors $(M, \omega)$ est symplectomorphe à un fibré en sphères symplectique au-dessus d'une surface, dans lequel $S$ est identifiée à une fibre.
\item Si $[S]^2 = 1$, alors $(M, \omega)$ est symplectomorphe à $(\mathbb{C}P^2, c \omega_{FS})$, avec $c >0$, et $S$ est identifiée à une droite.
\item Si $[S]^2 > 1$, on peut trouver une autre courbe rationnelle $S'$ symplectiquement plongée dans $M$ qui vérifie $[S']^2 \in \{0, 1 \}$. Autrement dit, on se ramène à l'un des deux cas précédents.
\end{enumerate}
\end{thm}

La démonstration du Théorème \ref{t:McDuff} se base sur les techniques de courbes pseudoholomorphes présentées dans la Section \ref{s:Jholo}. On fera usage de la proposition suivante, qui est une conséquence du théorème de transversalité automatique pour les courbes pseudoholomorphes rationnelles immergées (Théorème~\ref{t:transauto2}).

\begin{prop}[{\cite[Proposition~2.53]{Wendl}}] \label{p:feuilletage}
Soit $(M, J)$ une surface presque complexe et $u : \mathcal{S}^2 \rightarrow M$ une courbe rationnelle plongée $J$--holomorphe telle que $m : = [u]^2 \geq 0$. Alors pour tous points deux à deux distincts $p_1, \dots , p_m$, la courbe $u$ est Fredholm régulière pour le problème avec contraintes, et un voisinage $\mathcal{U} \subset \mathcal{M}_{0,m}(J;p_1, \dots, p_m)$ de $u$ admet la structure d'une variété lisse de dimension $2$ qui vérifie :
\begin{enumerate}
\item toute courbe $v \in \mathcal{U}$ est plongée,
\item toutes courbes $v,w \in \mathcal{U}$ s'intersectent seulement en les points $p_1, \dots, p_m$ et toutes les intersections sont transverses,
\item les images des courbes de $\mathcal{U}$ privées des points $p_1, \dots p_m$ forment un feuilletage d'un voisinage ouvert de $u(\mathcal{S}^2) \backslash \{ p_1, \dots, p_m \}$ dans $M \backslash \{p_1, \dots, p_m \}$.
\end{enumerate} 
\end{prop}
\begin{figure}[h]
	\centering
	\includegraphics[scale=0.5]{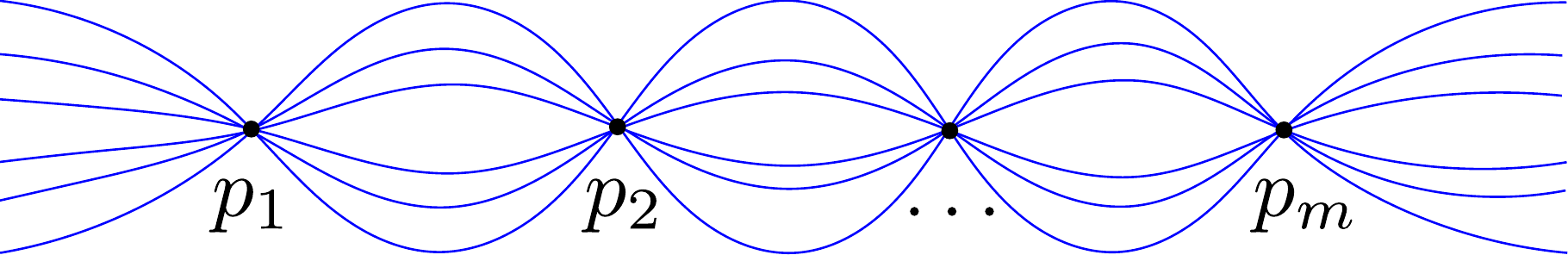}\\
	\caption{Feuilletage local par une famille de dimension $2$ de courbes pseudoholomorphes rationnelles d'auto-intersection $m$ passant par $m$ contraintes ponctuelles $p_1, \dots, p_m$.}
	\label{f:pencil}
\end{figure}

Dans un premier temps on se sert des propriétés des espaces de modules et de la Proposition \ref{p:feuilletage} pour construire sur $(M, \omega )$ un pinceau de Lefschetz dont $S$ est une fibre. Cette étape fait l'objet du théorème suivant. Plus de détails peuvent être trouvés dans \cite[Chapter~6]{Wendl}.

\begin{thm}\label{t:pinceaudeLefschetz}
Soit $S$ une courbe rationnelle symplectiquement plongée dans une surface symplectique $(M, \omega)$. On suppose $m := [S]^2 \geq 0$. Alors pour tous points deux à deux distincts $p_1, \dots p_m$ sur $S$, $(M, \omega)$ admet un pinceau de Lefschetz symplectique (ou une fibration de Lefschetz symplectique lorsque $m=0$) dont les points bases sont $p_1, \dots, p_m$, dont $S$ est une fibre lisse et qui contient au plus un point critique par fibre. De plus, l'ensemble des fibres singulières est vide si et seulement si $m \in \{0,1 \}$ et $(M , \omega)$ est relativement minimale par rapport à $S$.
\end{thm}

\begin{proof}[Esquisse de preuve.]
On commence par choisir une structure presque complexe $J$ dominée par $\omega$ telle que $S$ est l'image d'une courbe $J$--holomorphe plongée $u$. On peut alors calculer l'indice de $u$ à l'aide de la formule \eqref{eq:indexformula} : $\ind (u) = 2 + 2 [S]^2 \geq 2m$. Le théorème de transversalité automatique (Théorème~\ref{t:transauto1}) nous assure alors que la courbe $u$ est Fredholm régulière pour le problème avec contraintes $p_1, \dots, p_m$. Quitte à réaliser une perturbation $\mathcal{C}^\infty$ arbitrairement petite de la structure presque complexe $J$ et une perturbation $\mathcal{C}^\infty$ arbitrairement petite de la courbe $u$ passant par $p_1, \dots, p_m$, on peut supposer que $J$ est générique pour $p_1, \dots, p_m$ (c'est-à-dire que tout courbe $J$--holomorphe simple est Fredholm régulière pour le problème avec contraintes $p_1, \dots, p_m$).
Comme toutes les courbes considérées sont Fredholm régulières, l'espace de modules $\mathcal{M}_S (J) := \mathcal{M}^*_{0,m} ([S];J; p_1, \dots, p_m)$ est une variété lisse (a priori possiblement non compacte) de dimension $2$. L'intérêt d'imposer ces contraintes est double. D'une part, on impose à l'espace de modules d'avoir une dimension adaptée dans l'optique de construire un pinceau de Lefschetz sur $(M, \omega)$ dont les fibres régulières sont les images des courbes de $\mathcal{M}_S (J)$. D'autre part, il est bien plus facile de décrire les courbes nodales pouvant exister dans le compactifié de l'espace de modules lorsque l'indice (contraint) des courbes considérées est bas.

L'étape suivante consiste justement à comprendre les courbes nodales du compactifié $ \overline{\mathcal{M}}_S (J) := \overline{\mathcal{M}^*}_{0,m} ([S];J; p_1, \dots, p_m)$. Pour ce faire, on considère une courbe nodale $u_\infty \in \overline{\mathcal{M}}_S (J) \backslash \mathcal{M}_S (J)$. En utilisant les conditions de stabilité (pour éliminer l'existence de composante constante), de conservation du genre géométrique, ainsi que la généricité de $J$ (les indices contraints des courbes $J$--holomorphe simples qui pourraient éventuellement apparaître comme image de composantes de courbes nodales sont positifs ou nuls), on peut montrer que $u_\infty$ possède exactement deux composantes plongées, rationnelles, d'indices contraints respectifs nuls, qui s'intersectent exactement une fois, de manière transverse et positive, en un point qui est distinct de $p_1, \dots, p_m$ (voir~\cite[Chapter~4]{Wendl} ou l'Annexe~\ref{c:annexe} pour plus de détails). L'union de ces deux composantes passe évidemment par les contraintes ponctuelles. De par la généricité de $J$, et puisque leurs indices contraints sont nuls, les composantes des courbes nodales sont en nombre fini. Notons $\mathcal{N} \subset M$ l'union des images des courbes nodales.

On montre ensuite que l'ensemble $\mathcal{R}$ des points de $M \backslash \left( \{ p_1, \dots, p_m \} \cup \mathcal{N} \right)$ par lesquels passent les courbes de $\mathcal{M}_S (J)$ est :
\begin{itemize}
\item ouvert dans $M \backslash \left( \{ p_1, \dots, p_m \} \cup \mathcal{N} \right)$ grâce à la Proposition \ref{p:topfibrationLefschetz},
\item fermé dans $M \backslash \left( \{ p_1, \dots, p_m \} \cup \mathcal{N} \right)$ par compacité de $\overline{\mathcal{M}}_S(J)$.
\end{itemize}
Comme l'espace de modules $\mathcal{M}_S (J)$ est non vide (car il contient la courbe $u$) et l'ensemble $M \backslash \left( \{ p_1, \dots, p_m \} \cup \mathcal{N} \right)$ est connexe, on obtient $\mathcal{R} = M \backslash \left( \{ p_1, \dots, p_m \} \cup \mathcal{N} \right)$.

Pour finir, une étude minutieuse au voisinage des courbes nodales permet de munir $\overline{\mathcal{M}}_S(J)$ d'une structure de variété lisse compacte de dimension $2$. On considère alors l'application lisse $\pi : M \backslash \{p_1, \dots, p_m \}\rightarrow \overline{\mathcal{M}}_S(J)$ qui à un point $p \in  M \backslash \{p_1, \dots, p_m \}$ associe l'unique courbe de $\overline{\mathcal{M}}_S(J)$ passant par $p$. Dans le cas où $m=0$, l'application $\pi$ définit alors une fibration de Lefschetz symplectique sur $(M, \omega)$. Mais si $m >0$, pour tout $i \in \{ p_1, \dots , p_m \}$, on peut montrer que l'application $\overline{\mathcal{M}}_S(J) \rightarrow \mathbb{P}(T_{p_1} M,J) \simeq \mathbb{C}P^1$ qui à une courbe $u \in \overline{\mathcal{M}}_S(J)$ associe son espace tangent en $p_i$ est un difféomorphisme. L'application $\pi$ définit alors la fibration de Lefschetz symplectique voulue sur $(M, \omega)$.
\end{proof}

\begin{rk}
La preuve du Théorème~\ref{t:pinceaudeLefschetz} constitue un des rares cas où l'on sait déterminer le type de difféomorphisme d'un espace de modules de courbes pseudoholomorphes (ici $\overline{\mathcal{M}}_S(J)$) pour une structure presque complexe générique. On remarque également que dans le cas où $[S]^2 \in \{0,1 \}$ et où la surface symplectique $(M , \omega)$ est relativement minimale par rapport à $S$, l'espace de modules $\mathcal{M}_S(J)$ est compact pour $J$ générique (en effet, dans ce cas toute courbe nodale de $\overline{\mathcal{M}}_S(J)$ possède une composante qui est un diviseur exceptionnel disjoint de $S$).
\end{rk}

En considérant des espaces de modules à paramètre, on peut améliorer la preuve du Théorème \ref{t:pinceaudeLefschetz} et montrer que deux pinceaux construits de cette manière sont reliés par une isotopie. Cette amélioration ne nous servira pas pour démontrer le Théorème \ref{t:McDuff}, mais on en fera l'usage dans le Chapitre~3. On l'énonce ci-dessous.

\begin{thm}[{\cite[Theorem~G]{Wendl}}] \label{t:isotopiepinceaudeLefschetz}
Soit $S$ une courbe rationnelle symplectiquement plongée dans une surface symplectique $(M, \omega)$. On suppose $m := [S]^2 \geq 0$ et on choisit  $p_1, \dots p_m$ des points deux à deux distincts sur $S$. Soit une structure presque complexe générique $J$ dominée par $\omega$.
Alors :
\begin{enumerate}
\item Le pinceau (resp. fibration) de Lefschetz donné par le Théorème \ref{t:pinceaudeLefschetz} avec points bases $p_1, \dots, p_m$ est isotope (à points bases fixés) à un pinceau (resp. une fibration) de Lefschetz dont toutes les composantes des fibres sont des courbes $J$--holomorphes rationnelles plongées. De plus, toute courbe $J$--holomorphe qui est homologue à la fibre et qui passe par tous les points bases $p_1, \dots , p_m$ est une fibre.
\item Pour une famille lisse à paramètre $\{ \omega_s \}_{s \in [0,1]}$ de formes symplectiques avec $\omega_0 = \omega$, et pour une famille lisse à paramètre $\{ J_s \}_{s [0,1]}$ de structures presque complexes dominées par $\{ \omega_s \}_{s \in [0,1]}$ avec $J_0 = J$, il existe une isotopie lisse de pinceaux (resp. fibrations) de Lefschetz avec points bases $p_1, \dots, p_m$ dont les composantes des fibres sont des courbes $J_s$--holomorphes, pour $s \in [0,1]$.
\end{enumerate}
\end{thm}

\begin{proof}[Esquisse de preuve du Théorème \ref{t:McDuff}.]
D'après le Théorème \ref{t:pinceaudeLefschetz}, on peut construire un pinceau (ou une fibration) de Lefschetz symplectique sur $(M, \omega)$ dont $S$ est une fibre et dont le nombre de points bases $m$ est égal à $[S]^2$.

\begin{enumerate}
\item Lorsque $[S]^2 =0$, on obtient une fibration de Lefschetz symplectique dont les fibres singulières sont constituées de l'union de deux diviseurs exceptionnels qui s'intersectent exactement une fois, de manière transverse et positive. En contractant un des deux diviseurs exceptionnels dans chacune des fibres singulières, on obtient un fibré en sphères symplectique au-dessus d'une surface. Donc $(M, \omega)$ est un éclatement d'un fibré en sphères symplectique au-dessus d'une surface : c'est une surface symplectiquement réglée.

\item Lorsque $[S]^2 = 1$, le pinceau de Lefschetz possède un unique point base $p$. Chaque fibre singulière est l'union d'une courbe rationnelle plongée passant par $p$ d'auto-intersection $0$ et d'un diviseur exceptionnel qui l'intersecte exactement une fois, de manière positive et transverse. En contractant ces diviseurs exceptionnels, on obtient un pinceau avec un unique point base et sans fibre singulière. La Proposition \ref{p:pinceauCP2} nous assure que ce pinceau est isomorphe à un pinceau de droites complexes sur $\mathbb{C}P^2$. Enfin en utilisant le Théorème~\ref{t:GompfThurston}, le fait que $H_{dR}^2(\mathbb{C}P^2; \mathbb{R})$ est engendré par $[\omega_{FS}]$ et le théorème de stabilité de Moser (Théorème~\ref{t:Moser}) on a en fait l'existence d'une constante $c >0$ telle que la variété obtenue après contractions est symplectomorphe à $(\mathbb{C}P^2, c \omega_{FS})$. Par conséquent, $(M, \omega)$ est une surface rationnelle (donc une surface symplectiquement réglée).

\item Lorsque $[S]^2 >0$, le pinceau de Lefschetz possède au moins deux points bases. D'après la Proposition~\ref{p:pinceau2ptbases}, il existe alors au moins une fibre singulière. Cette fibre singulière possède deux composantes rationnelles plongées $S_1$ et $S_2$ qui s'intersectent exactement une fois, de manière transverse et positive. On note $m_1$ et $m_2$ le nombre de points bases du pinceau par lesquels $S_1$ et $S_2$ passent respectivement, de sorte que $m = m_1 + m_2$. Puisque la classe $[S_1]+[S_2]$ est homologue à la classe d'une fibre régulière, on a pour tout $i \in \{1,2\}$,
$$m_i = [S_i] \cdot ([S_1] +[S_2]) = 1+ [S_i]^2.$$
Donc pour tout $i \in \{1,2\}$, on a $[S_i]^2 = m_i -1 \geq -1$. De plus, puisque $([S_1] + [S_2])^2 =0$, on a également $m = [S_1]^2 + 2 +[S_2]^2$. Quitte à réindexer, on en déduit que $m > [S_1]^2 \geq 0$. On réitère alors le raisonnement avec $S_1$ jusqu'à trouver une courbe rationnelle symplectiquement plongée d'auto-intersection $0$ ou $1$. On conclut enfin en se ramenant aux cas précédents. 
\end{enumerate}
\end{proof}

\begin{rk}\label{r:sphereautointgeq2}
Dans le cas où $[S]^2 > 1$, on se ramène à trouver l'existence d'une courbe rationnelle symplectiquement plongée $S'$ d'auto-intersection $0$ ou $1$. Si $[S']^2 = 1$, on a montré que la surface symplectique $(M, \omega)$ était nécessairement rationnelle. Si en revanche $(M,\omega)$ ne contient pas de courbe rationnelle symplectiquement plongée d'auto-intersection $1$, on peut en fait aussi montrer que $(M, \omega)$ est une surface rationnelle. En effet, on peut utiliser certaines composantes apparaissant au cours de la récurrence pour trouver une autre courbe rationnelle symplectiquement plongée $S''$ d'auto-intersection $0$ qui intersecte $S'$ exactement une fois, de manière transverse et positive. On peut alors utiliser les courbes rationnelles $S'$ et $S''$ (et des techniques de courbes pseudoholomorphes) pour construire deux fibrations de Lefschetz symplectiques sur $(M, \omega)$, de sorte que chacune des fibres d'une des fibrations intersecte toutes les fibres de l'autre fibration transversalement positivement. Après contraction d'un nombre fini de diviseurs exceptionnels, on peut montrer qu'on obtient la surface symplectique $(\mathcal{S}^2 \times \mathcal{S}^2, \sigma_1 \oplus \sigma_2)$ où $\sigma_1$ et $\sigma_2$ sont deux formes volumes telles que $\int_{\mathcal{S}^2 \times \{* \}} \sigma_1 = \int_{S'} \omega$ et $\int_{\{* \} \times \mathcal{S}^2 } \sigma_1 = \int_{S''} \omega$ (voir~\cite[Theorem~6.8, Theorem~6.9]{Wendl}).
\end{rk}

\begin{figure}[h]
	\centering
	\includegraphics[scale=0.5]{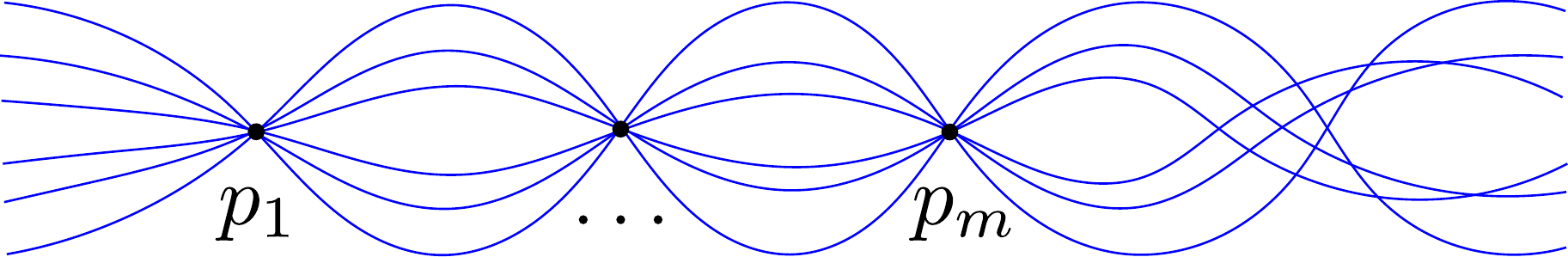}\\
	\caption{Représentation schématique d'une famille de dimension $2$ de courbes $J$--holomorphes de genre $1$ et d'auto-intersection $m+1$ passant par $m$ contraintes ponctuelles $p_1, \dots p_m$, pour une structure presque complexe $J$ générique.}
	\label{f:notapencil}
\end{figure}

Résumons les points centraux de la preuve du Théorème \ref{t:McDuff} :
\begin{itemize}
\item On utilise le théorème de transversalité automatique pour trouver une structure presque complexe $J$ générique qui préserve $TS$ : cette étape nous permet de munir les espaces de modules qui apparaissent dans la preuve d'une structure de variété lisse.
\item Comme on s'intéresse à l'espace de modules des courbes $J$--holomorphes \emph{rationnelles} simples homologues à $S$, une coïncidence numérique nous permet d'obtenir un espace de modules de dimension $2$ en imposant $m = [S]^2$ contraintes ponctuelles.
\item L'espace de modules avec contraintes $\mathcal{M}_S (J)$ est de dimension $2$, ce qui nous permet de comprendre aisément les courbes nodales du compactifié, puis de construire ainsi le pinceau (ou la fibration) de Lefschetz du Théorème \ref{t:pinceaudeLefschetz}.
\item Quand $m \geq 2$, des contraintes topologiques permettent d'affirmer l'existence de courbes nodales dans $\overline{\mathcal{M}}_S (J)$. On peut alors trouver une autre courbe rationnelle symplectiquement plongée d'auto-intersection strictement plus petite parmi les composantes d'une fibre singulière. \item On réitère le raisonnement jusqu'à trouver une courbe rationnelle symplectiquement plongée d'auto-intersection $0$ ou $1$, ce qui permet de conclure.
\end{itemize}
Dans quelle mesure peut-on étendre ce raisonnement aux courbes irrationnelles symplectiquement plongées dans la surface symplectique $(M, \omega)$ ?
Tout d'abord, pour que la condition de transversalité automatique soit respectée pour des courbes symplectiquement plongées $S$ de genre $g > 0$, on doit avoir $[S]^2  > 2g-2$. Ensuite, l'indice d'une courbe pseudoholomorphe de genre $g$ homologue à $S$ est égal à $\ind (u) = 2-2g + 2[S]^2$. Le nombre $m$ de contraintes ponctuelles $p_1, \dots, p_m$ à imposer sur les courbes pour obtenir un espace de modules $\mathcal{M}_S (J)$ de dimension $2$ est donc égal à $[S]^2 -g$. Cette fois les points $p_1, \dots, p_m$ ne permettent pas d'avoir un contrôle sur toutes les intersections entre les courbes de $\mathcal{M}_S (J)$ (pour deux courbes distinctes de $\mathcal{M}_S (J)$, il y a exactement $g$ points d'intersection sur lesquels on n'a aucun contrôle, voir la Figure~\ref{f:notapencil}). Cela ne permet pas de construire un pinceau de Lefschetz de la même façon que dans le Théorème~\ref{t:pinceaudeLefschetz}. L'Exemple~\ref{e:cubiquespar9pts} du pinceau de cubiques dans $\mathbb{C}P^2$ relève d'une situation non générique qui semble spécifique aux surfaces rationnelles complexes. Dans le Chapitre~3, on conserve cette idée de \og casser \fg{} la courbe $S$ pour trouver une composante plongée de genre strictement plus petit, mais d'intersection suffisamment grande par rapport à son genre, puis de continuer le raisonnement par récurrence sur le genre des courbes. Un des points cruciaux consiste à trouver un moyen de prouver l'existence de courbes nodales dans $\overline{\mathcal{M}}_S (J)$.

\begin{rk}
Une autre stratégie est présentée dans les démonstrations de~\cite[Theorem~9.4.1, Theorem~9.4.4]{mcduff2004j}. Elle consiste à montrer, sous l'hypothèse que l'espace de modules $\mathcal{M}^*_0 ([S];J)$ est compact, que l'application d'évaluation
$$\mathrm{ev} : \mathcal{M}^*_{0,m+1}([S];J) \rightarrow M^{m+1} \backslash \Delta$$
est un difféomorphisme, où $\Delta$ désigne la grande diagonale de $M^{m+1}$ (c'est-à-dire l'ensemble des $(m+1)$--uplets de $M^{m+1}$ qui possèdent au moins deux composantes identiques).
Si $m+1 \geq 3$, on peut alors se servir de l'application $\mathrm{ev}$ pour construire une suite de courbes de $\mathcal{M}^*_0 ([S];J)$ qui converge vers une courbe non plongée dans $\mathcal{M}^*_0 ([S];J)$ et ainsi aboutir à une contradiction.

Bien qu'il ne soit pas possible d'étendre telle quelle cette approche au cas des courbes de genre strictement plus grand pour des raisons évidentes concernant les dimensions des espaces considérés,  la stratégie qu'on emploie dans le Chapitre~3 s'en inspire.
\end{rk}

Le Théorème~\ref{t:McDuff} de McDuff et Gromov ne suffit malheureusement pas pour classifier les surfaces symplectiquement réglées. Il faut avant tout vérifier que le fait de contenir une sphère symplectiquement plongée d'auto-intersection positive est préservé par transformation birationnelle symplectique. L'opération d'éclatement ne pose aucun problème : quitte à perturber légèrement $S$, on peut toujours supposer que le point qu'on éclate est disjoint de $S$. En revanche, la propriété de contenir une sphère $S$ symplectiquement plongée d'auto-intersection positive n'est \emph{a priori} pas préservée par l'opération de contraction. Un éventuel diviseur exceptionnel $E$ pourrait intersecter strictement plus d'une fois $S$. L'image de $S$ après la contraction de $E$ ne serait alors pas plongée. Comme le genre géométrique des courbes est conservé par les transformations birationnelles, l'étape suivante consiste à étudier les courbes symplectiques singulières rationnelles dans les surfaces symplectiques. 

\section{Courbes rationnelles symplectiquement positivement immergées}\label{s:spheresimmergees}

On présente ici une caractérisation des surfaces symplectiquement réglées par les courbes symplectiques rationnelles (éventuellement singulières) qu'elle contiennent. Celle-ci nous permet ensuite d'obtenir la classification des surfaces symplectiquement réglées et des surfaces symplectiques rationnelles. Le point clé pour obtenir la caractérisation est un résultat démontré par McDuff dans~\cite{McDuffImmersed}. 
La preuve du résultat étant plutôt technique, on se contente ici d'en présenter les idées principales (les détails techniques seront mis sous silence ici). Une exposition complète et détaillée peut être trouvée dans~\cite{Wendl}. Le résultat en question s'intéresse au courbes (rationnelles) symplectiquement positivement immergées dans les surfaces symplectiques. Elles constituent un certain type de courbes symplectiques singulières dont on donne la définition ci-dessous.
  
\begin{df}
Soit $(M , \omega)$ une surface symplectique, $\Sigma$ une surface et $\iota :  \Sigma \looparrowright M$ une immersion. On dit que $\iota (\Sigma)$ est une \emph{courbe symplectiquement positivement immergée} dans $(M , \omega)$ si $\iota^* \omega$ est une forme symplectique sur $\Sigma$ et tous les points multiples de $\iota (\S)$ sont des points doubles (c'est-à-dire tout $p \in  M$ possède au plus deux antécédents par $\iota$) qui sont transverses et positifs.
\end{df}

On peut désormais énoncer la caractérisation des surfaces symplectiquement réglées.

\begin{thm}[Caractérisation des surfaces symplectiquement réglées]\label{t:caracterisationsurfacereglee}
Soit $(M, \omega )$ une surface symplectique. Les propositions suivantes sont équivalentes :
\begin{enumerate}[label=(\arabic*)]
\item \label{i:CSR1} $(M, \omega )$ est une surface symplectiquement réglée,
\item \label{i:CSR2} il existe une courbe rationnelle $S \subset M$ symplectiquement plongée telle que $[S]^2 \geq 0$,
\item \label{i:CSR3} il existe une courbe rationnelle $S \looparrowright M$ symplectiquement positivement immergée vérifiant $c_1([S]) \geq 2$,
\item \label{i:CSR4} il existe $J \in \mathcal{J}_{\tau} (M, \omega )$ et une courbe rationnelle $J$--holomorphe $u$ injective quelque part, Fredholm régulière, satisfaisant $\ind (u) \geq 2$,
\end{enumerate}
En particulier, chacune des propositions énoncées ci-dessus est invariante par transformation birationnelle symplectique.
\end{thm}

\begin{rk}
Les propriétés du Théorème~\ref{t:caracterisationsurfacereglee} peuvent également être interprétées en terme d'invariants de Gromov--Witten. En effet, elles sont vérifiées si et seulement si l'invariant $\mathrm{GW}^{(M, \omega)}_{0,m,A}$ n'est pas identiquement nul (on dit dans ce cas que la surface symplectique $(M, \omega)$ est uniréglée). Plus de détails à ce sujet peuvent être trouvés dans~\cite[Section~7]{Wendl}.
\end{rk}

La démonstration du Théorème \ref{t:caracterisationsurfacereglee} repose principalement sur le théorème suivant. On en esquisse une démonstration à la fin de la section. 

\begin{thm}[Théorème de McDuff pour les sphères symplectiquement positivement immergées, \cite{McDuffImmersed}]\label{t:McDuffimmersed}
Soit $S$ une courbe rationnelle symplectiquement positivement immergée dans une surface symplectique $(M, \omega)$ telle que $c_1([S]) \geq 2$. Alors $(M, \omega)$ contient également une courbe rationnelle symplectiquement plongée d'auto-intersection positive.
\end{thm}

\begin{rk}
Le Théorème~\ref{t:McDuffimmersed} peut aisément être généralisé au cas où $S$ est une courbe rationnelle symplectique singulière. En effet, on peut toujours perturber $S$ symplectiquement, de manière complexe au voisinage des singularités, de façon à obtenir une courbe rationnelle symplectiquement positivement immergée homologue à $S$.
\end{rk}

\begin{proof}[Démonstration du Théorème~\ref{t:caracterisationsurfacereglee}]
L'implication $\ref{i:CSR3} \Rightarrow \ref{i:CSR2}$ est donnée par le Théorème~\ref{t:McDuffimmersed} de McDuff pour les sphères symplectiquement positivement immergées. L'implication $\ref{i:CSR2} \Rightarrow \ref{i:CSR1}$ découle du Théorème~\ref{t:McDuff} de McDuff et Gromov.

On montre ensuite l'implication $\ref{i:CSR1} \Rightarrow \ref{i:CSR3}$. 
D'après la formule d'adjonction, une droite $\ell \subset \mathbb{C}P^2$ satisfait $c_1([\ell]) = 2+1 \geq 2$ et une fibre $F$ d'un fibré en sphères symplectique au-dessus d'une surface vérifie $ c_1([F])= 2+0 \geq 2$. Il suffit alors de montrer que la propriété de contenir une sphère $S \looparrowright M$ symplectiquement positivement immergée vérifiant $c_1([S]) \geq 2$ est invariante par éclatement, contraction et déformation symplectique. 
Lorsqu'on éclate en un point $p \in M$, quitte à réaliser une perturbation arbitrairement petite de $S$, on peut supposer que $p$ n'appartient pas à $S$, ce qui permet de conclure. 
Lorsqu'on contracte un diviseur exceptionnel $E$, on introduit une structure presque complexe auxiliaire $J \in \mathcal{J}_{\tau} (M, \omega )$ telle que $S$ est l'image d'une courbe $J$--holomorphe. D'après le Corollaire~\ref{c:divexisotopy}, quitte à réaliser une isotopie, on peut supposer sans perte de généralité que $E$ est l'image d'une courbe $J$--holomorphe. On a alors la propriété de positivité d'intersection entre $S$ et $E$. Quitte à perturber symplectiquement $S$, de manière complexe au voisinage des points d'intersection entre $E$ et $S$, on peut alors supposer que $E$ et $S$ s'intersectent transversalement positivement. De cette manière, l'image $S'$ de $S$ par la contraction de $E$ est symplectiquement positivement immergée. On note $k = [E] \cdot [S]$. On peut alors calculer la première classe de Chern de $S'$ en utilisant la formule d'adjonction $c_1([S']) = c_1([S]) + k c_1 ([E]) = c_1([S]) + k \geq 2$, ce qui conclut. 
En ce qui concerne les déformations symplectiques, on commence par appliquer le Théorème~\ref{t:McDuffimmersed} à $S$ pour trouver une courbe rationnelle $S_0$ symplectiquement plongée dans $(M,\omega)$ d'auto-intersection positive (qui vérifie donc $c_1([S_0]) \geq 2$ d'après la formule d'adjonction). Le Théorème~\ref{t:isotopiepinceaudeLefschetz} assure alors l'existence d'une autre courbe rationnelle symplectiquement plongée d'auto-intersection positive une fois la déformation symplectique effectuée.

Il ne reste plus qu'à montrer l'équivalence entre les deux dernières propositions. L'implication $\ref{i:CSR3} \Rightarrow \ref{i:CSR4}$ est immédiate. Pour démontrer sa réciproque $\ref{i:CSR4} \Rightarrow \ref{i:CSR3}$, on perturbe la structure presque complexe $J$ en une structure générique $J' \in \mathcal{J}_{\tau}^{\mathrm{reg}} (M, \omega) $ et on utilise la régularité Fredholm de la courbe $u$ pour trouver une courbe $J'$--holomorphe arbitrairement proche de $u$. En choisissant $J'$ suffisamment générique, on peut alors trouver une courbe symplectiquement positivement immergée proche de $u'$ homologue à $u$, comme vu au dernier paragraphe de la Section~\ref{s:contraintesdérivées}. Le fait que $\ind(u) = -2 + 2 c_1([u])$ permet de conclure.
\end{proof}

Grâce au Théorème~\ref{t:caracterisationsurfacereglee}, on peut désormais obtenir la classification des surfaces symplectiquement réglées et des surfaces symplectiques rationnelles.

\begin{thm}[Classification des surfaces symplectiquement réglées] \label{t:classificationsurfacesreglees}
Une surface symplectique $(M, \omega)$ est symplectiquement réglée si et seulement si elle peut être munie d'une fibration de Lefschetz symplectique avec fibre régulière de genre $0$ et au plus une singularité par fibre singulière. Autrement dit $(M, \omega)$ est symplectiquement réglée si et seulement si c'est un éclatement de $(\mathbb{C}P^2,c \omega_{FS})$ avec $c >0$, ou un éclatement d'un fibré en sphères symplectique au-dessus d'une surface.
\end{thm}

\begin{proof}
Si $(M, \omega)$ est symplectiquement réglée, le Théorème~\ref{t:caracterisationsurfacereglee} assure l'existence d'une sphère $S \subset M$ symplectiquement plongée d'auto-intersection positive. Le Théorème~\ref{t:McDuff} de McDuff et Gromov permet de conclure.
La réciproque est une conséquence immédiate de la définition des surfaces symplectiquement réglées.
\end{proof}

\begin{thm}\label{t:pinceaurationnel}
Une surface symplectique $(M, \omega)$ est rationnelle si et seulement si elle admet un pinceau de Lefschetz symplectique dont les fibres régulières sont des courbes rationnelles (éventuellement sans point base, mais dont la base est une sphère).
\end{thm}

\begin{proof}
On suppose que $(M, \omega)$ est rationnelle. D'après le Théorème~\ref{t:classificationsurfacesreglees}, $(M, \omega)$ est soit un éclatement de $(\mathbb{C}P^2,c \omega_{FS})$ avec $c >0$, soit un éclatement d'un fibré en sphères symplectique au-dessus d'une surface.
Dans le premier cas, on peut directement conclure. Dans le second cas il suffit de remarquer qu'une surface rationnelle symplectique est simplement connexe (voir la Remarque~\ref{r:blowuptop}). La Proposition~\ref{p:topfibrationLefschetz} nous assure alors que la base du fibré est de genre $0$.

La réciproque découle directement du Théorème~\ref{t:McDuff} de Gromov et McDuff, de la Remarque~\ref{r:sphereautointgeq2} et de la Remarque~\ref{r:blowuptop} lorsque le pinceau admet au moins un point base (car dans ce cas, $(M, \omega)$ contient une courbe rationnelle symplectiquement plongée d'auto-intersection strictement positive). On suppose que le pinceau de Lefschetz n'admet pas de point base. Puisque ses fibres régulières sont des courbes rationnelles symplectiquement plongées d'auto-intersection $0$, le Théorème~\ref{t:McDuff} de Gromov et McDuff nous assure que $(M, \omega)$ est l'éclatement d'un fibré en sphères symplectique au-dessus d'une surface. D'après la Proposition~\ref{p:topfibrationLefschetz}, $(M, \omega)$ est simplement connexe donc la base du fibré est une sphère. D'après les résultats de la Section~\ref{s:ssr,ssr}, les seuls modèles minimaux possibles pour $(M, \omega)$ sont alors $\mathcal{S}^2 \times \mathcal{S}^2 $ et $\mathbb{C}P^2$, ce qui conclut la preuve.
\end{proof}

\begin{rk}
Pour la réciproque dans la preuve du Théorème~\ref{t:pinceaurationnel}, lorsque le pinceau admet au moins un point base, on pouvait aussi éclater tous les points bases pour se ramener au cas sans point base.
\end{rk}

\begin{thm}[Classification des surfaces symplectiques rationnelles] \label{t:classificationsurfacesrationnelles}
Une surface symplectique $(M, \omega)$ est rationnelle si et seulement si elle peut être munie d'une fibration de Lefschetz symplectique dont les fibres régulières et la base sont des sphères. Autrement dit $(M, \omega)$ est symplectiquement réglée si et seulement si c'est un éclatement de $(\mathbb{C}P^2,c \omega_{FS})$  avec $c >0$ ou un éclatement d'un fibré en sphères symplectique au-dessus d'une sphère.
\end{thm}

\begin{proof}
On suppose que $(M, \omega)$ est une surface symplectique rationnelle. En particulier $(M, \omega)$ est symplectiquement réglée, donc d'après le Théorème~\ref{t:classificationsurfacesreglees}, c'est un éclatement de $(\mathbb{C}P^2,c \omega_{FS})$ avec $c >0$, ou un éclatement d'un fibré en sphères symplectique au-dessus d'une surface. Dans le second cas, comme $M$ est simplement connexe, la base du fibré est une sphère.

La réciproque provient du Théorème~\ref{t:pinceaurationnel}.
\end{proof}

La proposition suivante est un corollaire du Théorème~\ref{t:McDuffimmersed} de McDuff sur les sphères symplectiquement positivement immergées. Elle affirme que les surfaces symplectiquement réglées sont les seules surfaces symplectiques à posséder plusieurs modèles minimaux. Autrement dit, une surface symplectique non rationnelle admet un unique modèle minimal.

\begin{prop}
Soit $(M, \omega)$ une surface symplectique et $(N_1, \omega_1)$, $(N_2, \omega_2)$ deux surfaces minimales obtenues à partir de $(M, \omega)$ en contractant des collections maximales de diviseurs exceptionnels deux à deux disjoints. Si $(M, \omega)$ n'est pas une surface symplectiquement réglée alors $(N_1, \omega_1)$ et $(N_2, \omega_2)$ sont dans la même classe de déformation symplectique.
\end{prop}

\begin{proof}
Soit $\{ E_1, \dots, E_k \}$ et $\{ E_1', \dots, E_\ell' \}$ deux collections maximales de diviseurs exceptionnels deux à deux disjoints. On suppose que les surfaces symplectiques minimales $(N_1, \omega_1)$ et $(N_2, \omega_2)$ obtenues respectivement en contractant $\{ E_1, \dots, E_k \}$ et $\{ E_1', \dots, E_\ell' \}$ ne sont pas dans la même classe de déformation symplectique.  D'après le Corollaire~\ref{c:divexisotopy}, quitte à réaliser des isotopies symplectiques de diviseurs exceptionnels, on peut supposer que les $E_i$ et les $E_j'$ s'intersectent deux à deux localement positivement. Supposons sans perte de généralité qu'on a $k \geq \ell$. S'il existait $i \in \{1, \dots ,k \}$ tel que pour tout $j \in \{1, \dots , \ell \}$, $E_i \cdot E_j' =0$ (c'est-à-dire $E_i$ disjoint de $E_j$ par positivité d'intersection), cela contredirait l'hypothèse de minimalité sur $(N_2, \omega_2)$. Puisque $(N_1, \omega_1)$ et $(N_2, \omega_2)$ ne sont pas dans la même classe de déformation symplectique, on a également $\{ E_1, \dots, E_k \} \neq \{ E_1', \dots, E_\ell' \}$. Il existe donc $i \in \{1, \dots ,k \}$ et $j \in \{1, \dots , \ell \}$ tels que $[E_i] \cdot [E_j'] >0$. L'image de $E_j'$ par la contraction de $E_i$ est alors une courbe symplectique singulière rationnelle $S$ qui vérifie $c_1([S])= c_1([E_j']) + ([E_i] \cdot [E_j']) c_1([E_i]) = 1 +[E_i] \cdot [E_j']  \geq 2$. Le Théorème~\ref{t:McDuffimmersed} permet alors de conclure.
\end{proof}

On termine en présentant les grandes lignes de la preuve du théorème clé de cette section : le Théorème~\ref{t:McDuffimmersed} de McDuff sur les courbes rationnelles symplectiquement positivement immergées. Comme dans la Section \ref{s:GromovMcDuff}, la preuve se base sur des techniques de courbes pseudoholomorphes. \'Etant donné une surface symplectique $(M, \omega)$ et $S$ l'image d'une courbe rationnelle symplectiquement positivement immergée, on dispose d'une structure presque complexe $J$ dominée par $\omega$ telle que $S$ est l'image d'une courbe $J$--holomorphe simple. Cette fois-ci on ne peut pas espérer décrire la topologie de $(M, \omega)$ en construisant un pinceau de Lefschetz dont les fibres sont des courbes rationnelles homologues à $S$. En effet, d'après la formule d'adjonction, aucune telle courbe ne serait plongée. Une autre difficulté provient du fait que certaines des courbes de l'espace de modules $\mathcal{M}_{0}^*([S];J)$ peuvent ne pas être immergées. Ces courbes non immergées ne satisfont pas les hypothèses du théorème de transversalité automatique, ce qui rend la preuve assez technique par endroit (notamment lorsqu'on souhaite obtenir certaines propriétés globales sur l'espace de modules) et pousse à manipuler diverses conditions de généricité avec prudence. Les détails techniques ne seront pas évoqués ici, on renvoie le lecteur ou la lectrice à~\cite[Chapter~7]{Wendl} pour une exposition complète.

On va néanmoins montrer que l'espace de modules des courbes rationnelles simples homologues à $S$ n'est pas compact afin de prouver soit l'existence d'une courbe rationnelle symplectiquement plongée d'auto-intersection positive, soit l'existence d'une courbe rationnelle symplectiquement positivement immergée $S'$ vérifiant $c_1([S]) > c_1([S']) \geq 2$. On procède ensuite de manière récursive jusqu'à trouve la sphère symplectiquement plongée d'auto-intersection positive voulue.

Une solution pour remédier au problème évoqué plus tôt consiste à considérer la construction suivante qu'on appelle la \emph{courbe universelle}. Elle est définie pour tous entiers naturels $m$, $g$, pour toute classe d'homologie $A \in H_2(M; \mathbb{Z})$ et pour toute structure presque complexe $J \in \mathcal{J}_\tau(M, \omega)$, par 
$$\overline{\mathcal{U}}_{g,m}(A;J) = \overline{\mathcal{M}}_{g, m+1}(A;J).$$
La courbe universelle est accompagnée de deux applications lisses naturelles : l'application d'oubli du $(m+1)$--ième point marqué
$$\pi_m : \overline{\mathcal{U}}_{g,m}(A;J) \rightarrow \overline{\mathcal{M}}_{g, m}(A;J)$$
et l'application d'évaluation en le $(m+1)$--ième point marqué
$$\mathrm{ev}_{m+1} : \overline{\mathcal{U}}_{g,m}(A;J) \rightarrow M.$$
Pour tout représentant $(\Sigma ,j,u,(\zeta_1 , \dots , \zeta_m), \Delta)$ d'une courbe de $\overline{\mathcal{M}}_{g, m+1}(A;J)$, on définit l'espace topologique
$\check{\Sigma} = \Sigma  \diagup \sim$ où $\check{z} \sim \hat{z} $  pour toute paire de points nodaux $\{\check{z}, \hat{z} \} \in \Delta$.

\begin{prop}[{\cite[Proposition~7.40]{Wendl}}]\label{p:preimageev}
Soit $(S,j,u,(\zeta_1 , \dots , \zeta_m), \Delta)$ un représentant d'un élément de $\overline{\mathcal{M}}_{g,m} (A ; J)$ avec groupe d'automorphisme $\mathrm{Aut}(u)$. On dispose d'un homéomorphisme naturel entre $\pi_m^{-1} (u)$ et le quotient $\check{\Sigma} / \mathrm{Aut}(u)$ qui identifie l'application $\mathrm{ev}_{m+1}|_{\pi_m^{-1}(u)} : \pi_m^{-1} (u) \rightarrow M$ avec $u : \check{\Sigma} / \mathrm{Aut}(u) \rightarrow M$.
\end{prop}

\`A partir de maintenant, on fixe $m = c_1([S]) -2$ et $p_1, \dots , p_m \in S$ de sorte que l'espace de modules contraint $\overline{\mathcal{M}}_S (J):= \overline{\mathcal{M}}_{0,m}([S];J;p_1, \dots , p_m)$ a pour dimension virtuelle $-2 +2c_1([S]) -2m =2$. Comme $S$ est immergée, elle satisfait alors les hypothèses du théorème de transversalité automatique avec contraintes ponctuelles $p_1, \dots, p_m$. On peut alors supposer sans perte de généralité que $J$ est générique (dans un sens qu'on ne précisera pas ici). On note $\overline{\mathcal{U}}_S (J) = \overline{\mathcal{M}}_{0,m+1}([S];J;p_1, \dots , p_m)$ la courbe universelle correspondante. De la même façon que dans la Section \ref{s:GromovMcDuff}, le fait d'imposer ces contraintes permet d'avoir une meilleur compréhension des courbes nodales pouvant apparaître dans le compactifié de Gromov de l'espace de modules correspondant. En effet, il est possible de montrer que ces courbes nodales sont en nombre fini et que chacune d'entre elle est l'union de deux courbes pseudoholomorphes immergées qui s'intersectent exactement une fois et de manière transverse et positive (voir~\cite[Lemma~7.43]{Wendl}). Cette compréhension des courbes nodales de $\overline{\mathcal{M}}_S (J)$, couplée à la Proposition~\ref{p:preimageev}, permet de montrer le lemme suivant, qui fournit de précieuses informations sur la topologie de $\overline{\mathcal{U}}_S (J)$.

\begin{lem}[{\cite[Theorem~7.45]{Wendl}}] \label{l:fibrationcourbeuniverselle}
Les espaces $\overline{\mathcal{U}}_S (J)$ et $\overline{\mathcal{M}}_S (J)$ peuvent être munis de structures de variétés lisses de dimensions respectives $4$ et $2$ telles que $\pi_m : \overline{\mathcal{U}}_S (J) \rightarrow \overline{\mathcal{M}}_S (J)$ est une fibration de Lefschetz dont les fibres régulières sont de genre $0$ et avec un unique point critique par fibre singulière.
\end{lem}

Les images réciproques des points $p_1, \dots, p_m$ par l'application $\mathrm{ev}_{m+1}$ définissent des sections $\Sigma_1, \dots, \Sigma_N$ deux à deux disjointes de la fibration de Lefschetz donnée par le Lemme~\ref{l:fibrationcourbeuniverselle}.
La stratégie consiste alors à utiliser ces sections prouver l'existence de fibres singulières dans cette fibration de Lefschetz via des contraintes de nature topologique. L'application d'évaluation $\mathrm{ev}_{m+1}$ va jouer un rôle crucial à cet effet.

\begin{lem}[{\cite[Lemma~7.49]{Wendl}}]\label{l:evdiffeolocal}
Soit $u \in \mathcal{M}_S (J)$ une courbe immergée. Alors sur un voisinage du complémentaire des points marqués dans $\pi_m^{-1}(u)$, l'application $\mathrm{ev}_{m+1} : \overline{\mathcal{U}}_S (J) \rightarrow M$ est un difféomorphisme local préservant l'orientation.
\end{lem}

La démonstration du Lemme~\ref{l:evdiffeolocal} repose essentiellement sur le théorème de transversalité automatique avec $m+1 = c_1 ([S])-1$ contraintes. On ne peut donc pas espérer étendre aussi simplement ces résultats aux courbes de genre strictement positif. Le Lemme~\ref{l:evdiffeolocal} permet de calculer le degré de l'application d'évaluation  $\mathrm{ev}_{m+1}$.

\begin{lem}[{\cite[Lemma~7.51]{Wendl}}]
L'application $\mathrm{ev}_{m+1} : \overline{\mathcal{U}}_S (J) \rightarrow M$ vérifie $\deg ( \mathrm{ev}_{m+1}) \geq 1$, avec égalité si et seulement si $S$ est plongée.
\end{lem}

On note $N = \deg ( \mathrm{ev}_{m+1})$ et, si $m \geq 1$, on introduit pour tout $i \in \{1, \dots , m \}$, l'application 
$$\sigma_i : \overline{\mathcal{M}}_S (J) \rightarrow \mathbb{P}(T_{p_i} M) \cong \mathbb{C}P^1,$$
qui à une courbe $u$ associe l'espace tangent à $u$ en $p_i$. Le degré de $\mathrm{ev}_{m+1}$ permet de calculer le degré de $\sigma_i$.


\begin{lem}[{\cite[Lemma~7.54]{Wendl}}]
Pour tout $i \in \{1, \dots , m \}$, $\deg ( \sigma_i)= N$.
\end{lem}

Le degré des applications $\sigma_i$ permettent le calcul des nombres d'auto-intersection des sections $\Sigma_i$.

\begin{lem}[{\cite[Lemma~7.55]{Wendl}}]\label{l:autointsection}
Les sections $\Sigma_i \subset \overline{\mathcal{U}}_S (J)$ sont deux à deux disjointes et vérifient $[\Sigma_i]^2 = - N$.
\end{lem}

D'après le Lemme~\ref{l:autointsection}, les classes d'homologies des sections $\Sigma_1, \dots, \Sigma_m$ et la classe d'homologie d'une fibre régulière $F$ de la fibration de Lefschetz forment une famille libre de $H_2 (\overline{\mathcal{U}}_S (J);\mathbb{Z})$, donc $b_2 (\overline{\mathcal{U}}_S (J)) \geq m+1$. Si $m > 1$ (c'est-à-dire $c_1([S]) > 3$), la Proposition~\ref{p:topfibrationLefschetz} assure l'existence d'une fibre singulière. Cette fibre singulière assure en retour l'existence d'une courbe nodale dans $\overline{\mathcal{M}}_S (J)$, constituée de deux composantes symplectiquement positivement immergées $S'$ et $S''$ qui vérifient $[S]=[S']+[S'']$. Comme $S'$ et $S''$ sont des images de courbes $J$--holomorphes simples Fredholm régulières, la condition de positivité des indices assure que $c_1([S']) \geq 1$ et $c_1([S'']) \geq 1$. Quitte à renommer les composantes, obtient ainsi l'existence d'une autre sphère symplectiquement positivement immergée $S'$ qui vérifie $ c_1 ([S]) > c_1([S']) \geq 2$. En procédant de manière itérative, on peut donc se ramener au cas où $c_1([S]) \in \{2,3 \}$.

On suppose désormais $c_1([S]) \in \{2,3 \}$. La dernière étape consiste à montrer l'existence d'une sphère symplectiquement plongée d'auto-intersection positive parmi les composantes des courbes nodales de $\overline{\mathcal{M}}_S (J)$. Les outils utilisés sont les mêmes que pour le reste de la preuve. Néanmoins cette partie requiert plus d'astuce et nécessite potentiellement de changer l'immersion de départ. Cette partie de la preuve n'éclairant pas les raisonnements abordés dans le chapitre suivant, on renvoie à~\cite[Section~7.3.7]{Wendl} pour plus détails.

\begin{rk}
Si $m \geq 1$ (c'est-à-dire si $c_1([S]) \geq 3$), alors on peut voir \emph{a posteriori} que $(M, \omega)$  est une surface symplectique rationnelle. En effet, $H_1(\overline{\mathcal{U}}_S (J); \mathbb{Z})$ est engendré par les lacets sur $\Sigma_1$, et tous ces lacets sont envoyés sur un point par $\mathrm{ev}_{m+1}$. Puisque $\mathrm{ev}_{m+1}$ induit une application surjective en homologie rationnelle, on obtient $H_1(M; \mathbb{Q}) =0$. Les Théorèmes~\ref{t:classificationsurfacesreglees} et~\ref{t:classificationsurfacesrationnelles}, couplés à la Proposition~\ref{p:topfibrationLefschetz}, permettent de conclure.
\end{rk}

Lorsqu'on tente de généraliser la stratégie consistant à étudier la courbe universelle aux courbes de genre $g>0$, on se heurte à quelques difficultés.

Tout d'abord, l'application d'oubli $\pi_m : \overline{\mathcal{U}}_S (J) \rightarrow \overline{\mathcal{M}}_S (J)$ définit une fibration de Lefschetz sur $\overline{\mathcal{U}}_S (J)$, mais la topologie de tels objets (quand la fibre est de genre strictement positif) est encore assez mal comprise. Cette difficulté pourrait éventuellement être contournée en ayant plus d'informations sur la topologie de la base de la fibration de Lefschetz, par exemple en montrant que $\overline{\mathcal{M}}_S (J)$ est difféomorphe à une sphère. Un résultat allant dans ce sens pour les tores symplectiquement plongées d'auto-intersection $1$ dans certaines surfaces symplectiques est discuté dans la Section~\ref{s:isotopiebirationnelle} (voir la Remarque~\ref{r:moduledetorespherique}). Ce problème constitue une motivation supplémentaire pour étudier la topologie de tels espaces de modules.

Ensuite, le fait que l'espace de modules $\mathcal{M}_S (J)$ soit non vide et naturellement muni d'une structure de variété lisse est une conséquence du théorème de transversalité automatique avec $m$ contraintes ponctuelles. De même, le Lemme~\ref{l:evdiffeolocal} (qui sert à montrer $\deg (\mathrm{ev}_{m+1}) >0$) repose de manière cruciale sur le fait que les courbes de $\overline{\mathcal{M}}_S (J)$ sont automatiquement Fredolm régulières pour le problème avec $m+1$ contraintes ponctuelles. Or dans le cas où $S$ est une courbe de genre $g$, on doit poser $m = c_1([S])+g-2$ pour que la dimension virtuelle de l'espace de modules $\overline{\mathcal{M}}_S (J)$ soit égale à $2$. L'hypothèse du Théorème~\ref{t:transautoSiko} avec $m+1$ contraintes ponctuelles n'est alors plus respectée dès que $g \geq 1$ puisque la condition $2g-2 + 2c_1([S])-2(m+1) > 2g-2$ équivaut à $c_1 ([S]) >m+1$.


\clearemptydoublepage

\chapter{Courbes symplectiques irrationnelles de haute auto-intersection dans les surfaces symplectiques}
On discute maintenant des effets de la présence de courbes symplectiquement plongées et d'auto-intersection grande par rapport à leur genre sur la topologie des surfaces symplectiques qui les contiennent. Le cas des courbes rationnelles ayant été traité dans le Chapitre 2, on mettra l'emphase sur l'étude des courbes irrationnelles. Ce chapitre est en grande partie consacré à la démonstration du théorème suivant.

\begin{thm} \label{t:thmprincipal}
Soit $(M, \omega )$ une surface symplectique et $S$ une courbe symplectiquement plongée dans $(M, \omega)$ de genre $g$ telle que $(M, \w)$ est relativement minimale par rapport à $S$. Si $[S]^2 > 4g +5$, alors
la surface symplectique $(M,\omega)$ est un fibré en sphères symplectique au-dessus d'une surface de genre $g$ et $S$ est l'image d'une section de ce fibré.
De plus, l'énoncé est également vrai si on suppose seulement que $[S]^2 \geq 4g +5$, à l'unique exception près des cubiques symplectiques non singulières dans le plan projectif complexe.
\end{thm}

Remarquons tout d'abord que le Théorème~\ref{t:thmprincipal} permet de déterminer complètement le type de difféomorphisme de la variété symplectique $(M, \omega)$. 

\begin{cor}\label{c:corprincipal}
Soit $(M, \omega )$ une surface symplectique et $S$ une courbe symplectiquement plongée dans $(M, \omega)$ de genre $g$ telle que $(M, \w)$ est relativement minimale par rapport à $S$. Si $[S]^2 > 4g +5$ (ou $[S]^2 \geq 4g + 5$ si $g \neq 1$), alors :
\begin{itemize}
\item si $[S]^2$ est pair, $M$ est difféomorphe à $\Sigma_g \times \mathcal{S}^2$ ;
\item si $[S]^2$ est impair, $M$ est difféomorphe à $\Sigma_g \, \tilde \times \, \mathcal{S}^2$ ;
\end{itemize}
où $\Sigma_g$ désigne une surface de genre $g$.
\end{cor}

\begin{proof}
D'après le Théorème~\ref{t:thmprincipal}, $M$ est difféomorphe à un fibré en sphères sur $\Sigma_g$. On a vu dans la Section~\ref{s:ssr,ssr} qu'il n'y a que deux types de difféomorphismes possibles pour de tels objets : $\Sigma_g \times \mathcal{S}^2$ et $\Sigma_g \, \tilde \times \, \mathcal{S}^2$. D'après la Proposition~\ref{p:topfibrationLefschetz}, le groupe $H_2(M;\mathbb{Z})$ est engendré par $[S]$ et la classe de la fibre $[F]$. Si $[S]^2$ est impair, alors la forme d'intersection $\mathcal{Q}_M$ est impaire. On en déduit dans ce cas que $M$ est difféomorphe à $\Sigma_g \, \tilde \times \, \mathcal{S}^2$. Si $[S]^2$ est pair, alors pour tous entier $k,l$, l'entier $(k[S]+\ell [F])^2 = k^2[S]^2 +2k \ell$ est pair, donc la forme d'intersection $\mathcal{Q}_M$ est paire. Dans ce cas, $M$ est difféomorphe à $\Sigma_g \times \mathcal{S}^2$.
\end{proof}

\begin{rk}\label{r:déformationsymp}
Au vu de la construction présentée dans la Section~\ref{s:ssr,ssr}, on constate que pour deux fibrés en sphères symplectiques $(M_1, \omega_1)$ et $(M_2, \omega_2)$ au-dessus de surfaces tels que $M_1$ et $M_2$ sont difféomorphes, on peut choisir un difféomorphisme $\Phi : M_1 \rightarrow M_2$ de façon à ce que $\Phi$ préserve les fibres des deux fibrations. Le Théorème~\ref{t:GompfThurston} assure alors que les formes symplectiques $\omega_1$ et $\Phi^* \omega_2$ sont reliées par un chemin de formes symplectiques compatibles avec la structure de fibré en sphères sur $M_1$ (notons également que d'après le Théorème~\ref{t:isotopiepinceaudeLefschetz}, les déformations symplectiques préservent les structures de fibrés en sphères symplectiques, considérées à isotopie près). Autrement dit, le type de difféomorphisme d'un fibré en sphères symplectique $(M, \omega)$ au-dessus d'une surface détermine complètement sa classe de déformation symplectique (voir aussi~\cite[Example~3.6]{Sal12}). Le Corollaire~\ref{c:corprincipal} permet donc de déterminer complètement la classe de déformation symplectique de $(M, \omega)$.
\end{rk}

\begin{rk}
Avec les hypothèses du Théorème~\ref{t:thmprincipal}, dans le cas où $g>0$, le fait que $(M, \omega)$ soit relativement minimale par rapport à $S$ implique la minimalité de $(M, \omega)$ (ce qui n'est pas nécessairement vrai dans le cas où $g=0$).
\end{rk}

Notons ensuite qu'étant donné une surface symplectique $(M, \omega)$ et $S$ une courbe symplectiquement plongée dans $(M, \omega)$, on peut aisément ramener au cas où $(M, \omega)$ est relativement minimale par rapport à $S$. En effet, il suffit pour cela de contracter une famille maximale de diviseurs exceptionnels symplectiques disjoints de $S$ et deux à deux disjoints. On peut alors appliquer le Théorème~\ref{t:thmprincipal}, puis éclater les points correspondant aux diviseurs exceptionnels précédemment contractés pour obtenir de nouveau $(M, \omega)$. Quitte à perturber légèrement $S$ et la structure de fibré en sphères symplectique, on peut supposer que ces éclatements sont disjoints de $S$ et qu'il y a au plus un éclatement par fibre. On obtient alors le corollaire suivant.

\begin{cor}
Soit $(M, \omega )$ une surface symplectique et $S$ une courbe symplectiquement plongée dans $(M, \omega)$ de genre $g$. Si $[S]^2 > 4g +5$, alors la surface symplectique $(M,\omega)$ est munie d'une fibration de Lefschetz symplectique générique (i.e avec au plus un point singulier par fibre) dont la base est une surface de genre $g$, les fibres non singulières sont des sphères et telle que $S$ est l'image d'une section.
De plus, le théorème est également vrai si $[S]^2 \geq 4g +5$ et $g\neq 1$ à l'exception près des transformées propres des cubiques symplectiques non singulières dans les éclatés du plan projectif complexe.
\end{cor} 

Le Théorème~\ref{t:thmprincipal} est inspiré de son analogue algébrique (vérifié pour les corps algébriquement clos), démontré par Hartshorne dans~\cite{hartshorne}. Le cas se rapprochant le plus du cadre de ce manuscrit étant le cas complexe, on donne une présentation très succincte de ce théorème de Hartshorne du point de vue complexe dans la Section~\ref{s:hartshorne}. On présente ensuite deux démonstrations du Théorème~\ref{t:thmprincipal}. 
La première, esquissée dans la Section~\ref{s:marteaupilonSW}, est plus \og concise \fg{} (pour nuancer, elle nécessite quand même de reprendre la totalité des arguments utilisés par Hartshorne) mais possède l'inconvénient de faire appel à la théorie de Seiberg--Witten. Elle est par conséquent moins visuelle et géométrique.
La seconde démonstration, qui se base sur les propriétés des espaces de modules de courbes pseudoholomorphes, est détaillée dans la Section~\ref{s:démo}. Elle a l'avantage d'être plus directe, plus éclairante de faire intervenir des techniques intéressantes concernant les courbes pseudoholomorphes de genre strictement positif. Elle permet également de prendre un raccourci non négligeable dans les arguments utilisés par Hartshorne. 
On discute dans la Section~\ref{s:pistes} de l'optimalité des bornes du Théorème~\ref{t:thmprincipal}, ainsi que d'éventuelles pistes pour aller plus loin.

\section{Courbes de haute auto-intersection dans les surfaces complexes : un théorème de Hartshorne}\label{s:hartshorne}

Le Théorème~\ref{t:thmprincipal} s'inspire d'un théorème de Hartshorne dont on donne ici une présentation très succincte. Afin de rester au plus proche du cadre de ce manuscrit, on n'en présente ici qu'un cas particulier : celui de la géométrie algébrique complexe, qui permet un point de vue plus topologique que le cadre général donné par la géométrie algébrique sur les corps finis.

Commençons par quelques définitions préliminaires. 
Une \emph{surface réglée} est une surface birationnellement équivalente à un produit de la forme $\mathbb{C}P^1 \times C$, où $C$ désigne une courbe complexe non singulière. Une \emph{surface rationnelle} est une surface birationnellement équivalente $\mathbb{C}P^1 \times \mathbb{C}P^1$. Une \emph{surface géométriquement réglée} est un fibré holomorphe en $\mathbb{C}P^1$ au-dessus d'une courbe complexe non singulière $C$. \'Etant donné deux surfaces complexes $V_1$ et $V_2$, ainsi qu'une courbe complexe $X$, on dit que deux plongements $X \rightarrow V_1$ et $X \rightarrow V_2$ sont \emph{équivalents} s'il existe une application birationnelle $f : V_1 \rightarrow V_2$ qui est un isomorphisme sur un voisinage ouvert de $X$ dans $V_1$ et qui induit l'application identité sur $X$.

On peut maintenant énoncer le théorème de Hartshorne en question.

\begin{thm}[{\cite[Theorem~4.1]{hartshorne}}]\label{t:hartshorne}
Soit $X$ une courbe algébrique complexe non singulière de genre $g$ dans une surface algébrique complexe non singulière $V$. On suppose que $X^2 > 4g +5$. Alors $V$ est une surface réglée et le plongement de $X$ dans $V$ est équivalent à une section d'une surface géométriquement réglée. De plus, l'énoncé est également vrai si $X^2 \geq 4g +5$, à l'exception de la cubique non singulière de $\mathbb{C}P^2$ (ou d'un plongement équivalent).
\end{thm}

Soit $V$ une surface complexe non singulière. On désigne par $K_V$ la classe canonique de $V$, définie dans ce cadre comme la première classe de Chern du fibré en droites complexe des $2$--formes holomorphes sur $V$ (en particulier on a $K_V= - c_1(TV)$). La démonstration du Théorème~\ref{t:hartshorne} consiste essentiellement en une utilisation astucieuse des outils suivants :
\begin{itemize}
\item  Le critère d'Enriques pour les surfaces réglées (voir~\cite[Corollary VI.18]{beauville1996complex}). Celui-ci affirme qu'une surface complexe non singulière $V$ est réglée si et seulement si elle contient une courbe complexe irréductible $X$ qui n'est pas un diviseur exceptionnel et telle que $K_V \cdot X <0 $ (notons que c'est aussi équivalent au fait que la dimension de Kodaira de $V$ est égale à $-\infty$).
\item La formule d'adjonction, qui affirme que pour toute courbe complexe irréductible dans $V$, on a $K_V \cdot X + 2 - 2 p_a(X) + X^2 =0$, où $p_a(X)$ désigne le genre arithmétique de $X$.
\item La propriété de positivité d'intersection entre deux courbes complexes dans une surface complexe, qui affirme que deux telles courbes sans composante commune s'intersectent localement positivement.
\item Les transformations birationnelles, qui peuvent être décrites comme des successions d'éclatements, d'isomorphismes et de contractions (voir la Remarque~\ref{r:birationalequiv}).
\end{itemize}

Commençons par présenter brièvement l'idée de la preuve. Comme $X$ est non singulière, on a $p_a(X) =g$. Puisque $X^2 > 4g +5$, la formule d'adjonction nous indique alors que $K_V \cdot X < 0$. D'après la classification d'Enriques pour les surfaces réglées, on sait donc que $V$ est une surface réglée.
On étudie ensuite les courbes irréductibles $Y$, possiblement singulières, dans un modèle minimal de $V'$ (les modèles minimaux des surfaces réglées sont les surfaces géométriquement réglées et le plan projectif complexe $\mathbb{C}P^2$, voir~\cite[Theorem~2.1, Theorem~3.3]{hartshorne}) de façon à obtenir une borne supérieure sur le nombre d'auto-intersection $Y^2$ en fonction de $p_a(Y)$. Cette étape est obtenue en utilisant judicieusement la formule d'adjonction, certains nombres d'intersection entre différentes courbes, ainsi que la propriété de positivité d'intersection. Les cas des surfaces réglées non rationnelles et des surfaces rationnelles sont traités à part (ce dernier cas étant plus astucieux). 
La dernière étape consiste alors à examiner ce qu'il se passe au niveau de la borne quand on éclate des points de $V'$ pour obtenir de nouveau la surface complexe $V$.

L'objectif des deux démonstrations du Théorème~\ref{t:thmprincipal} (présentées respectivement dans la Section~\ref{s:marteaupilonSW} et la Section~\ref{s:démo}) est de se ramener à une situation dans laquelle on peut adapter au cadre symplectique les arguments et les outils utilisés par Hartshorne. On adaptera ces arguments de la manière suivante :

\begin{itemize}
\item Dans l'esquisse de preuve présentée dans la Section~\ref{s:marteaupilonSW}, un analogue du critère d'Enriques sera donné par le Théorème~\ref{t:TaubesSW}, issu de la théorie de Seiberg--Witten. Ce rôle sera joué par le Théorème~\ref{t:ruledFKversion2} dans la Section~\ref{s:demo}, démontré grâce à l'étude de certains espaces de modules de courbes pseuholomorphes.

\item La formule d'adjonction s'utilise sans problème dans le cadre symplectique. Quand on compare la formule énoncée plus haut dans le cas complexe avec la formule d'adjonction dans le cadre symplectique (voir le Théorème~\ref{t:adjformula}), appliquée à une courbe symplectique singulière $S$ de genre $g$ dans une surface symplectique : $c_1([S]) = 2-2g + [S]^2 -2 \delta(S)$, le terme $K_V \cdot X$ correspond à $- c_1([S])$, le terme $p_a(X)$ correspond à $g + \delta(S)$ et le terme $X^2$  correspond à $[S]^2$.

\item Pour la propriété de positivité d'intersection, on veut se retrouver dans une situation où la courbe symplectique de départ intersecte localement positivement toutes les fibres d'une fibration (ou pinceau) de Lefschetz sur la surface symplectiquement réglée ambiante qu'on considère. Ceci sera effectué via l'utilisation de structures presque complexes auxiliaires et de courbes pseudoholomorphes.

\item Les explications données dans la Section~\ref{s:eclat} décrivent les similarités des opérations d'éclatements et de contractions en symplectique ou en complexe.
\end{itemize}

\section{Esquisse de preuve du Théorème~\ref{t:thmprincipal} à l'aide de la théorie de Seiberg--Witten}\label{s:marteaupilonSW}

Pour cette preuve, on utilise le théorème suivant, obtenu à l'aide de la théorie de Seiberg--Witten. Les grandes lignes de sa démonstration sont exposées dans~\cite[Section~7.3.1]{Wendl}.

\begin{thm}[\cite{Taubes95,Taubes96a,Taubes96b}]\label{t:TaubesSW}
Soit $(M, \omega )$ une surface symplectique et $S$ une courbe symplectiquement plongée dans $M$ de genre $g$. Si $[S]^2 > 2g-2$ et $S$ n'est pas un diviseur exceptionnel, alors $(M,\omega)$ contient une courbe rationnelle symplectiquement plongée d'auto-intersection positive. 
\end{thm}

\begin{rk}
Le Théorème~\ref{t:TaubesSW}, combiné au Théorème~\ref{t:McDuff} de McDuff et Gromov est à rapprocher, dans la classification des surfaces algébriques complexes de Enriques-Kodaira, du critère évoqué dans la Section~\ref{s:hartshorne} permettant de caractériser les surfaces complexes réglées. En effet, d'après la formule d'adjonction, on a $c_1([S]) = 2-2g +[S]^2$. La condition $[S]^2 >2g-2$ est donc équivalente à $c_1([S]) >0$.
\end{rk}

Soit $(M, \omega )$ une surface symplectique et $S$ une courbe  symplectiquement plongée dans $M$ de genre $g$. On suppose que $[S]^2 > 4g +5$. On utilise tout d'abord le Théorème~\ref{t:TaubesSW} pour trouver une courbe rationnelle symplectiquement plongée $S_0$ dans $(M, \omega)$ d'auto-intersection positive, puis on se ramène au cas où $S_0$ est d'auto-intersection $0$ ou $1$ grâce au Théorème~\ref{t:McDuff}. On sait alors à ce stade que $(M, \omega)$ est une surface symplectiquement réglée. Plus particulièrement, si $[S_0]^2 = 0$, on sait que $(M, \omega)$ est un éclaté d'un fibré en sphères symplectique dont $S_0$ est une fibre et si $[S_0]^2 = 1$, on sait que $(M, \omega)$ est un éclaté de $\mathbb{C}P^2$.

Le but est maintenant de réaliser une isotopie symplectique pour trouver, à partir de $S_0$, une courbe rationnelle symplectiquement plongée homologue à $S_0$ qui intersecte $S$ localement positivement. Pour ce faire, on utilise des techniques de courbes pseudoholomorphes. On commence par choisir une structure presque complexe $J$ dominée par $\omega$ telle que $S$ est l'image d'une courbe pseudoholomorphe $u : \Sigma_g \rightarrow M$ plongée dans $M$. Grâce à la formule d'adjonction pour les courbes symplectiques, on a 
$$\ind (u) = -\chi (\Sigma_g) + 2 c_1 ( [S]) = \chi (\Sigma_g) + 2 [S]^2 > 6g + 12 > 2g -2.$$ 
D'après le Théorème~\ref{t:transauto1}, la courbe $u$ est donc automatiquement Fredholm régulière. Quitte à réaliser des perturbations $\mathcal{C}^\infty$ arbitrairement petites de $J$ et $S$, on peut supposer sans perte de généralité que $J$ est une structure presque complexe générique dominée par $\omega$, de sorte que toutes les courbes $J$--holomorphes injectives quelque part sont Fredholm régulières. En utilisant le Théorème~\ref{t:pinceaudeLefschetz}, on construit ensuite un pinceau de Lefschetz symplectique générique sur $M$ (possédant $[S_0]^2$ points bases) dont $S_0$ est une fibre lisse. Puis en appliquant le Théorème~\ref{t:isotopiepinceaudeLefschetz}, on obtient que ce pinceau de Lefschetz est isotope à un pinceau de Lefschetz générique dont toutes les composantes des fibres sont des courbes rationnelles $J$--holomorphes plongées. En particulier, $S$ intersecte localement positivement chacune des fibres.

On peut alors appliquer des arguments similaires à ceux utilisés dans~\cite{hartshorne}, en les adaptant au cadre symplectique, pour ainsi démontrer le résultat voulu. Une partie de ces arguments sera présentée dans la section suivante, plus précisément dans la deuxième moitié de la Sous-section~\ref{ss:plongementISsection}, on ne s'y attarde donc pas plus ici.

\section{Démonstration du Théorème~\ref{t:thmprincipal} basée sur des techniques pseudoholomorphes}\label{s:démo}

Soit $(M, \omega )$ une surface symplectique et $S$ une courbe symplectiquement plongée dans $(M, \omega)$ de genre $g$. On suppose $[S]^2 > 4g +5$ et $(M,\w)$ relativement minimale par rapport à $S$. On procède de la manière suivante. On commence par choisir une structure presque complexe $J$ dominée par $\omega$ telle que $S$ est l'image d'une courbe pseudoholomorphe $u : \Sigma \rightarrow M$ plongée dans $M$. Grâce à la formule d'adjonction pour les courbes symplectiques, on a 
$$\ind (u) = -\chi (\Sigma) + 2 c_1 ( [S]) = \chi (\Sigma) + 2 [S]^2 > 6g + 8 > 2g -2.$$
D'après le Théorème~\ref{t:transauto1}, la courbe $u$ est automatiquement Fredholm régulière. On peut alors perturber $J$ en une structure presque complexe $J'$ générique, de sorte que toutes les courbes $J'$--holomorphes injectives quelque part soient Fredholm régulières, et trouver une courbe $J'$--holomorphe $u' : \Sigma' \rightarrow M$ $ \mathcal{C}^1$--proche de $u$. On peut ensuite trouver une isotopie ambiante $(\Phi_t)_{t \in [0,1]}$ telle que $\Phi_0 = \mathrm{id}_M$ et $\Phi_1(u(\Sigma)) = u'(\Sigma')$. La structure presque complexe $J''=\Phi_1^*J'$ est à la fois générique et rend $J''$--holomorphe la courbe symplectique $S$. Par conséquent, on peut supposer sans perte de généralité que $J$ est générique.

On commence dans un premier temps par démontrer le théorème suivant.

\begin{thm}\label{t:ruledFKversion2}
Soit $S$ une courbe de genre $g$ symplectiquement plongée dans une surface symplectique $(M, \omega )$. Si $[S]^2 \geq 4g$, alors $(M, \omega )$ est une surface symplectiquement réglée.
\end{thm}

\begin{rk}\label{r:comparaisondeTh}
Le Théorème~\ref{t:ruledFKversion2} est à comparer avec le Théorème~\ref{t:TaubesSW} de Taubes, issu de la théorie de Seiberg--Witten. Les hypothèses du Théorème~\ref{t:ruledFKversion2} sont plus fortes mais la démonstration ne faisant appel qu'aux propriétés des courbes pseudoholomorphes, elle est à la fois plus géométrique et élémentaire. Dans la Section~\ref{s:pistes}, on discute d'éventuelles stratégies pour essayer d'améliorer la borne dans les hypothèses.
\end{rk}

Le Théorème~\ref{t:ruledFKversion2} sera une conséquence directe du Théorème~\ref{t:ruledFKversion} et du Théorème~\ref{t:McDuff} de McDuff et Gromov. Sa démonstration fera l'objet des Sous-sections~\ref{ss:smash}, \ref{ss:glue} et~\ref{ss:recurrence}. On commence par montrer dans la Sous-section~\ref{ss:smash} que l'espace de modules des courbes $J$--holomorphes simples de genre $g$ homologue à $S$ n'est pas compact. Dans la Sous-section~\ref{ss:glue}, on trouve grâce aux composantes des courbes nodales dans le compactifié de l'espace de modules, une courbe symplectiquement plongée de genre strictement plus petit que $g$, mais avec auto-intersection suffisamment grande par rapport à son genre. Enfin, on termine le raisonnement dans la Sous-section~\ref{ss:recurrence} en effectuant un récurrence sur le genre. \emph{Dans la démonstration et dans toute la suite de cette partie, on suppose sans perte de généralité que $(M, \omega)$ est relativement minimale par rapport à $S$.}

Dans la Sous-section~\ref{ss:plongementISsection}, on montre la deuxième partie du Théorème~\ref{t:thmprincipal}, c'est-à-dire que le plongement de la courbe de départ est une section d'un fibré en sphères symplectique sur la surface de départ.

\subsection{Trouver des courbes nodales}\label{ss:smash} 

\subsubsection{Surjectivité de l'application d'évaluation}

\begin{lem} \label{l:surjectev}
Soit $J \in \mathcal{J}^{\mathrm{reg}}_{\tau} (M, \omega)$, $u : \Sigma \rightarrow M$ une courbe $J$--holomorphe plongée de genre $g$ et $m$ un entier naturel. On suppose que $\ind (u) -2m > 2g -2$. Si l'espace de modules $\mathcal{M}_{g}^* ([u];J)$ est compact, alors l'application $\ev : \mathcal{M}_{g,m}^* ([u];J) \rightarrow M^m \backslash \Delta$ est surjective, où $\Delta$ désigne la grande diagonale (l'ensemble des $m$--uplets qui possèdent au moins deux composantes identiques).
\end{lem}

\begin{rk}
Les espaces de modules de la forme $\mathcal{M}_{g,m}^* ([u];J)$ ne sont pas compact dès que $m>1$. En effet, une suite de courbes convergente dont au moins deux points marqués convergent vers un même point $\zeta$ tend vers une courbe nodale possédant une bulle fantôme en $\zeta$, sur laquelle se trouve au moins deux points marqués. En revanche ce fait ne nous est pas d'une grande aide, car ce genre de courbes nodales ne nous permet pas de trouver d'autres courbes de genre strictement plus petit.
\end{rk}

\begin{proof}
On montre que l'image de $\ev$, qui est non vide, est à la fois ouverte et fermée dans $M^m \backslash \Delta$, ce qui permet de conclure par connexité de $M^m \backslash \Delta$.

Soit $\left( (p_1^{(n)}, \dots, p_m^{(n)}) \right)$ une suite de points de $\ev \left( \mathcal{M}_{g,m}^* ([u];J) \right)$ qui converge vers un point $(p_1, \dots, p_m) \in  M^m \backslash \Delta$. On dispose d'une suite de courbes $J$--holomorphes $\left((\Sigma_n, j_n,u_n,(\zeta_1^{(n)}, \dots, \zeta_m^{(n)})) \right)$ dans $\mathcal{M}_{g,m}^* ([u];J)$ telle que pour tout entier naturel $n$ et tout $i \in \{1, \dots, m \}$, on a $u_n \left(\zeta_i^{(n)}\right)=p_i^{(n)}$. On peut supposer, quitte à extraire, que la suite $(u_n)$ converge vers une courbe $u_\infty \in \mathcal{M}_{g}^* ([u];J)$ (qui est plongée, car $u$ est plongée) puisque $\mathcal{M}_{g}^* ([u];J)$ est compact. La courbe $u_\infty$ ainsi obtenue passe par $p_1, \dots, p_m$, et les pré-images $\zeta_i$ des points $p_i$ par la courbe $u_\infty$ permettent d'obtenir un antécédent $(\Sigma_\infty, j_\infty , u_\infty ,(\zeta_1, \dots, \zeta_m))$ de $(p_1, \dots, p_m)$ par $\ev$. Donc l'ensemble $\ev\left( \mathcal{M}_{g,m}^* ([u];J) \right)$ est fermé dans $M^m \backslash \Delta$.

Puisque $\ind (u) -2m > 2g -2$, le Théorème~\ref{t:transautoSiko} nous assure que toutes les courbes de l'espace de modules $\mathcal{M}_{g,m}^* ([u];J)$ sont automatiquement Fredholm régulières pour le problème avec $m$ contraintes ponctuelles. Ceci signifie précisément que l'application $\ev$ est transverse à tout point de $M^m \backslash \Delta$, c'est donc une submersion. Puisqu'une submersion est une application ouverte, l'ensemble $\ev\left( \mathcal{M}_{g,m}^* ([u];J) \right)$ est ouvert dans $M^m \backslash \Delta$.
\end{proof}

\begin{rk}\label{r:surjCC}
La démonstration du Lemme~\ref{l:surjectev} montre en fait plus généralement que toute composante connexe de l'espace de modules $\mathcal{M}_{g,m}^* ([u];J)$ s'évalue surjectivement sur $M \backslash \Delta$ (il suffit pour cela de raisonner indépendamment sur chacune des composantes connexes).
\end{rk}

\begin{rk}
La borne sur l'indice dans le Lemme~\ref{l:surjectev} est optimale. En effet, les droites dans $\mathbb{C}P^2$ sont d'indice $4$ et aucune droite ne passe par trois points en position générale.
\end{rk}

\begin{prop}\label{p:casserlescourbes}
Soit $J \in \mathcal{J}^{\mathrm{reg}}_{\tau} (M, \omega)$ et $u : \Sigma \rightarrow M$ une courbe $J$--holomorphe plongée de genre $g$. Si $\ind (u) > 2g +4$, alors l'espace de modules $\mathcal{M}_g^* ([u];J)$ n'est pas compact.
\end{prop}

Pour démontrer la Proposition~\ref{p:casserlescourbes}, on se sert du Lemme~\ref{l:surjectev} afin d'obtenir une application d'évaluation surjective pour trois points marqués. On se sert ensuite de cette application d'évaluation pour forcer l'existence d'une courbe nodale dans le compactifié de l'espace de modules considéré (voir~\cite[Theorem~9.4.1]{mcduff2004j} pour un argument similaire dans le cas des courbes rationnelles).

\begin{proof}
On suppose que $\mathcal{M}_g^* ([u];J)$ est compact. D'après le Lemme~\ref{l:surjectev}, l'application $\ev : \mathcal{M}_{g,3}^* ([u];J) \rightarrow M^3 \backslash \Delta$ est surjective. On choisit alors $p$ un point dans $M$, une application exponentielle $\mathrm{exp}_p$ d'un voisinage de $0$ dans $T_p M$ vers un voisinage de $p$ dans $M$ et deux vecteurs $v,w \in T_p M$ dans un voisinage de $0$ qui sont $J_p$--linéairement indépendants. On pose alors, pour tout entier naturel $n$, $q_n = \mathrm{exp}_p\left( 2^{-n}v \right)$ et $r_n =  \mathrm{exp_p}\left( 2^{-n}w \right)$. On construit ainsi deux suites points de $M\backslash \{p\}$ qui convergent vers $p$ dans des directions $J_p$--linéairement indépendantes. Par surjectivité de l'application d'évaluation $\ev : \mathcal{M}_{g,3}^* ([u];J) \rightarrow M^3 \backslash \Delta$, il existe alors une suite de courbes $J$--holomorphes $\left( (\S_n, j_n, u_n) \right)$ dans $\mathcal{M}_g^* ([u];J)$ telle que pour tout entier naturel $n$, la courbe $u_n$ passe par les points $r_n$, $q_n$ et $p$ (on note que dans cette étape, on oublie les points marqués). Par compacité de $\mathcal{M}_g^* ([u];J)$, on peut supposer, quitte à extraire, que la suite $(u_n)$ converge vers une courbe $(\S_\infty, j_\infty, u_\infty) \in \mathcal{M}_g^* ([u];J)$ qui passe nécessairement par $p$. On rappelle que puisque $u$ est plongée, la formule d'adjonction garantit que toutes les courbes dans l'espace de modules $\mathcal{M}_g^* ([u];J)$ sont également plongées.

Par définition de la convergence des courbes $J$--holomorphes, pour tout entier naturel $n$, on dispose d'un biholomorphisme $\varphi_n : (\S_\infty, j_n') \rightarrow (\S_n, j_n)$ tel que les suites $(u_n \circ \varphi_n)$ et $(j_n')$ convergent respectivement vers $u_\infty$ et $j_\infty$ pour la topologie $\mathcal{C}^\infty$. On note alors, pour tout entier naturel $n$, $\zeta^q_n = (u_n \circ \varphi_n)^{-1} (q_n)$ et $\zeta^r_n = (u_n \circ \varphi_n)^{-1} (r_n)$. Comme $\S$ est compact, quitte à extraire, on peut supposer que les suites de points marqués $(\zeta^q_n)$ et $(\zeta^r_n)$ convergent respectivement vers des points marqués $\zeta^q$ et $\zeta^r$ sur $\S_\infty$. On a alors $u_\infty (\zeta^q) = \lim u_n \circ \varphi_n (\zeta^q_n) = \lim q_n =p$, la première égalité provenant de la continuité de l'application d'évaluation en un point marqué $(u,\zeta) \mapsto u(\zeta)$. De même on a $u_\infty (\zeta^r) = p$. Puisque $u_\infty$ est plongée, on a alors $\zeta^r = \zeta^q$.

À partir de maintenant, on travaille dans une carte $j_\infty$--holomorphe au voisinage de $\zeta := u_\infty^{-1}(p) \in \S_\infty$ qui envoie $\zeta$ sur $0\in \mathbb{C}$ et dans la carte donnée par l'application exponentielle $\mathrm{exp}_p$ au voisinage de $p \in M$. 
On définit pour tout entier naturel $n$ et pour tout $z$ au voisinage de $0\in \mathbb{C}$, $f_n (z) = \frac{(u_n \circ \varphi_n) (z)}{z}$ si $z \neq 0$ et $f_n (0) = (u_n\circ \varphi_n)'(0)$, ainsi que $f_\infty (z) = \frac{u_\infty (z)}{z}$ si $z \neq 0$ et $f_\infty (0) = u_\infty'(0)$. Les fonctions $f_n$ et $f_\infty$ sont continues et la suite $(f_n)$ converge localement uniformément vers $f_\infty$ pour la topologie $\mathcal{C}^0$. Par continuité de l'application $(f , \zeta) \mapsto f (\zeta)$,
 on en déduit 
$$u_\infty '(0) = \lim\limits_n \frac{u_\infty (\zeta_n^q)}{\zeta_n^q} = \lim\limits_n \frac{(u_n \circ \varphi_n) (\zeta_n^q)}{\zeta_n^q} = \lim\limits_n \frac{1}{2^n \zeta_n^q}v.$$
Puisque $u_\infty$ est plongée, le côté gauche de l'égalité est un vecteur non nul de $T_p u$, donc $v$ est un vecteur de $T_p u$.
De même on obtient 
$$u_\infty '(0) = \lim\limits_n \frac{1}{2^n \zeta_n^p}w,$$
donc $w$ est également un vecteur de $T_p u$. C'est une contradiction puisque $v$ et $w$ sont des vecteurs $J_p$--linéairement indépendants.
\begin{figure}[h]
	\centering
	\includegraphics[scale=0.5]{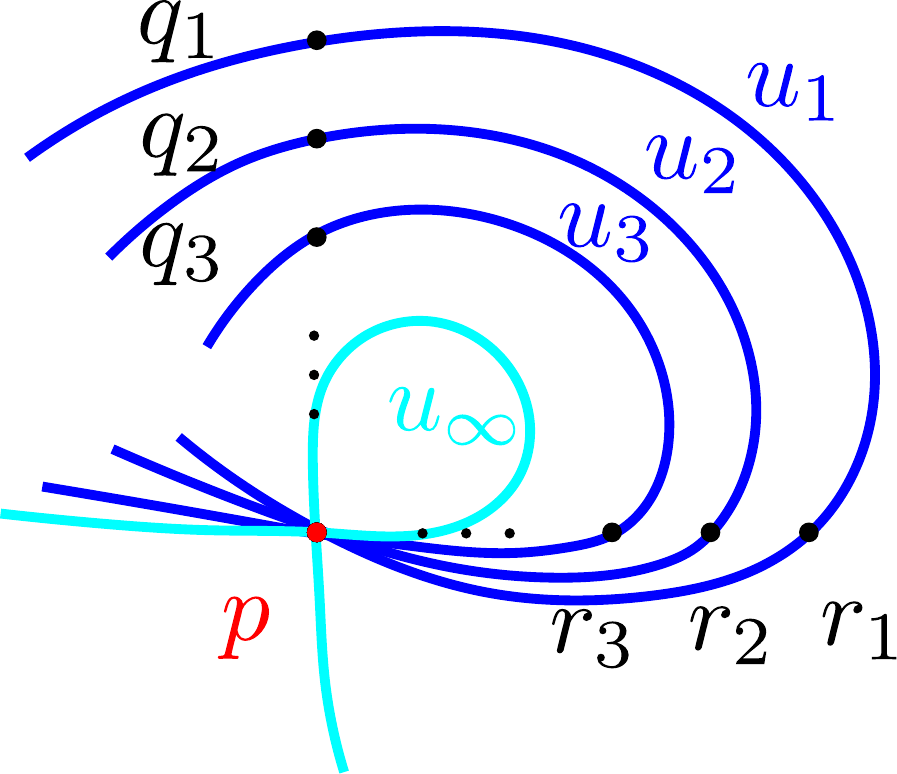}\\
	\caption{La limite $u_\infty$ de la suite de courbes pseudoholomorphes plongées $(u_n)$ dans la démonstration de la Proposition~\ref{p:casserlescourbes} ne peut pas être plongée.}
	\label{f:sequence1}
\end{figure}
\end{proof}

\begin{rk}\label{r:CCnoncompactes}
Dans la démonstration de la Proposition~\ref{p:casserlescourbes}, on peut raisonner indépendamment sur chacune des composantes connexes (voir également la Remarque~\ref{r:surjCC}). Ainsi, on a montré que sous les hypothèses de la Proposition~\ref{p:casserlescourbes}, aucune des composantes connexes de l'espace de modules  $\mathcal{M}_g^* ([u];J)$ n'est compacte.
\end{rk}

\begin{rk}\label{r:casseràpoint}
\'Evoquons une légère modification de la preuve de la Proposition~\ref{p:casserlescourbes} qui permet d'obtenir un résultat plus fort dont on aura besoin par la suite. En supposant que l'espace de modules $\mathcal{M}_g^* ([u];J;p)$ est compact, avec $p$ générique fixé au départ, on aboutit également à la même contradiction. Donc pour tout point $p$ générique, sous les hypothèses de la Proposition~\ref{p:casserlescourbes}, l'espace de modules $\mathcal{M}_g^* ([u];J;p)$ n'est pas compact. On souligne aussi qu'il suffit de montrer que $\mathcal{M}_{g,1}^* ([u];J;p)$ n'est pas compact pour montrer que $\mathcal{M}_g^* ([u];J)$ n'est pas compact (grâce à l'application d'oubli du point marqué et à la Remarque~\ref{r:limitenonplongée}).
\end{rk}

\begin{rk}\label{r:casserlescourbes}
Dans la démonstration de la Proposition~\ref{p:casserlescourbes}, on pouvait aussi procéder de la manière suivante. Une fois la suite $(u_n)$ construite, on oublie la structure presque complexe $J$ et on éclate en $p$. On note $E$ le diviseur exceptionnel apparu suite à l'éclatement. En considérant les transformées propres des éléments de la suite $(u_n)$ et de $u_\infty$, on obtient une suite de courbes symplectiquement plongées $(S_n)$ dans $M \# \overline{\mathbb{C}P}^2$, qui converge vers une courbe symplectiquement plongée $S_\infty$, telle que pour tout entier naturel $n$, $S_n$ intersecte $E$ exactement une fois, de manière transverse et positive. Comme les suites $(q_n)$ et $(r_n)$ dans $M$ tendent vers $p$ dans des directions $J_p$--linéairement indépendantes, elles tendent dans $M \# \overline{\mathbb{C}P}^2$ respectivement vers des points $q,r \in E$ distincts. La courbe symplectique $S$ doit alors nécessairement intersecter $E$ en $q$ et $r$. L'image de $S$ par la contraction de $E$ possède par conséquent une singularité en $p$. Comme $S$ est la transformée propre de $u_\infty$, qui est plongée, on aboutit à une contradiction. 
\begin{figure}[h]
	\centering
	\includegraphics[scale=0.5]{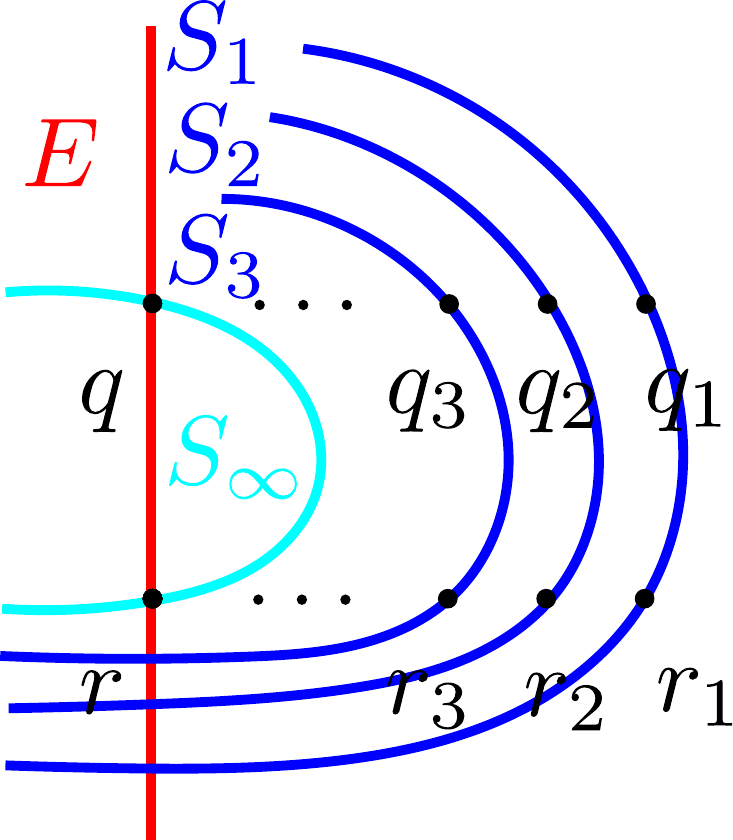}\\
	\caption{La limite $S_\infty$ de la suite de courbes symplectiquement plongées $(S_n)$ dans la Remarque~\ref{r:casserlescourbes} doit intersecter le diviseur exceptionnel $E$ en deux points distincts.}
	\label{f:sequence2}
\end{figure}
\end{rk}

\subsubsection{Le cas des tores d'auto-intersection 3}

Dans la suite, on aura besoin de montrer que l'espace de modules de courbes plongées $\mathcal{M}_g^* ([u];J)$ n'est pas compact dès que $[u]^2 \geq 4g-1$. Le seul cas que la Proposition~\ref{p:casserlescourbes} ne permet pas de traiter est celui où $g=1$ et $[u]^2 = 3$ (ce point est détaillé dans la Sous-section~\ref{ss:recurrence}). La proposition suivante nous permettra de remédier à ce problème.

\begin{prop} \label{p:autoint3}
Soit $J \in \mathcal{J}^{\mathrm{reg}}_{\tau} (M, \omega)$ et $u : \Sigma \rightarrow M$ une courbe $J$--holomorphe plongée de genre $g$. On suppose qu'on est dans une des deux situations suivantes :
\begin{enumerate}
\item $\ind (u) > 2g +2$ et il existe une courbe $J$--holomorphe $v$ dont l'image est distincte de celle de $u$ telle que $[u] \cdot [v] = 1$, ou bien
\item $\ind (u) > 2g $ et il existe une courbe $J$--holomorphe $v$ dont l'image est distincte de celle de $u$ telle que $[u] \cdot [v] = 0$.
\end{enumerate}
Alors l'espace de modules $\mathcal{M}_g^* ([u];J)$ n'est pas compact.
\end{prop}

\begin{proof}
On suppose que l'espace de modules $\mathcal{M}_g^* ([u];J)$ est compact. 

Dans la première situation, le Lemme~\ref{l:surjectev} nous donne la surjectivité l'application $\ev : \mathcal{M}_{g,2}^* ([u];J) \rightarrow M^2 \backslash \Delta$. En particulier, par deux points distincts dans l'image de $v$, il passe une courbe de $\mathcal{M}_g^* ([u];J)$. La positivité d'intersection nous donne alors une contradiction avec $[u]\cdot[v] = 1$. 

Dans la seconde situation, le Lemme~\ref{l:surjectev} nous donne la surjectivité l'application $\ev : \mathcal{M}_{g,1}^* ([u];J) \rightarrow M$. En particulier, par tout point dans l'image de $v$, il passe une courbe de $\mathcal{M}_g^* ([u];J)$. La positivité d'intersection nous donne alors une contradiction avec $[u] \cdot [v] = 0$. 
\end{proof}

\subsection{Trouver des courbes nodales avec une unique paire de points nodaux}\label{ss:glue}

Dans la suite, on note $\mathcal{M}_S(J)$ la composante connexe de l'espace de modules $\mathcal{M}_g^* ([u];J)$ qui contient la courbe $u$ d'image $S$.

Maintenant qu'on a montré l'existence de courbes nodales dans $\moduleadh$ pour une structure presque complexe générique $J$ (voir la Remarque~\ref{r:CCnoncompactes}), on s'intéresse à trouver des courbes symplectiques singulières homologues à $S$ d'une forme spécifique : 
\begin{itemize}
\item deux composantes plongées qui s'intersectent exactement une fois, de manière transverse, ou bien
\item une composante simple immergée dont l'unique point d'auto-intersection est un point double transverse.
\end{itemize}
\begin{figure}[htbp]
	\centering
	\subfloat[Courbe symplectique singulière consituée de deux composantes plongées qui s'intersectent exactement une fois, de manière transverse.]{\label{f:typeAbis}\includegraphics[width=0.45\linewidth]{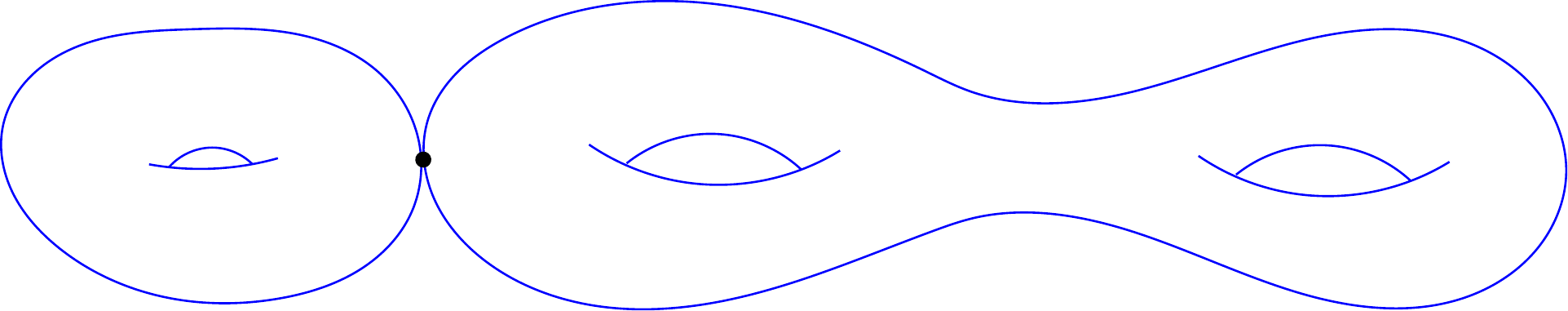}}\qquad
	\subfloat[Courbe symplectique singulière constituée d'une composante simple immergée dont l'unique point d'auto-intersection est un point double transverse.]{\label{f:typeBbis}\includegraphics[width=0.45\linewidth]{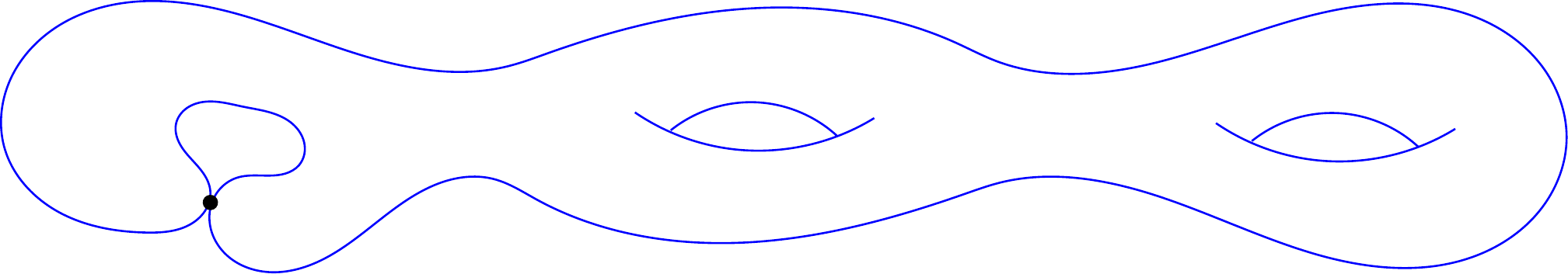}}
\caption{Exemples de courbes nodales avec une unique paire de points nodaux.}\label{f:typebis}
\end{figure}

\begin{rk}
Quand $\ind (u) =2$, l'existence d'une courbe symplectique singulière d'une de ces deux formes est immédiate par la Proposition~\ref{p:morceauxcourbesind2}. 
\end{rk}

Ces éléments sont ceux qui vont permettre d'effectuer l'hérédité dans la récurrence servant à démontrer le théorème principal. Pour montrer l'existence de telles courbes symplectiques singulières, on les construit explicitement à partir d'un autre élément $\ulim$ de $\moduleadh$. On distingue trois cas :
\begin{itemize}
\item le cas où toutes les composantes non constantes de $\ulim$ sont simples et d'images distinctes, qui fait l'objet du Lemme~\ref{l:lissage0},
\item le cas où les composantes non constantes de $\ulim$ ont toutes la même image (c'est le cas -- entre autres -- lorsque $\ulim$ possède une unique composante qui est un revêtement multiple par exemple), qui fait l'objet du Lemme~\ref{l:lissage1},
\item le cas où $\ulim$ possède au moins deux composantes non constantes d'images distinctes, qui fait l'objet du Lemme~\ref{l:lissage2}.
\end{itemize}
Le premier cas est couvert par l'union des deux autres cas, mais il est intéressant de l'énoncer à part puisqu'on s'en servira pour montrer les deux cas restants.

Dans les démonstrations de ces différents lemmes, on sera souvent amené à lisser des points doubles transverses positifs. Cependant les points singuliers qui apparaissent peuvent être de natures différentes. On utilisera alors la régularité Fredholm des courbes simples qui nous intéressent pour pouvoir perturber ces courbes et ainsi remplacer ces points singuliers par des points doubles transverses positifs.


\begin{lem} \label{l:lissage0}
Soit $\ulim \in \moduleadh$. On suppose que les composantes non constantes de $\ulim$ sont des composantes simples d'images distinctes. Alors $(M, \omega)$ contient une courbe symplectique singulière d'une des deux formes suivantes :
\begin{itemize}
\item deux composantes plongées de genres respectifs $g_1$ et $g_2$, avec $g=g_1+g_2$, qui s'intersectent exactement une fois, de manière transverse,
\item une composante symplectiquement positivement immergée de genre $g-1$, avec un unique point double transverse.
\end{itemize}
De plus, dans chacun des cas, la somme des classes d'homologie des composantes est égale à $[S]$.
\end{lem}

\begin{proof}
On note $v_1, \dots , v_\ell$ les composantes de $\ulim$. Les courbes s'intersectent deux à deux en des points isolés, mais pas nécessairement de manière transverse. Puisque les courbes $v_1, \dots, v_\ell$ sont simples, la généricité de $J$ assure que les courbes $v_1, \dots, v_\ell$ sont Fredholm régulières. Par les propriétés de généricité évoquées au dernier paragraphe de la Sous-section~\ref{s:contraintesdérivées}, on peut alors supposer, quitte à réaliser une perturbation arbitrairement petite de $J$ et de $v_1, \dots , v_\ell$, que :
\begin{itemize}
\item pour tout $i$, la courbe $v_i$ est symplectiquement positivement immergée,
\item pour tous $i$,$j$, les courbes $v_i$ et $v_j$ s'intersectent transversalement positivement en des points injectifs,
\item pour tous $i$, $j$, $k$, la courbe $v_i$ ne passe pas par les points d'intersection entre $v_j$ et $v_k$.
\end{itemize}
La configuration de courbes $J$--holomorphes $\mathcal{C}$ constituée de l'union des courbes $v_1, \dots , v_\ell$ possède au moins un point double. En effet, si $\ell \geq 2$ cela provient de la connexité de $\mathcal{C}$ et si $\ell = 1$, la courbe $v_1$, qui est simple, ne peut pas être une courbe plongée car $u_\infty$ appartiendrait sinon à $\mathcal{M}_S(J)$ (voir la Remarque~\ref{r:limitenonplongée}).

On choisit alors un point double transverse positif de la configuration de courbes $\mathcal{C}$ et on lisse un à un tous les autres points doubles transverses positifs de $\mathcal{C}$ (voir la Remarque~\ref{r:lissageptdoubleJholo}). Puisque la configuration $\mathcal{C}'$ ainsi construite ne possède qu'un seul point double, on obtient alors une des deux possibilités suivantes : 
\begin{itemize}
\item deux courbes symplectiquement plongées qui s'intersectent transversalement positivement exactement une fois,
\item une courbe symplectiquement positivement immergée avec un unique point double transverse.
\end{itemize}
Les conditions sur les genres des différentes composantes sont alors aisément obtenues via la formule d'adjonction pour les courbes symplectiques singulières puisque $[S]=[S']$.
\end{proof}

Avant de continuer, on rappelle que d'après la Remarque~\ref{r:casseràpoint}, pour un point $p \in M$ fixé à l'avance et pour une structure presque complexe $J \in \mathcal{J}_\tau (M,\omega)$ générique pour $p$, l'espace de modules $\overline{\M}_S (J;p)$ n'est pas compact. 

\begin{lem} \label{l:lissage1}
Soit $\ulim \in \moduleadhcont$. On suppose que les composantes non constantes de $\ulim$ ont toutes la même image. Alors $(M, \omega)$ contient une courbe symplectique singulière d'une des deux formes suivantes :
\begin{itemize}
\item deux composantes plongées de genres respectifs $g_1$ et $g_2$, avec $g=g_1+g_2$, qui s'intersectent exactement une fois, de manière transverse,
\item une composante symplectiquement positivement immergée de genre $g-1$, avec un unique point double transverse.
\end{itemize}
De plus, dans chacun des cas, la somme des classes d'homologie des composantes est égale à $[S]$.
\end{lem}

\begin{proof}
Par soucis de clarté, on suppose tout au long de la démonstration que $\ulim$ ne possède qu'une composante, qui est un revêtement multiple d'ordre $k \geq 1$ d'une courbe $J$--holomorphe simple $v$. Les autres cas se traitent de manière strictement similaire. Par généricité de la structure presque complexe, et puisque la courbe $J$--holomorphe simple $v$ passe par le point $p$, l'indice contraint de $v$ est positif. Autrement dit, on a $\ind (v) -2 \geq 0$. On note $g'$ le genre de $v$. L'espace de modules des courbes simples sans contrainte $\M_{g'}^*([v]; J)$ est par conséquent une variété lisse de dimension au moins $2$. On choisit désormais $k$ éléments deux à deux distincts $v_1, \dots, v_k$ de $\M_{g'}^*([v]; J)$. Les courbes $v_i$ ne sont pas deux à deux disjointes puisque $[v]^2 = \frac{1}{k^2} [u]^2 >0$. L'union de ces $k$ courbes forme donc un élément $\ulim '$ de $\moduleadh$ qui vérifie les hypothèses du Lemme~\ref{l:lissage0}, ce qui permet de conclure.
\end{proof}

\begin{lem} \label{l:lissage2}
Soit $\ulim \in \moduleadh$. On suppose que $\ulim$ possède au moins deux composantes d'images distinctes. Alors $(M, \omega)$ contient une courbe symplectique singulière d'une des deux formes suivantes :
\begin{itemize}
\item deux composantes plongées de genres respectifs $g_1$ et $g_2$, avec $g=g_1+g_2$, qui s'intersectent exactement une fois, de manière transverse,
\item une composante symplectiquement positivement immergée de genre $g-1$, avec un unique point double transverse.
\end{itemize}
De plus, dans chacun des cas, la somme des classes d'homologie des composantes est égale à $[S]$.
\end{lem}

Avant de démontrer le Lemme~\ref{l:lissage2}, on rappelle que la surface symplectique $(M,\omega)$ est supposée relativement minimale par rapport à $S$ et que la courbe $J$--holomorphe $u$ plongée d'image $S$ vérifie $[u]^2 \geq 4g$ (voir l'énoncé du Théorème~\ref{t:ruledFKversion2} et le paragraphe qui suit la Remarque~\ref{r:comparaisondeTh}).

\begin{proof}
Notons $\ell$ le nombre de composantes de $\ulim$ et pour tout $i \in \{1, \dots, \ell \}$, $u_i : \S_i \rightarrow M$ les différentes composantes de $\ulim$. Pour tout $i \in \{1, \dots, \ell \}$ tel que la courbe $u_i$ est une composante non constante, $u_i$ est un revêtement multiple d'ordre $k_i \geq 1$ d'une courbe $J$--holomorphe simple $v_i : \S_i' \rightarrow M$. L'objectif est de se ramener au cas où $\ulim$ possède exactement deux composantes simples. Pour ce faire, on raisonne par récurrence sur $N = \sum_{i=1}^\ell k_i$. Si $N=2$, puisque $\ulim$ possède au moins deux composantes d'images distinctes, les deux composantes sont simples et on conclut en appliquant le Lemme~\ref{l:lissage0}. Si $N > 2$, quitte à ré-indexer, on peut supposer que $v_1$ et $v_2$ sont d'images distinctes et $[v_1] \cdot [v_2] >0$. Par les propriétés de généricité évoquées au dernier paragraphe de la Sous-section~\ref{s:contraintesdérivées}, on peut supposer, quitte à réaliser une perturbation arbitrairement petite de $J$, $v_1$ et $v_2$, que les courbes $v_1$ et $v_2$ s'intersectent seulement en des points réguliers et de manière transverse (comme $J$ est générique, toutes les composantes simples $v_i$ survivent à la perturbation et restent Fredholm régulières). On lisse un point double (voir la Remarque~\ref{r:lissageptdoubleJholo}) entre $v_1$ et $v_2$ pour obtenir une courbe $J$--holomorphe simple $v_0 : \S_1' \# \S_2' \rightarrow M$, qu'on peut supposer distincte des autres composantes, telles que $[v_0] = [v_1] + [v_2]$. On considère désormais l'union des courbes suivantes : $v_0$ avec multiplicité $1$, $v_1$ avec multiplicité $k_1 -1$, $v_2$ avec multiplicité $k_2-1$ et pour $i \geq 3$, $v_i$ avec multiplicité $k_i$.

Il reste encore à montrer que la configuration de courbes ainsi obtenue est connexe. Autrement dit, on doit prouver que $v_0$ (ainsi que $v_1$ et $v_2$ si elles existent encore dans la nouvelle configuration) intersecte l'union des autres composantes de la configuration.

Dans ce paragraphe, on suppose que $\ulim$ possède au moins trois composantes. Il existe alors une composante $v_3$ (quitte à ré-indexer) qui intersecte l'union de l'image de $v_1$ et de l'image de $v_2$ (voir Figure~\ref{f:gluev}). 
\begin{figure}[h]
	\centering
	\includegraphics[scale=0.6]{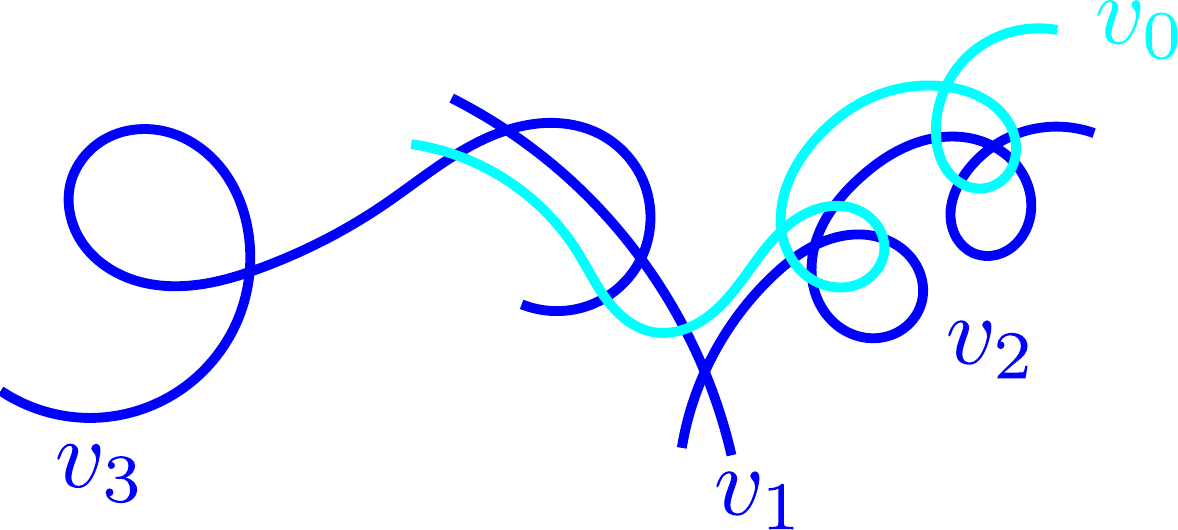}\\
	\caption{La courbe $v_0$ est obtenue en lissant un point d'intersection transverse entre $v_1$ et $v_2$.}
	\label{f:gluev}
\end{figure}
Par positivité d'intersection, on a alors $[v_0] \cdot [v_3] = \left( [v_1] + [v_2] \right) \cdot [v_3] > 0$. Il nous reste à montrer que dans la nouvelle configuration de courbes, $v_1$ et $v_2$ (si elles existent encore) intersectent l'union des autres composantes. Le seul cas (à ré-indexation près) à considérer est celui où $\ulim$ possède au moins trois composantes et où $k_1 >1$, $k_2 = 1$ et $[v_0] \cdot[v_1] =0$ (la courbe $v_1$ pourrait alors a priori devenir disjointe des autres composantes suite à la disparition de $v_2$ à l'étape suivante de la récurrence). Montrons que $v_1$ intersecte au moins une autre composante. Puisque $[v_0] \cdot [v_1] = 0$, on a  $[v_1]^2 = - [v_1] \cdot [v_2]$. Or par généricité de $J$ et puisque $v_1$ est simple, on a $[v_1]^2  \geq -1$. Donc $[v_1] \cdot[v_2] \leq 1$. Puisque les composantes $v_1$ et $v_2$ s'intersectent au moins une fois, on a nécessairement $[v_1] \cdot [v_2] = 1$, d'où $[v_1]^2 = -1$ (notons que par la formule d'adjonction, on a $\ind (v_1) = \chi (\Sigma_1') + 2[v_1]^2 -4 \delta (v_1) \geq 0$ ; ainsi $v_1$ est nécessairement une courbe rationnelle plongée, donc un diviseur exceptionnel).
Par positivité d'intersection on a $[v_1] \cdot [u] \geq 0$, donc $[v_1] \cdot \sum\limits_{i=1}^\ell k_i [v_i] \geq 0$, d'où $\sum\limits_{i=3}^\ell k_i [v_1] \cdot [v_i] \geq k_1-1$. Comme $k_1 \geq 2$, on obtient alors $\sum\limits_{i=3}^\ell k_i [v_1] \cdot [v_i] > 0$. Ainsi, la courbe $v_1$ n'est pas disjointe de l'union des autres composantes de la nouvelle configuration, ce qui prouve bien que cette nouvelle configuration de courbes est connexe.

Dans ce paragraphe, on suppose que $\ulim$ possède exactement deux composantes. Si $[v_0] \cdot [v_1] = 0$, alors $[v_1]^2 = - [v_1] \cdot [v_2]$. Or par généricité de $J$ et puisque $v_1$ est simple, on obtient de la manière manière qu'au paragraphe précédent $[v_1]^2  \geq -1$. On en déduit  alors que $[v_1] \cdot [v_2] \leq 1$. Puisque les composantes $v_1$ et $v_2$ s'intersectent au moins une fois, on a nécessairement $[v_1] \cdot [v_2] = 1$, d'où $[v_1]^2 = -1$ (notons que, comme précédemment, par la formule d'adjonction, $v_1$ est nécessairement une courbe rationnelle plongée, donc un diviseur exceptionnel). Par hypothèse de relative minimalité, on a $[v_1] \cdot [u] >0$, donc $[v_1] \cdot \left( k_1 [v_1] +k_2 [v_2] \right) > 0$, d'où $k_2 > k_1$. La courbe $v_2$ était donc une composante comptée avec multiplicité au moins $2$ de $\ulim$ et sera par conséquent présente à l'étape suivante de la récurrence. On a finalement $[v_2] \cdot [v_0] = [v_2] \cdot \left( [v_1] + [v_2] \right) = [v_2]^2 +1$.
On peut remarquer que $[v_2]^2 \geq 0$. En effet, si on avait $[v_2]^2 = -1$, on aurait $[u]^2 = \left( k_1 [v_1] +k_2 [v_2] \right)^2 = -k_1^2 +2k_1 k_2 -k_2^2 < 0$, ce qui contredirait les hypothèses sur $u$. Donc $[v_2] \cdot [v_0] = [v_2] \cdot \left([v_1]+[v_2]\right) >0$. Autrement dit, $v_0$ intersecte bien l'union des autres composantes de la configuration, ce qui prouve bien que cette nouvelle configuration de courbes est connexe.

Pour l'étape suivante de la récurrence, l'indice de récurrence est égal à $N' = N + 1 -1 -1 < N$. On conclut alors en appliquant le principe de récurrence, puis en utilisant finalement le Lemme~\ref{l:lissage0}.
\end{proof}

On résume à présent les différents résultats obtenus dans la proposition suivante.

\begin{prop}\label{p:existeunidonale}
Si l'ensemble $\moduleadhcont$ est non vide, alors $(M,\omega)$ contient une courbe symplectique singulière d'une des deux formes suivantes :
\begin{itemize}
\item deux composantes plongées de genres respectifs $g_1$ et $g_2$, avec $g=g_1+g_2$, qui s'intersectent exactement une fois, de manière transverse et positive,
\item une composante symplectiquement positivement immergée de genre $g-1$, avec un unique point double transverse.
\end{itemize}
De plus, dans chacun des cas, la somme des classes d'homologie des composantes est égale à $[S]$.
\end{prop}

\begin{proof}
Les Lemmes~\ref{l:lissage1} et~\ref{l:lissage2} permettent de conclure.
\end{proof}

\begin{rk}\label{r:imagesdecourbesFR}
En reprenant les preuves des Lemmes~\ref{l:lissage0},\ref{l:lissage1} et~\ref{l:lissage2}, on remarque que chacune des composantes de la courbe symplectique singulière obtenue dans la Proposition~\ref{p:existeunidonale} est en fait l'image d'une courbe pseudoholomorphe simple Fredholm régulière.
\end{rk}

\subsection{La récurrence} \label{ss:recurrence}

Dans cette sous-section, on démontre le théorème suivant en effectuant une récurrence sur le genre. 

\begin{thm}\label{t:ruledFKversion}
Soit $S$ une courbe de genre $g$ symplectiquement plongée dans une surface symplectique $(M, \omega )$. Si $[S]^2 > 4g -1$, alors il existe une courbe rationnelle symplectiquement plongée dans $(M,\w)$ d'auto-intersection $0$ ou $1$ qui intersecte $S$ localement positivement.
\end{thm}

\begin{rk}\label{r:ruledFKversion}
Les Théorèmes~\ref{t:ruledFKversion} et~\ref{t:McDuff} impliquent que $(M, \omega )$ est une surface symplectiquement réglée, ce qui démontre le Théorème~\ref{t:ruledFKversion2}. Il est alors important de noter que $(M, \w)$ admet ou bien une fibration de Lefschetz symplectique (générique) où les fibres régulières sont des sphères, ou bien un pinceau symplectique avec un unique point base où les fibres régulières sont des sphères. De plus, on peut choisir ce pinceau ou cette fibration (en s'aidant d'une structure presque complexe $J$ générique dominée par $\omega$ auxiliaire bien choisie et en appliquant le Théorème~\ref{t:isotopiepinceaudeLefschetz}) de sorte que $S$ intersecte toutes les fibres  localement positivement.
\end{rk}

Quitte à contracter une collection maximale de diviseurs exceptionnels disjoints de $S$ et deux à deux disjoints, on peut supposer sans perte de généralité tout au long de la preuve que $(M, \w)$ est relativement minimale par rapport à $S$ (il suffira ensuite d'effectuer les éclatements adéquats à la fin de la preuve). Cette hypothèse va nous être indispensable pour trouver des courbes de genre strictement plus petit que $g$ avec auto-intersection suffisamment grande, comme le montre l'exemple suivant. 

\begin{ex}\label{e:nonminrel}
Soit $J$ une structure presque complexe générique (en particulier les indices des courbes $J$--holomorphes simples sont positifs). Une suite $(u_n)$ de courbes $J$--holomorphes plongées de genre $g$ peut éventuellement converger vers une courbe nodale composée de deux courbes plongées $u_1$ et $u_2$ qui s'intersectent exactement une fois, de manière transverse. Sans l'hypothèse de minimalité relative, $u_1$ pourrait être un diviseur exceptionnel disjoint de $u$ et $u_2$ une courbe de genre $g$ d'auto-intersection $[u]^2 -1$. Cette courbe nodale ne permet pas dans ce cas de garantir à coup sûr l'existence d'une composante plongée de genre strictement plus petit que $g$ avec auto-intersection suffisamment grande.
\end{ex} 

Le lemme suivant permet de s'assurer que la situation décrite dans l'Exemple~\ref{e:nonminrel} n'arrive pas quand on a l'hypothèse de minimalité relative.

\begin{lem}\label{l:noexdiv}
Soit $(M,J)$ une surface presque complexe, $u$ une courbe $J$--holomorphe plongée et $v,w$ des courbes $J$--holomorphes plongées Fredholm régulières telles que $[u] = [v]+[w]$ et $[v] \cdot [w] =1$. Si $(M, J)$ est relativement minimale par rapport à l'image de $u$, alors on a $[v]^2 \geq0$ et $[w]^2 \geq 0$.
\end{lem}
\begin{proof}
Comme les rôles de $v$ et $w$ sont symétriques, il suffit de montrer que $[v]^2 \geq 0$. On note $g$ le genre de $v$. Puisque $v$ est une courbe Fredholm régulière, on a $\ind (v) \geq 0$. \`A l'aide de la formule d'adjonction, on obtient alors $[v]^2 \geq g -1$. Si $g \geq 1$, le lemme est vérifié. On suppose désormais $g=0$ et $[v]^2 <0$, c'est-à-dire $[v]^2 =-1$. Comme la courbe $J$--holomorphe $v$ est plongée et de genre $0$, on en déduit que c'est un diviseur exceptionnel. Puisque $[w] \cdot [v] =1$, on a  alors $[u]\cdot [v] = [v]^2 +[w] \cdot [v] =0$. Comme $u$ et $v$ sont des courbes d'images distinctes (sinon on aurait eu $[u]\cdot [v] = -1$), on en déduit par positivité d'intersection que les courbes $u$ et $v$ sont disjointes. Mais ceci contredit la minimalité relative de $M$ par rapport à l'image de $u$. 
\end{proof}

Le lemme suivant fait office d'hérédité pour le raisonnement par récurrence présenté plus loin dans cette sous-section.

\begin{lem}\label{l:hérédité}
Soit $S$ une courbe de genre $g>0$ symplectiquement plongée dans une surface symplectique $(M, \omega)$. Si $[S]^2 \geq 4g -1$ et $S$ n'est pas un tore d'auto-intersection $3$, alors  $(M, \omega)$ ou l'éclaté de $(M, \omega)$ en un point disjoint de $S$ contient une des deux courbes suivantes :
\begin{itemize}
\item une courbe rationnelle symplectiquement plongée d'auto-intersection positive, ou
\item une courbe symplectiquement plongée $C$ de genre $\g$ qui vérifie $0 < \g < g$ et $[C]^2 \geq 4\g-1$.
\end{itemize}
\end{lem}

\begin{proof}
On commence par se ramener au cas où la surface symplectique ambiante est relativement minimale par rapport à $S$ en contractant une collection maximale de diviseurs exceptionnels $E_1, \dots, E_\ell$ disjoints de $S$ et deux à deux disjoints.
Comme expliqué au début de la Section~\ref{s:démo}, on peut trouver une structure presque complexe générique $J$ dominée par la forme symplectique telle que $S$ est l'image d'une courbe $J$--holomorphe plongée $u$ de genre $g$.

On commence par prouver l'existence d'une courbe nodale dans $\moduleadhcont$. D'après la formule d'adjonction, et puisque $[S]^2 \geq 4$, on a $\ind (u) =  -\chi( S ) + 2 c_1([S]) = \chi( S ) +2 [S]^2 \geq 2g +6$. Comme $g >0$, cela signifie que $\ind (u) > 6$. D'après la Proposition~\ref{p:casserlescourbes}, l'espace de modules $\mathcal{M}_S (J) = \mathcal{M}^*_g([u];J)$ n'est alors pas compact. La Proposition~\ref{p:existeunidonale} nous donne donc l'existence d'une courbe symplectique singulière d'une des deux formes suivantes :
\begin{itemize}
\item deux composantes $S_1$ et $S_2$ symplectiquement plongées, qui s'intersectent exactement une fois, de manière transverse,
\item une composante symplectiquement positivement immergée $S'$, avec un unique point double transverse.
\end{itemize}

De plus, dans le premier cas, on a $[S] = [S_1]+[S_2]$. Notons $g_1$ le genre de $S_1$ et $g_2$ le genre de $S_2$. Si une des courbes $S_i$ est de genre $0$, elle est d'auto-intersection positive d'après le Lemme~\ref{l:noexdiv} (en effet, d'après la Remarque~\ref{r:imagesdecourbesFR}, $S_1$ et $S_2$ sont obtenues comme images de courbes pseuholomorphes plongées Freholm régulières) et on a terminé. Supposons que les courbes $S_i$ sont de genre strictement positif. Puisque $g = g_1 + g_2$, chacune des courbes $S_i$ est de genre strictement plus petit que $g$. Supposons $[S_1]^2  \leq 4g_1 -2$ et $[S_2]^2 \leq 4g_2 -2$. L'égalité $[S]^2 = [S_1]^2 + 2 +[S_2]^2$ nous donne alors $[S]^2 \leq 4g_1 +4g_2 - 2$, c'est-à-dire $[S]^2 \leq 4g - 2$, ce qui est impossible. Ainsi, une des composantes, disons $S_1$ quitte à échanger les indices, vérifie $[S_1]^2 \geq 4g_1 -1$. 

Dans le second cas, on éclate le point double de $S'$ (qu'on peut supposer disjoint de $S$, quitte à perturber légèrement $S'$). On obtient alors une courbe $C$ plongée de genre $\g = g-1$, et on a $[C]^2 = [S']^2 - 4 = [S]^2 -4 \geq 4 \g -1$. Comme $S$ n'est pas un tore d'auto-intersection $3$, on a $[S]^2 \geq 4$. Par conséquent si $\g =0$, $C$ est une courbe rationnelle symplectiquement plongée d'auto-intersection positive. 


Pour finir, on éclate les points qui correspondent aux images de $E_1, \dots, E_\ell$ par la contraction effectuée plus tôt afin de retrouver $(M, \omega)$ (dans le premier cas) ou $(M, \omega)$ éclaté en un point disjoint de $S$ (dans le second cas).
\end{proof}

\begin{rk}
La seule raison pour laquelle l'inégalité est stricte dans l'énoncé du Théorème~\ref{t:ruledFKversion} tient à la gestion des tores symplectiquement plongés d'auto-intersection $3$.
\end{rk}

Dans la démonstration qui suit, la stratégie concernant les tores d'auto-intersection $3$ est la suivante : on peut selon les cas soit montrer qu'ils n'apparaissent pas au cours de la récurrence, soit les \og casser \fg{} à l'aide de la Proposition~\ref{p:autoint3}.

\begin{proof}[Démonstration du Théorème~\ref{t:ruledFKversion}]
On raisonne par récurrence sur le genre de $S$. 
Si $g=0$, on conclut en appliquant le Théorème~\ref{t:McDuff} de McDuff. On suppose désormais que $g>0$.

D'après le Lemme~\ref{l:hérédité}, dès que la surface symplectique ambiante considérée (un éclaté de $(M, \omega)$ en des points disjoints de $S$) à une étape donnée contient une courbe $C$ de genre $\g >0$ symplectiquement plongée, qui n'est pas un tore d'auto-intersection $3$, telle que $[C]^2 \geq 4 \g -1$, on peut faire apparaître (quitte à éclater en un point disjoint de $S$) une courbe rationnelle symplectiquement plongée d'auto-intersection positive ou une courbe symplectiquement plongée $C'$ de genre $\g '$, avec $0 < \g' < \g$, telle que $[C']^2 \geq 4 \g' -1$. Dans le premier cas, la récurrence se termine. Dans le second cas, si la courbe $C'$ n'est pas un tore d'auto-intersection $3$, elle satisfait les hypothèses du Lemme~\ref{l:hérédité} et on peut continuer la récurrence.

Comme le genre de courbes considérées décroît strictement à chaque étape, au bout d'un nombre fini d'étapes on trouve, dans un éclaté de $(M, \omega)$ en des points disjoints de $S$, une courbe rationnelle symplectiquement plongée d'auto-intersection positive ou un tore symplectiquement plongé d'auto-intersection $3$.

Avant de continuer, on rappelle qu'à chaque étape la courbe $C'$ est obtenue en appliquant la Proposition~\ref{p:existeunidonale} à $C$. Cette proposition affirme qu'on trouve une courbe symplectique singulière d'une des deux formes suivantes :
\begin{itemize}
\item Deux composantes $S_1$ et $S_2$ symplectiquement plongées, qui s'intersectent exactement une fois, de manière transverse. On dira dans ce cas que cette étape de la récurrence est de type A.
\item Une composante symplectiquement positivement immergée $S'$, avec un unique point double transverse. On dira dans ce cas que cette étape de la récurrence est de type B.
\end{itemize}
La composante $C'$ est obtenue comme une des deux composantes $S_1$ ou $S_2$ si l'étape est de type A. Elle est en revanche obtenue comme la transformée propre de $S'$ après l'éclatement de son point double si l'étape est de type B. 
\begin{figure}[htbp]
	\centering
	\subfloat[\'Etape de type A.]{\label{f:typeA}\includegraphics[width=0.45\linewidth]{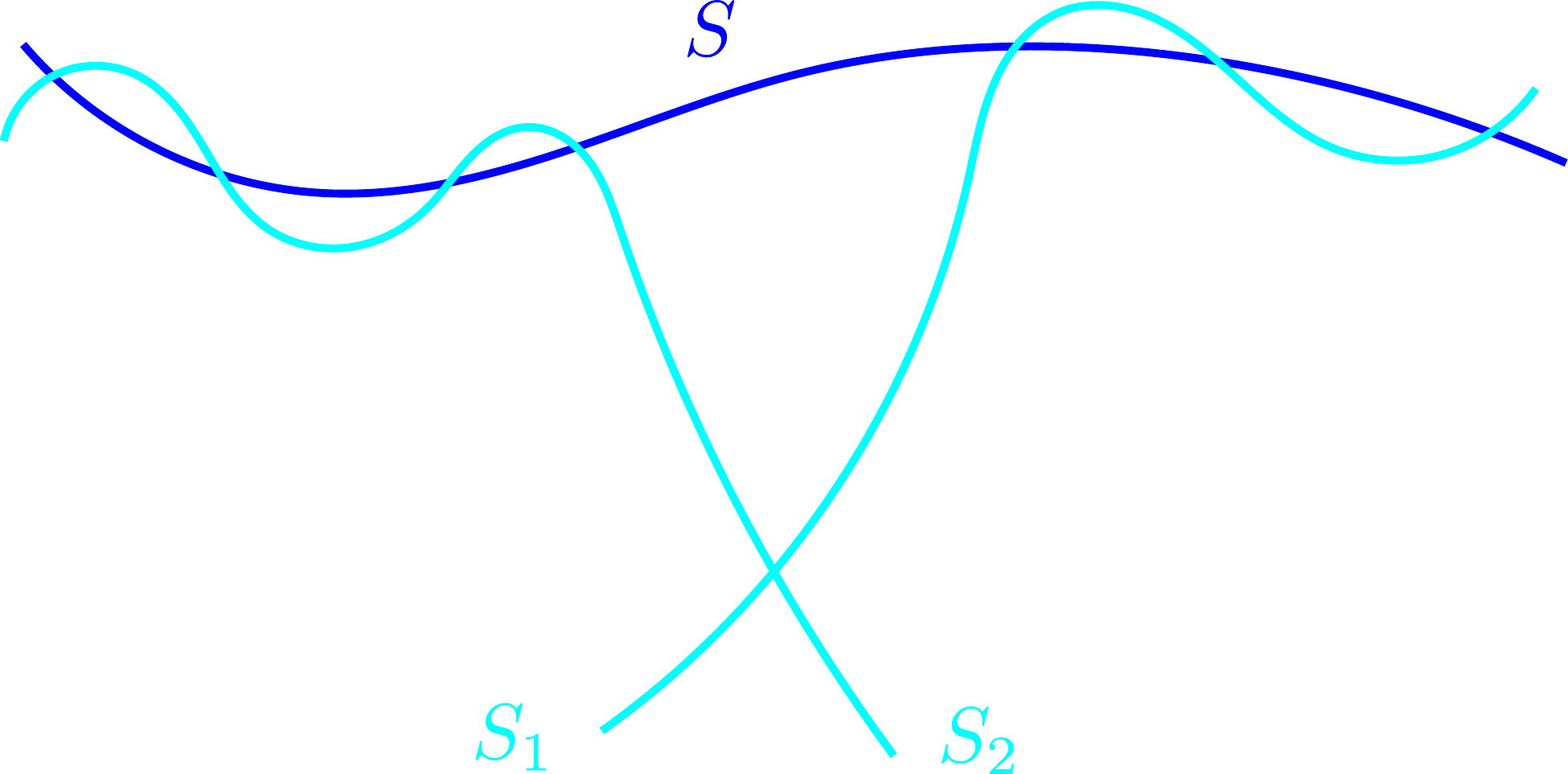}}\qquad
	\subfloat[\'Etape de type B.]{\label{f:typeB}\includegraphics[width=0.45\linewidth]{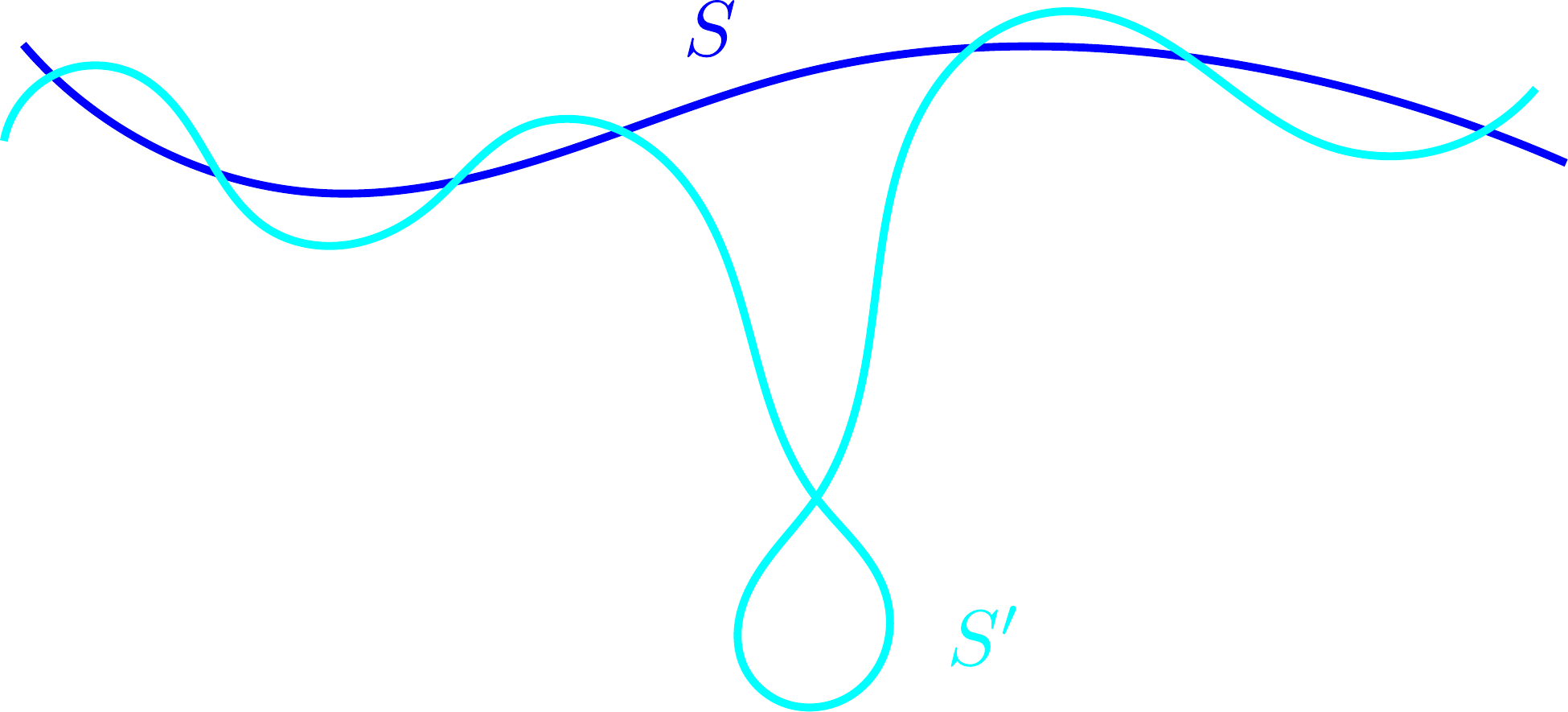}}
\caption{Les deux types d'étapes possibles au cours de la récurrence dans la démonstration du Théorème~\ref{t:ruledFKversion}.}\label{f:type}
\end{figure}

Avec ces rappels en tête, discutons les conditions d'apparition d'un tore $T$ symplectiquement plongé d'auto-intersection $3$ au cours de la récurrence. Si la récurrence ne faisait intervenir que des étapes de type B, alors on aurait $[S]^2= 4(g-1) +3 = 4g-1$, ce qui est impossible par hypothèse. (Il suffit de contracter les diviseurs exceptionnels apparus à chacune des étapes pour que l'image de $T$ par ces contractions soit homologue à $S$. Le genre des courbes considérées diminuant de $1$ à chaque étape de type B, ces diviseurs exceptionnels sont au nombre de $g-1$. Ils intersectent tous $T$ exactement deux fois, de manière transverse et positive. Voir la Figure~\ref{f:torusint3}.)
\begin{figure}[h]
	\centering
	\includegraphics[scale=0.5]{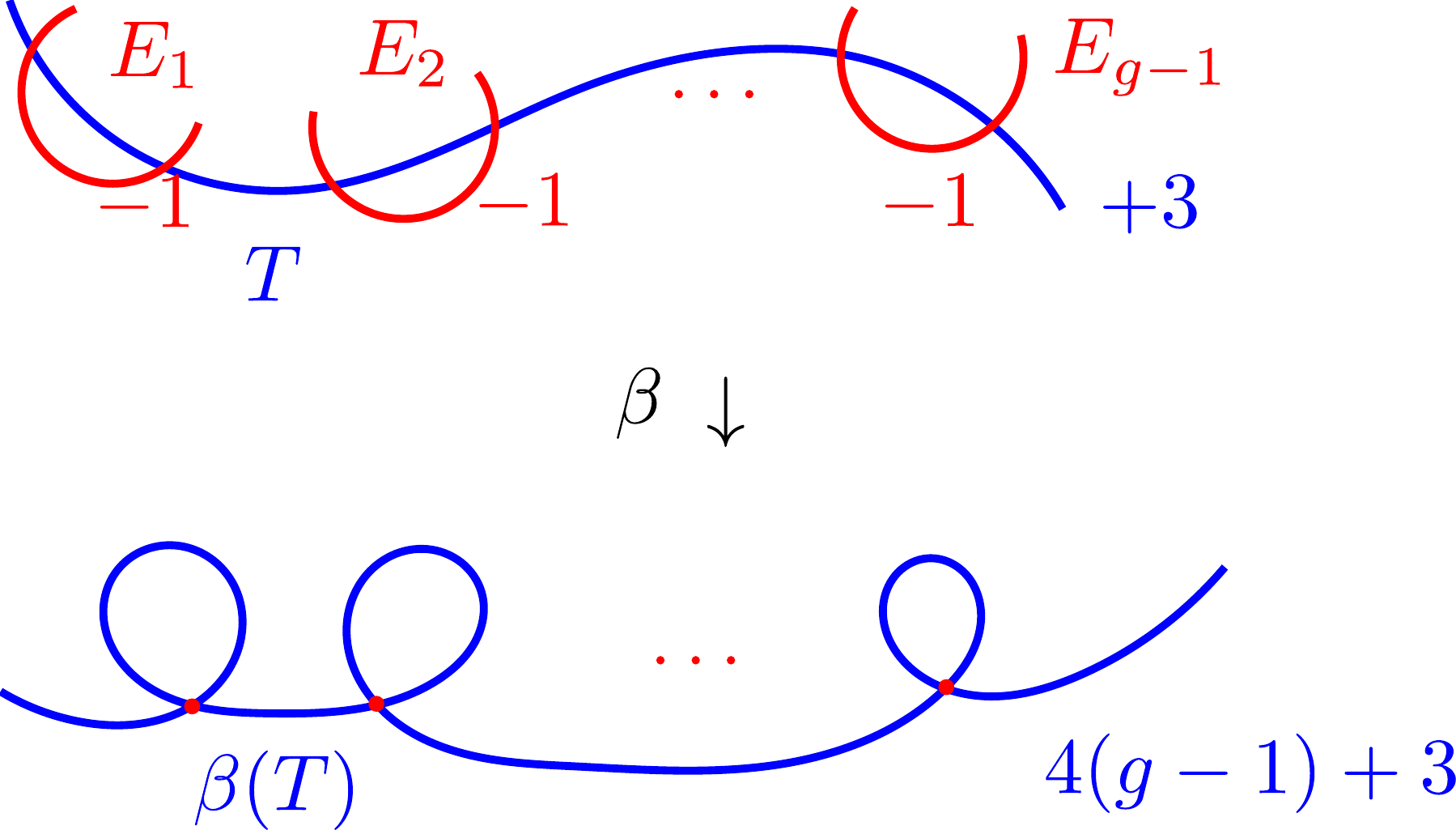}\\
	\caption{Dans le cas où la récurrence ne fait intervenir que des étapes de type B et où un tore $T$ d'auto-intersection $3$ apparaît, on peut contracter les diviseurs exceptionnels $E_1, \dots, E_{g-1}$ apparus au cours des étapes précédentes. L'image de $T$ par cette contraction $\beta$ est une courbe immergée d'auto-intersection $4g-1$ qui est homologue à $S$, ce qui fournit une contradiction.}
	\label{f:torusint3}
\end{figure}

Il y a donc au moins eu une étape de type A au cours de la récurrence. Considérons la dernière étape de type A avant l'apparition de $T$. On note $C$ la courbe de genre $\g$ à laquelle est appliquée cette étape, et on note $S_1$ et $S_2$ les deux composantes symplectiquement plongées, de genres respectifs $g_1$ et $g_2$ avec $\g = g_1 + g_2$, apparues suite à cette étape. La courbe $T$ apparaît ensuite après avoir effectué un certain nombre d'étapes de type B à partir d'une des deux composantes, disons $S_1$ par exemple. On a alors $$[S_1]^2 = 4 g_1 -1 \quad \text{et} \quad [S_2]^2 = [C]^2 - 2 -[S_1]^2 \geq 4g_2 - 2.$$
On choisit alors une structure presque complexe $J$ dominée par la forme symplectique telle que $S_2$ et $T$ sont des images de courbes $J$--holomorphes plongées. Puisque $g_2 >0$ (sinon la récurrence est déjà terminée), on a $[S_2]^2 > 2g_2-2$, donc la courbe $J$--holomorphe paramétrant $S_2$ est automatiquement Fredholm régulière. Il en est de même pour la courbe paramétrant $T$, on peut donc supposer, quitte à perturber légèrement $J$ et les deux courbes $J$--holomorphes, que $J$ est générique. Comme $[S_2] \cdot [T] = 1$, on peut appliquer successivement la Proposition~\ref{p:autoint3} et la Proposition~\ref{p:existeunidonale} à $T$ pour trouver soit une courbe rationnelle symplectiquement plongée d'auto-intersection positive (dans le cas d'une étape de type A), soit une courbe rationnelle symplectiquement positivement immergée d'auto-intersection $3$ avec un unique point double transverse (dans le cas d'une étape de type B).

À la fin de la récurrence, on trouve une courbe symplectique rationnelle $S_0$ qui est soit plongée d'auto-intersection positive, soit positivement immergée avec un unique point double (si $S_0$ est obtenue suite à une étape de type B appliquée à un tore d'auto-intersection $3$). Dans le premier cas, on a d'après la formule d'adjonction $c_1([S_0]) =2 + [S_0]^2 \geq 2$. Dans le second cas, on a $[S_0]^2 =3$ et la formule d'adjonction pour les courbes symplectiques singulières nous donne alors $c_1([S_0]) = 2 +  [S_0]^2 -2 = 3$. Pour terminer, on contracte les diviseurs exceptionnels disjoints de $S$ qui sont apparus au cours de la récurrence (qui correspondent aux éventuels éclatements de points doubles effectués à chaque étape  de type B de la récurrence). De cette manière, on retrouve bien la surface $(M, \omega)$. Quitte à perturber légèrement $S_0$ et à réaliser des isotopies de diviseurs exceptionnels (en s'aidant d'une structure presque complexe qui rend $S_0$ pseudoholomorphe et du Théorème~\ref{t:divexisotopy2}), on peut supposer que tous ces diviseurs exceptionnels intersectent $S_0$ au plus deux fois, de manière transverse et positive en des points non singuliers. L'image $S_0'$ de $S_0$ par la contraction de ces diviseurs exceptionnels est alors une courbe rationnelle symplectiquement positivement immergée satisfaisant $c_1([S_0']) \geq 2$ (la contraction ne peut qu'augmenter la classe de Chern car d'après la formule d'adjonction, tout diviseur exceptionnel $E$ vérifie $c_1(E) = 2 -1 =1 >0$). On conclut en appliquant le Théorème~\ref{t:McDuffimmersed} de McDuff sur les courbes rationnelles symplectiquement positivement immergées.
\end{proof}

\subsection{Le plongement est une section}\label{ss:plongementISsection}

Maintenant qu'on sait que $(M, \w )$ est une surface symplectiquement réglée, il nous reste à montrer que $S$ est une section d'un fibré en sphères symplectique sur $(M, \omega)$ au-dessus d'une surface de genre $g$ (ou une section d'une fibration de Lefschetz symplectique avec fibres sphériques sur $(M, \omega)$ au-dessus d'une base de genre $g$ si on enlève l'hypothèse de minimalité relative). Les arguments présentés dans cette sous-section sont inspirés de ceux utilisés dans le cas algébrique dans~\cite{hartshorne}. Le principe consiste à étudier en premier lieu les cas minimaux et les courbes symplectiques singulières qu'ils contiennent. Pour une telle courbe $C$ de genre arithmétique $p_a(C)$, on trouve une borne supérieure sur le nombre d'auto-intersection de $C$ en fonction du degré de l'application de projection du fibré en sphères symplectique correspondant et de $p_a(C)$ (le cas de $\mathbb{C}P^2$ étant traité à part). Puis en second lieu, on vérifie que les inégalités ainsi obtenues sont préservées par les opérations d'éclatement. On montrera ainsi -- entre autres -- que si le degré de l'application de projection est supérieur ou égal à $2$, l'auto-intersection de $C$ est inférieure ou égale $4 p_a (C) +4$, où $p_a (C)$ désigne le genre arithmétique de $C$. Ces arguments s'adaptent bien au cadre symplectique grâce aux techniques pseudoholomorphes, puisqu'ils reposent essentiellement sur les propriétés suivantes : la positivité d'intersection entre certaines courbes privilégiées dans la surface, la formule d'adjonction et les opérations d'éclatement ou de contraction.

Dans le cas où le raisonnement de la partie précédente fait apparaître une courbe rationnelle $S_0$ d'auto-intersection $0$, une étude une plus fine sur le nombre d'intersection entre $S_0$ et $S$ nous permet de prendre un raccourci dans les arguments de Hartshorne (lorsqu'on étudie les courbes dans les surfaces rationnelles), ce qui simplifie grandement la preuve.

\subsubsection{\'{E}tude du nombre d'intersection entre $S$ et $S_0$} \label{ss:SinterS0} 

On réitère ici notre argument de \og cassage \fg{} des courbes, en imposant une condition légèrement plus forte sur le nombre d'auto-intersection, en sachant désormais \emph{a posteriori} que $(M, \w )$ est une surface symplectiquement réglée et en examinant attentivement les intersections entre les différentes composantes qui apparaissent au cours de la récurrence. L'objectif de cette sous-section est d'aboutir à la proposition suivante.

\begin{prop} \label{p:S0interS}
Soit $S$ une courbe symplectiquement plongée de genre $g$ dans une surface symplectiquement réglée $(M, \w)$. On suppose $(M, \w)$ relativement minimale par rapport à $S$ et $[S]^2 \geq 4g +5$. Si $M$ n'est pas difféomorphe à $\mathbb{C}P^2$, alors il existe une courbe rationnelle symplectiquement plongée $S_0$ d'auto-intersection $0$ qui intersecte localement positivement $S$ telle que $[S] \cdot [S_0] \leq g +1$.
\end{prop}

\begin{rk}
L'hypothèse sur la surface symplectique $(M, \w)$ d'être symplectiquement réglée est redondante d'après le Théorème~\ref{t:ruledFKversion2}. On choisit néanmoins de conserver cette hypothèse dans l'énoncé afin de souligner que la Proposition~\ref{p:S0interS} est obtenue \emph{a posteriori}.
\end{rk}

Pour démontrer cette proposition, on utilisera le lemme suivant. Celui-ci nous permettra de comprendre de manière plus précise les composantes qui apparaissent lors du \og cassage \fg{} des courbes rationnelles d'auto-intersection positive.

\begin{lem} \label{l:detneg}
Soit $M$ une variété de dimension $4$ vérifiant $b_2^+ (M)=1$. Alors pour tous $A,B \in H_2(M; \mathbb{Z})$ tels que $A^2 > 0$, on a $A^2 B^2 \leq (A \cdot B)^2$.
\end{lem}

\begin{proof}
La matrice de la forme d'intersection de $M$ restreinte au sous-espace de $H_2(M;\mathbb{R})$ engendré par $A$ et $B$ est une sous-matrice de la matrice
$\begin{pmatrix}
A^2 & A \cdot B \\
A \cdot B & B^2 \\
\end{pmatrix}$. 
Puisque $b_2^+ (M)=1$ et $A^2>0$, le déterminant de cette matrice est inférieur ou égal à $0$.
\end{proof}

\begin{rk}
Toutes les surfaces symplectiquement réglées vérifient les hypothèses du Lemme~\ref{l:detneg}.
\end{rk}

\begin{lem} \label{l:casseenunefoisrationnelle}
Soit $S$ une courbe rationnelle symplectiquement plongée dans une surface symplectique $(M, \omega )$. On suppose $(M, \omega)$ relativement minimale par rapport à $S$. Si $[S]^2 \geq 5$, alors il existe une courbe rationnelle symplectiquement plongée $S_0$ d'auto-intersection $0$ qui intersecte localement positivement $S$ telle que $[S] \cdot [S_0] = 1$.
\end{lem}

\begin{proof}
D'après le Théorème~\ref{t:pinceaudeLefschetz}, on peut trouver deux courbes rationnelles symplectiquement plongées $S_1$ et $S_2$ telle que l'union de $S$, $S_1$ et $S_2$ forme une courbe symplectique singulière, $[S] =[S_1] + [S_2]$ et $[S_1] \cdot [S_2] = 1$. Comme les courbes $S_1$ et $S_2$ sont obtenues comme images de courbes pseudoholomorphes plongées Fredholm régulières, le Lemme~\ref{l:noexdiv} nous assure que $[S_1]^2 \geq 0$ et $[S_2]^2 \geq 0$.
Mais comme $5 \leq [S]^2 = [S_1]^2 +2 + [S_2]^2 $, on a $[S_1]^2 \geq 2$ ou $[S_2]^2 \geq 2$.
Le Lemme~\ref{l:detneg} nous donne alors $[S_1]^2 [S_2]^2 \leq 1$, d'où $[S_1]^2 =0$ ou $[S_2]^2 = 0$. Si on note $S_0$ la courbe parmi $S_1$ et $S_2$ qui est d'auto-intersection $0$, on obtient $[S] \cdot [S_0] = ([S_1] + [S_2]) \cdot [S_0] = 1$.
\end{proof}

\begin{rk}
Si, dans les hypothèses du Lemme~\ref{l:casseenunefoisrationnelle}, on suppose seulement que $[S]^2 \geq 2$, alors on a $2 \leq [S]^2 = [S_1]^2 +2 + [S_2]^2 $. Si $[S_1]^2 \neq 0$ et $[S_2]^2 \neq 0$, alors le Lemme~\ref{l:detneg} nous assure que $[S_1]^2 [S_2]^2 \leq 1$, d'où $[S_1]^2 = [S_2]^2 =1$. Dans tous les cas, une des deux courbes $S_1$ et $S_2$ est une courbe rationnelle symplectiquement plongée d'auto-intersection $0$ ou $1$. On a alors démontré \emph{a posteriori} que la récurrence effectuée lors de la démonstration du Théorème~\ref{t:McDuff} de McDuff et Gromov se termine en réalité dès la première étape.
\end{rk}

\begin{proof}[Démonstration de la Proposition~\ref{p:S0interS}]
On procède par récurrence sur le genre de $S$. L'initialisation fait l'objet du Lemme~\ref{l:casseenunefoisrationnelle}. On suppose désormais $g \geq 1$ et la proposition vérifiée pour toutes les courbes de genre strictement inférieur à $g$.

On commence par réexaminer la démonstration par récurrence présentée dans la Sous-section~\ref{ss:recurrence}.
On suppose dans un premier temps qu'on effectue une étape de type A à un certain moment au cours de la récurrence présentée dans la Sous-section~\ref{ss:recurrence}. Une fois la première étape de type $A$ effectuée, on obtient une courbe symplectique singulière formée de l'union de $S$, de deux courbes symplectiquement plongées d'auto-intersection positive $S_1$ et $S_2$ satisfaisant $[S_1] \cdot [S_2] = 1$, ainsi que des diviseurs exceptionnels deux à deux disjoints $E_1, \dots , E_\ell$ qui apparaissent lorsqu'on éclate un point double au cours de chaque étape de type B précédant cette étape de type A. Notons que $\ell \leq g-1$ car l'entier $\ell+1$ correspond au nombre d'étapes effectuées dans la récurrence jusqu'à cette étape de type A (incluse). Chaque diviseur exceptionnel $E_i$ est disjoint de $S$ et intersecte l'union de $S_1$ et $S_2$ en exactement deux points, de manière transverse et positive (voir la Figure~\ref{f:findS0a}). 
\begin{figure}[h]
	\centering
	\includegraphics[scale=0.5]{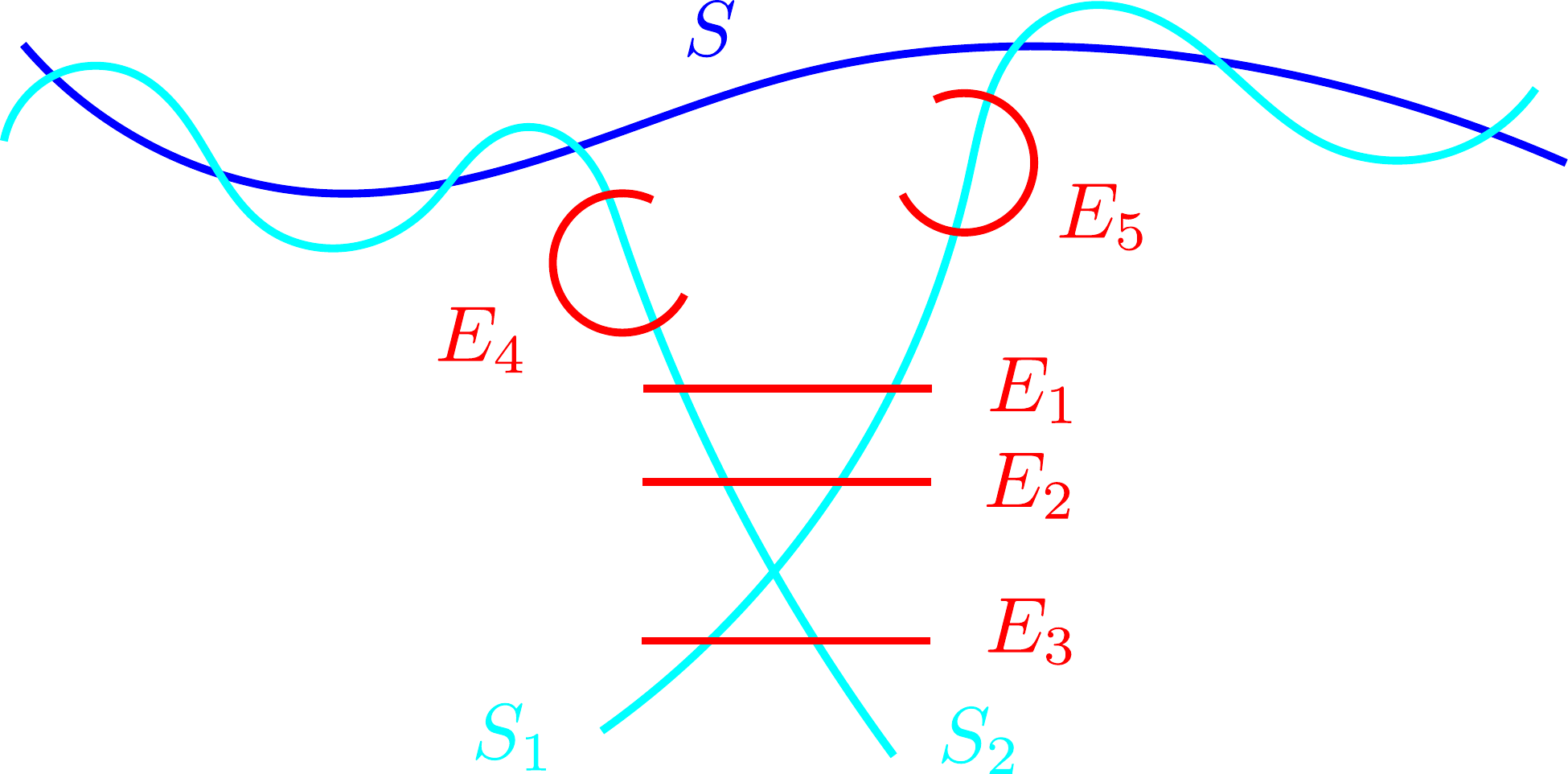}\\
	\caption{Exemple possible de configuration de courbes obtenue au cours de la récurrence après la première étape de type A.}
	\label{f:findS0a}
\end{figure}
Parmi les diviseurs exceptionnels $E_1, \dots , E_\ell$, on contracte ceux qui n'intersectent qu'une seule des deux courbes $S_1$ ou $S_2$ (et par conséquent intersectent la courbe en question exactement deux fois, de manière transverse et positive). On lisse ensuite les points doubles des images de $S_1$ et $S_2$ par ces contractions afin d'obtenir des courbes symplectiquement plongées $\tilde S_1$ et $\tilde S_2$ d'auto-intersection positive (en effet, on a pour tout $i \in \{1,2 \}$, $[\tilde S_i]^2 \geq [S_i]^2 \geq 0$). Quitte à réindexer, on peut supposer que $E_1, \dots,E_k$ sont les diviseurs exceptionnels restants qu'on n'a pas contractés. On obtient alors dans $M \# k\overline{\mathbb{C}P}^2$ une courbe symplectique singulière formée de l'union de $S$, de $ E_1, \dots, E_k$ et de $\tilde S_1$, $\tilde S_2$ telle que $[S]= [\tilde S_1] + [\tilde S_2] + 2( [E_1] + \dots + [E_k]) $ et $[ \tilde S_1] \cdot [ \tilde S_2] = 1$. De plus, chaque $E_j$ intersecte chaque courbe $\tilde S_i$ exactement une fois, de manière transverse et positive. 
On en déduit que $[S]^2 = [\tilde S_1]^2 +[\tilde S_2]^2 -4 k +2 + 8 k$, donc $[\tilde S_1]^2 + [\tilde S_2]^2 \geq 4(g-k)+3 \geq 7$. Quand on contracte les diviseurs exceptionnels $E_1, \dots, E_k$ et qu'on applique ensuite le Lemme~\ref{l:detneg} aux classes d'homologies des images $\tilde S_1'$ et $\tilde S_2'$ de $\tilde S_1$ et $\tilde S_2$ par ces contractions (voir la Figure~\ref{f:findS0c}), on obtient 
\begin{figure}[h]
	\centering
	\includegraphics[scale=0.5]{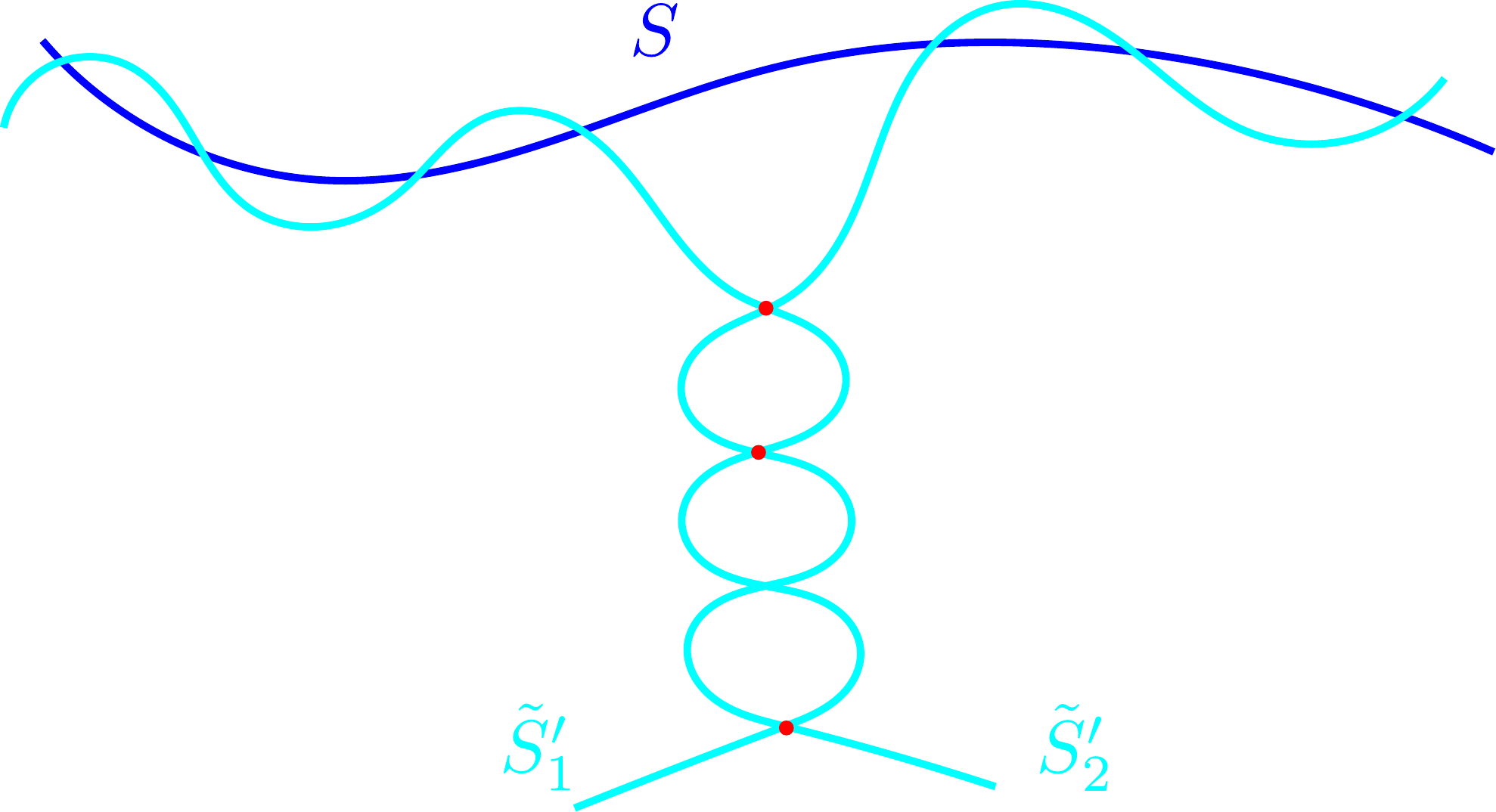}\\
	\caption{Configuration de courbes obtenue à partir de la configuration de courbes présentée dans la Figure~\ref{f:findS0a} en contractant les diviseurs exceptionnels, puis en lissant les points doubles de chaque composante.}
	\label{f:findS0c}
\end{figure}
$$[\tilde S_1']^2 [\tilde S_2']^2= \left( [\tilde S_1]^2+k \right) \left( [ \tilde S_2]^2+k \right) \leq (k+1)^2 = \left( [\tilde S_1'] \cdot [\tilde S_2'] \right)^2.$$
Ainsi on a $[\tilde S_1]^2 [\tilde S_2]^2 +k ([\tilde S_1]^2+[\tilde S_2]^2-2) \leq 1$. Comme $[\tilde S_1]^2 [\tilde S_2]^2 \geq 0$ et $[\tilde S_1]^2 + [\tilde S_2]^2 -2 \geq 5$, on a nécessairement $k=0$ et une des courbes $\tilde S_i$, disons $\tilde S_1$, est d'auto-intersection $0$. Notons que $\tilde S_1$ étant obtenue comme l'image d'une courbe pseudoholomorphe plongée Fredholm régulière, elle satisfait la condition de positivité des indices, c'est-à-dire $\chi (\tilde S_1) \geq - 2 [\tilde S_1]^2 =0$. Donc $\tilde S_1$ est de genre $0$ ou $1$. Si $\tilde S_1$ est une courbe rationnelle, on a $[S] \cdot [\tilde S_1] = [\tilde S_1] ^2 +[\tilde S_1] \cdot [\tilde S_2] =1$, ce qui permet de conclure. Si $\tilde S_1$ est de genre $1$, le genre de $\tilde S_2$ est égal à $g-1$ et on a alors $[\tilde S_2]^2 = [S]^2-2 \geq 4 (g-1) + 7$. Par hypothèse de récurrence (pour garantir l'hypothèse de minimalité relative, on contracte une collection maximale de diviseurs exceptionnels deux à deux disjoints et disjoints de $\tilde S_2$, puis on inverse cette procédure avec des éclatements une fois l'hypothèse de récurrence appliquée), on dispose d'une courbe $S_0$ rationnelle symplectiquement plongée dans $(M, \omega)$ d'auto-intersection $0$ qui vérifie $[\tilde S_2] \cdot [S_0] \leq g$. La courbe $S_0$ étant obtenue comme composante d'une courbe nodale homologue à $\tilde S_2$ via des techniques pseudoholomorphes, elles satisfait la propriété de positivité d'intersection avec les autres courbes considérées.
Finalement, on a $[\tilde S_1] \cdot [S_0] \leq [\tilde S_1] \cdot [\tilde S_2] =1$ et par conséquent, $[S] \cdot [S_0] = [\tilde S_1] \cdot [S_0] + [\tilde S_2] \cdot [S_0] \leq g+1$.

On suppose désormais qu'on n'effectue que des étapes de type $B$ au cours de la récurrence présentée dans la Sous-section~\ref{ss:recurrence}. On obtient alors au bout de $g$ étapes une courbe symplectique singulière donnée par l'union de $S$, d'une courbe rationnelle symplectiquement plongée $S'$ et des diviseurs exceptionnels deux à deux disjoints $E_1, \dots , E_g$ qui apparaissent lorsqu'on éclate un point double au cours de chaque étape de type B. Chaque diviseur exceptionnel $E_i$ est disjoint de $S$ et intersecte $S'$ en exactement deux points, de manière transverse et positive. On a de plus $[S] = [S'] + 2( [E_1] + \dots + [E_{g}])$ et $[S']^2 = [S]^2 -4g \geq 5$. Comme la surface symplectique ambiante est relativement minimale par rapport à $S'$, le Théorème~\ref{t:McDuff} et le Lemme~\ref{l:noexdiv} nous assurent l'existence de deux courbes rationnelles symplectiquement plongées d'auto-intersection positive $S_1'$ et $S_2'$, telle que l'union de $S$, $S'$, $S_1'$, $S_2'$ et des $E_i$ forme une courbe symplectique singulière, $[S'] =  [S_1'] + [S_2']$ et $[S_1 '] \cdot [S_2'] = 1$.  Notons qu'on peut supposer sans perte de généralité que chaque diviseur exceptionnel $E_i$ intersecte l'union des $S_i'$ en exactement deux points, de manière transverse et positive. Autrement dit, on s'est ramené à une situation similaire à celle traitée précédemment. En reprenant les mêmes notations, la seule différence avec le cas précédent est que cette fois-ci on a seulement $[\tilde S_1]^2 +[ \tilde S_2]^2 =[S']^2-2 \geq 3$. Ainsi, l'inégalité $[\tilde S_1]^2 [\tilde S_2]^2 +k ([\tilde S_1]^2+[\tilde S_2]^2-2) \leq 1$ nous indique que $k=0$, ou bien $k=1$, $[\tilde S_1]^2 + [\tilde S_2]^2 =3$ et une des deux courbes $\tilde S_i$ est d'auto-intersection $0$. Si $k=0$, on conclut comme précédemment. Si on se trouve dans le deuxième cas de figure, remarquons qu'on a $k=g$, $\tilde S_1 =S_1'$ et $\tilde S_2 = S_2'$ (en effet la contraction d'un diviseur exceptionnel qui intersecterait une courbe $S_i'$ transversalement deux fois contribuerait à augmenter de $4$ l'auto-intersection de $\tilde S_i$). La contraction de l'unique diviseur exceptionnel $E_1$, qui intersecte chacune des courbes $S_i'$ transversalement exactement une fois, fait alors apparaître une courbe rationnelle symplectiquement plongée $L$ d'auto-intersection $1$ et une courbe rationnelle symplectiquement plongée $C$ d'auto-intersection $4$. De plus, les courbes $L$ et $C$ s'intersectent en exactement deux points, de manière transverse et positive. Par le Théorème~\ref{t:McDuff} de McDuff et Gromov, $(M, \omega)$ est donc un éclaté d'une déformation symplectique de $(\mathbb{C}P^2, \omega_{FS})$ et $L$ représente la classe d'homologie d'une droite. Comme $g=1$, $[S]^2 = [S']^2 +4 = 9$, $[S] \cdot [L] = [L]^2 +[L] \cdot [C] =3$ et $(M, \omega)$ est relativement minimale par rapport à $S$, on obtient alors que la surface symplectique $(M, \omega)$ est minimale. Par conséquent $M$ est difféomorphe à $\mathbb{C}P^2$, ce qui est exclu par hypothèse.
\end{proof}

\subsubsection{Fibrés en sphères symplectiques au-dessus de surfaces}

On considère une surface symplectique $(N, \w)$ munie d'une structure de fibré en sphères symplectique $\pi : N \rightarrow B$, d'une section symplectique d'auto-intersection positive $\mathcal{S}$, ainsi que d'une structure presque complexe $J$ dominée par $\omega$ qui rend $\mathcal{S}$ et toutes les fibres $J$--holomorphes.

\begin{rk}
On a vu lors de la Remarque~\ref{r:paritéréglée} que $(N, \w)$ admet une section lisse d'auto-intersection positive $\mathcal{S}$. On choisit une forme d'aire $\sigma$ sur $\Sigma$. D'après~\cite[Lemma~3.41, Remark 3.42]{Wendl}, pour tout réel $K \geq 0$, la $2$--forme différentielle
$$\omega_K := \omega +K \pi^* \sigma$$
est une forme symplectique compatible avec $\pi$
et, pour tout $K\geq0$ suffisamment grand, la section $\mathcal{S}$ est une sous-variété symplectique de $(N, \w_K)$. Par conséquent, quitte à effectuer une déformation symplectique, on peut toujours supposer qu'une telle section $\mathcal{S}$ est symplectique.
\end{rk}

On rappelle que le second groupe d'homologie $H_2 (N; \mathbb{Z})$ est engendré par $[\mathcal{S}]$ et la classe d'homologie d'une fibre $[F]$. On rappelle aussi que d'après le Théorème~\ref{t:classificationsurfacesreglees}, mis à part celles difféomorphes à $\mathbb{C}P^2$, toutes les surfaces symplectiquement réglées sont des éclatés de fibrés en sphères symplectiques au-dessus de surfaces.

On s'intéresse aux propriétés des courbes $J$--holomorphes dans $(N, \w)$. Notons $C$ l'image d'une telle courbe et $\g$ le genre de $B$. On dispose d'entiers $k$ et $\ell$ tels que 
$$[C] = k [ \mathcal{S}] + \ell [F],$$ 
d'où 
$$c_1([C]) = k c_1([ \mathcal{S}]) + \ell c_1([F]) \quad \text{et}\quad [C]^2 = k^2 [ \mathcal{S}]^2 +2k\ell.$$ En appliquant la formule d'adjonction de chaque côté de l'égalité, on obtient 
$$[C]^2 + 2 - 2 p_a (C) = k([\mathcal{S}]^2 + 2 - 2 \g) + 2 \ell,$$ c'est-à-dire 
$$k^2[\mathcal{S}]^2 + 2k\ell + 2 - 2 p_a (C) = k[\mathcal{S}]^2 +k( 2 - 2 \g) + 2 \ell.$$ En regroupant tous les termes du même côté, on obtient l'égalité suivante 
$$k(k-1) [\mathcal{S}]^2 +2(k-1)\ell + 2- 2 p_a (C) - k(2-2 \g) = 0,$$
puis en multipliant par $k$ et en reconnaissant l'expression de $[C]^2$, on trouve
\begin{equation} \label{eq:courbedansfibré}
(k-1) [C]^2  + k \left(2- 2 p_a (C) - k(2-2 \g) \right) = 0.
\end{equation}

\begin{rk}\label{r:projofdegk}
Puisque $[F]^2=0$ et $[\mathcal{S}] \cdot [F]=1$, on a $k = [C] \cdot [F]$. Par positivité d'intersection, on obtient alors $k \geq 0$. Notons également que $k$ correspond au nombre de points d'intersection, comptés avec multiplicité, entre $C$ et chacune des fibres. Autrement dit, $k$ correspond au nombre d'antécédents, comptés avec multiplicité, de chaque point de $B$ par la restriction de l'application de projection $\pi$ à $C$. Ainsi, on a $k = \deg ( \pi_{\vert C}) \geq 0$.
\end{rk}

On va maintenant distinguer les cas selon que $\g = 0$ (i.e. $(N, \w)$ est une surface symplectique rationnelle, d'après le Théorème~\ref{t:classificationsurfacesrationnelles}) ou $\g > 0$ (i.e. $(N, \w)$ est une surface symplectique réglée non rationnelle, d'après le Théorème~\ref{t:classificationsurfacesrationnelles}) pour obtenir des contraintes sur $[C]^2$.

\subsubsection{Courbes symplectiques singulières dans les surfaces symplectiquement réglées non rationnelles}

Dans cette sous-section, on utilise des techniques pseudoholomorphes afin de reprendre, en les simplifiant légèrement, les arguments utilisés dans le cadre algébrique par Hartshorne dans \cite{hartshorne}.

\begin{prop}\label{p:courbesdanssurfacesNR}
Soit $(N, \w)$ une surface symplectiquement réglée non rationnelle (i.e. une fibration de Lefschetz symplectique générique avec fibres sphériques et base de genre strictement positif) et $C$ une courbe symplectique singulière irréductible dans $(N, \w)$ de genre arithmétique $p_a(C)$. On note $k \in \mathbb{N}$ le degré de l'application de projection $\pi$ restreinte à $C$. Si $k \geq 2$, alors $$ [C]^2 \leq \frac{2k}{k-1} (p_a(C) -1).$$
\end{prop}

\begin{proof}
On choisit tout d'abord une structure presque complexe $J$ dominée par $\omega$, générique en dehors d'un voisinage de $C$ (voir la Remarque~\ref{rk:généricitéhorsrelcompact}), telle que $C$ est l'image d'une courbe $J$--holomorphe. Comme $C$ n'est pas homologue à une composante irréductible d'une fibre (sinon on aurait $k=0$), on peut choisir ce voisinage de manière à ce qu'il ne contienne aucune courbe $J$--holomorphe homologue à une composante d'une fibre. Grâce au Théorème~\ref{t:isotopiepinceaudeLefschetz}, on peut alors supposer sans perte de généralité (quitte à réaliser une isotopie de fibration de Lefschetz) que chaque composante des fibres est $J$--holomorphe. De cette façon, on a la propriété de positivité d'intersection entre $C$ et chacune des fibres.

On commence par traiter le cas où $(N,\w)$ est minimale, c'est-à-dire quand $(N,\w)$ est un fibré en sphères symplectique. Dans ce cas, on a d'après l'égalité \eqref{eq:courbedansfibré} (en reprenant les mêmes notations),
$$(k-1) [C]^2  + k \left(2- 2 p_a (C) - k(2-2 \g) \right) = 0,$$
donc 
$$(k-1) [C]^2  = k \left( 2 p_a (C)-2 \right) - k^2(2 \g-2).$$
Puisque $\g \geq 1$, on a alors
$$(k-1) [C]^2  \leq k \left( 2 p_a (C)-2 \right).$$
On en conclut que 
$$[C]^2 \leq \frac{k}{k-1}(2p_a(C)-2),$$
ce qui est l'inégalité voulue.


Pour le cas général, on raisonne par récurrence sur le nombre de contractions à effectuer pour aboutir à un modèle minimal. Notons $\beta_p : N \rightarrow N'$ une contraction, où $p$ correspond à l'image du diviseur exceptionnel contracté. On note $C' = \beta_p (C)$. On suppose que l'inégalité est vraie pour $C'$, c'est-à-dire 
$$[C']^2 \leq \frac{2k}{k-1} (p_a(C') -1).$$
Notons $r \geq 0$ la multiplicité de $p$ en tant que point de $C'$. On a $[C]^2 = [C']^2 -r^2$ et la formule d'adjonction nous donne $p_a(C) = p_a (C') - \frac{1}{2}r (r-1)$. On obtient alors
$$[C]^2 \leq \frac{2k}{k-1} (p_a(C) -1) + \frac{k}{k-1} r(r-1) -r^2.$$
Par positivité d'intersection, la multiplicité de tout point de $C'$ est bornée par le degré de la projection $\pi$ restreinte à $C'$, qui est égal à $k$ (on rappelle que $k=[C'] \cdot [F]$, voir la Remarque~\ref{r:projofdegk}). On a donc $r \leq k$ et 
$$\frac{k}{k-1} r(r-1) -r^2 = \frac{r}{k-1} (r-k) \leq 0,$$
ce qui permet de conclure.
\end{proof}

\begin{cor}\label{c:courbesdanssurfacesNR}
Soit $(N, \w)$ une surface symplectiquement réglée non rationnelle (i.e. une fibration de Lefschetz symplectique générique avec fibres sphériques et base $B$ de genre strictement positif), munie d'une structure presque complexe $J$ dominée par $\omega$ qui rend chaque composante des fibres $J$--holomorphe.
Soit $v : \S \rightarrow N$ une courbe $J$--holomorphe de genre $g>0$ plongée. On a alors une des deux possibilités suivantes :
\begin{itemize}
\item $\Sigma$ est difféomorphe à $B$ et $v$ est une section de la fibration de Lefschetz symplectique, ou bien
\item $[v]^2 \leq 4g -4$.
\end{itemize}
\end{cor}

\begin{proof}
On note $C$ l'image de $v$, $F$ une fibre générique et $k$ le degré de l'application de projection $\pi$ restreinte à $C$. 

Si $k=0$, alors on a $[F] \cdot [C] = 0$. Par positivité d'intersection, $C$ est alors une composante d'une fibre, ce qui est impossible puisque toutes les composantes des fibres sont de genre $0$.

Si $k=1$, la positivité d'intersection nous assure que $v$ intersecte chaque fibre transversalement exactement une fois. L'application $\pi_{| C} \circ v : \Sigma \rightarrow B$ est donc un difféomorphisme et $v$ est une section de la fibration de Lefschetz.

Si $k \geq 2$, on a $\frac{k}{k-1} \leq 2$, et la Proposition~\ref{p:courbesdanssurfacesNR} nous donne $[v]^2 \leq 4g - 4.$
\end{proof}

\subsubsection{Courbes symplectiques singulières dans les surfaces symplectiques rationnelles}

On étudie désormais les nombres d'auto-intersection possibles pour les courbes symplectiques singulières dans les surfaces symplectiques rationnelles. Cette étape se trouve être un peu plus délicate que l'étape précédente. On procède de la même manière que précédemment : on s'intéresse en premier lieu aux cas minimaux, puis on traite le cas général en examinant les effets des éclatements et des contractions. On s'inspire des arguments utilisés par Hartshorne dans~\cite{hartshorne}, mais on utilise les résultats de la Sous-section~\ref{ss:SinterS0} pour prendre des raccourcis et simplifier significativement la preuve.

\begin{prop}\label{p:courbesdansCP2}
Soit $C$ une courbe symplectique singulière irréductible de genre arithmétique $p_a(C)$ dans une déformation symplectique de $(\mathbb{C} P^2, \omega_{FS})$. Alors on a $[C]^2 \leq 4 p_a(C)+5$. De plus, on a $[C]^2 = 4 p_a(C)+5$ si et seulement si $C$ est une cubique (i.e. $p_a(C)=1$ et $[C]^2 = 9$).
\end{prop}

\begin{proof}
Notons $h \in H_2 (\mathbb{C} P^2; \mathbb{Z})$ la classe d'homologie de la droite et $d = h \cdot [C]$ (c'est le degré de la courbe). On a alors $[C] = dh$ et $[C]^2 = d^2$. En appliquant la formule d'adjonction à $C$, on obtient
$d c_1(h) = 2 - 2 p_a(C) + d^2$. Puis en appliquant la formule d'adjonction à $h$, on obtient $d(2 +1) = 2 - 2p_a(C) +d^2$, ce qui donne $$p_a(C) = \frac{d^2 -3d +2}{2} = \frac{(d-1)(d-2)}{2}.$$
On a alors les équivalences suivantes 
\begin{align*}
 d^2 \leq 4p_a(C)+5
 &\Leftrightarrow d^2 \leq 2(d-1)(d-2) +5,\\
 &\Leftrightarrow 0 \leq d^2 -6d +9,\\
 &\Leftrightarrow 0  \leq (d-3)^2,
\end{align*}
avec $d^2 = 4 p_a(C) +5$ si et seulement si $d=3$.
\end{proof}

\begin{prop}\label{p:courbesdanssurfacesRNCP2}
Soit $(N, \w)$ une surface symplectique rationnelle qui n'est pas difféomorphe à $\mathbb{C}P^2$ (i.e. une fibration de Lefschetz symplectique générique avec fibres sphériques et base de genre $0$) et $C$ une courbe symplectique singulière irréductible dans $(N, \w)$ de genre arithmétique $p_a(C)$. On note $k \in \mathbb{N}$ le degré de l'application de projection restreinte à $C$. Si $k \geq 2$, alors $$ [C]^2 \leq \frac{2k}{k-1}(p_a(C) -1 +k).$$
\end{prop}

\begin{proof}
Comme dans la démonstration de la Proposition~\ref{p:courbesdanssurfacesNR}, on peut supposer sans perte de généralité que $C$ et chacune des composantes irréductibles des fibres sont des images de courbes $J$--holomorphes pour une certaine structure presque complexe $J$ dominée par $\omega$.

On commence par traiter le cas où $(N, \w)$ est un fibré en sphères symplectique. D'après l'équation \eqref{eq:courbedansfibré} (avec $\g =0$), on a 
$$(k-1) [C]^2 + k (2-2p_a(C) -2k) =0.$$
Donc $$[C]^2 = \frac{2k}{k-1} (p_a(C) -1 +k).$$

Pour le cas général, on raisonne par récurrence sur le nombre de contractions à effectuer pour aboutir à un modèle minimal. Notons $\beta_p : N \rightarrow N'$ une contraction, où $p$ correspond à l'image du diviseur exceptionnel contracté. On note $C' = \beta_p (C)$. On suppose que l'inégalité est vraie pour $C'$, c'est-à-dire 
$$[C']^2 \leq \frac{2k}{k-1} (p_a(C') -1 + k).$$

Notons $r \geq 0$ la multiplicité de $p$ en tant que point de $C'$. On a $[C]^2 = [C']^2 -r^2$ et la formule d'adjonction nous donne $p_a(C) = p_a (C') - \frac{1}{2}r (r-1)$. On obtient alors
$$[C]^2 \leq \frac{2k}{k-1} (p_a(C) -1 +k) + \frac{k}{k-1} r(r-1) -r^2.$$
De la même façon que dans la démonstration de la Proposition~\ref{p:courbesdanssurfacesNR}, on a $r \leq k$, donc 
$$\frac{k}{k-1} r(r-1) -r^2 = \frac{r}{k-1} (r-k) \leq 0,$$
ce qui permet de conclure.
\end{proof}

\begin{thm}\label{t:courbesdanssurfacesR}
Soit $(N, \w)$ une surface symplectique rationnelle (i.e une fibration de Lefschetz symplectique générique avec fibres sphériques et base de genre $0$ ou une déformation symplectique de $(\mathbb{C}P^2, \omega_{FS})$, c'est-à-dire un pinceau de Lefschetz symplectique avec fibres sphériques possédant un unique point base et aucune fibre singulière), munie d'une structure presque complexe $J$ générique dominée par $\omega$ qui rend chaque composante des fibres $J$--holomorphe. Soit $v : \S \rightarrow N$ une courbe $J$--holomorphe plongée de genre $g$. On a alors une des trois possibilités suivantes :
\begin{itemize}
\item $\Sigma$ est difféomorphe à $\mathbb{C}P^1$ et $v$ est une section d'une fibration de Lefschetz $J$--holomorphe générique avec fibres sphériques et base de genre $0$ sur $(N,\w)$ (pas nécessairement celle de départ), ou bien
\item $g=1$ et, quitte à contracter des diviseurs exceptionnels disjoints de $v$, $v$ est une cubique (symplectique) de $\mathbb{C}P^2$, auquel cas $[v]^2 =9$, ou bien
\item $[v]^2 \leq 4g + 4$.
\end{itemize}
\end{thm}

\begin{proof}
On suppose sans perte de généralité que $(N,\omega)$ est relativement minimale par rapport à l'image de $v$.

Si $(N,\w)$ est une déformation symplectique de $(\mathbb{C}P^2, \omega_{FS})$, la Proposition~\ref{p:courbesdansCP2} nous assure que $[v]^2 \leq 4g+4$, avec pour seule exception la cubique lisse.

On suppose dorénavant que la surface symplectique $(N,\w)$ n'est pas une déformation symplectique de $(\mathbb{C}P^2, \omega_{FS})$. C'est donc une fibration de Lefschetz $J$--holomorphe avec fibres sphériques et base de genre $0$. Notons $C$ l'image de $v$, $F$ une fibre générique et $k$ le degré de l'application de projection $\pi$ restreinte à $C$. 

Si $k=0$, alors on a $[F] \cdot [C] = 0$. Par positivité d'intersection, $C$ est une composante d'une fibre, donc $[C]^2 \leq 0$.

Si $k=1$, la positivité d'intersection nous assure que $v$ intersecte chaque fibre transversalement exactement une fois. L'application $\pi_{| C} \circ v : \Sigma \rightarrow B$ est donc un difféomorphisme et $v$ est une section de la fibration de Lefschetz.

On suppose désormais $k \geq 2$ et $[v]^2 > 4g +4$. D'après la Proposition~\ref{p:S0interS}, il existe une courbe rationnelle symplectiquement plongée $S_0$ d'auto-intersection $0$ qui intersecte localement positivement $v$, telle que $[v] \cdot [S_0] \leq g +1$ (l'hypothèse de minimalité relative nous sert ici). La courbe rationnelle $S_0$ étant obtenue comme image d'une courbe $J$--holomorphe, on peut, en considérant l'espace de modules approprié, construire une fibration de Lefschetz $J$--holomorphe (potentiellement différente de la fibration de Lefschetz de départ ; même si ce genre de situation n'arrive que dans les éclatés de $\mathbb{C}P^1 \times \mathbb{C}P^1$, voir la Remarque~\ref{r:sphereautointgeq2}) dont $S_0$ est une fibre générique. Notons $h = [v] \cdot[S_0]$ le degré de la nouvelle application de projection $\pi'$ restreinte à $C$. Si $h=0$ ou $h=1$ on conclut comme ci-dessus (dans le cas où $h=1$, $v$ peut être une section d'une fibration de Lefschetz $J$--holomorphe différente de celle de départ). Si $h \geq 2$, la Proposition~\ref{p:courbesdanssurfacesRNCP2} nous donne $$[v]^2 \leq \frac{2h}{h-1}(g -1 +h).$$
On a alors $2hg -2h +2h^2 > (h-1)(4g+4)$, donc  $2h^2-2h -4h+4 > g(4h -4 -2h)$, d'où $h^2-3h+2 > g(h-2)$. Au final on obtient $(h-1)(h-2)>g(h-2)$, donc $h \neq2$ et $[v] \cdot [S_0] = h > 1+g$, ce qui est une contradiction.
\end{proof}

\begin{rk}
L'utilisation de la Proposition~\ref{p:S0interS}, qui affirme l'existence d'une courbe rationnelle symplectiquement plongée $S_0$ d'auto-intersection $0$ qui intersecte localement positivement $v$ et qui vérifie $[v] \cdot [S_0] \leq g +1$, constitue un raccourci substantiel dans les arguments de~\cite{hartshorne}.
\end{rk}

\subsubsection{Conclusion}

\begin{proof}[Démonstration du Théorème~\ref{t:thmprincipal}]
Commençons par remarquer que $M$ n'est pas difféomorphe à $\mathbb{C}P^2$ par la Proposition~\ref{p:courbesdansCP2}. D'après le Théorème~\ref{t:ruledFKversion} et la Remarque~\ref{r:ruledFKversion}, il existe alors une structure presque complexe générique $J$ dominée par $\omega$ telle que $S$ est l'image d'une courbe $J$--holomorphe plongée et $(M, \w)$ admet une fibration de Lefschetz $J$--holomorphe générique avec fibres sphériques (c'est donc une surface réglée symplectique). 

Si $(M, \w)$ est une surface symplectique non rationnelle, le Corollaire~\ref{c:courbesdanssurfacesNR} nous permet d'affirmer que $S$ est l'image d'une section de la fibration de Lefschetz. 
 
Si $(M, \w)$ est une surface symplectique rationnelle, le Théorème~\ref{t:courbesdanssurfacesR} nous assure que $S$ est l'image d'une section d'une fibration de Lefschetz $J$--holomorphe générique avec fibres sphériques (possiblement différente de celle de départ).

Dans les deux cas, il ne reste plus qu'à s'assurer que la fibration de Lefschetz n'admet pas de fibre singulière. Chaque fibre singulière de la fibration est l'union de deux diviseurs exceptionnels qui s'intersectent transversalement exactement une fois. Par positivité d'intersection et puisque $S$ est une section, on sait également que chaque fibre singulière possède une composante disjointe de $S$. Il ne peut donc pas exister de fibre singulière car $(M,\w)$ est relativement minimale par rapport à $S$. Autrement dit, $(M, \w)$ est un fibré en sphères symplectique au-dessus d'une surface de genre $g$.
\end{proof}

\section{Commentaires sur le Théorème~\ref{t:thmprincipal} et sa démonstration}\label{s:pistes}

On discute dans cette section de l'optimalité ou non des bornes dans le Théorème~\ref{t:thmprincipal} et le Théorème~\ref{t:ruledFKversion2}.

On constate dans un premier temps que les bornes du Théorème~\ref{t:thmprincipal} sont optimales. En effet, dans $\mathbb{C}P^1 \times \mathbb{C}P^1$, on peut considérer la configuration de courbes constituée de l'union de deux sections disjointes dans la classe d'homologie $[\mathbb{C}P^1 \times \{*\}]$ et de $n \in \mathbb{N}^*$ fibres deux à deux distinctes. En lissant symplectiquement tous les points singuliers de cette configuration, qui sont des points doubles transverses positifs, on obtient une courbe $C$  symplectiquement plongée de genre $g= n-1$ et dans la classe d'homologie $2[\mathbb{C}P^1 \times \{*\}] + n[ \{*\} \times \mathbb{C}P^1]$, qui vérifie donc $[C]^2 = 4n = 4g +4$. Ainsi, pour tout entier naturel $g$, les bornes du Théorème~\ref{t:thmprincipal} sont optimales.

Pour le Théorème~\ref{t:ruledFKversion2}, qui concerne les courbes symplectiquement plongées de genre $g$ et d'auto-intersection supérieure ou égale à $4g$, la borne n'est clairement pas optimale. Le Théorème~\ref{t:TaubesSW} de Taubes montre, grâce à la théorie de Seiberg--Witten, que le résultat est également valable pour les courbes symplectiquement plongées de genre $g$ et d'auto-intersection strictement supérieure à $2g-2$. Même si le Théorème~\ref{t:ruledFKversion2} n'est utilisé qu'en tant qu'étape intermédiaire dans la démonstration du Théorème~\ref{t:thmprincipal}, qui lui est optimal, il serait néanmoins intéressant d'obtenir une version du Théorème~\ref{t:ruledFKversion2} ne faisant appel qu'aux techniques de courbes pseudoholomorphes.

Le coefficient $4$ dans la borne provient des éclatements effectués sur les points doubles des composantes symplectiquement positivement immergées (en effet, chaque tel éclatement diminue le genre arithmétique de $1$ et le nombre d'auto-intersection de $4$). Pour se débarrasser de cette constante, il suffirait d'avoir une compréhension suffisante (régularité, compacité ou non compacité, etc) de certains espaces de modules de courbes pseudoholomorphes simples non plongées, ce qui pourrait éventuellement permettre de se passer des éclatements. Plusieurs problèmes se dessinent alors. D'un côté, de tels espaces de modules pourraient contenir des courbes non immergées, qui ne satisfont pas les hypothèses des théorèmes de transversalité automatique évoqués dans ce manuscrit (notamment celui concernant le problème avec contraintes ponctuelles), ce qui pose un problème pour obtenir des propriétés de surjectivité de certaines applications d'évaluation comme dans la Section~\ref{ss:smash} (plus précisément, pour montrer que les images de ces applications sont ouvertes). Ce genre de problème pourrait être contourné grâce à des théorèmes de transversalité automatique qui s'appliqueraient également aux courbes non immergées. D'un autre côté, dans le cas des courbes de genre $1$ (parmi lesquelles les courbes d'indice contraint égal à $2$ satisfont les hypothèses du théorème de transversalité automatique avec contraintes), on pourrait être tenté d'adopter une stratégie faisant intervenir la courbe universelle au-dessus de certains espaces de modules de dimension virtuelle $2$. Mais comme expliqué dans la Section~\ref{s:spheresimmergees}, cela nécessiterait soit d'avoir une compréhension avancée de la topologie des fibrés en tores sur les surfaces, soit de montrer que les espaces de modules de dimension virtuelle $2$ considérés sont homéomorphes à des sphères. Il faudrait également trouver un moyen d'obtenir des informations sur le degré des applications d'évaluation associées sans faire appel aux théorèmes de transversalité automatique avec une contrainte ponctuelle supplémentaire (les courbes de genre $1$ et d'indice contraint nul n'étant pas automatiquement Fredholm régulières pour le problème avec contraintes ponctuelles).

Le terme constant de la borne est quant à lui déterminé pour ne pas avoir à s'occuper des courbes symplectiquement plongées de genre $g$ et d'auto-intersection appartenant à $\{2g-1, 2g, 2g +1 \}$ (ce qui correspond aux courbes pseudoholomorphes plongées de genre $g$ d'indice strictement supérieur à $2g-2$ et inférieur ou égal à $2g+4$), pour lesquelles la Proposition~\ref{p:casserlescourbes} ne s'applique pas. L'idéal serait de trouver un moyen de montrer que les espaces de modules des courbes pseudoholomorphes de ce type ne sont pas compacts afin d'assurer l'existence de courbes nodales appropriées. On verra cependant dans le Chapitre $5$ que ce n'est pas toujours le cas, puisqu'on y montre notamment que pour une structure presque complexe générique, les espaces de modules des courbes pseudoholomorphes de genre $1$ et d'auto-intersection $1$ plongées dans les surfaces symplectiques \emph{minimales} sont compacts. La question reste ouverte pour les autres courbes plongées de genre $g \geq 1$ et d'auto-intersection appartenant à $\{2g-1, 2g, 2g +1 \}$.
\clearemptydoublepage

\chapter{Application à la classification de remplissages symplectiques}
À toute courbe symplectiquement plongée d'auto-intersection strictement positive dans une surface symplectique, on peut associer canoniquement une variété de contact de dimension $3$. Dans ce chapitre, on utilise le Théorème~\ref{t:thmprincipal} pour classifier les remplissages symplectiques forts des variétés de contact associées aux courbes qui possèdent une auto-intersection haute par rapport à leur genre. On commence par quelques rappels concernant les variétés de contact et leurs interactions avec les variétés symplectiques. Une présentation plus détaillée, orientée vers les basses dimensions, peut être trouvée dans~\cite{ozbagci2004surgery}.

Soit $X$ une variété différentielle de dimension $2n+1$ munie d'une $1$--forme différentielle $\alpha$. On dit que $\alpha$ est une \emph{forme de contact} si $\alpha \wedge (d \alpha)^n \neq 0$. Une distribution d'hyperplan $\xi \subset T X$ est une \emph{structure de contact} sur $X$ si elle est localement définie comme le noyau d'une forme de contact. On dit dans ce cas que $(X, \xi)$ est une \emph{variété de contact}. 

\begin{ex}
\leavevmode
\begin{enumerate}
\item Pour l'espace vectoriel $\mathbb{R}^{2n+1}$ muni des coordonnées $(x_1, y_1, \dots, x_n, y_n, z)$, la structure de contact standard $\xi_{st}$ est donnée par $\ker \left(dz + \sum_{i=1}^n x_i dy_i \right)$.
\item L'ensemble des tangentes complexes à $\mathcal{S}^{2n+1} \subset \mathbb{C}^{n+1}$ (c'est-à-dire les espaces de la forme $T_x \mathcal{S}^{2n+1} \cap i(T_x \mathcal{S}^{2n+1})$ pour $x \in \mathcal{S}^{2n+1}$) définit une structure de contact sur $\mathcal{S}^{2n+1}$.
\end{enumerate}
\end{ex}

Deux variétés de contact $(X, \xi)$ et $(X', \xi')$ sont dites \emph{contactomorphes} s'il existe un difféomorphisme $\Phi :X \rightarrow X'$ tel que $\Phi_* \xi = \xi'$, on dit alors que $\Phi$ est un \textit{contactomorphisme}. Cette notion donne lieu à une relation d'équivalence sur l'ensemble des variétés de contact. 
De la même manière que pour les variétés symplectiques, les variétés de contact ne possèdent pas d'invariant local. 

\begin{thm}[Théorème de Darboux pour les structures de contact]
Soit $(X, \xi)$ une variété de contact. Pour tout point $p \in X$, il existe un voisinage ouvert $U$ de $p$ dans $X$ tel que $(U, \xi_{ \vert U})$ est contactomorphe à un ouvert $V \subset \mathbb{R}^{2n+1}$ muni de la forme symplectique standard $\xi_{st \vert V}$.
\end{thm}

Un champ de vecteurs $v$ sur une variété symplectique $(M, \omega)$ est un \emph{champ de vecteur de Liouville} si $\mathcal{L}_v \omega = \omega$ (où $\mathcal{L}$ désigne la dérivée de Lie). Puisque $d \omega = 0$, on a d'après la formule de Cartan $\mathcal{L}_v \omega = d \iota_v \omega + \iota_v d \omega = d \iota_v \omega$. La condition pour que $v$ soit un champ de vecteurs de Liouville est donc équivalente à $d \iota_v \omega =0$. On dit qu'une sous-variété $X$ de codimension $1$ de $(M, \omega)$ est de \emph{type contact} s'il existe un champ de vecteur de Liouville $v$ défini au voisinage de $X$ qui est transverse à $X$. 

\begin{rk}
\'Etant donné un champ de vecteur de Liouville $v$, son opposé $-v$ n'est pas un champ de vecteur de Liouville puisque $\mathcal{L}_{-v} \omega = -\omega$.
\end{rk}

Une sous-variété de type contact est une variété de contact, comme le montre le théorème suivant.

\begin{thm}[\cite{WEINSTEIN1381413841}]
Soit $(M, \omega)$ une variété symplectique. Une sous-variété $X$ de codimension $1$ de $(M, \omega)$ est de type contact si et seulement si on dispose d'une $1$--forme différentielle $\alpha$ définie au voisinage de $X$ telle que $d \alpha = \omega$ et $\alpha$ est non nulle sur le champ de droite donné par $L_X = \{ v \in T M  \mid  \text{pour tout} ~ w \in T X, ~ \omega ( v, w) = 0 \}$.
Dans ce cas, la restriction de $\alpha$ à $X$ définit une forme de contact sur $X$.
\end{thm}

On dit qu'une sous-variété à bord $U$ de codimension $0$ de $(M, \omega)$ possède un bord  $\omega$--convexe (resp. $\omega$--concave) si $\partial U$ est de type contact et le champ de vecteur $v$ pointe à l'extérieur (resp. à l'intérieur) de $U$.

Soit $(X, \xi)$ une variété de contact de dimension $3$ et $(W, \omega)$ une surface symplectique \emph{compacte} à bord $\omega$--convexe. On dit que $(W, \omega)$ est un \emph{remplissage symplectique fort} de $(X, \xi)$ si son bord est contactomorphe à $(X, \xi)$. En d'autres termes, la forme symplectique $\omega$ est exacte au voisinage du bord de $(W, \omega)$ et on peut choisir une $1$--forme $\alpha$ au voisinage du bord de $(W, \w)$ telle que $d \alpha =\omega$ et $\ker ( \alpha_{\vert \partial W} ) = \xi$. De manière équivalente, la condition d'être un remplissage symplectique fort peut se reformuler en l'existence d'un champ de vecteur de Liouville défini au voisinage de $\partial W$, transverse à $\partial W$ et pointant vers l'\emph{extérieur}.

De même, étant donné une surface symplectique \emph{compacte} $(N, \omega)$ à bord $\omega$--concave, on dit que $(N, \omega)$ est un \emph{bouchon symplectique} (ou \emph{remplissage concave}) de $(X, \xi)$ si son bord est contactomorphe à $(X, \xi)$. Cette fois cette condition est équivalente à l'existence d'un champ de vecteur de Liouville, défini au voisinage de $\partial N$, transverse à $\partial N$ et pointant vers l'\emph{intérieur}. Ainsi, une surface symplectique $(N , \omega)$ compacte à bord est un bouchon symplectique si son bord est $\omega$--concave.

\begin{ex}
Le champ de vecteurs sur $\mathbb{R}^4$ défini dans les coordonnées $(x_1,y_1,x_2,y_2)$ par
$$v(x_1,y_1,x_2,y_2) = x_1 \frac{\partial}{\partial x_1}+y_1 \frac{\partial}{\partial y_1} + x_2 \frac{\partial}{\partial x_2}+y_2 \frac{\partial}{\partial y_2},$$
est un champ de vecteurs de Liouville sur $(\mathbb{R}^4, \omega_{st})$. La sphère unité $\mathcal{S}^3 \subset \mathbb{R}^4$ hérite donc d'une structure de contact $\xi_{st}$, qu'on appelle la structure de contact standard. Ceci montre que la boule unité $(B^4, \omega_{st})$ est un remplissage symplectique fort de la variété de contact $(\mathcal{S}^3,\xi_{st})$.
\end{ex}

Une variété de contact peut admettre beaucoup de bouchons symplectiques. En effet, le complémentaire de toute boule dans une surface symplectique est un bouchon symplectique de $(\mathcal{S}^3,\xi_{st})$. Les remplissages symplectiques peuvent en revanche être plus rares. On dispose en effet de certains résultats d'unicité concernant les remplissages symplectiques forts (voir les résultats évoqués dans l'introduction).

Le théorème suivant permet de recoller de manière symplectique un bouchon symplectique d'une variété de contact de dimension $3$ avec un remplissage symplectique fort de cette même variété de contact.

\begin{thm}[\cite{ETNYRE19983}]\label{t:collesymplectique}
Soit $(W, \omega)$ une surface symplectique avec une composante de bord $\omega$--convexe $(X, \xi)$ et $(N, \omega')$ une surface symplectique avec une composante de bord $\omega'$--concave $(X', \xi')$. On suppose qu'il existe un contactomorphisme $\Phi : ( X, \xi) \rightarrow ( X', \xi')$. Alors la variété obtenue en identifiant $( X, \xi)$ et $(X', \xi')$ via $\Phi$ admet une structure symplectique.
\end{thm}

\begin{rk}
Notons qu'il est possible de définir les opérations d'éclatement et de contraction dans le cadre symplectique à l'aide du Théorème~\ref{t:collesymplectique} (mais ce n'est pas la seule façon de le faire).
\end{rk}

On s'intéresse à présent aux structures de contact définies à partir des courbes symplectiquement plongées dans les surfaces symplectiques.

\begin{thm}[{\cite{Gay11}}]\label{t:voisinageconcave} 
Soit $C$ une courbe symplectiquement plongée dans une surface symplectique $(M, \omega)$, avec genre et nombre d'auto-intersection fixés. Si $[C]^2 >0$, alors il existe un voisinage tubulaire $(N, \omega)$ de $C$ dont le bord est $\omega$--concave. De plus, tout plongement symplectique de $C$ (avec le même nombre d'auto-intersection) dans une surface symplectique $(M', \omega')$ possède un voisinage à bord $\omega'$--concave qui est une déformation symplectique de $(N,\omega)$ et qui induit la même structure de contact sur le bord.
\end{thm}

En particulier, à un type de courbe $C$ donné (c'est-à-dire un genre et un nombre d'auto-intersection strictement positif), on peut associer une variété de contact $(X_C, \xi_C)$ de dimension $3$, où $X_C = - \partial N$ est le bord de $N$ muni de l'orientation opposée et $\xi_C$ est induit par un champ de vecteurs de Liouville pointant vers l'intérieur défini sur un voisinage de $\partial N$.

On peut aisément obtenir d'autres remplissages symplectiques à partir d'un remplissage symplectique donné en éclatant des points dans son intérieur. En effet, cette opération ne change pas la structure symplectique au voisinage du bord, et par conséquent ne change pas non plus la structure de contact sur le bord. Pour classifier les remplissages symplectiques, il est donc nécessaire et suffisant de considérer les cas minimaux.

Le théorème suivant est une application du Théorème~\ref{t:thmprincipal}.


\begin{thm}
Soit $S$ une courbe de genre $g$ symplectiquement plongée dans une surface symplectique. Si on a $S^2 > 4g+5$ ou bien $S^2 \geq 4g+ 5$ et $g \neq 1$, alors la variété de contact $(X_S, \xi_S)$ admet un unique remplissage symplectique fort minimal à difféomorphisme près, qui est un fibré en disques au-dessus d'une surface de genre $g$. Si en revanche $g=1$ et $S^2 = 9$, la variété de contact $(X_S, \xi_S)$ admet deux remplissages symplectiques forts minimaux à difféomorphisme près.
\end{thm}

\begin{rk}
En fait, dans le cas où le remplissage symplectique fort est difféomorphe un fibré en disques au-dessus d'une surface de genre $g$, on peut choisir le difféomorphisme de façon à ce que la forme symplectique se restreigne à une forme symplectique sur chaque fibre. La classe de déformation symplectique du remplissage est alors complètement déterminée (voir la Remarque~\ref{r:déformationsymp}).
\end{rk}

\begin{proof}
On considère un remplissage symplectique fort minimal $(W, \omega)$ de la variété de contact $(X_S, \xi_S)$. D'après le Théorème~\ref{t:voisinageconcave}, la courbe $S$ admet un voisinage tubulaire $(N, \omega')$ qui est un bouchon symplectique de la variété de contact $(X_S, \xi_S)$ qui lui est associée. Grâce au Théorème~\ref{t:collesymplectique}, on peut recoller de manière symplectique $(W, \omega)$ et $(N,\omega')$ selon leurs bords afin d'obtenir une surface symplectique qui contient $S$ et qui est relativement minimale par rapport à $S$. Le Théorème~\ref{t:thmprincipal} permet alors d'affirmer que la surface symplectique ainsi obtenue est un fibré en sphères symplectique dont $S$ est une section.
\begin{figure}[h]
	\centering
	\includegraphics[scale=0.6]{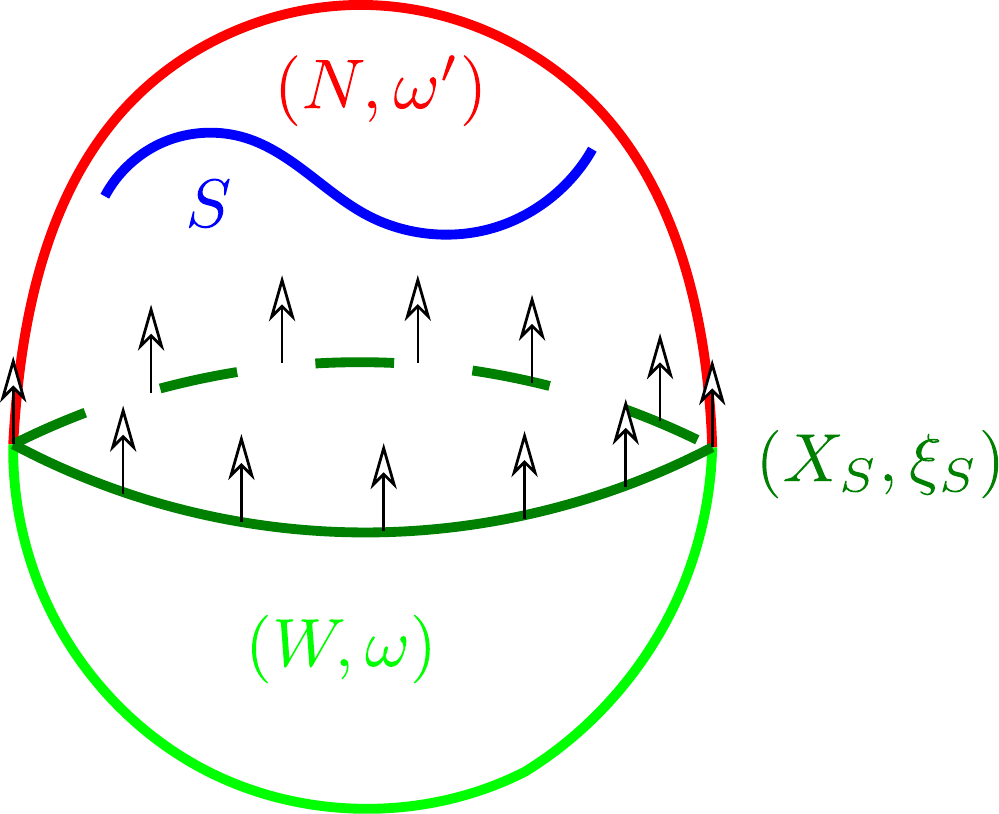}\\
	\caption{Représentation schématique de la surface symplectique obtenue en recollant selon leurs bords le bouchon symplectique $(N,\omega')$ contenant la courbe $S$ avec un remplissage symplectique fort $(W, \omega)$ de la variété de contact $(X_S, \xi_S)$.}
	\label{f:symplecticfeelings}
\end{figure}

La variété $X_S$ est donc un fibré en cercles (orienté) au-dessus d'une surface $B$ de genre $g$ ; le voisinage tubulaire $N$ de $S$ est un fibré en disques (orienté) au-dessus de $B$ qui possède $S$ pour section ; et son complémentaire $(W, \omega)$ est également un fibré en disques (orienté) au-dessus de $B$. Remarquons que chaque fibré en cercles orienté au-dessus de $B$ correspond de manière unique à un fibré en disques orienté au-dessus de $B$ (chacun d'entre eux étant également associé de manière unique à un fibré en plans orienté au-dessus de $B$). Rappelons également que les fibrés en disques orientés au-dessus d'une base donnée sont complètement classifiés par leurs premières classes de Chern. Rappelons enfin que la première classe de Chern d'un fibré en disques orienté peut être vue comme le nombre de zéros d'une section générique comptés avec signes. La première classe de Chern du fibré en disques sur $N$ est donc égale à $[S]^2$. Par définition, la variété $X_S$ est difféomorphe au bord de ce fibré en disques muni de l'orientation opposée. Le fibré en cercles orienté sur $X_S$ est donc le bord du fibré en disques orienté dont la première classe de Chern est égale à $-[S]^2$. Ceci permet de déterminer complètement le type de difféomorphisme du remplissage $(W, \omega)$.
\end{proof}
\clearemptydoublepage

\chapter{Isotopie symplectique pour les sections des surfaces géométriquement réglées au-dessus de courbes elliptiques}
On s'intéresse maintenant au problème d'isotopie symplectique dans certaines surfaces kählériennes réglées non rationnelles. Une \emph{isotopie symplectique} entre deux courbes symplectiquement plongées $S$ et $S'$ dans une surface symplectique donnée est une isotopie lisse $(S_t)_{t \in [0,1]}$ telle que $S_0 = S$, $S_1 = S'$ et pour tout $t \in [0,1]$, $S_t$ est une courbe symplectiquement plongée. Une \emph{surface kählérienne} est une surface symplectique munie d'une structure complexe intégrable compatible avec la forme symplectique. Rappelons que toutes les surfaces projectives complexes peuvent être munies d'une structure kählérienne via la restriction de la forme de Fubini--Study sur $\mathbb{C}P^n$ à la surface projective en question. Le problème d'isotopie symplectique consiste à montrer qu'une courbe symplectique d'un type donné (c'est-à-dire un genre, un nombre d'auto-intersection et éventuellement le nombre et le type de singularités) dans une surface kählérienne est symplectiquement isotope à une courbe complexe (dans le cas des courbes singulières, on demande en plus que l'isotopie soit \emph{équisingulière}, c'est-à-dire que les courbes possèdent les mêmes types de singularités au cours de l'isotopie).

On appelle \emph{courbe elliptique} toute courbe (complexe ou symplectique selon le contexte) non singulière de genre  $1$. On résout dans ce chapitre le problème d'isotopie symplectique pour les sections des surfaces géométriquement réglées au-dessus de courbes elliptiques.  

\begin{thm}\label{t:isotopiesectionell}
Soit $T$ une courbe elliptique symplectiquement plongée dans une surface complexe géométriquement réglée au-dessus d'une courbe elliptique. Si $T$ est une section, alors elle est symplectiquement isotope à une section complexe. 
\end{thm}

On commence par s'intéresser dans la Section~\ref{s:courbellplongées} aux plongements symplectiques possibles des courbes elliptiques d'auto-intersection strictement positive dans les surfaces symplectiques minimales. Dans la Section~\ref{s:toreautoint1}, on s'intéresse aux courbes elliptiques symplectiquement plongées d'auto-intersection $1$ dans les surfaces symplectiques minimales. On montre dans un premier temps, sans utiliser la théorie de Seiberg--Witten, que pour une structure presque complexe générique, les espaces de modules de courbes pseudoholomorphes associés à ce type de courbes sont compacts. La Section~\ref{s:surfaceregsurell} présente un résultat dû à Suwa concernant la classification à biholomorphisme près des surfaces complexes minimales réglées au-dessus de courbes elliptiques. On clôt ce chapitre par la Section~\ref{s:isotopiebirationnelle}. Dans cette section, on utilise des techniques pseudoholomorphes pour montrer, grâce au résultat de compacité obtenu dans la Section~\ref{s:toreautoint1}, que les sections symplectiques d'auto-intersection $1$ des surfaces géométriquement réglées stables au-dessus de courbes elliptiques sont symplectiquement isotopes à des courbes complexes. On utilise ensuite les techniques d'isotopie symplectique via des configurations de courbes birationnellement équivalentes, introduites par Golla et Starkston dans~\cite{GS}, pour démontrer le Théorème~\ref{t:isotopiesectionell}.

\section{Courbes elliptiques symplectiquement plongées dans les surfaces symplectiques minimales}\label{s:courbellplongées}

\begin{prop}\label{p:Tdansfibré}
Soit $T$ une courbe elliptique symplectiquement plongée dans une surface symplectique minimale $(M, \omega)$.

Si $[T]^2 \in \mathbb{N}^* \backslash \{ 8 ,9\}$, alors $(M, \omega)$ est un fibré en sphères symplectique au-dessus d'un tore dont $T$ est une section.

Si $[T]^2 =8$, alors $(M, \omega)$ est un fibré en sphères symplectique au-dessus d'un tore dont $T$ est une section ou $(M, \omega)$ est un fibré en sphères symplectique au-dessus d'une sphère dont la restriction de l'application de projection à $T$ est de degré $2$.

Si $[T]^2 =9$, alors $(M, \omega)$ est un fibré en sphères symplectique au-dessus d'un tore dont $T$ est une section ou $T$ est une cubique symplectique non singulière dans $\mathbb{C}P^2$.
\end{prop}

\begin{rk}
La Proposition~\ref{p:Tdansfibré}, contrairement au Théorème~\ref{t:thmprincipal}, ne permet pas de classifier les remplissages symplectiques des variétés de contact associées aux courbes considérées. En effet, les méthodes employées dans le Chapitre 4 se basent sur la connaissance de tous les plongements symplectiques \emph{relativement} minimaux dans les surfaces symplectiques, et pas seulement sur celle des plongements symplectiques dans les surfaces symplectiques minimales.

Pour les courbes elliptiques d'auto-intersection comprise entre $1$ et $8$, il existe d'autres plongements relativement minimaux. Par exemple, une cubique lisse de $\mathbb{C}P^2$ éclatée en $1 \leq k \leq 8$ points fournit une courbe elliptique symplectiquement plongée $T$ d'auto-intersection $[T]^2 = 9-k \in \{1, \dots,8 \}$ telle que $\mathbb{C}P^2 \# k\overline{\mathbb{C}P}^2$ est relativement minimale par rapport à $T$.
\end{rk}

\begin{proof}
Si $[T]^2 \geq 9$, le résultat est une conséquence immédiate du Théorème~\ref{t:thmprincipal}. Il ne reste donc que les cas où $1 \leq [T]^2 \leq 8$ à traiter. D'après le Théorème~\ref{t:TaubesSW} de Taubes, $(M, \omega)$ est une surface symplectiquement réglée. Par la formule du genre, tous les tores symplectiquement plongés dans $\mathbb{C}P^2$ vérifient $[T]^2 = 9$, donc d'après le Théorème~\ref{t:classificationsurfacesreglees}, la surface symplectiquement réglée $(M, \omega)$ est un fibré en sphères symplectique au-dessus d'une surface. En choisissant une structure presque complexe $J$ dominée par $\omega$ telle que $T$ est l'image d'une courbe pseudoholomorphe simple (on peut supposer $J$ générique car $T$ est d'auto-intersection strictement positive et, par conséquent, est automatiquement Fredholm régulière d'après le Théorème~\ref{t:transauto1}) et en utilisant le Théorème~\ref{t:isotopiepinceaudeLefschetz}, on peut supposer que $T$ intersecte localement positivement chacune des fibres.

On note $\S$ une section du fibré et $F$ une fibre. On dispose d'entiers $k,\ell$ tels que $[T] = k[\S]+\ell[F]$, avec $k = [T] \cdot [F] >0$ (on a $k \neq 0$ car $T$ n'est pas une fibre et intersecte au moins une fibre). On a alors
\begin{equation}\label{e:Tcarré1}
[T]^2 = k^2 [\S]^2 +2k \ell,
\end{equation}
et la formule d'adjonction nous donne 
\begin{equation}\label{e:Tcarré2}
\begin{split}
[T]^2 
&= c_1([T]) - \chi (T) \\
&= k c_1([\S]) + \ell c_1([F]) \\
&= k c_1([\S]) + 2\ell.
\end{split}
\end{equation}
D'après l'égalité~\eqref{e:Tcarré1}, l'entier $k$ divise $[T]^2$. Par conséquent l'égalité~\eqref{e:Tcarré2} implique que $k$ divise $\ell$. En reprenant l'égalité~\eqref{e:Tcarré1}, on obtient que $k^2$ divise $[T]^2$. 

Si $[T]^2$ n'admet pas de facteur carré (ce qui est le cas pour $[T]^2 \in \{1,2,3,5,6,7 \}$), alors $k=1$, ce qui permet de conclure.

%

Si $[T]^2= 4$, alors $k \in \{1,2 \}$.
Si $k=2$, on a $1 = [\S]^2 + \ell$ d'après l'égalité~\eqref{e:Tcarré1} et $2 = c_1([\S]) + \ell$ d'après l'égalité~\eqref{e:Tcarré2}. On a alors $\ell \equiv [\S]^2 +1 \mod 2$ et $\ell \equiv c_1([\S]) \mod 2$. Mais d'après la formule d'adjonction, on a $\chi (\S) = c_1([\S])-[\S]^2$, donc $c_1([\S])$ et $[\S]^2$ ont la même parité, ce qui est une contradiction. On a alors $k=1$, donc $T$ est une section.

Si $[T]^2 = 8$, alors $k \in \{1,2 \}$. Si $k =2$, l'équation~\eqref{e:Tcarré2} nous donne $[T]^2 =  2c_1([\S]) + 2\ell$. En appliquant la formule d'adjonction, on obtient $8=  2 \chi(\S) +2 [\S]^2  + 2\ell$. De plus, on a d'après l'équation~\eqref{e:Tcarré1}, $8 = 4 [\S]^2 + 4\ell$. Donc $\chi(\S) = 2$, autrement dit, la base du fibré est une sphère. Si $k=1$, $T$ est une section.
\end{proof}

\begin{ex}
Les courbes elliptiques symplectiquement plongées d'auto-intersection $0$, $8$ ou $9$ dans les surfaces symplectiques minimales ne sont pas forcément des sections de fibrés en sphères symplectiques. En effet, la surface symplectique $\mathcal{T}^2 \times \mathcal{T}^2 $ contient des tores symplectiquement plongés d'auto-intersection $0$. De même, dans $\mathcal{S}^2 \times \mathcal{S}^2$ on considère la configuration de courbes donnée par l'union de deux sections d'auto-intersection $0$ distinctes et de deux fibres distinctes. En lissant symplectiquement les quatre points singuliers de la configuration, qui sont des points doubles transverses positifs, on obtient une courbe symplectiquement plongée de genre $1$ et d'auto-intersection $8$.
\end{ex}

\section{Courbes elliptiques symplectiquement plongées d'auto-intersection $1$}\label{s:toreautoint1}

Dans l'optique de redémontrer dans sa totalité le Théorème~\ref{t:TaubesSW} de Taubes sans faire appel à la théorie de Seiberg--Witten, il est naturel de s'intéresser aux courbes symplectiquement plongées de genre $g \geq 1$ dont le nombre d'auto-intersection $s$ vérifie $2g-2 < s <4g$ (le cas $s \geq 4g$ faisant l'objet du Théorème~\ref{t:ruledFKversion2}). Le premier cas dans l'ordre lexicographique est celui des courbes symplectiquement plongées de genre $1$ et d'auto-intersection $1$. Dans cette section on montre, sans utiliser la théorie de Seiberg--Witten, que pour une structure presque complexe générique, les espaces de modules associés à de telles courbes sont compacts.

\begin{lem}\label{l:compacitéTautoint1réglée}
Soit $T$ un tore symplectiquement plongé d'auto-intersection $1$ dans une surface symplectiquement réglée minimale $(M, \omega)$. Alors $T$ est une section d'un fibré en sphères symplectique au-dessus d'un tore et pour toute structure presque complexe $J$ \emph{générique} dominée par $\omega$, l'espace de modules $\mathcal{M}_{1}([T];J)$ est une variété lisse \emph{compacte} de dimension $2$.
\end{lem}

\begin{rk}
L'hypothèse de généricité sur la structure presque complexe $J$ est nécessaire. Sans cette hypothèse, l'espace de modules $\overline{\mathcal{M}}_{1}([T];J)$ pourrait contenir des courbes nodales constituées de l'union d'une section d'auto-intersection $-1$ avec une fibre (voir la Section~\ref{s:surfaceregsurell} pour un exemple). Comme l'indice des courbes pseudholomorphes plongées de genre $1$ et d'auto-intersection $-1$ est égal à $-2$, la généricité de $J$ permet de s'assurer que de telles courbes n'existent pas dans $(M,J)$.
\end{rk}

\begin{proof}
D'après la Proposition~\ref{p:Tdansfibré}, $T$ est une section d'un fibré en sphères symplectique au-dessus d'un tore. On peut supposer chacune des fibres $J$--holomorphes grâce au Théorème~\ref{t:isotopiepinceaudeLefschetz}. D'après la Proposition~\ref{p:morceauxcourbesind2}, énoncée dans l'Annexe~\ref{c:annexe}, on a deux possibilités pour les courbes de $\overline{\mathcal{M}}_{1}([T];J)$ qui ne sont pas plongées : soit une courbe nodale avec deux composantes symplectiquement plongées d'indice $0$ qui s'intersectent exactement une fois, de manière transverse et positive ; soit une courbe nodale avec une unique composante symplectiquement positivement immergée qui possède un unique point double.

Dans le premier cas, une des composantes plongées est rationnelle. Comme elle est aussi plongée et d'indice $0$, son auto-intersection est égale à $-1$, ce qui contredit la minimalité de $(M, \omega)$.

Dans le second cas, en considérant la fibre $F$ qui passe par l'unique point double, on a par positivité d'intersection $[T] \cdot [F] > 1$, ce qui est absurde puisque $T$ est une section.
\end{proof}

\begin{rk}~\label{r:morceausection}
Considérons une surface presque complexe $(M,J)$ géométriquement réglée dont les fibres sont des images de courbes $J$--holomorphes, $F$ une fibre de ce fibré et $\Sigma$ une section $J$--holomorphe de ce fibré. On rappelle que le groupe $H_2(M; \mathbb{Z})$ est engendré par $[\S]$ et $[F]$. On considère une courbe nodale $\Sigma_\infty$ homologue à $\Sigma$. Même si la structure presque complexe n'est pas générique, par positivité d'intersection, une composante $\Sigma'$ de $\Sigma_\infty$ doit vérifier $[\Sigma'] \cdot [F] =1$. La courbe $\Sigma'$ est donc également une section $J$--holomorphe de $(M,J)$ et on dispose d'un entier $k$ tel que 
$$[\Sigma] = [\Sigma'] + k [F].$$ Comme les autres composantes $\Sigma''$ vérifient $[\Sigma''] \cdot [F] =0$, elles ont toutes pour images des fibres. 
La classe d'homologie $[\Sigma] - [\Sigma'] = k [F]$ représente la somme des classes d'homologie des composantes de $\Sigma_\infty$ différentes de $\Sigma'$. Par positivité d'intersection, on a alors $k =  ([\Sigma] - [\Sigma']) \cdot [\Sigma] \geq 0$ avec égalité si et seulement si $\Sigma_\infty$ est plongée.
\end{rk}

On pourrait directement en déduire la compacité de l'espace de modules $\mathcal{M}_{1}([T];J)$ dans le cas où on suppose seulement que $T$ est symplectiquement plongé dans une surface symplectique minimale grâce au Théorème~\ref{t:TaubesSW} de Taubes. On préfère néanmoins démontrer la proposition suivante sans recourir à la théorie de Seiberg--Witten.

\begin{prop}\label{p:compacitéTautoint1}
Soit $T$ un tore symplectiquement plongé d'auto-intersection $1$ dans une surface symplectique minimale $(M, \omega)$. Alors pour toute structure presque complexe $J$ générique dominée par $\omega$, l'espace de modules $\mathcal{M}_{1}([T];J)$ est une variété lisse \emph{compacte} de dimension $2$.

De plus, pour tout chemin $(\omega_t)_{t \in [0,1]}$ de formes symplectiques sur $M$ et tout chemin générique de structures presque complexes $(J_t)_{t \in [0,1]} \in \mathcal{J}_{\tau}^\mathrm{reg} (M, (\omega_s)_{s \in [0,1]})$, l'espace de modules à paramètre $\mathcal{M}_{1}([T];(J_t)_{t \in [0,1]})$ est une variété lisse orientée, \emph{compacte}, à bord et de dimension $3$, avec 
$$\partial \mathcal{M}_{1}([T]; (J_t)_{t \in [0,1]}) = -(\{ 0 \} \times \mathcal{M}_{1}([T]; J_0)) \cup (\{ 1 \} \times \mathcal{M}_{1}([T] ; J_1)).$$
De plus, la projection
$$
\begin{array}{ccc}
\mathcal{M}_{1} ([T]; (J_t)_{t \in [0,1]})  & \longrightarrow & [0,1] \\
(t,u) & \longmapsto & t 
\end{array}
$$
est une submersion.
\end{prop}

On utilisera le résultat suivant, démontré par Korkmaz et Özba\u{g}c\i , concernant les nombres minimaux de fibres singulières dans les fibrations de Lefschetz. 
 
\begin{prop}[{\cite[Lemma~2.2]{KO}}]\label{p:victoireparKO}
Une fibration de Lefschetz topologique sur une variété de dimension $4$, dont les fibres régulières sont des tores et qui possède au moins une fibre singulière en possède au moins douze.
\end{prop}

\begin{rk}
La borne de la Proposition~\ref{p:victoireparKO} est atteinte par le pinceau de cubiques dans $\mathbb{C}P^2$, qui possède exactement douze fibres singulières.
\end{rk}

\begin{proof}[Démonstration de la Proposition~\ref{p:compacitéTautoint1}]
Puisque l'entier $[T]^2=1$ n'a pas de facteur carré, toutes les courbes de l'espace de modules $\mathcal{M}_{1}([T];J)$ sont simples. Comme la courbe $T$ est plongée, la formule d'adjonction nous assure que toutes les courbes de $\mathcal{M}_{1}([T];J)$ sont plongées. D'après le Théorème~\ref{t:transauto1}, toutes les courbes elliptiques symplectiquement plongées d'auto-intersection positive sont automatiquement Fredholm régulières. L'espace de modules $\mathcal{M}_{1}([T];J)$ est donc une variété lisse de dimension égale à $\chi(T) + 2 [T]^2 = 2$.

On suppose qu'il existe une courbe $u_\infty \in \overline{\mathcal{M}}_{1}([T];J)$ qui n'est pas plongée. D'après la Proposition~\ref{p:morceauxcourbesind2}, énoncée dans l'Annexe~\ref{c:annexe}, on a deux possibilités :
\begin{enumerate}
\item la courbe $u_\infty$ est constituée de deux composantes plongées $u_1$ et $u_2$, de genres respectifs $0$ et $1$, et d'indice $0$ qui s'intersectent transversalement en un unique point, ou bien
\item la courbe $u_\infty$ est constituée d'une unique composante rationnelle immergée, dont l'unique point multiple est un point double transverse.
\end{enumerate}

Dans le premier cas, puisque $u_1$ est une courbe pseudoholomorphe rationnelle plongée d'indice $0$, la formule d'adjonction nous donne $[u_1]^2 = -1$, ce qui contredit la minimalité de $(M, \omega)$.

Donc toutes les courbes nodales de $\overline{\mathcal{M}}_{1}([T];J)$ sont irréductibles. On considère alors la courbe universelle $\overline{\mathcal{U}}_T(J) =\overline{\mathcal{M}}_{1,1}([T];J)$ au-dessus de $\overline{\mathcal{M}}_{1}([T];J)$. D'après la Proposition~\ref{p:preimageev}, $\overline{\mathcal{U}}_T (J)$ est une variété lisse de dimension $4$, et l'application naturelle de projection $\pi : \overline{\mathcal{U}}_T (J) \rightarrow \overline{\mathcal{M}}_{1}([T];J)$ définit une fibration de Lefschetz avec fibres génériques de genre $1$ et au plus une singularité par fibre. De plus, toutes les fibres de la fibration de Lefschetz sont irréductibles et chacune d'entre elles correspond à une unique courbe nodale de $\overline{\mathcal{M}}_{1}([T];J)$.

D'après la Proposition~\ref{p:victoireparKO}, puisqu'on dispose d'une courbe nodale dans $\overline{\mathcal{M}}_{1}([T];J)$, il en existe au moins onze autres. Considérons une courbe nodale $u_\infty' \in \overline{\mathcal{M}}_{1}([T];J)$ distincte de $u_\infty$. On a $[u_\infty] \cdot [u_\infty'] = [T]^2 = 1$, donc par positivité d'intersection, les deux courbes $u_\infty$ et $u_\infty'$ s'intersectent exactement une fois, en des points réguliers, de manière transverse et positive (voir la Figure~\ref{f:2rational1doublepoint1}).
\begin{figure}[h]
	\centering
	\includegraphics[scale=0.5]{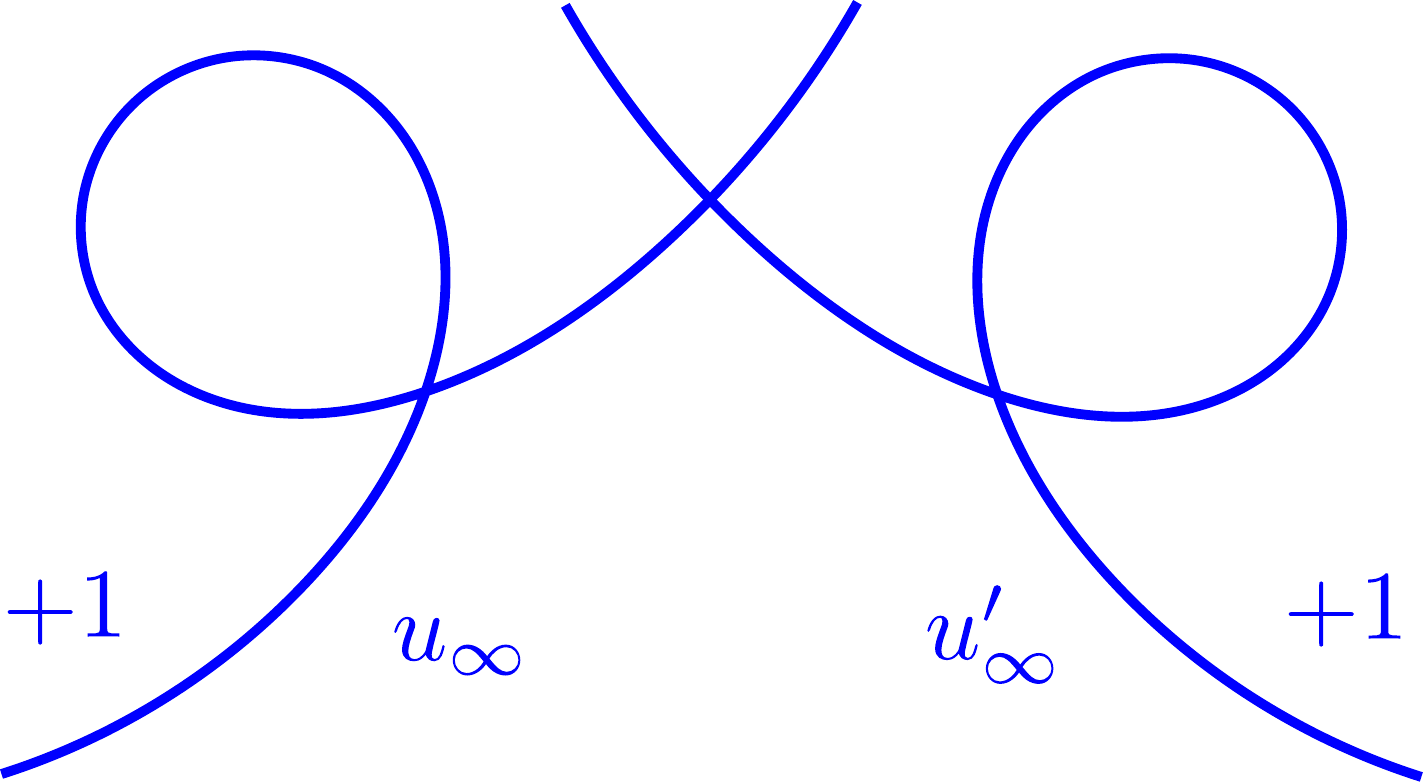}\\
	\caption{Deux copies d'une courbe pseudoholomorphe rationnelle immergée d'auto-intersection $1$ avec un unique point double transverse.}
	\label{f:2rational1doublepoint1}
\end{figure}
On peut alors lisser symplectiquement l'unique point d'intersection entre $u_\infty$ et $u_\infty'$ pour obtenir une courbe $S$ rationnelle symplectiquement positivement immergée vérifiant $c_1([S]) = 2 c_1([T])= 2[T]^2 =2$. D'après le Théorème~\ref{t:McDuffimmersed}, $(M, \omega)$ est une surface symplectiquement réglée. Comme $(M, \omega)$ est minimale, la Proposition~\ref{l:compacitéTautoint1réglée} permet d'aboutir à une contradiction. Donc l'espace de modules $\mathcal{M}_{1}([T];J)$ est compact.

En ce qui concerne la deuxième partie de la proposition, la compacité de l'espace de modules $\mathcal{M}_{1}([T];(J_t)_{t \in [0,1]})$ se démontre de manière strictement similaire grâce à la Remarque~\ref{r:morceauxparam} de l'Annexe~\ref{c:annexe}. Le reste découle directement du Théorème~\ref{t:moduliparamreg}.
\end{proof}

\begin{rk}
On pourrait être tenté d'étendre la stratégie employée ci-dessus aux tores symplectiquement plongés d'auto-intersection strictement plus grande que $1$ ou aux courbes symplectiquement plongées de genre $g \geq2$ et d'auto-intersection $s > 2g-2$. Pour une telle courbe $S$, on pose $m = s -g >0$ et on choisit $m$ contraintes ponctuelles $p_1, \dots, p_m \in M$. Pour une structure presque complexe $J$ générique dominée par $\omega$, l'espace de modules $\mathcal{M}_{1,m}([S];J; p_1, \dots, p_m)$ est une variété lisse de dimension $2$. En revanche, le compactifié $\overline{\mathcal{M}}_{1,m}([S];J; p_1, \dots, p_m)$ pourrait contenir une courbe nodale constituée de deux composantes plongées d'indice contraint nul qui s'intersectent exactement une fois, de manière positive et transverse : une courbe rationnelle d'auto-intersection $0$ avec une contrainte ponctuelle et une courbe de genre $g$ d'auto-intersection $[S]^2-2$ passant par les $m-1$ contraintes ponctuelles restantes (de telles configurations de courbes existent dans certains fibrés en sphères symplectiques au-dessus de surfaces de genre $g$).
\end{rk}

La Proposition~\ref{p:compacitéTautoint1} va à l'encontre de ce qu'on aurait souhaité pour redémontrer le Théorème~\ref{t:TaubesSW} sans utiliser la théorie de Seiberg--Witten. En effet, un résultat de non compacité aurait permis de trouver des courbes de genre strictement plus petit parmi les composantes des courbes nodales provenant du compactifié de Gromov. Il faudra donc trouver d'autres astuces pour arriver au résultat souhaité. Néanmoins la situation est idéale pour réaliser des isotopies symplectiques via des techniques de courbes pseudoholomorphes. Dans la sous-section suivante, on présente les surfaces kählériennes dans lesquelles on va réaliser ces isotopies symplectiques.

\section{Classification des surfaces complexes réglées au-dessus de courbes elliptiques}\label{s:surfaceregsurell}

On commence par rappeler des propriétés générales des surfaces complexes géométriquement réglées au-dessus de courbes complexes non singulières. On met notamment l'accent sur les propriétés des sections qu'elles possèdent. On présente ensuite un résultat de classification à biholomorphisme près des surfaces géométriquement réglées au-dessus de courbes elliptiques dû à Suwa~\cite{suwa}. Ce dernier résultat se base sur la classification des fibrés vectoriels complexes en plans au-dessus de courbes elliptiques, effectuée par Atiyah dans~\cite{atiyah}.

Soit $C$ une courbe complexe non singulière. Chaque fibré vectoriel complexe en plans $E$ au-dessus de $C$ définit une surface géométriquement réglée $\mathbb{P} (E)$ au-dessus de $C$, obtenue en projectivisant fibre à fibre. En réalité, chaque surface géométriquement réglée peut être obtenue de cette manière, comme le montre la proposition suivante, dont la démonstration peut être trouvée dans~\cite[Chapter~V]{hartshorne1977algebraic}.

\begin{prop}
Soit $\pi : V \rightarrow C$ une surface géométriquement réglée au-dessus de $C$, alors 
\begin{enumerate}
\item il existe un fibré vectoriel complexe en plans $E$ au dessus de $C$ tel que $V = \mathbb{P} (E)$,
\item $\pi : V \rightarrow C$ admet une section complexe,
\item si $E$ et $E'$ sont des fibrés en plans complexes au-dessus de $C$, alors $\mathbb{P(E)}$ et $\mathbb{P(E')}$ sont isomorphes si et seulement si on dispose d'un fibré en droites $L$ au-dessus de $C$ tel que $E$ est isomorphe à $E' \otimes L$.
\end{enumerate}
\end{prop}

Le lemme suivant, démontré dans~\cite[Lemma~1.14]{maruyama}, montre le lien entre un fibré vectoriel complexe en plans $E$ et son projectivisé $\mathbb{P} (E)$.

\begin{lem}\label{l:autointsectionruled}
Soit $E$ un fibré vectoriel complexe en plans au-dessus d'une courbe complexe non singulière $C$. Les sous-fibrés en droites complexes de $E$ sont en correspondance bijective avec les sections complexes de $\mathbb{P} (E)$. De plus, étant donné un tel sous-fibré $L$, la section $\Sigma_L$ qui lui est associée vérifie 
$$ \Sigma_L^2 = c_1(E) - 2 c_1(L).$$
\end{lem}

\'Etant donné un fibré en droites complexe $L$, on appelle parfois $c_1(L)$ le \emph{degré} de $L$.

Lorsqu'un fibré en plans $E$ s'écrit comme la somme directe de deux fibrés en droites, on dit que $E$ est \emph{décomposable}. Dans ce cas, on dit que $\mathbb{P} (E)$ est une surface géométriquement réglée \emph{décomposable}. Sinon on dit que $E$ est \emph{indécomposable} et que $\mathbb{P} (E)$ est une surface géométriquement réglée \emph{indécomposable}. Lorsque $E$ est décomposable, on dispose de deux fibrés en droites complexes $L_1$ et $L_2$ tels que $E = L_1 \oplus L_2$. On a alors $c_1(E) = c_1(L_1) + c_1(L_2)$. Les sections $\Sigma_{L_1}$ et $\Sigma_{L_2}$ de $\mathbb{P} (E)$ provenant de la projectivisation des fibrés en droites $L_1$ et $L_2$ sont disjointes et vérifient d'après le Lemme~\ref{l:autointsectionruled}, $\Sigma_{L_1}^2 = - c_1(L_1) + c_1(L_2) = - \Sigma_{L_2}^2 $. On rappelle que le produit tensoriel munit l'ensemble des fibrés en droites complexes au-dessus d'une même courbe complexe, considérés à isomorphisme près, d'une structure de groupe. Par ailleurs, pour deux fibrés en droites complexes $L_1$ et $L_2$ au-dessus d'une même courbe complexe, on a $c_1(L_1 \otimes L_2) =c_1(L_1) + c_1(L_2)$. Par conséquent, les fibrés $L_1 \oplus L_2$ et $(L_1 \oplus L_2) \otimes L_1^{-1} \simeq L_0 \oplus L_2 \otimes L_1^{-1}$, où $L_0$ désigne le fibré en droites complexe trivial, définissent la même surface géométriquement réglée décomposable.

\'Etant donnée une surface géométriquement réglée, on dit qu'elle est \emph{stable} si toutes ses sections sont d'auto-intersection strictement positive et \emph{semi-stable} si toutes ses sections sont d'auto-intersection positive.

On peut passer d'une surface géométriquement réglée au-dessus d'une courbe complexe $C$ à une autre grâce à des transformations élémentaires de Nagata.
On rappelle qu'une transformation élémentaire $\mathrm{elm}_p$ de Nagata en un point $p$ d'une surface géométriquement réglée $V$ consiste à éclater $p$, puis à contracter la transformée propre $\tilde F$ de la fibre $F$ passant par $p$.
\'Etant donné une courbe $S$ symplectiquement plongée dans $V$, on appelle image de $S$ par $\mathrm{elm}_p$ l'image de la transformée propre $\tilde S$ de $S$ par la contraction de $\tilde F$.

\begin{rk}
Soit $\Sigma$ une section complexe d'une surface complexe géométriquement réglée $V$ et $p$ un point de $V$. On note $\S '$ l'image de $\S$ par la transformation élémentaire de Nagata au point $p$. Si $p$ appartient à $\S$, alors on a $[\S ']^2 = [\S ]^2 -1$, sinon on a $[\S ']^2 = [\S ]^2 +1$.
\end{rk}

On présente maintenant quelques propriétés générales concernant les sections complexes des surfaces complexes géométriquement réglées.

\begin{lem}\label{l:sectionsnomansland}
Soit $V$ une surface complexe géométriquement réglée et $\S_-$ une section complexe d'auto-intersection négative. Alors toute section complexe $\S_+$ différente de $\S_-$ vérifie $[\S_+]^2 \geq - [\S_-]^2$.
\end{lem}

\begin{proof}
Le groupe $H_2(M ; \mathbb{Z})$ est engendré par $[\S_-]$ et la classe d'une fibre $[F]$. Comme $[F] \cdot [\S_+] =1$, on dispose d'un entier $k$ tel que 
$$[\S_+] = [\S_-] + k [F].$$
Par positivité d'intersection, on a $[\S_+] \cdot [\S_-] \geq 0$, donc $[\S_-]^2 + k \geq 0$, d'où $k \geq - [\S_-]^2$.
Au final, on a bien $[\S_+]^2 = [\S_-]^2 + 2k \geq - [\S_-]^2$.
\end{proof}

\begin{lem}\label{l:paritésection}
Soit $M$ un fibré en sphères au-dessus d'une surface. Tous les nombres d'auto-intersection des sections de $M$ ont la même parité. 
\end{lem}

\begin{proof}
Soit $\S$ une section de $M$ et $F$ une fibre de $M$. Le groupe $H_2 (M ; \mathbb{Z})$ est engendré par $[\S ]$ et $[F]$. Soit $\S'$ une autre section de $M$. Comme $[\S ] \cdot [F] =1$, on dispose d'un entier $k$ tel que $[\S' ] = [\S ] + k [F]$. 
On a alors $[\S' ]^2 = [\S ]^2 +2k$, ce qui conclut.
\end{proof}

\begin{rk}
La forme d'intersection de $M$ a la même parité que les nombres d'auto-intersection des sections de $M$.  
\end{rk}

Mentionnons aussi le théorème suivant, dû à Nagata. 

\begin{thm}[\cite{nagata1970}] 
Toute surface complexe géométriquement réglée au-dessus d'une courbe de genre $g$ admet une section complexe $\S$ qui vérifie $[\S ]^2 \leq g$.
\end{thm}

On note désormais $C$ une courbe elliptique complexe. Suwa~\cite[Section~1]{suwa} a classifié les surfaces géométriquement réglées au-dessus de $C$ à biholomorphisme près. Elles se décomposent en deux types : les décomposables et les indécomposables. On peut distinguer la plupart d'entre elles par les propriétés de leurs sections.

L'ensemble des surfaces géométriquement réglées décomposables au-dessus de $C$ comporte :
\begin{itemize}
\item Une famille infinie $\mathcal{S}_0$ de surfaces géométriquement réglées décomposables au-dessus de $C$, qui forme un espace de modules isomorphe à $\mathbb{C}P^1$. Chaque surface complexe de $\mathcal{S}_0$ est le projectivisé de la somme de deux fibrés en droites de degré $0$ au-dessus de $C$. Le fibré trivial $S_0$ (provenant de la projectivisation de $C \times \mathbb{C}^2$), qui possède une famille infinie de sections d'auto-intersection $0$, en fait partie. Toute autre surface complexe $S \in \mathcal{S}_0 \backslash \{ S_0 \}$ (provenant de la projectivisation de la somme du fibré trivial avec un fibré non trivial de degré $0$) possède exactement deux sections d'auto-intersection $0$, qui sont disjointes (voir le Lemme~\ref{l:autointsectionruled}). En particulier, d'après le Lemme~\ref{l:sectionsnomansland}, toutes les surfaces complexes géométriquement réglées de $\mathcal{S}_0$ sont semi-stables et d'après le Lemme~\ref{l:paritésection}, elles ont toutes une forme d'intersection paire.

\item Pour tout $n \geq 1$, une unique surface complexe $S_n$, provenant de la projectivisation du fibré trivial avec un fibré de degré $n$. Elle possède deux sections disjointes dont les nombres d'auto-intersection respectifs sont $n$ et $-n$ d'après le Lemme~\ref{l:autointsectionruled}. La forme d'intersection de $S_n$ possède la même parité que $n$ d'après le Lemme~\ref{l:paritésection}.
\end{itemize}

Il existe seulement deux surfaces géométriquement réglées indécomposables au-dessus de $C$ : 
\begin{itemize}
\item La surface complexe géométriquement réglée stable $A_{-1}$. Elle possède une infinité de sections d'auto-intersection $1$, qui forment un espace de modules lisse de dimension $2$ d'après le théorème de transversalité automatique. Cet espace de modules est compact d'après la Remarque~\ref{r:morceausection}. La forme d'intersection de $A_{-1}$ est impaire. Pour une étude plus détaillée de la géométrie de $A_{-1}$, on renvoit à~\cite{diaw}.

\item La surface $A_0$, qui est semi-stable. Elle possède une unique section d'auto-intersection $0$. D'après le Lemme~\ref{l:paritésection}, sa forme d'intersection est paire.
\end{itemize}

\begin{ex}
\leavevmode
\begin{enumerate}
\item Lorsqu'on effectue, $n \geq 1$ transformations élémentaires de Nagata sur une même section du fibré trivial au-dessus de $C$, on obtient une surface complexe géométriquement réglée qui possède une section d'auto-intersection $-n$. On obtient de cette manière la surface $S_n$.
\item On peut passer de la surface complexe $S_1$ à la surface complexe $A_{-1}$ en effectuant deux transformations élémentaires de Nagata respectivement en des points $p, q \in S_1$ tels que $p$ et $q$ n'appartiennent ni à la section d'auto-intersection $-1$ ni à une même section d'auto-intersection $1$ (car on obtient de cette façon une surface réglée stable). Notons que par tout point de $A_{-1}$ passe une section d'auto-intersection $1$.
\item Lorsqu'on effectue une transformation élémentaire de Nagata sur un point de $A_0$ qui n'appartient pas à l'unique section d'auto-intersection $0$, on obtient une surface complexe stable, c'est-à-dire la surface complexe $A_{-1}$.
\end{enumerate}
\end{ex}

On rappelle que la forme d'intersection d'une surface complexe géométriquement réglée au-dessus d'une courbe complexe non singulière de genre donné permet de déterminer complètement son type de difféomorphisme. Le théorème suivant permet de caractériser l'ensemble des surfaces réglées au-dessus de courbes elliptiques par les types d'homéomorphismes possibles pour ces surfaces complexes.

\begin{thm}[{\cite[Theorem~2]{suwa}}]\label{t:homeoreglee}
Si une surface complexe $V$ est homéomorphe à un fibré en sphères au-dessus d'un tore, alors $V$ est une surface géométriquement réglée au-dessus d'une courbe elliptique.
\end{thm}

\section{Isotopie symplectique via des transformations birationnelles}\label{s:isotopiebirationnelle}

On rappelle dans cette section les techniques introduites dans~\cite{GS} par Golla et Starkston pour relier des classes d'isotopie symplectique équisingulière de configurations de courbes symplectiques entre elles via des transformations birationnelles. On résout le problème d'isotopie symplectique pour les sections d'auto-intersection $1$ de la surface complexe géométriquement réglée $A_{-1}$, puis on utilise ces techniques birationnelles pour en déduire le Théorème~\ref{t:isotopiesectionell}.

On rappelle que les surfaces réglées complexes sont des surfaces projectives complexes. N'importe quel plongement complexe d'une surface réglée $V$ dans un espace projectif complexe permet de munir $V$ d'une forme symplectique qui domine la structure complexe, donnée par la restriction de la forme de Fubini--Study à $V$.

\begin{df}
Soit $(M_1, \omega_1)$ et $(M_2, \omega_2)$ deux surfaces symplectiques. On dit qu'une configuration de courbes symplectiques $\mathcal{C}_2$ dans $(M_2, \omega_2)$ est \emph{birationnellement dérivée} d'une configuration de courbes symplectiques $\mathcal{C}_1$ dans $(M_1, \omega_1)$ s'il existe une succession d'éclatements de la paire $(M_1,\mathcal{C}_1)$ vers la transformée totale $(M_1 \# N \overline{\mathbb{C}P}^2,\overline{\mathcal{C}}_1)$, suivie d'une succession de contractions $\pi : M_1 \# N \overline{\mathbb{C}P}^2 \rightarrow M_2$ telle que $\mathcal{C}_2 = \pi (\overline{\mathcal{C}}_1)$.
\end{df}

Notons que cette relation n'est pas symétrique.

\begin{ex}
Soit $C$ une courbe complexe non singulière. Dans $C \times \mathbb{C}P^1$, on considère la configuration $\mathcal{C}_1$ constituée d'une unique courbe de la forme $C \times \{* \}$. On éclate un point $p \in \mathcal{C}_1$. La transformée totale $\overline{\mathcal{C}}_1$ de $\mathcal{C}_1$ est constituée de l'union de la transformée propre de $\mathcal{C}_1$ et du diviseur exceptionnel provenant de l'éclatement du point $p$. On contracte ensuite la transformée propre $\tilde F$ de la fibre $F$ de $C \times \mathbb{C}P^1$ passant par $p$. On note $\mathcal{C}_2$ la configuration de courbes obtenue comme l'image de $\overline{\mathcal{C}}_1$ par cette contraction (voir la Figure~\ref{f:cexbirderiv}). La configuration $\mathcal{C}_2$ est birationnellement dérivée de la configuration $\mathcal{C}_1$, mais la configuration $\mathcal{C}_1$ n'est pas birationnellement dérivée de la configuration $\mathcal{C}_2$.
\begin{figure}[h]
	\centering
	\includegraphics[scale=0.5]{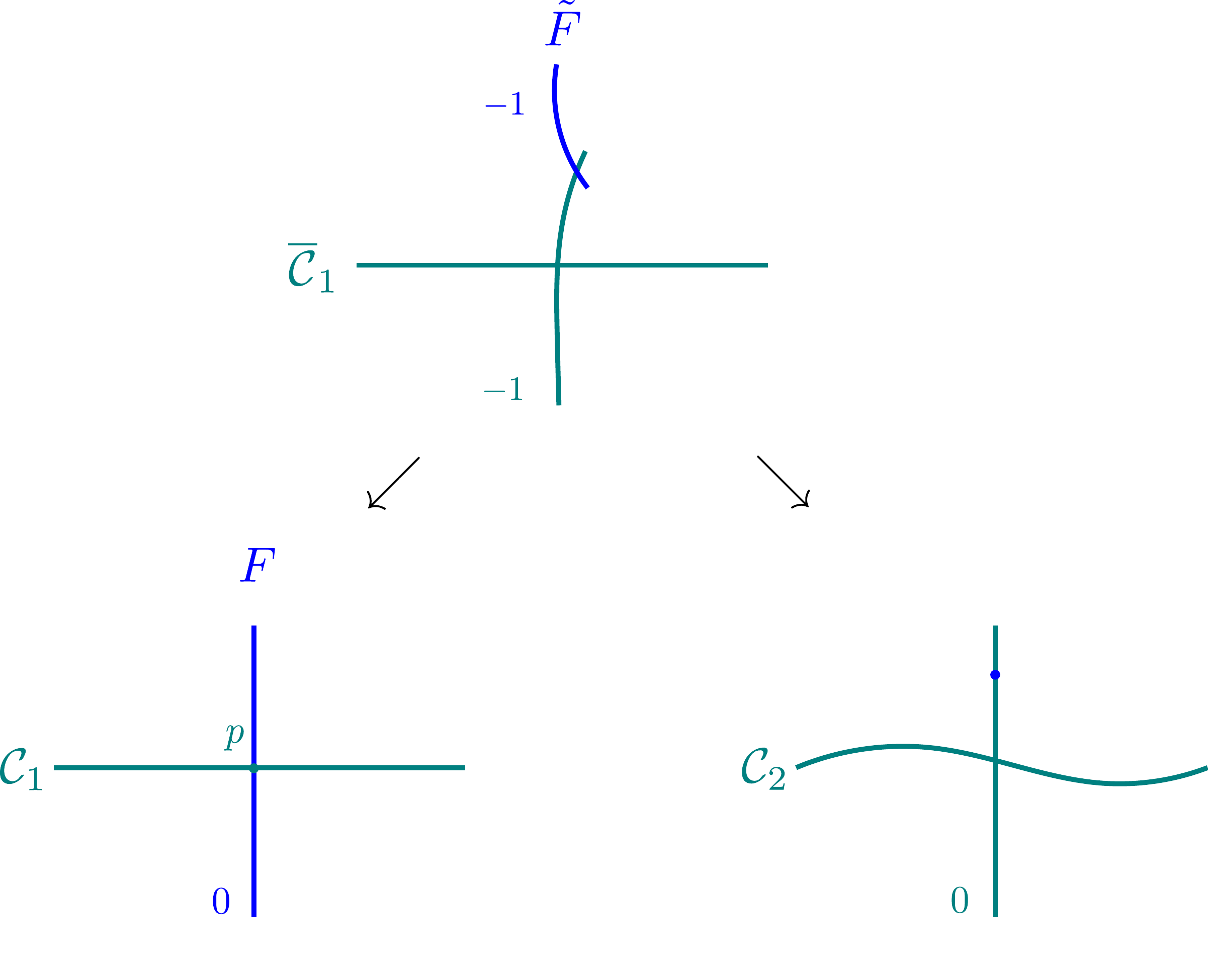}\\
	\caption{La relation de dérivation birationnelle n'est pas symétrique.}
	\label{f:cexbirderiv}
\end{figure}
\end{ex}

\begin{prop}[{\cite[Proposition~3.18]{GS}}]\label{p:birationalequivalence}
Soit $(M_1, \omega_1)$, $(M_2, \omega_2)$ deux surfaces symplectiques, $\mathcal{C}_1$, $\mathcal{C}_1'$ deux configurations de courbes symplectiques dans $(M_1, \omega_1)$ et $\mathcal{C}_2$, $\mathcal{C}_2'$ deux configurations de courbes symplectiques dans $(M_2, \omega_2)$. On suppose que $\mathcal{C}_2$ et $\mathcal{C}_2'$ sont respectivement birationnellement dérivées de $\mathcal{C}_1$ et $\mathcal{C}_1'$ (via les mêmes opérations d'éclatements et de contractions). Si $\mathcal{C}_2$ est équisingulièrement isotope à $\mathcal{C}_2'$ dans $(M_2, \omega_2)$, alors  $\mathcal{C}_1$ est équisingulièrement isotope à $\mathcal{C}_1'$ dans $(M_1, \omega_1)$.

En particulier, si $(M_1, \omega_1)$ et $(M_2, \omega_2)$ sont munies de structures kählériennes, si les éclatements et contractions sont réalisés de manière \emph{complexe} et si $\mathcal{C}_2'$ est une courbe complexe, alors $\mathcal{C}_1'$ est également une courbe complexe.
\end{prop}

On commence dans un premier temps par résoudre le problème d'isotopie symplectique pour les tores symplectiquement plongés d'auto-intersection $1$ dans la surface complexe $A_{-1}$. Une stratégie pour résoudre ce genre de problème repose sur des techniques de courbes pseudoholomorphes : on choisit une structure presque complexe $J$ qui rend la courbe $C$ en question pseudoholomorphe et on considère un chemin générique de structure presque complexe $(J_t)_{t \in [0,1]}$ qui relie $J$ à une structure complexe intégrable. L'objectif est alors de prouver, pour tout $t \in [0,1]$, l'existence d'une courbe $J_t$--holomorphe $C_t$ qui remplit les conditions souhaitées.

\begin{lem}\label{l:A-1isotopiepbTautoint1}
Soit $T$ un tore symplectiquement plongé d'auto-intersection $1$ dans la surface $A_{-1}$ (munie d'une forme symplectique $\omega$ qui domine la structure complexe). Alors $T$ est symplectiquement isotope à une courbe complexe.
\end{lem}

\begin{proof}
On choisit une structure presque complexe $J$ dominée par $\omega$ telle que $T$ est l'image d'une courbe pseudoholomorphe plongée $u$. Comme $[T]^2 >0$, la courbe $u$ est automatiquement Fredholm régulière et on peut supposer sans perte de généralité que $J$ est générique. On choisit alors un chemin générique $(J_t)_{t \in [0,1]}$ de structures presque complexes dominées par $\omega$ tel que $J_0 = J$ et $J_1$ est la structure complexe sur $A_{-1}$. Notons que l'extrémité $J_1$ n'est pas a priori pas générique, mais comme la surface $A_{-1}$ ne contient pas de section complexe d'auto-intersection négative, les arguments de la Remarque~\ref{r:morceausection} nous montrent que l'espace de modules $\mathcal{M}_{1}([T];J_1)$ est une variété lisse compacte de dimension $2$.

D'après la Proposition~\ref{p:compacitéTautoint1}, l'espace de modules à paramètre $\mathcal{M}_{1}([T];(J_t)_{t \in [0,1]})$ est alors une variété lisse orientée compacte à bord de dimension $3$, telle que la projection
$$
\begin{array}{ccc}
\mathcal{M}_{1} ([T]; (J_t)_{t \in [0,1]})  & \longrightarrow & [0,1] \\
(t,u) & \longmapsto & t 
\end{array}
$$
est une submersion.
L'espace de modules $\mathcal{M}_{1} ([T]; (J_t)_{t \in [0,1]})$ est donc difféomorphe  à $\mathcal{M}_{1} ([T]; J) \times [0,1]$ (par la théorie de Morse). Par conséquent, on peut trouver un chemin de courbes pseudoholomorphes $(u_t)_{t \in [0,1]} \in \mathcal{M}_{1} ([T]; (J_t)_{t \in [0,1]})$ tel que pour tout $t \in [0,1]$, $u_t$ est une courbe $J_t$--holomorphe et $u_0 = u$. En particulier, $u_1$ est une courbe complexe, ce qui conclut.
\end{proof}

\begin{rk}\label{r:moduledetorespherique}
Pour toute structure presque complexe $J$ générique dominée par $\omega$, les composantes connexes de l'espace de modules $\mathcal{M}_{1}(J;[T])$ sont des variétés lisses difféomorphes à $\mathbb{C}P^1$.
En effet, on a montré que pour tout chemin générique $(J_t)_{t \in [0,1]}$ tel que $J_0 = J$ et $J_1$ est la structure complexe sur $A_{-1}$, l'espace de modules à paramètre $\mathcal{M}_{1}([T];(J_t)_{t \in [0,1]})$ est une variété lisse orientée compacte à bord de dimension $3$, telle que la projection
$$
\begin{array}{ccc}
\mathcal{M}_{1} ([T]; (J_t)_{t \in [0,1]})  & \longrightarrow & [0,1] \\
(t,u) & \longmapsto & t 
\end{array}
$$
est une submersion. L'espace de modules $\mathcal{M}_{1}(J;[T])$ est  donc difféomorphe à l'espace de modules $\mathcal{M}_{1}(J_1;[T])$. Or l'espace $\mathcal{M}_{1}(J_1;[T])$ correspond exactement à une union disjointe de systèmes linéaires complets de diviseurs (on rappelle que tout diviseur irréductible $D \subset V$ est naturellement associé à un fibré en droites complexe $\mathcal{L}_D$ au-dessus de $V$ et que les diviseurs linéairement équivalents à $D$ sont les intersections des sections holomorphes non nulles de $\mathcal{L}_D$ avec la section nulle ; l'ensemble des diviseurs linéairement équivalents à un même diviseur irréductible forme un système linéaire complet de diviseurs). Chacun de ces ensembles est donc isomorphe au projectivisé d'un espace vectoriel complexe (plus précisément d'un espace de sections d'un fibré en droites complexe au-dessus de $V$, correspondant au système linéaire complet de diviseurs en question).
\end{rk}

Pour traiter le cas des autres sections des surfaces géométriquement réglées au-dessus de courbes elliptiques, on va passer par des configurations de courbes auxiliaires. 

\begin{df}
Soit $S$ une courbe symplectiquement plongée dans une surface symplectique. On appelle configuration en \emph{myriapode} de \emph{corps} $S$ et avec $k$ \emph{paires de pattes} toute configuration de courbes symplectiques constituée de l'union de $S$ et de $k \in \mathbb{N}^*$ courbes rationnelles symplectiquement plongées d'auto-intersection $0$ deux à deux distinctes et qui intersectent chacune $S$ exactement une fois, de manière transverse et positive. 
\end{df}
\begin{figure}[h]
	\centering
	\includegraphics[scale=0.4]{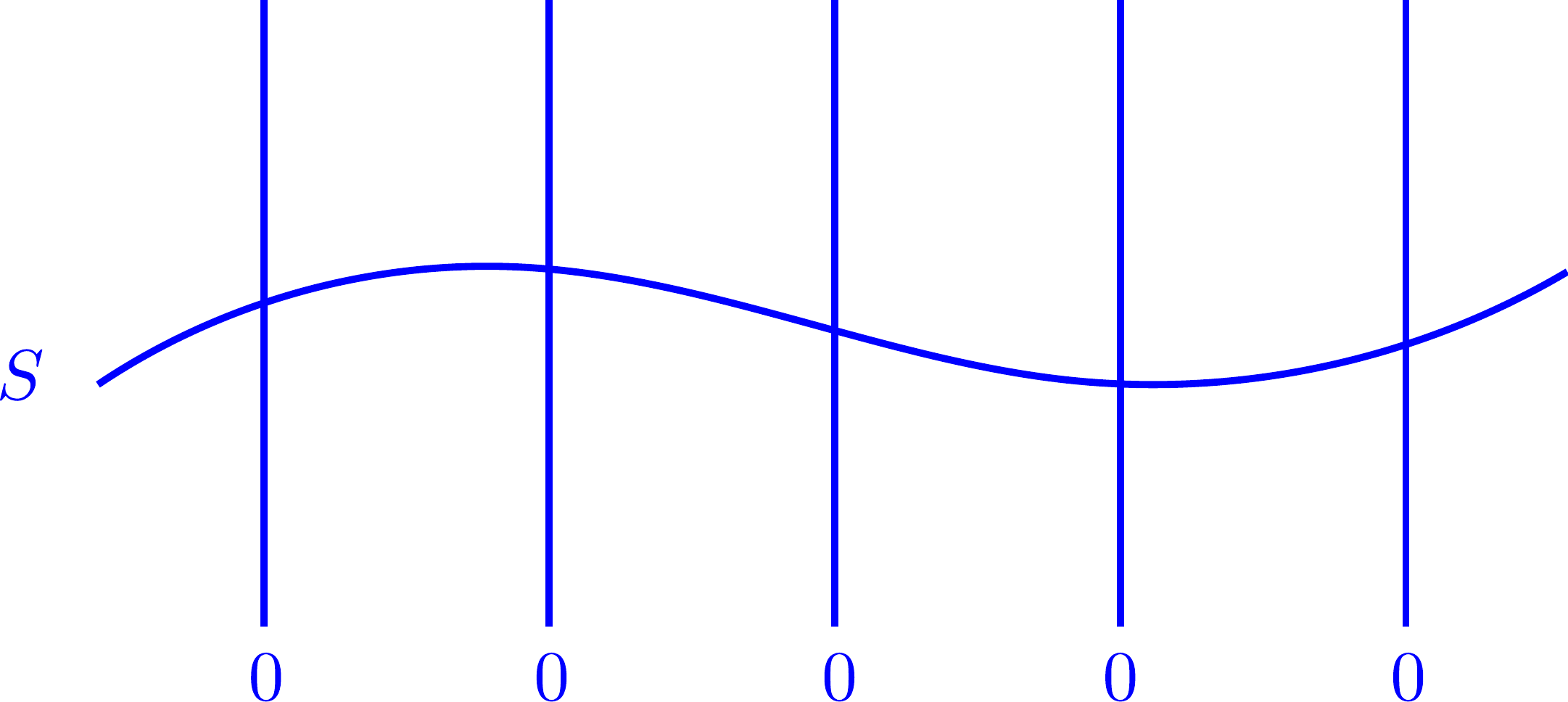}\\
	\caption{Une configuration en myriapode de corps $S$ avec $5$ paires de pattes.}
	\label{f:myriapod}
\end{figure}

\begin{rk}\label{r:ruledmyriapoda}
D'après le Théorème~\ref{t:McDuff}, toute surface symplectique qui contient une configuration en myriapode est une surface symplectiquement réglée $(M, \omega)$. Notons que les paires de pattes sont toutes homologues à une même fibre d'une fibration de Lefschetz symplectique avec fibres sphériques sur $(M, \omega)$ et que le corps est homologue à une section de cette fibration de Lefschetz.
\end{rk}

Soit $S$ une section symplectique d'une surface kählérienne géométriquement réglée. Pour tout entier $k$ strictement positif, l'union de $S$ et de $k$ fibres deux à deux distinctes est une configuration en myriapode. L'avantage des techniques de courbes pseudoholomorphes réside dans le fait de garder le contrôle sur les intersections géométriques des configurations considérées.
%

\begin{lem}\label{l:isotopieconfigTautoint1}
Soit $V$ une surface projective complexe minimale (munie d'une forme symplectique $\omega$ qui domine la structure complexe), $S$ une courbe symplectiquement plongée dans $(V, \omega)$ et $\mathcal{C}_{S,k} \subset V$ une configuration de courbes symplectiques en myriapode de corps $S$ et avec $k$ paires de pattes. Si $S$ est symplectiquement isotope à une courbe complexe via un chemin de courbes pseudoholomorphes, alors la configuration $\mathcal{C}_{S,k}$ est équisingulièrement symplectiquement isotope à une configuration de courbes complexes.
\end{lem}

\begin{proof}
D'après la Remarque~\ref{r:ruledmyriapoda} et le Théorème~\ref{t:homeoreglee}, la surface complexe $V$ est réglée, les paires de pattes de $\mathcal{C}_{S,k}$ sont homologues à des fibres et $S$ est homologue à une section. On note $F_1, \dots, F_k$ les paires de pattes de $\mathcal{C}_{S,k}$, et on fixe des points $p_1, \dots, p_k$ tels que pour tout $i \in \{1, \dots , k \}$, $p_i$ appartient à $F_i$. On choisit une structure presque complexe $J$ dominée par $\omega$ qui rend $J$--holomorphe la configuration $\mathcal{C}_{S,k}$. Quitte à perturber la structure presque complexe $J$ en dehors d'un voisinage de $S$, comme les courbes rationnelles $J$--holomorphes plongées d'auto-intersection $0$ sont automatiquement transverses pour le problème avec une contrainte ponctuelle, on peut supposer sans perte de généralité que $J$ est générique pour $p_1, \dots , p_k$. Soit $i \in \{1, \dots , k \}$, l'espace de modules $\mathcal{M}_{0,1}([F];J; p_i)$ est un point (sa dimension virtuelle est $0$ et la positivité d'intersection assure que deux éléments de cet espace de modules ne peuvent passer par le même point). On choisit désormais un chemin générique $(J_t)_{t \in [0,1]}$ de structures presque complexes dominées par $\omega$ tel que $J_0 = J$, $J_1$ est la structure complexe sur $V$ et $S$ est isotope à une courbe complexe à travers un chemin de courbes $(J_t)$--holomorphes. Comme l'espace de modules $\mathcal{M}_{0,1}([F];J_1; p_i)$ est aussi un point, d'après le Théorème~\ref{t:moduliparamreg}, l'espace de modules à paramètre $\mathcal{M}_{0,1}([F];(J_t)_{t \in [0,1]};p_i)$ est alors une variété lisse orientée compacte à bord de dimension $1$, telle que la projection
$$
\begin{array}{ccc}
\mathcal{M}_{0,1} ([F]; (J_t)_{t \in [0,1]}; p_i)  & \longrightarrow & [0,1] \\
(t,u) & \longmapsto & t 
\end{array}
$$
est une submersion.
L'espace $\mathcal{M}_{0,1} ([F]; (J_t)_{t \in [0,1]}; p_i)$ est donc difféomorphe à $[0,1]$ et on peut trouver un chemin de courbes pseudoholomorphes $(u_t)_{t \in [0,1]} \in \mathcal{M}_{0,1} ([F]; (J_t)_{t \in [0,1]}; p_i)$ tel que pour tout $t \in [0,1]$, $u_t$ est une courbe $J_t$--holomorphe et $u_0 = u$. En particulier, $u_1$ est une courbe complexe. Comme les isotopies de chacune des courbes de la configuration $\mathcal{C}_{S,k}$ se réalisent à travers des chemins de courbes $(J_t)_{t \in [0,1]}$--holomorphes, la propriété de positivité d'intersection est vérifiée à chaque instant au cours de l'isotopie. Enfin, pour tous $i,j \in \{1, \dots , k \}$, on a $[F_i] \cdot [F_j] =0$ (le fait qu'aucune courbe $J_t$--holomorphe simple homologue à $[F]$ ne passe par deux points distincts parmi $p_1, \dots, p_k$ assure que les paires de pattes de la configuration en myriapode restent deux à deux disjointes au cours de l'isotopie) et $[F_i] \cdot [S] = 1$, ce qui assure que la configuration $\mathcal{C}_{S,k}$ est \emph{équisingulièrement} symplectiquement isotope à une configuration de courbes complexes.
\end{proof}

Dans tout ce qui suit, on effectuera les opérations d'éclatement et de contraction de manière complexe. On rappelle que les éclatements et les contractions complexes peuvent également être interprétés comme des éclatements et des contractions symplectiques. On renvoit la lectrice ou le lecteur à~\cite[Section~3.1]{GS} pour de plus amples détails sur les classes d'isotopies des transformées propres de courbes symplectiques.

\begin{prop}\label{p:isotopiepbTautoint1}
Soit $T$ une courbe elliptique symplectiquement plongée dans une surface complexe géométriquement réglée $V$ au-dessus d'une courbe elliptique. Si $T$ est une section d'auto-intersection $1$, alors elle est symplectiquement isotope à une section complexe. 
\end{prop}

\begin{proof}
Comme $[T]^2=1$, la forme d'intersection de $V$ est impaire. Donc $V$ est biholomorphe à $A_{-1}$ ou à $S_{2k+1}$ pour $k \in \mathbb{N}$. On souhaite se ramener au cas où $V \simeq A_{-1}$ en utilisant des transformations birationnelles. On suppose dans un premier temps que $V \simeq S_{2k+1} $. On note $\Sigma_+$ une section complexe d'auto-intersection $2k+1$, $\Sigma_-$ la section complexe d'auto-intersection $-2k-1$ et $F$ une fibre complexe (on rappelle que $\Sigma_+$ et $\Sigma_-$ sont disjointes). 
Comme le groupe $H_2 (V ;\mathbb{Z})$ est engendré par $[T]$ et $[F]$, et puisque $[\Sigma_+] \cdot [F] =1$ et $[\Sigma_+]^2 = 2k+1$, on a 
$$[\Sigma_+] = [T] + k [F ].$$
On en déduit que $[\Sigma_+] \cdot [T] = 1 +k$. Quitte à réaliser une perturbation $\mathcal{C}^\infty$ arbitrairement petite de $T$ (d'après les résultats de~\cite{HI10}, deux courbes symplectiques arbitrairement proches pour la topologie $\mathcal{C}^0$ sont symplectiquement isotopes), on peut supposer que $T$ et $\Sigma_+$ s'intersectent transversalement. Comme $[\Sigma_+] \cdot [T] = 1 +k$, il existe au moins $k+1$ points d'intersection transverses entre $\Sigma_+$ et $T$. On choisit $k$ de ces points d'intersection $p_1, \dots, p_k$ et $k$ autres points génériques $p_{k+1}, \dots, p_{2k}$ (n'appartenant pas à $T$, $\Sigma_+$, $\Sigma_-$ ou aux fibres passant par les autres $p_i$). On note pour tout $i \in \{1, \dots, 2k \}$, $F_i$ la fibre complexe passant par $p_i$. On effectue une transformation élémentaire de Nagata (on éclate le point concerné et on contracte la transformée propre de la fibre passant par ce point) en chacun des points $p_1, \dots, p_{2k}$.
On note $V'$ la surface complexe géométriquement réglée ainsi obtenue et $T'$, $\Sigma_+'$, $\Sigma_-'$, $F_i'$ les images respectives de $T$, $\Sigma_+$, $\Sigma_-$, $F_i$ par cette transformation birationnelle.
Comme les $2k$ transformations de Nagata ont été réalisées en des points distincts de $\Sigma_-$, on a $[\Sigma_-']^2 = [\Sigma_-]^2 +2k = -1$. De même, comme autant de transformations élémentaires de Nagata ont été effectuées sur $T$ et en dehors de $T$, on a $[T']^2 = [T]^2 = 1$. Comme la surface complexe $V'$ contient une section complexe d'auto-intersection $-1$, elle est biholomorphe à la surface complexe $S_1$ d'après la classification des surfaces réglées au-dessus de courbes elliptiques. La configuration en myriapode $\mathcal{C}_{T',2k}$ dans $V'$ de corps $T'$ avec paires de pattes $F_1', \dots, F_{2k}'$ est birationnellement dérivée de la configuration en myriapode $\mathcal{C}_{T,2k}$ dans $V$ de corps $T$ avec paires de pattes $F_1, \dots, F_{2k}$. Par la Proposition~\ref{p:birationalequivalence}, la configuration $\mathcal{C}_{T',2k}$ est équisingulièrement isotope à une configuration de courbes complexes si et seulement si la configuration $\mathcal{C}_{T,2k}$ est équisingulièrement isotope à une configuration de courbes complexes. On s'est donc ramené au cas où la surface complexe est biholomorphe à $S_1$.

On choisit désormais un point générique $p$ sur $T'$ et un point générique $q$ en dehors de $T'$. En particulier aucun de ces deux points n'appartient à $\Sigma_-'$ et aucune section complexe d'auto-intersection $1$ ne passe à la fois par $p$ et $q$ (en effet, au plus deux sections d'auto-intersection $1$ passent par un même point dans $A_{-1}$, voir~\cite{diaw}). On effectue les transformations élémentaires de Nagata en $p$ et $q$. On obtient alors une surface complexe géométriquement réglée $V''$, et on note $T''$, $S_- ''$ les images respectives de $T'$, $S_-'$ par cette transformation birationnelle. Pour montrer que $V''$ est biholomorphe à $A_{-1}$, il suffit de montrer qu'elle est stable d'après le résultat de classification de Suwa. 

Par le Lemme~\ref{l:sectionsnomansland}, toute section complexe $\Sigma$ de $V'$ différente de $\Sigma_-'$ vérifie $[\Sigma]^2 \geq 1$. De plus, si $[\Sigma]^2 > 1$, le Lemme~\ref{l:paritésection} nous assure que $[\Sigma]^2 \geq 3$. 

Comme $p$ et $q$ sont disjoints de $\Sigma_-'$, on a  $[\Sigma_-'']^2 = [\Sigma_-']^2 +2 = 1$. Soit $\Sigma$ une section de $V'$ différente de $\Sigma_-'$ et notons $\Sigma'$ son image par la transformation birationnelle vers $V''$. Si $[\Sigma]^2 = 1$, alors au plus un des deux points $p$ et $q$ appartient à $\Sigma$, donc $[\Sigma']^2 \geq [\Sigma]^2  =1$. Si $[\Sigma]^2 > 1$, alors on a $[\Sigma']^2 \geq [\Sigma]^2 -2 >0$. Ainsi, la surface complexe $V''$ est stable, donc biholomorphe à $A_{-1}$. De plus, on a $[T'']^2=[T']^2=1$.

D'après les Lemmes~\ref{l:A-1isotopiepbTautoint1} et~\ref{l:isotopieconfigTautoint1}, la configuration en myriapode $\mathcal{C}_{T'',2}$ dans $V''$, constituée de l'union de $T''$ et des deux fibres obtenues suites aux transformations élémentaires de Nagata en $p$ et $q$, est équisingulièrement symplectiquement isotope à une configuration de courbes complexes.

De plus, la configuration $\mathcal{C}_{T'',2}$ dans $V''$ est birationnellement dérivée de la configuration $\mathcal{C}_{T',2}$ dans $V'$ définie comme l'union de $T'$ et des deux fibres passant respectivement par les points $p$ et $q$. La Proposition~\ref{p:birationalequivalence} nous permet alors d'affirmer que la configuration $\mathcal{C}_{T',2}$ est équisingulièrement isotope à une configuration de courbes complexes.
\end{proof}

\begin{rk}
D'après la Proposition~\ref{p:isotopiepbTautoint1}, pour $k \neq 1$, les surfaces complexes $S_{k}$ ne contiennent pas de sections symplectiques d'auto-intersection $1$.
\end{rk}

On peut maintenant démontrer le théorème central de ce chapitre.

\begin{proof}[Démonstration du Théorème~\ref{t:isotopiesectionell}]
On note $m=[T]^2$. 
Si $m=1$, on conclut avec la Proposition~\ref{p:isotopiepbTautoint1}. 
Si $m > 1$ on effectue $m-1$ transformations élémentaires de Nagata  sur $T$, sur des fibres deux à deux distinctes. On note $T'$ l'image de $T$ par cette transformation birationnelle. On a $[T']^2 = [T]^2 -(m-1) =1$. On obtient de cette manière une configuration en myriapode $\mathcal{C}_{T',m-1}$ de corps $T'$ avec $m-1$ paires de pattes qui est birationnellement dérivée d'une configuration en myriapode $\mathcal{C}_{T,m-1}$ de corps $T$ avec $m-1$ paires de pattes. D'après la Proposition~\ref{p:isotopiepbTautoint1} et le Lemme~\ref{l:isotopieconfigTautoint1}, la configuration $\mathcal{C}_{T',m-1}$ est équisingulièrement isotope à une configuration de courbes complexes. Par conséquent, la Proposition~\ref{p:birationalequivalence} nous assure que la configuration $\mathcal{C}_{T,m-1}$ est également équisingulièrement symplectiquement isotope à une configuration de courbes complexes, ce qui conclut.
Si $m < 1$, on procède de manière similaire en effectuant cette fois $-m+1$ transformations élémentaires de Nagata en dehors de $T$, sur des fibres deux à deux distinctes. Cette fois encore, l'image $T'$ de $T$ par cette transformation birationnelle vérifie $[T']^2 = [T]^2 -m+1 =1$. On obtient de cette manière une configuration en myriapode $\mathcal{C}_{T',-m+1}$ de corps $T'$ avec $-m+1$ paires de pattes qui est birationnellement dérivée d'une configuration en myriapode $\mathcal{C}_{T,-m+1}$ de corps $T$ avec $-m+1$ paires de pattes. On conclut de la même manière que précédemment.
\end{proof}

\begin{cor}
Soit $V$ une surface complexe et $T$ une courbe elliptique symplectiquement plongée dans $V$. On suppose que $V$ admet une fibration de Lefschetz complexe générique, dont la fibre générique est $\mathbb{C}P^1$, au-dessus d'une courbe elliptique. Si $T$ est une section, alors elle est symplectiquement isotope à une section complexe. 
\end{cor}

\begin{proof}
Comme la fibration de Lefschetz sur $V$ est générique, les fibres singulières sont constituées de l'union de deux diviseurs exceptionnels qui s'intersectent exactement une fois. On contracte le diviseur exceptionnel qui intersecte $T$ dans chaque fibre singulière de façon à obtenir une surface géométriquement réglée $V'$. L'image par ces contractions de la configuration de courbes symplectiques $\mathcal{C}$, formée de l'union de $T$ et des diviseurs exceptionnels contractés, est une section $T'$ de $V'$. D'après le Théorème~\ref{t:isotopiesectionell}, la courbe $T'$ est symplectiquement isotope à une courbe complexe. Comme $T'$ est birationnellement dérivée de la configuration $\mathcal{C}$, la Proposition~\ref{p:birationalequivalence} nous assure que $T$ est également symplectiquement isotope à une courbe complexe.
\end{proof}


\begin{rk}
Deux courbes symplectiquement plongées, symplectiquement isotopes à des courbes complexes et qui appartiennent à un même système linéaire de diviseurs sont symplectiquement isotopes. La connaissance de ces systèmes linéaires de diviseurs (tous isomorphes à des espaces projectifs complexes) permet de donner une borne supérieure sur le nombre de classes d'isotopie symplectique d'un type (genre et auto-intersection) de courbe donné.
\end{rk}

%
%
%

Il serait souhaitable de trouver un moyen d'améliorer le Théorème~\ref{t:isotopiesectionell} en affaiblissant l'hypothèse sur $T$ d'être une section en le fait de juste être homologue à une section. Une telle amélioration, combinée à la Proposition~\ref{p:Tdansfibré}, permettrait de résoudre complètement le problème d'isotopie symplectique dans les surfaces kählériennes minimales pour les courbes elliptiques symplectiquement plongées dont l'auto-intersection est strictement positive. On rappelle en effet que d'après le Théorème~\ref{t:thmprincipal}, les courbes elliptiques symplectiquement plongées d'auto-intersection $9$ dans les surfaces kählériennes minimales qui ne sont pas homologues à des sections de surfaces géométriquement réglées sont des cubiques non singulières du plan projectif complexe (qui sont symplectiquement isotopes à des courbes complexes d'après~\cite{Sikorav}) ; et les courbes elliptiques symplectiquement plongées d'auto-intersection $8$ dans les surfaces kählériennes minimales qui ne sont pas homologues à des sections de surfaces géométriquement réglées sont des courbes dans des surfaces géométriquement réglées \emph{rationnelles} pour lesquelles la restriction de l'application de projection est de degré $2$ (qui sont symplectiquement isotopes à des courbes complexes d'après~\cite{SiebertTian}). En ce qui concerne les courbes elliptiques symplectiquement plongées d'auto-intersection $0$, Fintushel et Stern~\cite{FS99} ont construit des familles infinies de telles courbes deux à deux non isotopes (de manière lisse).

\clearemptydoublepage

\appendix
\chapter{Classification des dégénérescences possibles des courbes pseudoholomorphes plongées d'indice contraint égal à $2$}\label{c:annexe}
On s'intéresse aux courbes nodales pouvant apparaître dans les compactifiés de Gromov des espaces de modules de courbes pseudoholomorphes plongées dans les surfaces presque complexes, satisfaisant des contraintes ponctuelles et dont l'indice contraint est égal à $2$. On rappelle qu'une courbe nodale possède des points nodaux qui viennent par paires, que ces points nodaux sont répartis sur les domaines des différentes composantes de la courbe nodale et que les points nodaux d'une même paire sont envoyés par la courbe sur un même point (voir la Sous-section~\ref{s:CompactificationGromov}). On fera usage de la proposition suivante, dont la démonstration repose essentiellement sur la formule de l'indice et sur la définition de la compactification de Gromov.

\begin{prop}[{\cite[Proposition~4.8]{wendl2014lectures}}]\label{p:morceauxwendl}
Soit $(M,J)$ une surface presque complexe et $u_\infty$ une courbe $J$--holomorphe nodale non constante dans $(M,J)$. On note  $u^1_\infty, \dots, u^\ell_\infty$ les composantes irréductibles de $u_\infty$ et $N_i$ le nombre de points nodaux sur chaque composante $u^i_\infty$. Alors 
$$\ind (u_\infty) \geq \sum_{\{i \vert u_\infty^i\ \neq c^{ste}\}} \left( \ind(u_\infty^i) +N_i \right),$$
avec égalité si et seulement si la courbe $u_{\infty}$ n'a pas de composante constante. En particulier, dès que $u_\infty$ possède des points nodaux, on a
$$ \ind (u_\infty) \geq 2 + \sum_{\{i \vert u_\infty^i\ \neq c^{ste}\}} \ind(u_\infty^i).$$
\end{prop}

On s'en sert pour démontrer la proposition suivante.

\begin{prop}\label{p:morceauxcourbesind2}
Soit $(M,J)$ une surface presque complexe et $u$ une courbe $J$--holomorphe de genre $g$ plongée dans $(M,J)$. On pose $m= [u]^2-g$ et $p_1, \dots, p_m$ des points deux à deux distincts dans l'image de $u$ de sorte que la courbe $u$ est d'indice contraint égal à $\ind (u)-2m = 2$. On suppose que $J$ est générique pour les points $p_1, \dots, p_m$. Si $g \neq 4$ ou $[u]^2 \neq 4$, alors toute courbe $u_\infty \in \overline{\M}_{g,m}([u];J;p_1, \dots, p_m)$ qui n'est pas plongée est d'une des formes suivantes :
\begin{enumerate}
\item la courbe $u_\infty$ est constituée d'une unique composante  immergée de genre $g-1$, dont l'unique point multiple est un point double transverse disjoint des points $p_1, \dots, p_m$, ou bien 
\item la courbe $u_\infty$ est constituée de deux composantes plongées $u_1$ et $u_2$, de genres respectifs $g_1$ et $g_2$ avec $g=g_1+ g_2$, qui s'intersectent transversalement en un unique point distinct des points $p_1, \dots, p_m$. De plus, chacune des courbes $u_i$ est d'indice égal à $2m_i$, où $m_i$ correspond au nombre de contraintes ponctuelles dans l'image de $u_i$.
\end{enumerate}
Si $g=4$ et $[u]^2=4$ (donc $m =0$), en plus des possibilités évoquées ci-dessus, une courbe nodale $u_\infty$ peut également être d'une des formes suivantes :
\begin{enumerate}
\item[3.] la courbe $u_\infty$ est constituée d'une unique composante  qui est le revêtement double d'une courbe plongée de genre $2$ et d'auto-intersection $1$, ou bien 
\item[4.] la courbe $u_\infty$ est constituée de deux composantes plongées $u_1$ et $u_2$ de genre $2$ et d'auto-intersection $1$ ayant la même image. 
\end{enumerate}
\end{prop}

\begin{rk}
Il est en pratique difficile d'appliquer la Proposition~\ref{p:morceauxcourbesind2} lorsque $g \geq 2$. En effet, le meilleur moyen de satisfaire les hypothèses de généricité consiste à appliquer le Théorème~\ref{t:transautoSiko} de transversalité automatique avec contraintes afin de pouvoir perturber la structure presque complexe. Sans hypothèse de transversalité sur $u$ pour le problème avec contraintes ponctuelles $p_1, \dots, p_m$, rien ne garantit l'existence d'une courbe pseudoholomorphe adéquate.
\end{rk}

Pour démontrer la Proposition~\ref{p:morceauxcourbesind2}, on aura besoin des lemmes suivants.

\begin{lem}\label{l:morceaux1}
Soit $(M,J)$ une surface presque complexe. Soit $u : \S \rightarrow M$ une courbe $J$--holomorphe plongée de genre $g$ telle que $\ind(u) \geq 2$ et $u_\infty \in \M_{g}([u];J)$. On suppose que $u_\infty$ est un revêtement de degré $k \geq 1$ d'une courbe $J$--holomorphe simple $v : \L \rightarrow M$. Si le revêtement n'admet pas de point de ramification (c'est-à-dire si $\chi (\S) = k \chi (\L)$), alors on a $k=1$ et $u_\infty$ est plongée.
\end{lem}

\begin{proof}
On calcule l'indice de $u$ de deux manières différentes.
On a 
\begin{align*}
\ind (u) 
&= -\chi(\S) + 2 c_1([u]) \\
&= -\chi(\S) +2k c_1([v]) \\
&= -\chi(\S) +2k \chi(\L) +2k[v]^2-4k \delta (v) \quad \text{par la formule d'adjonction,} \\
&= \chi(\S) +2k [v]^2-4k \delta (v) \quad \text{par hypothèse.}
\end{align*}
Puisque la courbe $u$ est plongée, on obtient aussi les égalités suivantes en utilisant la formule d'adjonction 
\begin{align*}
\ind (u) 
&= \chi(\S) + 2 [u]^2 \\
&= \chi(\S) + 2 k^2 [v]^2.
\end{align*}
Remarquons que puisque $\ind (u) \geq 2$ et $\chi(\S)\leq 2$, on a $[u]^2 \geq 0$. Si $[u]^2 = 0$, alors $\chi(\S)= 2$. Comme $\chi (\S) = k \chi (\L)$, on a alors $\chi (\L)=2$ et $k=1$. On suppose désormais qu'on a $[u]^2 > 0$. En utilisant les deux façons de calculer l'indice de $u$, on obtient $2 k^2 [v]^2 = 2k [v]^2-4k \delta (v)$, autrement dit $ (k-1) [v]^2 + 2 \delta (v) =0$.
Comme tous les termes dans la somme sont positifs et puisque $[v]^2 =\frac{1}{k^2} [u]^2 > 0$, on obtient $k=1$ et $\delta (v) =0$.
\end{proof}

\begin{lem}\label{l:morceauxrelou}
Soit $(M,J)$ une surface presque complexe. Soit $u : \S \rightarrow M$ une courbe $J$--holomorphe plongée de genre $g$ telle que $\ind(u) \geq 2$ et $u_\infty \in \M_{g}([u];J)$. On suppose que $u_\infty$ est un revêtement de degré $k \geq 1$ d'une courbe $J$--holomorphe simple $v : \L \rightarrow M$. Si $k \chi (\L) -\chi (\S) =2 $, alors on a $k=2$, $[u]^2=4$ et $v$ est plongée.
\end{lem}

\begin{proof}
Si $k=1$, alors $u_\infty$ est une reparamétrisation de $v$, ce qui est absurde puisque $\chi (\L) -\chi (\S) =2$. On a donc $k \geq 2$.
On procède maintenant de manière similaire à la démonstration du Lemme~\ref{l:morceaux1}. On calcule l'indice de $u$ de deux manières différentes.
On a d'un côté
\begin{align*}
\ind (u) 
&= -\chi(\S) +2k \chi(\L) +2k[v]^2-4k \delta (v) \quad \text{par la formule d'adjonction,} \\
&= \chi(\S) +2k [v]^2-4k \delta (v) +4 \quad \text{par hypothèse.}
\end{align*}
Puisque la courbe $u$ est plongée, on obtient aussi en utilisant la formule d'adjonction 
\begin{align*}
\ind (u) 
&= \chi(\S) + 2 [u]^2 \\
&= \chi(\S) + 2 k^2 [v]^2.
\end{align*}
Si $[u]^2 = 0$, alors $\chi(\S)= 2$ et $\ind (u) =2$. Comme $k \chi (\L) -\chi (\S) =2$, on a alors $k\chi (\L)=4$, d'où $\chi (\L)=2$ et $k=2$. Mais $v$ est alors une courbe rationnelle d'auto-intersection $0$, donc $\ind(v) =2$, ce qui contredit le fait que d'après la formule de Riemann--Hurwitz pour les indices (Proposition~\ref{p:RHindiceversion}), on a $\ind(u) - 2 \ind(v) =2$. On a donc $[u]^2 >0$ et $[v]^2 = \frac{1}{k^2}[u]^2 >0$.
En utilisant les deux façons de calculer l'indice de $u$, on obtient $2 k^2 [v]^2 = 2k [v]^2-4k \delta (v) +4$, autrement dit $k(k-1) [v]^2 + 2k \delta (v) =2$.
Puisque $k(k-1) [v]^2 >0$ et $2k\delta (v) \geq 0$, on a $k(k-1) [v]^2 =2$. On en conclut que $k=2$, $[v]^2 =1$, $[u]^2 = k^2 [v]^2 =4$ et $\delta (v) =0$.
\end{proof}

\begin{lem}\label{l:morceauxmegarelou}
Soit $(M,J)$ une surface presque complexe. Soit $u : \S \rightarrow M$ une courbe $J$--holomorphe plongée de genre $g$ telle que $\ind(u) \geq 2$ et $u' \in \M_{g-1}([u];J)$. On suppose que $u': \Sigma' \rightarrow M$ est un revêtement de degré $k \geq 1$ d'une courbe $J$--holomorphe simple $v : \L \rightarrow M$. Si le revêtement n'admet pas de point de ramification (c'est-à-dire si $\chi(\S)+2 = \chi (\S') = k \chi (\L)$), alors on a $k=1$ et $\delta (v)=1$, ou bien $k=2$, $[u]^2 =4$ et $v$ est plongée.
\end{lem}

\begin{proof}
On procède de manière similaire à la démonstration du Lemme~\ref{l:morceaux1}. On calcule l'indice de $u$ de deux manières différentes.
On a d'un côté, en appliquant la formule d'adjonction à $v$,
\begin{align*}
\ind (u) 
&= -\chi(\S) +2k \chi(\L) +2k[v]^2-4k \delta (v)\\ 
&= \chi(\S) +2k [v]^2-4k \delta (v)+4\quad \text{par hypothèse.}
\end{align*}
De l'autre côté, puisque la courbe $u$ est plongée, on obtient grâce à la formule d'adjonction appliquée à $u$,
$$
\ind (u) = \chi(\S) + 2 k^2 [v]^2.
$$
Remarquons que puisque $\ind (u) \geq 2$ et $\chi( \S ) \leq 2$, on a $[u]^2 \geq 0$. Si $[u]^2 = 0$, alors $\chi ( \S )= 2$, d'où $\chi (\Sigma' ) = 4$, ce qui est absurde. On a donc $[u]^2 > 0$. En utilisant les deux façons de calculer l'indice de $u$, on obtient 
$k(k-1) [v]^2 + 2k \delta (v) =2$.
Comme tous les termes dans la somme sont positifs, et puisque $[v]^2 =\frac{1}{k^2} [u]^2 > 0$, on obtient $k=1$ et $\delta (v)=1$, ou bien $k=2$, $[v]^2=1$ (c'est-à-dire $[u]^2 =4$) et $\delta (v) =0$.
\end{proof}

%

\begin{lem}\label{l:lediviseurdelexception}
Soit $(M,J)$ une surface presque complexe et $u : \S \rightarrow M$ une courbe $J$--holomorphe simple non constante. Si $\ind(u) \geq 0$ et $u$ n'est pas un diviseur exceptionnel, alors on a $[u]^2 \geq 0$.
\end{lem}

\begin{proof}
D'après la formule d'adjonction, on a 
\begin{align*}
\ind (u) 
&= -\chi(\S) + 2 c_1([u]) \\
&= \chi(\S) + 2 [u]^2 - 4 \delta (u).
\end{align*}
On obtient donc $[u]^2 \geq - \frac{1}{2} \chi(\S) +2 \delta (u)  \geq -1$, avec $[u]^2 = -1$ si et seulement si $\chi(\S) = 2 $, $\delta (u)=0$ et $\ind(u)=0$, c'est-à-dire si et seulement si $u$ est un diviseur exceptionnel.
\end{proof}

\begin{lem}\label{l:morceaux2}
Soit $(M,J)$ une surface presque complexe. Soit $u : \S \rightarrow M$ une courbe $J$--holomorphe plongée telle que $\ind(u) \geq 2$ et $u_\infty \in \overline{\M}_{g}([u];J)$ une courbe nodale. On suppose que $u_\infty$ est constituée de deux composantes $u_1$ et $u_2$ non constantes reliées par une unique paire de points nodaux. Si les courbes $u_1$ et $u_2$ sont des revêtements non ramifiés de degrés respectifs $k_1 \geq 1$ et $k_2 \geq 1$ de courbes $J$--holomorphes simples $v_1 : \L_1 \rightarrow M$ et $v_2 : \L_2 \rightarrow M$ telles que $\ind(v_1) \geq 0$ et $\ind(v_2) \geq 0$, alors $k_1 = k_2 =1$. De plus, $u_1$ et $u_2$ sont deux composantes plongées telles que $[u_1] \cdot [u_2] =1$, dont les images sont distinctes dès que $[u]^2 \neq 4$.
\end{lem}

\begin{proof}
Puisque $u_\infty$ n'a qu'une paire de points nodaux, on a $ \chi(\S) = \chi(\S_1)+\chi(\S_2) -2$. Par hypothèse on a également pour tout $i\in \{1,2 \}$, $\chi (\S_i) = k_i \chi (\L_i)$. En s'aidant de la formule d'adjonction, appliquée aux courbes $J$--holomorphes simples $v_1$ et $v_2$, on obtient alors les égalités suivantes :
\begin{align*}
\ind (u) 
&= \ind (u_1) +\ind (u_2) +2\\
&= -\chi(\S_1) + 2 c_1([u_1]) -\chi(\S_2) + 2 c_1([u_2]) +2 \\
&= -\chi(\S_1) +2k_1 c_1([v_1]) -\chi(\S_2) +2k_2 c_1([v_2])+2 \\
&= -\chi(\S_1) +2k_1 \chi(\L_1) +2k_1[v_1]^2-4k_1 \delta (v_1) -\chi(\S_2) +2k_2 \chi(\L_2) +2k_2[v_2]^2
\\ & \quad -4k_2 \delta (v_2) +2 \\ 
&= \chi(\S_1)+\chi(\S_2) +2k_1 [v_1]^2-4k_1 \delta (v_1) +2k_2 [v_2]^2-4k_2 \delta (v_2) +2 \\
&= \chi(\S)+2k_1 [v_1]^2-4k_1 \delta (v_1) +2k_2 [v_2]^2-4k_2 \delta (v_2) +4.
\end{align*}
Puisque la courbe $u$ est plongée, on obtient aussi les égalités suivantes en utilisant la formule d'adjonction :
\begin{align*}
\ind (u) 
&= \chi(\S) + 2 [u]^2 \\
&= \chi(\S) + 2 k_1^2 [v_1]^2 + 4k_1k_2 [v_1] \cdot [v_2] +2 k_2^2 [v_2]^2.
\end{align*}
En combinant les deux égalités précédentes on a 
$$2k_1 [v_1]^2-4k_1 \delta (v_1) +2k_2 [v_2]^2-4k_2 \delta (v_2) +4 =   2 k_1^2 [v_1]^2 + 4 k_1k_2[v_1] \cdot [v_2] +2 k_2^2 [v_2]^2,$$
donc
$$k_1 (k_1 -1)[v_1]^2  + k_2 (k_2 -1)[v_2]^2 + 2k_1 \delta(v_1) + 2k_2 \delta(v_2) +2(k_1 k_2 [v_1] \cdot [v_2] -1) =0.$$
D'après le Lemme~\ref{l:lediviseurdelexception}, pour tout $i \in \{1,2\}$, si $[v_i]^2 < 0$, alors $v_i$ est un diviseur exceptionnel, et on a en particulier $\chi (\S_i) =2 k_i$, donc $k_i =1$. Les termes dans la somme décrite ci-dessus sont donc tous positifs, et on a par conséquent $k_1 =k_2=1$, $\delta(v_1) = \delta(v_2) = 0$ et $[v_1] \cdot [v_2] =1$. Si les deux courbes $u_1$ et $u_2$ ont la même image, alors on a nécessairement $[u_1]^2 =[u_2]^2=[v_1] \cdot [v_2] =1$, donc $[u]^2 =4$. 
\end{proof}

\begin{proof}[Démonstration de la Proposition~\ref{p:morceauxcourbesind2}]

On note $u_{\infty}^1, \dots, u_{\infty}^\ell$ les composantes irréductibles de $u_{\infty}$, $g_i$ le genre de la composante $u_{\infty}^i$ et $N_i$ le nombre de points nodaux sur la composante $u_{\infty}^i$.
Chaque contrainte ponctuelle appartient à l'image d'une composante non constante. En considérant chacune des composantes individuellement, on peut donc répartir des points marqués sur les composantes non constantes de façon à ce que 
\begin{itemize}
\item chaque composante $u_\infty^i$ possède $m_i$ points marqués qui sont envoyés sur des points deux à deux distincts de $\{p_1, \dots, p_m\}$,
\item $m_i = 0$ si la composante $u_{\infty}^i$ est constante,
\item $m = \sum_{i=1}^{\ell} m_i$.
\end{itemize}

Il est alors possible de voir chaque composante non constante $u_\infty^i$ comme un élément de l'espace de modules $\M_{{g_i},{m_i}}([u_{\infty}^i];J)$ qui envoie ses $m_i$ points marqués sur des points deux à deux distincts parmi $p_1,\dots, p_m$. Puisque la courbe $u_\infty^i$ appartient à un espace de modules avec $m_i$ points contraintes, son indice contraint est égal à $\ind(u_\infty^i) -2m_i$. Rappelons que par la généricité de $J$ par rapport aux points $p_1, \dots, p_m$, l'indice contraint $\ind(u_\infty^i) -2m_i$ est positif dès que la courbe $u_\infty^i$ est injective quelque part. Par la Proposition~\ref{p:morceauxwendl}, on a 
$$\ind(u_\infty) -2m \geq \sum_{\{i \vert u_\infty^i\ \neq c^{ste}\}} \left( \ind(u_\infty^i) -2m_i +N_i \right),$$
avec égalité si et seulement si la courbe $u_{\infty}$ n'a pas de composante constante. 
Pour tout $i$ tel que la courbe $u_\infty^i$ n'est pas constante, on écrit $u_\infty^i= v^i \circ \varphi_i$, où $v^i$ est une courbe $J$--holomorphe simple de genre inférieur ou égal à $g_i$, et $\varphi_i : \S_i \ra \L_i$ est une application holomorphe de degré $k_i \geq 1$ entre des surfaces de Riemann. Par la formule de Riemann-Hurwitz, et puisque $[u^i]= k_i [v^i]$, on obtient $\mathrm{Z}(d\varphi_i) = k_i \chi (\L_i) -\chi(\S_i) = \ind(u_\infty^i) - k_i \ind(v^i) \geq 0$.
Par conséquent, on a 
\begin{align*}
\ind (u_\infty^i) - 2 m_i 
&= k_i \ind(v^i) + \mathrm{Z}(d\varphi_i) -2m_i \\
&= k_i (\ind(v^i) - 2m_i) + k_i \chi(\L_i) - \chi(\S_i) +2m_i(k_i -1).
\end{align*}
En utilisant l'inégalité précédente, on obtient
\begin{align*}
2 &= \ind(u_\infty) -2m \\
&\geq \sum_{\{i \vert u_\infty^i\ \neq c^{ste}\}} k_i ( \ind(v^i) -2m_i) + k_i \chi(\L_i) -\chi(\S_i) +2m_i(k_i-1) +N_i \\
&\geq 0.
\end{align*}
Remarquons que par généricité de $J$ par rapport aux points $p_1, \dots, p_m$, pour tout $i$ tel que $u^i_\infty$ est non constante, on a $\ind(v^i) -2m_i \geq 0$.

\begin{flushleft}
\emph{Premier cas} : $\sum\limits_{\{i \vert u_\infty^i\ \neq c^{ste}\}} N_i =2$.
\end{flushleft}

Si le nombre de points nodaux sur les composantes non constantes de $u_\infty$ est égal à $2$, alors la première inégalité est une égalité. Par conséquent, la courbe $u_\infty$ ne possède pas de composante constante et il existe au plus deux composantes non constantes.

On suppose dans un premier temps que la courbe $u_\infty$ ne possède qu'une seule composante non constante $u_\infty^1$, qui est nécessairement de genre $g-1$. En regroupant les termes dans la somme de manière à ce que chacun d'entre eux soit positif (en utilisant la généricité de $J$ par rapport aux points $p_1,\dots, p_m$), on obtient 
$$\left\{ 
\begin{array}{cc}
k_1( \ind(v^1) -2m) &= 0\\
k_1 \chi( \L_1) -\chi(\S_1) &=0 \\
2m(k_1-1) &=0
\end{array}
\right..$$

D'après le Lemme~\ref{l:morceauxmegarelou}, on a deux possibilités : $k_1 = 1$ et $\delta(v^1)=1$, ou bien $k_1=2$, $[u]^2 =4$ et $v^1$ est plongée (auquel cas on a également $m=0$).
Pour la première possibilité, $u_\infty$ est une courbe $J$--holomorphe simple avec une unique composante. De plus, comme l'indice contraint de $v^1$ est nul, si $J$ est choisie suffisamment générique au départ, la courbe $v^1$ est immergée (plus précisément, pour une structure presque complexe générique, toutes les courbes pseudoholomorphes simples d'indice contraint nul par rapport aux contraintes ponctuelles $p_1, \dots,p_m$ sont immergées, voir le Théorème~\ref{t:crit} et le Théorème~\ref{t:critk}), et son unique point double est distinct des points $p_1, \dots, p_m$. Ainsi, $u_\infty$ est une courbe $J$--holomorphe immergée de genre $g-1$, avec un unique point multiple qui est un point double transverse distinct des points $p_1, \dots, p_m$. Pour la seconde possibilité, on a nécessairement $m=0$. Comme $u$ est une courbe plongée, on a alors $2 = \ind (u) = \chi(\Sigma) +2[u]^2$, donc $g=4$. De même, comme $\ind(v^1) = 0$ et $v^1$ est plongée, on a $\chi(v^1) = -2$.

Supposons à présent que $u_\infty$ possède deux composantes. En regroupant comme précédemment les termes dans la somme de manière à ce que chacun d'entre eux soit positif, on obtient pour $i \in \{1,2 \}$,
$$\left\{ 
\begin{array}{cc}
k_i( \ind(v^i) -2m_i) &= 0\\
k_i \chi( \L_i) -\chi(\S_i) &=0 \\
2m_i(k_i-1) &=0\\
\end{array}
\right..$$
Remarquons qu'on a également $m= m_1 +m_2$ et, puisque le genre de $u$ est égal à $g$, on a $\chi(\S_1) +\chi(\S_2) = \chi(\S) +2$.
D'après le Lemme~\ref{l:morceaux2}, on a $k_1 = k_2 =1$. De plus, $u^1_\infty$ et $u^2_\infty$ sont deux composantes plongées telles que $[u^1_\infty] \cdot [u^2_\infty] =1$, qui sont distinctes si $[u]^2 \neq 4$. Remarquons finalement que $u^1_\infty$ et $u^2_\infty$ sont d'indice contraint nul par rapport à leurs contraintes respectives. Elles ne peuvent donc pas passer par d'autres contraintes ponctuelles parmi $p_1, \dots, p_m$, sous peine de contredire la généricité de $J$. Par conséquent, si $m \neq 0$, les courbes $u^1_\infty$ et $u^2_\infty$ ont des images distinctes. Si $m=0$ et $[u]^2=4$, on a $2=\ind (u) =\chi(\Sigma)+ 2[u]^2$, donc $g=4$.
Comme les indices contraints des courbes $u^1_\infty$ et $u^2_\infty$ par rapport à leurs contraintes respectives sont nuls, si les images de ces deux courbes sont distinctes, leur unique point d'intersection est distinct des points $p_1, \dots,p_m$.

\begin{flushleft}
\emph{Deuxième cas} : $\sum\limits_{\{i \vert u_\infty^i\ \neq c^{ste}\}} N_i =1$.
\end{flushleft}

Si le nombre de points nodaux sur les composantes non constantes de $u_\infty$ est égal à $1$, alors $u_\infty$ possède exactement une composante non constante et une composante constante. Puisque deux points marqués distincts sont envoyés sur des points distincts, il y a au plus un point marqué sur la composante constante. Comme la courbe $u_\infty$ est stable (voir la Définition~\ref{d:stable}), la composante constante est donc de genre strictement positif. Sans perte de généralité, on suppose que $u^1_\infty$ est la composante non constante. Comme les termes considérés sont pairs, on a 
$$0 = k_1 ( \ind(v^1) -2m) + k_1 \chi(\L_1) -\chi(\S_1) +2m(k_1-1).$$
Puisque c'est une somme de termes positifs, on a $k_1 \chi(\L_1) -\chi(\S_1)=0$ et $\ind(v^1) -2m =0$. On a également $2m = 
\ind(u_\infty^1) = - \chi(\Sigma_1)+ 2k_1 c_1 ([v^1])$ et $2m+2 = 
\ind(u) = - \chi(\Sigma)+ 2k_1 c_1 ([v^1])$. On a alors $\chi(\Sigma_1) = \chi(\Sigma) +2$, donc d'après le Lemme~\ref{l:morceauxmegarelou}, on a $k_1=1$, ou bien $k_1=2$, $[u]^2 =4$ et $v^1$ est plongée (auquel cas on a $g=4$ et on peut conclure). Mais si $k_1 = 1$, puisque $J$ est générique pour les points $p_1, \dots, p_m$ et $v^1$ est une courbe $J$--holomorphe d'indice contraint nul, la courbe $v^1$ est immergée. Ceci contredit l'existence de la singularité cuspidale provenant de la composante constante de genre strictement positif (voir la Remarque~\ref{r:cuspisgenus}).

\begin{flushleft}
\emph{Troisième cas} : $\sum\limits_{\{i \vert u_\infty^i\ \neq c^{ste}\}} N_i =0$.
\end{flushleft}

Si $u_\infty$ ne possède pas de point nodal, alors $u_\infty$ possède une unique composante, qui est nécessairement non constante et de genre $g$. On a dans ce cas 
$$2 = k_1 ( \ind(v^1) -2m) + k_1 \chi(\L_1) -\chi(\S_1) +2m(k_1-1).$$
Si $k_1 \chi(\L_1) -\chi(\S_1) =0$, on a  $k_1 =1$ d'après le Lemme~\ref{l:morceaux1}. On en conclut que $u_\infty$ est une courbe de genre $g$ simple homologue à $u$, donc plongée, qui appartient par conséquent à $\M_{g,m}^*([u];J;p_1, \dots, p_m)$, ce qui contredit les hypothèses. 
Si $k_1 \chi(\L_1) -\chi(\S_1) =2$ et $m \neq 0$, alors $k_1=1$. Ceci signifie que $\varphi_1$ est un biholomorphisme, et on a donc $\chi(\L_1) =\chi(\S_1)$, ce qui est absurde. Si $k_1 \chi(\L_1) -\chi(\S_1) =2$ et $m =0$, on a d'après le Lemme~\ref{l:morceauxrelou}, $k=2$ et $[u]^2=4$. Comme $u$ est une courbe plongée, on a $2 = \ind (u) = \chi(\Sigma) +2[u]^2$, donc $g=4$.
\end{proof}

\begin{rk}\label{r:morceauxparam}
Les arguments invoqués dans cette annexe sont également valables pour des chemins génériques de structures presque complexes car, par parité des indices, pour tout courbe pseudoholomorphe $u$, $\ind (u) \geq-1$ implique $\ind (u) \geq 0$ (voir le Corollaire~\ref{c:indgeq-1} et Remarque~\ref{r:indgeq-1à0}). La Proposition~\ref{p:morceauxcourbesind2} est donc également valable pour des espaces de modules à paramètre définis à partir de chemins génériques de structures presque complexes.
\end{rk}

%
\clearemptydoublepage

\backmatter

\addcontentsline{toc}{chapter}{Bibliographie}
\nocite{*} 
\printbibliography

\clearemptydoublepage

\markboth{}{}
 \AddToShipoutPicture{\put(0,0){\includegraphics[width=\paperwidth,height=\paperheight]{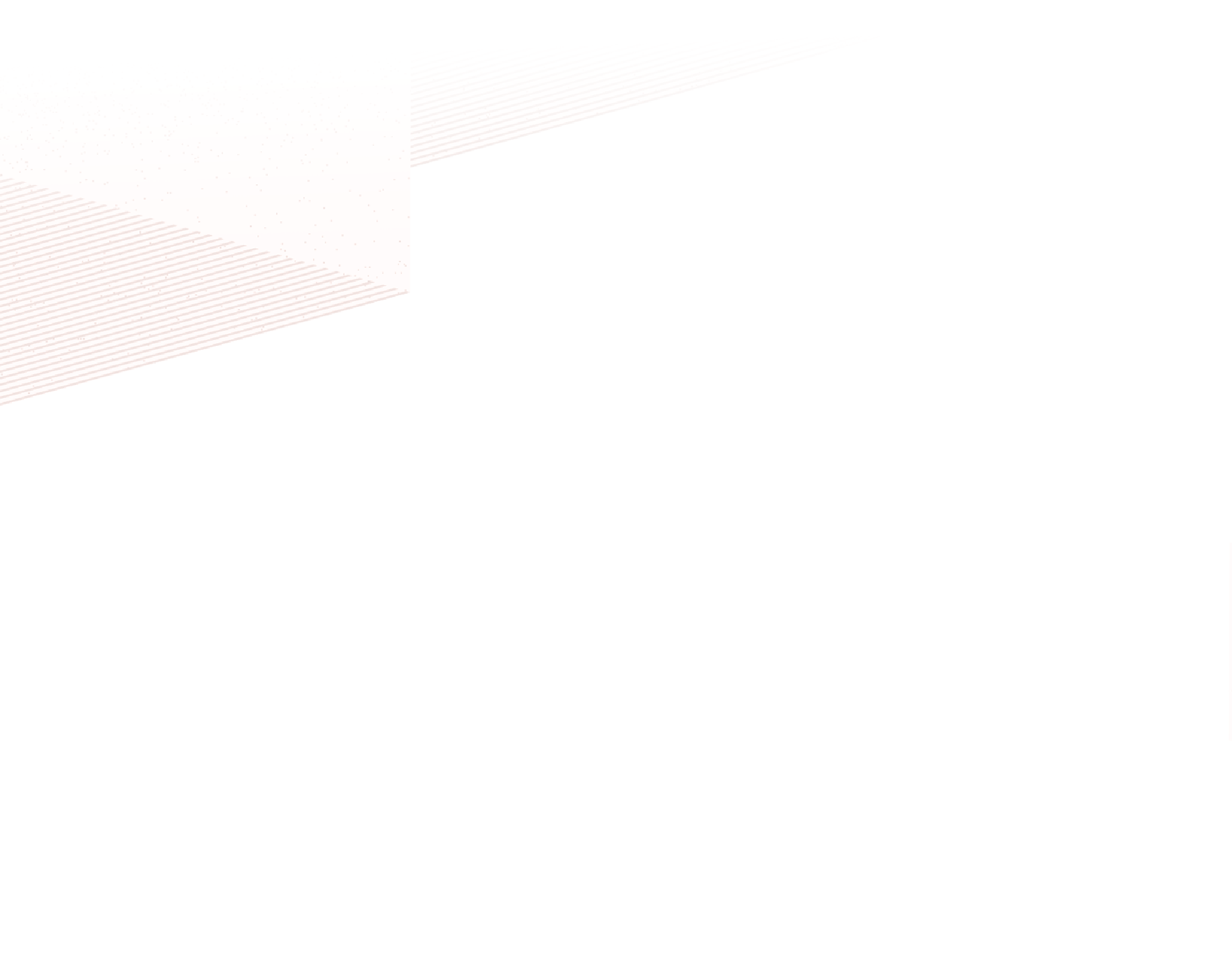}}}
\pagestyle{empty}
\vspace{-2cm}
\hspace{0.05cm}
\logouniversite{./Couverture-these/MathSTIC/logo-mathSTIC} 
\hspace{1cm}\logoetablissement{./Couverture-these/MathSTIC/logo-etablissements/logo-UN_noir}
\begin{tikzpicture}[remember picture,overlay,line width=0mm]
  \draw [draw=white,fill=white]
    (current page.north west) rectangle (\paperwidth,1);
  \node[xshift=0\paperwidth,yshift=2cm,text=white,font=\bf\Large] {
  \includegraphics[scale=1.2]{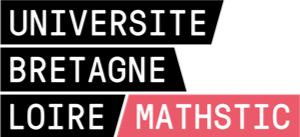}
  };
  \node[xshift=.45\paperwidth,yshift=2cm,text=white,font=\bf\Large] {
  \includegraphics[scale=0.5]{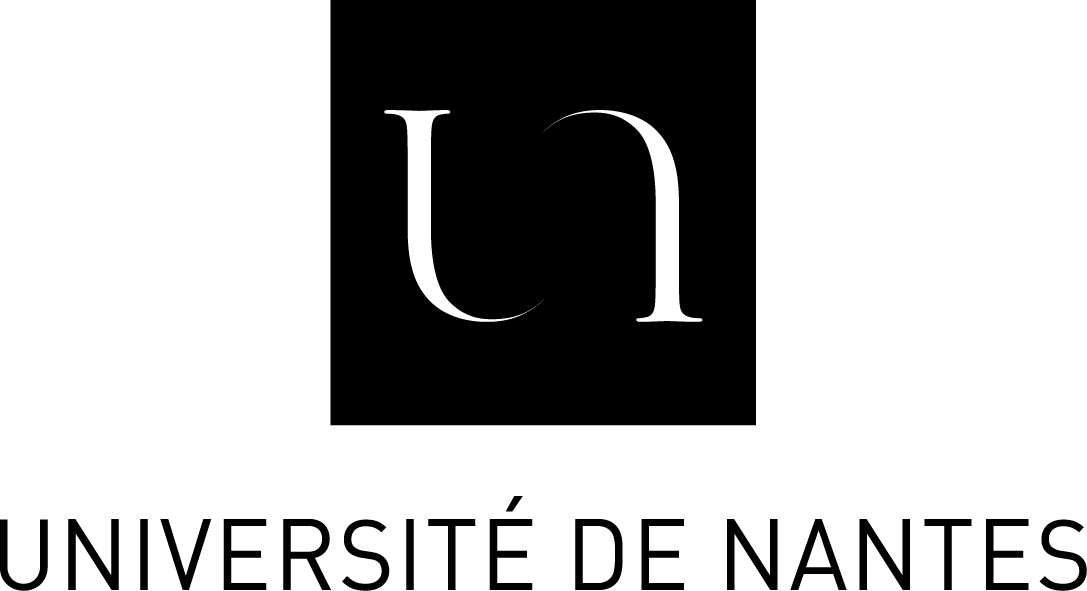}
  };
\end{tikzpicture}
\par\nobreak
\hspace{- 1.75cm}\noindent \textcolor{mathSTIC-Color}{\rule{\textwidth }{0.2cm}}  
\section*{\textcolor{mathSTIC-Color}{Titre}: Courbes symplectiques de haute auto-intersection dans les surfaces symplectiques}
\noindent \keywordsF{Variétés symplectiques de dimension $4$, Surfaces réglées, Courbes pseudoholomorphes, Remplissages symplectiques, Isotopies symplectiques.}

\begin{multicols}{2}
\noindent \textbf{Resum\'{e} : }
On étudie dans un premier temps les courbes symplectiquement plongées dans les surfaces symplectiques dont les nombres d'auto-intersection sont suffisamment grands par rapport leurs genres. On montre de deux manières différentes qu'une telle courbe détermine à la fois la classe de difféomorphisme de la surface symplectique qui la contient et la manière dont elle est plongée dans cette surface. La première démonstration fait appel à la théorie de Seiberg--Witten, alors que la seconde se restreint aux techniques pseudoholomorphes. On déduit de ce résultat l'unicité à difféomorphisme près des remplissages symplectiques forts des variétés de contact de dimension $3$ naturellement associées à ce type de courbes.

Dans un second temps, on s'intéresse aux sections symplectiques des surfaces complexes géométriquement réglées au-dessus de courbes elliptiques. On montre qu'une telle section est symplectiquement isotope à une section complexe.

\end{multicols}
\hspace{- 1.75cm}\noindent \textcolor{mathSTIC-Color}{\rule{\linewidth}{0.2cm}}
\section*{\textcolor{mathSTIC-Color}{Title}: Symplectic curves with high self-intersection in symplectic surfaces}
\noindent \keywordsE{Symplectic $4$-manifolds, Ruled surfaces, Pseudoholomorphic curves, Symplectic fillings, Symplectic isotopy.}

\begin{multicols}{2}
\noindent \textbf{Abstract : }
We first study symplectically embedded curves in symplectic surfaces with high self-intersection numbers compared to their genus. We prove in two different ways that such a curve completely determines both the diffeomorphism type of the surface in which it is embedded and the embedding itself. The first proof uses Seiberg--Witten theory whereas the second one only involves pseudoholomorphic techniques. We deduce from this result that the contact $3$--manifolds naturally associated with those curves admit a unique strong symplectic filling up to diffeomorphism.

We next examine symplectic sections of geometrically ruled complex surfaces over elliptic curves. We show that such a section is symplectically isotopic to a complex section.
\end{multicols}

\end{document}